\documentclass[12pt]{article}

\usepackage{latexsym,amsmath,amssymb,amsfonts,amscd,multirow,dsfont}
\usepackage{epsfig,verbatim,epstopdf,graphics}
\usepackage{bm,color}
\usepackage{subfigure}

\extrafloats{100}

\newtheorem{thm}{Theorem}[section]

\newtheorem{exa}[thm]{Example}

\newcommand{\bk}{\mathbf{k}}
\newcommand{\bx}{\mathbf{x}}

\newcommand{\eps}{\varepsilon}
\newcommand{\bE}{\mathbf{E}}

\newfont{\iams}{msbm9}

\newcommand{\commentbis}[1]{}
\newcommand{\be}{\begin{eqnarray}}
\newcommand{\ee}{\end{eqnarray}}
\newcommand{\beno}{\begin{eqnarray*}}
\newcommand{\eeno}{\end{eqnarray*}}

\newcommand{\barr}[1]{\begin{array}{#1}}
\newcommand{\earr}{\end{array}}

\newcommand{\beq}{\begin{equation}}
\newcommand{\eeq}{\end{equation}}
\newcommand{\beqa}{\begin{eqnarray}}
\newcommand{\eeqa}{\end{eqnarray}}

\newcommand{\df}{\partial}

\newcommand{\bv}{{\bf v}}
\newcommand{\bB}{{\bf B}}
\newcommand{\bJ}{{\bf J}}

\newcommand{\bU}{{\bf U}}
\newcommand{\bV}{{\bf V}}
\newcommand{\bn}{{\bf n}}

\newcommand{\Ox}{{\Omega_x}}
\newcommand{\Oxi}{{\Omega_\xi}}

\newcommand{\mE}{{\mathcal E}}

\newcommand{\mG}{{\mathcal G}}

\newcommand{\mU}{{\mathcal U}}

\newcommand{\bzero}{\mathbf{0}}
\newcommand{\bone}{\mathbf{1}}
\newcommand{\bl}{\mathbf{l}}
\newcommand{\bi}{\mathbf{i}}

\newcommand{\bj}{\mathbf{j}}
\newcommand{\bW}{\mathbf{W}}

\newcommand{\bal}{{\bm{\alpha}}}
\newcommand{\bb}{{\bm{\beta}}}

\title
{Sparse Grid Discontinuous Galerkin Methods for the  Vlasov-Maxwell System}

\author{
Zhanjing Tao
\thanks{Department of Mathematics, Michigan State University,
East Lansing, MI 48824 U.S.A.
 {\tt tzjnchy555@math.msu.edu}}
\and
 Wei Guo
\thanks{Department of Mathematics and Statistics, Texas Tech University, Lubbock, TX, 70409 U.S.A. {\tt weimath.guo@ttu.edu}. Research
	is supported by NSF grants DMS-1620047, DMS-1830838.}
\and
Yingda Cheng
\thanks{Department of Mathematics, Department of  Computational Mathematics, Science and Engineering, Michigan State University,
East Lansing, MI 48824 U.S.A.
 {\tt ycheng@msu.edu}. Research is supported by NSF grants  DMS-1453661, DMS-1720023 and the Simons Foundation under award number 558704.}
}

\date{\today}

\begin{document}

\maketitle

\begin{abstract}

In this paper, we develop  sparse grid discontinuous Galerkin (DG) schemes for  the Vlasov-Maxwell (VM) equations. The VM system is   a fundamental kinetic model in plasma physics, and its numerical computations are quite demanding, due to its intrinsic high-dimensionality and the need to retain many properties of the physical solutions.   To break the curse of dimensionality, we consider the sparse grid DG methods that were recently developed in \cite{guo2016sparse,guo2017adaptive} for transport equations. 
Such methods are based on multiwavelets  on tensorized nested grids and can significantly reduce the numbers of degrees of freedom.  We formulate two versions of the schemes: sparse grid DG and adaptive sparse grid DG methods for the VM system. Their key properties and implementation details are discussed. Accuracy and robustness  are demonstrated by  numerical tests,  with emphasis on   comparison of the performance of the two methods, as well as with their full grid counterparts. 

\end{abstract}


{\bf Keywords:}
discontinuous Galerkin methods;   sparse grids;  Vlasov-Maxwell system; streaming Weibel instability; Landau damping. 

\section{Introduction}

The Vlasov-Maxwell (VM) system is   a fundamental model in plasma physics for describing the dynamics of collisionless magnetized plasmas, which finds diverse applications in science and engineering,
including thermo-nuclear fusion, satellite amplifiers, high-power microwave generation, to name a few.
In this paper, we study the VM system that describes the evolution of a single species of nonrelativistic electrons under the self-consistent electromagnetic field while the ions are treated as uniform fixed background. Under the scaling of the characteristic time by the inverse of the plasma frequency $\omega_p^{-1}$, length by the Debye length $\lambda_D$, and  electric and magnetic fields by  $-m c \omega_p/e$ (with $m$ the electron mass, $c$ the speed of light, and $e$ the electron charge),  the  dimensionless form of the VM system is
\begin{subequations}
\begin{align}
&\partial_t f  + \xi \cdot \nabla_{\bx} f  +  (\bE + \xi \times \bB) \cdot \nabla_\xi f = 0~, \label{eq:vlasov}\\
&\frac{\df \bE}{\df t} = \nabla_{\bx}  \times \bB -  \bJ, \qquad
\frac{\df \bB}{\df t} = -  \nabla_{\bx}   \times \bE~,  \quad\label{eq:max:2} \\
&\nabla_\bx \cdot \bE = \rho-\rho_i, \qquad
\nabla_\bx \cdot \bB = 0~, \quad
\label{eq:max:4}\\
&f(\bx, \xi, 0)=f_0(\bx, \xi), \quad \bE(\bx,0)=\bE_0(\bx), \quad \bB(\bx,0)=\bB_0(\bx),
\end{align}
\end{subequations}
with
\begin{equation*}
\rho(\bx, t)= \int_\Oxi f(\bx, \xi, t)d\xi,\qquad \bJ(\bx, t)=  \int_\Oxi f(\bx, \xi, t)\xi d\xi~,
\end{equation*}
where the equations are defined on  $\Omega=\Ox\times\Oxi.$  $\mathbf{x} \in \Ox$ denotes position in  physical space, and $\xi \in \Oxi,$ which is the  velocity space. Here
$f(\bx, \xi, t)\geq 0$ is the distribution function of electrons at position $\bx$ with velocity $\xi$ at time $t$,  $\bE(\bx, t)$ is the electric field, $\bB(\bx, t)$ is the magnetic field, $\rho(\bx, t)$ is the electron charge density, and $\bJ(\bx, t)$ is the current density.  The charge  density of background ions is denoted by $\rho_i$, which is chosen to satisfy total charge neutrality,  $\int_\Ox \left(\rho(\bx, t)-\rho_i \right)\,d\bx=0$.   Ideally $\Oxi=\mathbb{R}^{d_\xi},$ however numerical computation requires a truncation of the space $\Oxi$ and the assumption that $f$ is compactly supported on $\Oxi.$  In this paper, for simplicity, we will assume  $\Ox, \Oxi$ to be box-shaped domains.  


The simulations of VM systems are quite challenging.  Particle-in-cell (PIC) methods  \cite{Birdsall_book1991, Hockney_book1981} have long been very popular  numerical tools, 
mainly because they can generate reasonable results with relatively low computational cost. 
However, as a Monte-Carlo type approach, the PIC methods are known to suffer from the statistical noise, which is $O(N^{-\frac12})$ with $N$ being the number of sampling particles. Such an inherent low order error of PIC methods prevents accurate description of physics of interest, when, for instance, the tail of the distribution function needs to be resolved accurately.  In recent years, there has been growing interest in deterministic simulations of the Vlasov equation. In the deterministic framework, the schemes are free of statistical noise and hence able to generate highly accurate results in phase space.  In the context of VM simulations, Califano \emph{et al.} have    used a semi-Lagrangian approach to compute the  streaming Weibel (SW) instability \cite{califano1998ksw}, current filamentation instability \cite{mangeney2002nsi}, magnetic vortices \cite{califano1965ikp}, magnetic reconnection \cite{califano2001ffm}. Also, various methods have  been  proposed for the relativistic VM system \cite{Sircombe20094773, Besse20087889, Suzuki20101643, Huot2003512}.
 In this paper, we are interested in a class of successful deterministic Vlasov solvers based on the discontinuous Galerkin (DG) finite element discretization, because of not only  their provable convergence
 and accommodation for adaptivity and parallel implementations, but also their excellent conservation property and superior performance in long time wave-like simulations. Those distinguishing properties of DG methods  are very much desired for the VM simulations, and they have been previously employed to solve VM system \cite{cheng2014discontinuous,cheng2014energy} and the relativistic VM system \cite{yang2017discontinuous}. However, due to the curse of dimensionality, traditional deterministic approaches including the DG methods are not competitive for high dimensional Vlasov simulations, even  
 with the aid of high performance computing systems.

 
 To break the curse of dimensionality, this paper will focus on the sparse grid approach. 
The sparse grid method \cite{bungartz2004sparse, garcke2013sparse} has long been an effective numerical tool to  reduce the degrees
of freedom for high-dimensional grid based methods.   In the context of wavelets or sparse grid methods for kinetic transport equations,  we  mention the work of using wavelet-MRA methods for Vlasov equations \cite{besse2008wavelet}, the combination technique for linear gyrokinetics \cite{kowitz2013combination},  sparse adaptive finite element method \cite{widmer2008sparse}, sparse discrete ordinates method
\cite{grella2011sparseo} and sparse tensor spherical harmonics \cite{grella2011sparse} for radiative transfer, among many others. In \cite{sparsedgelliptic, guo2016sparse},   a class of sparse grid DG schemes were proposed for solving high-dimensional partial differential equations (PDEs)   based on a novel sparse DG finite element approximation space. Such a sparse grid space can be regarded as a proper truncation of the standard tenor approximation space, which reduces degrees of freedom of from $O(h^{-d})$ to $O(h^{-1}|\log_2h|^{d-1})$, where $h$ is  the uniform mesh size in each direction and $d$ is the dimension of the problem. Coupling the DG weak formulation with the sparse grid space, the resulting sparse grid DG method is demonstrated to save computational and storage cost without deteriorating the approximation quality by much. In particular, when applied to a $d$-dimensional linear transport equation, the scheme is proven to be $L^2$ stable and convergent with rate $O(|\log_2h|^{d}h^{{k+1/2}})$ in $L^2$ norm for a smooth enough solution, where $k$ is the degree of polynomials \cite{guo2016sparse}.  Motivated by the development of sparse grid DG method \cite{guo2016sparse} and the adaptive multiresolution DG method   \cite{guo2017adaptive}  for transport equations, it is of interest to this paper to develop  sparse grid DG methods for solving the  VM system. 
The proposed methods are well suited for VM simulations, due to their ability to handle   high dimensional convection dominated problems, the ability to capture the main structures of the solution with feasible computational resource and the overall good performance in conservation of    physical quantities in long time simulations.  

  The rest of the paper is  organized as follows: in Section \ref{sec:method}, we describe the numerical algorithms, outlining the schemes as well as their key properties.   Section \ref{sec:numerical}  is devoted to discussions of simulation results. We conclude the paper in Section \ref{sec:conclusion}.

\section{Numerical methods}
\label{sec:method}

In this section, we  describe two sparse grid DG methods for
the VM system: the standard sparse grid DG method and the adaptive sparse grid DG method.
We first review the finite element space on sparse grid introduced in \cite{sparsedgelliptic}, and then describe the details of the schemes when applied to the VM system. 

\subsection{DG finite element space on sparse grid}

In this subsection, we   review the notations of DG finite element space on sparse grid. 
  First, we introduce the hierarchical decomposition of piecewise polynomial space in one dimension. Without   loss of generality, consider the  interval $[0,1]$, we define a set of nested grids, where the $l$-th level grid $\Omega_l$ consists of $2^l$ uniform cells $I_{l}^j=(2^{-l}j, 2^{-l}(j+1)]$, $j=0, \ldots, \max(0, 2^l-1),$ for any $l \ge 0.$
The nested grids result in the nested piecewise polynomial spaces. In particular, let
$$V_l^k:=\{v: v \in P^k(I_{l}^j),\, \forall \,j=0, \ldots, 2^l-1\}$$
  be the usual piecewise polynomials of degree at most $k$ on the $l$-th level grid $\Omega_l$. Then, we have $$V_0^k \subset V_1^k \subset V_2^k \subset V_3^k \subset  \cdots$$
We can now define the multiwavelet subspace $W_l^k$, $l=1, 2, \ldots $ as the orthogonal complement of $V_{l-1}^k$ in $V_{l}^k$ with respect to the $L^2$ inner product on $[0,1]$, i.e.,
\begin{equation*}
V_{l-1}^k \oplus W_l^k=V_{l}^k, \quad W_l^k \perp V_{l-1}^k.
\end{equation*}
Note that the space $W_l^k$, for $l>0$, represents the finer level details when the grid is refined and this is the key to the reduction of  the degrees of freedom in higher dimensions.
We further let
 $W_0^k:=V_0^k$, which is standard piecewise polynomial space of degree $k$ on $[0,1]$.  Therefore, we have found a hierarchical representation of  the standard piecewise polynomial space  $V_l^k$ on $\Omega_l$ as $V_l^k=\bigoplus_{0 \leq j\leq l} W_j^k$. 
In \cite{sparsedgelliptic}, we used  the orthonormal basis functions for space $W_l^k,$ which are constructed based on the one-dimensional orthonormal multiwavelet bases first introduced in \cite{alpert1993class}. We refer readers to \cite{alpert1993class} and \cite{sparsedgelliptic} for more details. We now denote the basis functions in level $l$ as 
$$v^j_{i,l}(x),\quad i=1,\ldots,k+1,\quad j=0,\ldots,2^{l-1}-1,$$
and they are orthonormal, i.e. 
 $\int_0^1 v^j_{i,l}(x)v^{j'}_{i',l'}(x)\,dx=\delta_{ii'}\delta_{ll'}\delta_{jj'}.$
 
 Below we give an example of the basis functions. We denote the bases for $W_1^k$   as 
$$h_i(x)= 2^{{1}/{2}} f_i(2x-1), \quad i = 1, \ldots, k+1,$$
where    $ \{f_i(x), \, i=1, \ldots, k+1 \}$ are functions supported on $[-1,1]$ depending on $k$.
For instance, when $k=3$, the formulas for $f_i(x), \, i=1, \ldots, 4$ restricted to the interval $(0, 1)$ are
\begin{align*}
f_1(x) =& \sqrt{\frac{15}{34}} (1 + 4x -30x^2 + 28x^3) \\
f_2(x) =& \sqrt{\frac{1}{42}} (-4 + 105x -300x^2 + 210x^3) \\
f_3(x) =& \frac{1}{2}\sqrt{\frac{35}{34}} (-5 + 48x - 105x^2 + 64x^3) \\
f_4(x) =& \frac{1}{2}\sqrt{\frac{5}{42}} (-16 + 105x - 192x^2 + 105x^3).
\end{align*}
The functions $f_i$ are extended to $(-1, 0)$ as   even or odd functions according to the parity of $i+k$:
$$f_i (x) = (-1)^{i+k} f_i (-x), \, x\in(-1,0),$$ and are zero outside $(-1,1)$.
Then, the bases for $W_l^k, \,l \geq 1$ are defined as
$$  v_{i,l}^j (x) = 2^{(l-1)/2} \,h_i(2^{l-1}x - j), \quad i = 1, \ldots, k+1,\, j=0, \ldots, 2^{l-1}-1.$$

 \medskip
Now we are ready to prescribe the sparse finite element space in $d$-dimensions on $\Omega=[0,1]^d$. Notice similar discussions apply for any box-shaped domain in $d$ dimensions. 
First we recall some basic notations about multi-indices. For a multi-index $\mathbf{\alpha}=(\alpha_1,\cdots,\alpha_d)\in\mathbb{N}_0^d$, where $\mathbb{N}_0$  denotes the set of nonnegative integers, the $l^1$ and $l^\infty$ norms are defined as 
$$
|\bal|_1:=\sum_{m=1}^d \alpha_m, \qquad   |\bal|_\infty:=\max_{1\leq m \leq d} \alpha_m.
$$
The component-wise arithmetic operations and relational operations are defined as
$$
\bal \cdot \bb :=(\alpha_1 \beta_1, \ldots, \alpha_d \beta_d), \qquad c \cdot \bal:=(c \alpha_1, \ldots, c \alpha_d), \qquad 2^\bal:=(2^{\alpha_1}, \ldots, 2^{\alpha_d}),
$$
$$
\bal \leq \bb \Leftrightarrow \alpha_m \leq \beta_m, \, \forall m,\quad
\bal<\bb \Leftrightarrow \bal \leq \bb \textrm{  and  } \bal \neq \bb.
$$

By making use of the multi-index notation, we indicate by $\bl=(l_1,\cdots,l_d)\in\mathbb{N}_0^d$ the mesh level in a multivariate sense, where $\mathbb{N}_0$ denotes non-negative integers. We  consider the tensor-product rectangular grid $\Omega_\bl=\Omega_{l_1}\otimes\cdots\otimes\Omega_{l_d}$ with mesh size $h_\bl=(h_{l_1},\cdots,h_{l_d}).$ Based on the grid $\Omega_\bl$, an elementary cell is denoted by $I_\bl^\bj=\{\bx:x_m\in(h_mj_m,h_m(j_{m}+1)),m=1,\cdots,d\}$, and 
$$\bV_\bl^k(\Omega):=\{\bv: \bv(\bx) \in Q^k(I^{\bj}_{\bl}), \,\,  \bzero \leq \bj  \leq 2^{\bl}-\bone \}= V_{l_1,x_1}^k\times\cdots\times  V_{l_d,x_d}^k$$
is the tensor-product piecewise polynomial space, where $Q^k(I^{\bj}_{\bl})$ denotes  polynomials of degree up to $k$ in each dimension on cell $I^{\bj}_{\bl}$. 
If $l_1=\ldots=l_d=N$, the  grid and space will be denoted by $\Omega_N$ and $\bV_N^k$, respectively.  
For the increment space $\bW_\bl^k=W_{l_1,x_1}^k\times\cdots\times  W_{l_d,x_d}^k,$
the orthonormal basis functions  can be defined as
$$v^\bj_{\bi,\bl}(\bx)\doteq\prod_{m=1}^d v^{j_m}_{i_m,l_m}(x_m),\quad i_m=1,\ldots,k+1,\,j_m=0,\ldots,\max(0,2^{l_m-1}-1),$$
where $v^{j_m}_{i_m,l_m}(x_m)$ denote orthonormal bases in m-th dimension defined in one-dimensional case.  
With the one-dimensional hierarchical decomposition,  we have 
$$
\bV_N^k(\Omega)=\bigoplus_{\substack{ |\bl|_\infty \leq N}} \bW_\bl^k.
$$
The sparse finite element approximation space on mesh $\Omega_N$ we use in this paper, on the other hand,  is  defined by \cite{sparsedgelliptic, guo2016sparse}
$$ \hat{\bV}_{N}^k(\Omega):=\bigoplus_{\substack{ |\bl|_1 \leq N}}\bW_\bl^k.$$
This is a subset of   $\bV_N^k(\Omega)$, and   its number of degrees of freedom  scales as $O((k+1)^d2^NN^{d-1})$ or $O((k+1)^d h_N^{-1} (\log h_N)^{d-1})$, where $h_N=2^{-N}$ denotes the finest mesh size in each direction \cite{sparsedgelliptic}. This is significantly less than that of $\bV_N^k(\Omega)$ with $O((k+1)^d h_N^{-d})$ number of elements when $d$ is large.
The sparse grid DG scheme  in Section \ref{subsec:formulation}  is defined using this space. The adaptive sparse grid DG scheme, however, relies on an adaptive choice of the space, and will be discussed in details in Section \ref{subsec:adaptive}.

\subsection{The  sparse grid DG scheme}
\label{subsec:formulation}

Below, we formulate a DG scheme in the sparse finite element space   for the VM system  \eqref{eq:vlasov}-\eqref{eq:max:2} inspired by \cite{cheng2014discontinuous}. On the PDE level, the two equations in \eqref{eq:max:4} involving the divergence of the magnetic and electric fields can be derived from the remaining part of the VM system.   However, how to impose \eqref{eq:max:4} numerically is an important and nontrivial issue \cite{munz2000divergence, barth2006role,jacobs2009implicit}. This will be studied in our future work.
Using the notations introduced in the previous subsections, and let $d_x$ and $d_\xi$ be  the dimension of $\Ox$ and $\Oxi$, respectively, we consider the partitions of the domain $\Omega$ into mesh level $N$ in all directions.
We distinguish between the $x$- and $\xi$-directions. Let $I^{\bj}_{x,N},  0 \leq j_m  \leq 2^{N}-1, \forall \,m=1, \ldots, d_x$ and $I^{\bj}_{\xi,N},  0 \leq j_m  \leq 2^{N}-1, \forall \,m=1, \ldots, d_\xi$ be the collection of all elements in $\Ox$ and $\Oxi$, respectively.
Let $\mE_x$ be the union of the edges for all elements in $I^{\bj}_{x,N}$, similarly let $\mE_\xi$ be the union of the edges for all elements in $I^{\bj}_{\xi,N}$. 
Furthermore, $\mE_\xi=\mE_\xi^i\cup\mE_\xi^b$ with $\mE_\xi^i$ and $\mE_\xi^b$ being the union of interior and boundary edges of $I^{\bj}_{\xi,N}$, respectively.

	For piecewise  functions, we further introduce the jumps and averages as follows. Suppose $T^+$ and $T^-$ are two elements in $I^{\bj}_{x,N}$.  For any edge $e=\{T^+\cap T^-\}\in\mE_x$, with $\bn_x^\pm$ as the outward unit normal to $\partial T^\pm$,
	$g^\pm=g|_{\partial T^\pm}$, and $\bU^\pm=\bU|_{\partial T^\pm}$, the jumps across $e$ are defined  as
	\begin{equation*}
	[g]_x={g^+}{\bn_x^+}+{g^-}{\bn_x^-},\qquad [\bU]_x={\bU^+}\cdot{\bn_x^+}+{\bU^-}\cdot{\bn_x^-},\qquad [\bU]_\tau=\bU^+\times\bn_x^++\bU^-\times\bn_x^-
	\end{equation*}
	and the averages are
	\begin{equation*}
	\{g\}_x=\frac{1}{2}({g^+}+{g^-}),\qquad \{\bU\}_x=\frac{1}{2}({\bU^+}+{\bU^-}).
	\end{equation*}
	When considering periodic boundary conditions, the jumps and averages can be naturally defined on the boundary edges.
	
	By replacing the subscript $x$ with $\xi$, the jumps and averages can be defined similarly for the $\xi$-direction. For a boundary edge $e\in\mE^b_\xi$ with $\bn_\xi$ being the outward unit normal, we use
	\begin{equation}
	[g]_\xi={g}{\bn_\xi},\qquad
	\{g\}_\xi=\frac{1}{2}g~.
	\label{eq:jump:ave:b}
	\end{equation}
	This is consistent with the fact that the exact solution $f$ is assumed to be compactly supported in $\xi$-domain.

We are now ready to describe the scheme. The    sparse discrete spaces on $\Omega$ and $\Omega_x$ we use are defined as
$$\hat \mG_h^{k} = \hat{\bV}_{N}^k(\Omega), \quad \hat \mU_h^k = [\hat{\bV}_{N}^k(\Ox)]^{d_x}.$$ 
Following \cite{guo2016sparse}, the semi-discrete DG methods for the VM system are: to find $f_h\in \hat \mG_h^k$, $\bE_h, \bB_h\in \hat \mU_h^k$, such that for any $g\in \hat \mG_h^k$, $\bU, \bV\in \hat \mU_h^k$,
\begin{subequations}
	\begin{align}
	\int_\Omega\df_t f_h g d\bx d\xi
	&- \int_\Omega f_h\xi \cdot \nabla_\bx g d\bx d\xi 
	- \int_\Omega f_h(\bE_h+\xi\times\bB_h)\cdot\nabla_\xi  g d\bx d\xi\notag\\
	&+ \int_{\Oxi}\int_{\mE_x} \widehat{f_h \xi} \cdot [g]_x ds_x d\xi + \int_{\Ox} \int_{\mE_\xi} \widehat{f_h (\bE_h+\xi\times \bB_h) } \cdot [g]_\xi ds_\xi dx
	=0~,\label{eq:scheme:1}\\
	\int_\Ox\df_t\bE_h\cdot\bU d\bx
	&=\int_\Ox \bB_h\cdot \nabla_\bx\times\bU d\bx+\int_{\mE_x}\widehat{\bB_h}\cdot [\bU]_\tau ds_x
	-\int_\Ox\bJ_h\cdot\bU d\bx~,\label{eq:scheme:2}\\
	\int_\Ox\df_t\bB_h\cdot\bV d\bx
	&=-\int_\Ox\bE_h\cdot \nabla_\bx\times\bV d\bx-\int_{\mE_x}\widehat{\bE_h}\cdot [\bV]_\tau ds_x~,\label{eq:scheme:3}
	\end{align}
\end{subequations}
with
\beq
\label{eq:J}
{\bJ}_h(\bx, t)=\int_{\Oxi} f_h(\bx, \xi, t)\xi d\xi~\in   \hat \mU_h^k.
\eeq
%
All  ``hat"  functions are numerical fluxes. For the Vlasov part,  we adopt the global Lax-Friedrichs flux:
\begin{subequations}
	\begin{align}
	\widehat{f_h \xi}:& = \{f_h\xi\}_x+\frac{\alpha_1}{2}[f_h]_x~,\label{eq:flux:1}\\
	\widehat{f_h (\bE_h+\xi\times \bB_h)}:& = \{f_h(\bE_h+\xi\times\bB_h)\}_\xi+\frac{\alpha_2}{2}[f_h]_\xi~,\label{eq:flux:3}
	\end{align}
\end{subequations}
where $\alpha_1=\max|\xi\cdot\bn_x|$ and  $\alpha_2=\max |(\bE_h+\xi\times\bB_h)\cdot\bn_\xi|$, where the maximum is taken for all possible $\bx$ and $\xi$ at time $t$ in the computational domain.
For the Maxwell part, we use the upwind flux \begin{subequations}
	\begin{align}
	\widehat{\bE_h}:&
	= \{\bE_h\}_x+\frac{1}{2}[\bB_h]_\tau~,\label{eq:flux:4}\\
	\widehat{\bB_h}:& 
	=\{\bB_h\}_x-\frac{1}{2}[\bE_h]_\tau~,
	\label{eq:flux:5}
	\end{align}
\end{subequations}
or the   alternating fluxes 
\begin{equation}
 \widehat{\bE_h}:=\bE_h^+, \; \widehat{\bB_h}:=\bB_h^-, \;\mbox{or}\;
\widehat{\bE_h}:=\bE_h^-,\; \widehat{\bB_h}:=\bB_h^+~.\label{eq:flux:7}
\end{equation}
Previous studies in \cite{cheng2014energy,cheng2014discontinuous} have demonstrated the
importance of numerical flux on the quality of DG simulations. It was understood that a dissipative flux choice is desired for the Vlasov equation. The alternating and upwind fluxes are both optimal in order for the Maxwell solver, but alternating flux is energy-conserving for the Maxwell's equation, while upwind flux is not. 
We will not consider the central flux based on the recommendations made in \cite{cheng2014discontinuous}.

 \medskip

Below, we will summarize  the conservation and stability properties of the semi-discrete scheme. The proof is very similar to those in \cite{cheng2014discontinuous} and is thus omitted.

\begin{thm}[{Mass conservation}]
The numerical solution $f_h\in\hat{\mG}_h^k$ with $k\geq 0$ satisfies
\begin{equation}
\frac{d}{dt}\int_{\Omega} f_h d\bx d\xi+\Theta_{h,1}(t)=0~,
\label{eq:MassC}
\end{equation}
where
\begin{equation*}
\Theta_{h,1}(t)=\int_{\Ox}\int_{\mE_\xi^b} f_h \max((\bE_h+\xi\times\bB_h)\cdot\bn_\xi, 0)ds_\xi d\bx~.
\end{equation*}
Equivalently, with $\rho_h(\bx, t)=\int_{\Oxi} f_h(\bx,\xi,t) d\xi$, for any $T>0$, the following holds:
\begin{align}
\int_{\Ox} \rho_h(\bx, T)d\bx+\int_0^T\Theta_{h,1}(t)dt
=\int_{\Ox} \rho_h(\bx, 0) d\bx~.
\label{eq:MassC:1}
\end{align}
\label{lem:MassC}
\end{thm}

\begin{thm}[{Energy conservation}]
\label{lem:EneC}
For $k\geq 2$,   the numerical solution $f_h\in\hat{\mG}_h^k$, $\bE_h, \bB_h\in\hat{\mU}_h^k$ with the upwind numerical flux \eqref{eq:flux:4}-\eqref{eq:flux:5} for the Maxwell part satisfies
\begin{equation*}
\frac{d}{dt}\left(\int_{\Omega} f_h|\xi|^2 d\bx d\xi+\int_{\Ox} (|\bE_h|^2+|\bB_h|^2) d\bx\right)
+\Theta_{h,2}(t)+\Theta_{h,3}(t)=0~,
\label{eq:EneC}
\end{equation*}
with
\begin{equation*}
\Theta_{h,2}(t)=\int_{\mE_x}\left( |[\bE_h]_\tau|^2+|[\bB_h]_\tau|^2 \right) ds_x~, \qquad \Theta_{h,3}(t)
=\int_{\Ox}\int_{\mE_\xi^b}f_h|\xi|^2 \max((\bE_h+\xi\times\bB_h)\cdot\bn_\xi, 0)  ds_\xi d\bx~.
\end{equation*}
While for the scheme with alternating fluxes \eqref{eq:flux:7} for the Maxwell part, we have
\begin{equation*}
\frac{d}{dt}\left(\int_{\Omega} f_h|\xi|^2 d\bx d\xi+\int_{\Ox} (|\bE_h|^2+|\bB_h|^2) d\bx\right)+\Theta_{h,3}(t)=0~.
\end{equation*}
%
%
\end{thm}

In the theorems above, terms like $\Theta_{h,1}(t), \Theta_{h,3}(t)$ come  from numerical errors of the velocity boundary truncation from an infinite to a finite domain. If those terms are negligible, i.e. when $\Oxi$ is chosen sufficiently large, we can conclude that the scheme is mass-conservative, and energy-conservative if the alternating flux is used for the Maxwell solver.  

\begin{thm}[{$L^2$-stability of $f_h$}]
For $k\geq 0$, the numerical solution $f_h\in\hat{\mG}_h^k$ satisfies
\begin{equation}
\label{eq:L2Stab}
\frac{d}{dt}\left(\int_{\Omega} |f_h|^2 d\bx d\xi\right)\le0~.\notag
\end{equation}
\label{lem:stability}
\end{thm}

As for time, we use the  total variation diminishing  Runge-Kutta (TVD-RK) methods \cite{Shu_1988_JCP_NonOscill} to solve the ordinary differential  equations resulting from the   discretization \eqref{eq:scheme:1}-\eqref{eq:scheme:3}, denoted by $(u_h)_t = R(u_h).$ In particular, we use the following third-order TVD-RK method in this paper  
\begin{align}
u_h^{(1)} &= u^{n} + \Delta t R(u^n_h), \notag \\
u_h^{(2)} &= \frac{3}{4}u^{n} + \frac14 u_h^{(1)} +\frac14 \Delta t R(u_h^{(1)}),\label{eq:tvd} \\
u_h^{n+1} &= \frac{1}{3}u^{n} + \frac23 u_h^{(2)} +\frac23 \Delta t R(u_h^{(2)}), \notag
\end{align}
where $u_h^{n}=\left[f_h^n, \bE_h^n, \bB_h^n \right]^T$ represents the numerical solutions at time level $t=t^n$.

\bigskip

Finally, we mention some details about implementation. A key to  computational efficiency is the efficient evaluation of multidimensional integrations. To compute multidimensional integration, we apply the \emph{unidirectional principle}. For example, if we want to evaluate $\int_\Omega \phi(\bx) \, d\bx$ with $\phi(\bx)=\phi_1(x_1) \ldots \phi_d(x_d)$, it is equivalent to multiplication of one-dimensional integrals $\int_{[0,1]} \phi_1 \, dx_1 \cdots \int_{[0,1]} \phi_d \, dx_d.$   Based on the hierarchical   structure of the basis functions, we only need some small overhead to compute one-dimensional integrals and assemble them to obtain the multi-dimensional integrations.  
In \eqref{eq:scheme:1}-\eqref{eq:scheme:3}, the numerical integrations with coefficients $\xi$ and $\bE_h, \bB_h$ (which belong to  $\hat{\mU}_h^k$) can be done very efficiently because they can all be computed using this trick.  This procedure is performed one time before time evolution starts, and the sparsity of the matrices due to the multiwavelet basis structures is utilized to accelerate the computation  \cite{sparsedgelliptic}.

\subsection{The adaptive sparse grid DG scheme}
\label{subsec:adaptive}

In this subsection, we will describe the adaptive sparse grid DG (also called adaptive multiresolution DG) method   \cite{guo2017adaptive}  for the VM system.
The motivation to study the adaptive scheme is to offer better numerical resolutions for the probability density function. It is known that the solution to Vlasov equation develops filamentation, therefore the standard sparse grid method can not offer enough resolution when $t$ becomes large (see \cite{guo2017adaptive} for  comparison in the case of the Vlasov-Poisson system). Therefore,  in this paper, we also consider the adaptive scheme.

The main idea of the algorithm is not to use $ \hat{\bV}_{N}^k(\Omega)$ in a pre-determined fashion, but rather to choose a subspace of $ \bV_N^k(\Omega)$ adaptively. In this work, adaptivity is only implemented for $f_h,$ not for the lower-dimensional variables $\bE_h, \bB_h,$ which are computed using the full finite element space.
 This will not cause much additional cost because the Maxwell's equations are in lower dimension than the Vlasov equation. If the computational and storage cost is a big concern, then the adaptive strategy can be potentially applied  to the Maxwell's part as well.
 
The algorithm is described in details below.
Given a maximum mesh level $N$ and an accuracy threshold $\varepsilon>0$,  we first use the adaptive multiresolution projection algorithm \cite{guo2017adaptive}  to get the numerical initial condition $f_h(\bx, \xi)$ for the DG scheme.  For completeness of the paper, the details of the method is listed below in Algorithm 1.

\medskip

\noindent\rule{16.5cm}{1pt}\\
\noindent{\bf Algorithm 1: Adaptive multiresolution projection }\\
\noindent\rule{16.5cm}{0.4pt}

{\bf Input:} Function $f(\bx, \xi)$.

{\bf Parameters:} Maximum level $N,$ polynomial degree $k,$  error threshold   $\varepsilon.$

{\bf Output:} Hash table $H$, leaf table $L$ and projected  solution $f_h(\bx, \xi) \in \bV_{N,H}^k.$

\begin{enumerate}
	\item Project $f(\bx, \xi)$ onto the coarsest level of mesh, e.g., level 0. Add all elements to the hash table $H$ (active list). Define an element without children as a leaf element, and add all the leaf elements to the leaf table $L$ (a smaller hash table). 
	\item For each leaf element $V_\bl^\bj$ in the leaf table, if the refinement criteria 
	\begin{align}
	&\left(\sum_{\mathbf{1}\leq\bi\leq\bk+\mathbf{1}}|f^\bj_{\bi,\bl}|^2\right)^{\frac12}>\varepsilon,\label{eq:l2}
	\end{align}
	holds, then we consider its child  elements: for a child element $V_{\bl'}^{\bj'}$, if it has not been added to the table $H$, then  compute the detail coefficients $ \{f^{\bj'}_{\bi,\bl'}, \mathbf{1}\leq\bi\leq\bk+\mathbf{1}\}$  and add $V_{\bl'}^{\bj'}$ to both table $H$ and table $L$.  For its parent elements in $H$, we increase the number of children by one. 
	\item Remove the parent elements from table $L$ for all the newly added elements. 
	\item Repeat step 2 - step 3, until no element can be further added.
\end{enumerate}
\noindent\rule{16.5cm}{1pt}

\medskip

When the adaptive projection algorithm completes, it will generate a   hash table $H$, leaf table $L$ and   $f_h(\bx, \xi)=\sum_{v^\bj_{\bi,\bl} \in H}f^\bj_{\bi,\bl} v^\bj_{\bi,\bl}(\bx, \xi)$.  We denote the approximation space  $\bV^k_{N,H}(\Omega)=\textrm{span}\{v^\bj_{\bi,\bl} \in H\}$ and it is a subspace of $\bV^k_N(\Omega).$ 
On the other hand, we compute the initializations of $\bE_h(\bx), \bB_h(\bx)$  by a simple $L^2$ projection of $\bE(\bx), \bB(\bx)$ onto $\bV^k_N(\Ox)$.

Then we begin the time evolution algorithm which consists of several key steps.
The first step is the prediction step, which means that, 
given the hash table $H$ that stores the numerical solution $f_h$ at time step $t^n$ and the associated leaf table $L$, we need to predict the location where the details becomes significant at the next time step $t^{n+1}$, then add more elements in order  to capture the fine structures.  We  solve for $f_h\in \bV_{N,H}^k(\Omega)$ from $t^n$ to $t^{n+1}$
based on the solutions $\bE_h^n$ and $\bB_h^n$ at time $t^n$.
The forward Euler discretization is used as the time integrator in this step and we denote the predicted solution at $t^{n+1}$  by $f_h^{(p)}.$  

The second step is the refinement step according to $f_h^{(p)}$. We traverse the hash table $H$ and if an element $I^\bj_{\bl}=\{v^\bj_{\bi,\bl}, \mathbf{1}\leq\bi\leq\bk+\mathbf{1}\}$ satisfies the refinement criteria  \eqref{eq:l2},
indicating that such an element becomes significant at the next time step,
then we need to add its children elements
to $H$ and $L$  provided they are not added yet, and set the associated detail coefficients to zero. We also need to make sure that all the parent elements of the newly added element are in $H$ (i.e., no ``hole" is allowed in the hash table) and increase the number of children for all its parent elements by one. This step  generates the updated hash table $H^{(p)}$ and leaf table $L^{(p)}$.
Note that    \eqref{eq:l2} corresponds to the  $L^2$-norm based refinement criteria in  \cite{guo2017adaptive}. 

Based on the updated hash table $H^{(p)}$, the third step is to evolve the numerical solution $f_h$ by the DG weak formulation   with space $\bV_{N,H^{(p)}}^k(\Omega)$. Namely, we solve for  $\bV_{N,H^{(p)}}^k(\Omega)$ from $t^n$ to $t^{n+1}$, 
by the TVD-RK scheme \eqref{eq:tvd} to generate the pre-coarsened numerical solution $\tilde{f}_h^{n+1}$. Meanwhile, we also evolve the numerical solutions $\bE_h$ and $\bB_h$ from $t^n$ to $t^{n+1}$   with space $\bV_{N}^k(\Ox)$.

The last step   is to coarsen the mesh by removing elements that become insignificant at time level $t^{n+1}.$  The hash table $H^{(p)}$ that stores the numerical solution $\tilde{f}_h^{n+1}$  is recursively coarsened by the following procedure. 
The leaf table $L^{(p)}$ is traversed, and if an element $I_\bl^\bj\in L^{(p)}$ satisfies the coarsening criterion
\begin{align}
&\left(\sum_{\mathbf{1}\leq\bi\leq\bk+\mathbf{1}}|f^\bj_{\bi,\bl}|^2\right)^{\frac12}<\eta, \label{eq:l2_c}
\end{align}
where $\eta$ is a prescribed error constant, then we remove the element from both table $L^{(p)}$ and $H^{(p)}$, and set the associated coefficients $f^{\bj}_{\bi,\bl}=0,\,\mathbf{1}\leq\bi\leq\bk+\mathbf{1}$. For each of its parent elements in table $H^{(p)}$, we decrease the number of children by one. If the number becomes zero, i.e, the element has no child any more, then it is added to the leaf table $L^{(p)}$ accordingly. Repeat the coarsening procedure until no element can be removed from the table $L^{(p)}$. The output of this coarsening procedure are the updated hash table and leaf table, denoted by $H$ and $L$ respectively, and the compressed numerical solution $f_h^{n+1} \in \bV_{N,H}^k$. In practice, $\eta$ is chosen to be smaller than $\varepsilon$ for safety. In the simulations presented in this paper, we use $\eta = \varepsilon/10$.    
For the adaptive sparse grid DG scheme, the properties of  conservation of moments are not as clear as the sparse grid DG schemes, and will be subject to the error thresholds $\eta, \varepsilon$, see \cite{guo2017adaptive}.

\bigskip

\section{Numerical results}
\label{sec:numerical}

In this section, we consider two numerical tests to benchmark the performance of the  proposed   schemes. 

\subsection{1D2V streaming Weibel instability}
\label{subsec:SW}
Here, we consider the streaming Weibel (SW) instability, which has been previously considered both analytically and numerically  \cite{califano1998ksw, cheng2014discontinuous,cheng2014energy}. 
%
\begin{exa}
	\label{ex:SW}
In this case, a reduced version of the single species VM system with one spatial variable $x_2$, and two velocity variables $\xi_1, \xi_2$ is considered:
\begin{align}
f_t &+ \xi_2 f_{x_2} + (E_1 + \xi_2 B_3)f_{\xi_1} + (E_2 -
\xi_1 B_3 )f_{\xi_2} = 0~, \\
\frac{\df B_3}{\df t} &=  \frac{\df E_1}{\df x_2}, \quad
\frac{\df E_1}{\df t} =  \frac{\df
	B_3}{\df x_2} -  j_1, \quad \frac{\df E_2}{\df t} =  -  j_2~,
\end{align}
where
\beq j_1=\int_{-\infty}^{\infty} \int_{-\infty}^{\infty}
f(x_2, \xi_1, \xi_2, t) \xi_1 \,d\xi_1 d\xi_2,\quad j_2=\int_{-\infty}^{\infty} \int_{-\infty}^{\infty} f(x_2, \xi_1,
\xi_2, t) \xi_2 \,d\xi_1 d\xi_2~.
\eeq
Here, $f=f(x_2, \xi_1, \xi_2, t)$ is the distribution function of electrons, 
$\textbf{E}=(E_1(x_2, t),E_2(x_2, t), 0)$ is the 2D electric field, and $\textbf{B}=(0, 0, B_3(x_2, t))$ is the 1D magnetic field. Ions are assumed to form a constant background. 
\end{exa}

The initial condition is given by
\begin{align}
f(x_2, \xi_1, \xi_2, 0)&=\frac{1}{\pi \beta}
e^{- \xi_2^2 /\beta} [\delta e^{- (\xi_1-v_{0,1})^2 /\beta}
+(1-\delta) e^{- (\xi_1+v_{0,2})^2 /\beta} ],\\
E_1(x_2, \xi_1, \xi_2, 0)&=E_2(x_2, \xi_1, \xi_2, 0)=0, \qquad B_3(x_2, \xi_1, \xi_2, 0)=b \sin(k_0 x_2)~.
\end{align}
where $b$ denotes the amplitude of the initial perturbation to the magnetic field, $\beta^{1/2}$ is the thermal velocity, which is take to be $\beta=0.01$, and $\delta$ is a parameter measuring the symmetry of the electron beams.
Note that when $b=0$, this initial condition is an equilibrium state representing counter-streaming beams propagating perpendicular to the direction of inhomogeneity.    As in \cite{califano1998ksw}, the instability can be triggered by taking   $b=0.001$ as a perturbation of the initial magnetic field. The computational domain is chosen to be $\Ox=[0, L_{y}]$, where $L_y=2 \pi/ k_0$ and  $\Oxi=[-1.2, 1.2]^2$. We consider two sets of parameters as in \cite{califano1998ksw}
\begin{eqnarray}
{\rm \underline{choice\  1}:} \  \ \delta&=&0.5, v_{0,1}=v_{0,2}=0.3, k_0=0.2;
\nonumber\\
{\rm  \underline{choice\  2}:} \ \ \delta&=&1/6,  v_{0,1}=0.5, v_{0,2}=0.1,  k_0=0.2,
\nonumber
\end{eqnarray}
which lead to initially symmetric and  strongly non-symmetric  counter-streaming electron beams, respectively.
For all numerical simulations, unless otherwise noted, the time step $\Delta t$  is chosen  according to the mesh on the most refined level, i.e.
$$\displaystyle\Delta t = \frac{\text{CFL}}{\displaystyle\sum_{m=1}^d \frac{c_m}{h_N}},$$
where $c_m$ is the maximum wave propagation speed in $m$-th direction, and we take $\text{CFL}=0.1$.  

\textbf{Accuracy test and comparisons.}
 It is well-known that the VM system is time reversible, and such property provides  practical means to test accuracy for VM solvers. To elaborate, let $f(\bx, \xi, 0), \bE(\bx, 0), \bB(\bx, 0) $ be the initial condition and $f(\bx, \xi, T),  \bE(\bx, T), \bB(\bx, T)$ be the solution  at $t=T$ for the VM system. When we reverse the velocity field of the solution and the magnetic filed, yielding $f(\bx, -\xi, T), \bE(\bx, T), -\bB(\bx, T) $, and evolve the VM system again to $t = 2T$,
then we can recover $f(\bx, -\xi, 0), \bE(\bx, 0), -\bB(\bx, 0)$, which is the initial condition with the reverse velocity field and the reverse magnetic field. For the accuracy test, we  compute the solutions to $T=1$ and then reversely back to $T=2$, and compare the numerical solutions with the initial condition.   
The  $L^2$ errors of  $f,\, \bE,\, \bB $  are defined as 
\begin{align}
	&\text{for $f$:  $\sqrt{\frac{1}{|\Omega|} \int_\Omega (f(x_2,\xi_1,\xi_2)-f_{h}(x_2,\xi_1,\xi_2))^2\, dx_2 d\xi_1 d\xi_2}\, ,$}\notag\\
	&\text{for $\bE$:  $ \sqrt{\frac{1}{L_y} \int_\Ox (E_1(x_2)-E_{1,h}(x_2))^2+(E_2(x_2)-E_{2,h}(x_2))^2\, dx_2}\, ,$ }\label{eq:l2error}\\
	&\text{ for $\bB$: $\sqrt{\frac{1}{L_y} \int_\Ox (B_3(x_2)-B_{3,h}(x_2))^2\, dx_2},$}\notag
\end{align}
where $f_{h}, E_{1,h}, E_{2,h}, B_{3,h}$ are the numerical approximations to the exact solutions $f, E_1,$ $E_2, B_3$.

 We test accuracy for the sparse grid DG method with $k=1,\,2,\,3$ on different levels of meshes. 
For $k=3$, we take $\displaystyle\Delta t=O( h_N^{4/3})$ to match the temporal and spatial orders in  the convergence study.
The  $L^2$ errors and orders of the numerical solutions with  upwind and alternating fluxes for parameter choice 1   are reported in Tables \ref{errorupwind_set_1} and \ref{erroralt_set_1}.  Similar to \cite{guo2016sparse}, we observe at least $(k+\frac12)$-th order accuracy for both fluxes.  
When comparing the results of the two tables, we find that the errors are  identical for $f$, while there are some slight differences in the errors of the electromagnetic fields. This is expected because the solvers only differ in the discretization of Maxwell's equations.

\begin{table}[htp]
	\caption{$L^2$ errors and orders for the sparse grid DG method in Example \ref{ex:SW} with parameter choice 1. Run to $T=1.0$ and back to $T=2.0$. Upwind flux for Maxwell's equations.  The orders are measured with respect to $h_N,$ which is the size of the smallest mesh in each direction. 
	}
	\centering
	\begin{tabular}{|c|c| c c|c c|c c|}
		\hline
		
		\multirow{2}{*}{} & \multirow{2}{*}{$N$}  &	\multicolumn{2}{|c|}{$ f$} &	\multicolumn{2}{|c|}{$ \bB$} &	\multicolumn{2}{|c|}{$ \bE$}   \\
		\cline{3-8}
		&  & $L^2$ error & order & $L^2$ error &  order & $L^2$ error & order\\
		
		\hline 
		
		\multirow{5}{*}{$ k=1$}& 7	&	1.25E-01	&	& 7.49E-08	&	& 9.82E-08	& \\
		& 8	&	4.44E-02	&	1.49	& 2.86E-08	& 1.39	& 2.32E-08	& 2.08\\
		& 9	&	2.01E-02	&	1.14	& 5.01E-09	& 2.51	& 5.61E-09	& 2.05\\
		&10	&	5.74E-03	&	1.81	& 1.50E-09	& 1.74	& 1.32E-09	& 2.09\\
		
		\hline
		\multirow{5}{*}{$ k=2$}& 7	&	3.09E-02	&	& 1.37E-09	&	& 9.25E-10	& \\
		& 8	&	6.65E-03	&	2.22	& 1.40E-10	& 3.29	& 9.09E-11	& 3.35\\
		& 9	&	2.56E-03	&	1.39	& 2.57E-11 	& 2.45	& 9.13E-12	& 3.32\\
		&10	&	4.07E-04	&	2.65	& 4.67E-12	& 2.46	& 1.24E-12	& 2.88 \\
		
		\hline
		\multirow{5}{*}{$ k=3$}& 7	&	8.47E-03	&	& 1.45E-11	&	& 7.82E-12	& \\
		& 8	&	9.37E-04	&	3.18	& 1.18E-12	& 3.62	& 3.48E-13	& 4.49\\
		& 9	&	9.48E-05	&	3.31	& 7.09E-14	& 4.06	& 3.87E-14	& 3.17\\
		&10	&	8.48E-06	&	3.48	& 9.49E-15	& 2.90	& 1.40E-15	& 4.79 \\
		\hline
		
	\end{tabular}
	\label{errorupwind_set_1}
\end{table}

\begin{table}[htp]
	\caption{ $L^2$ errors and orders for the sparse grid DG method in example \ref{ex:SW} with parameter choice 1. Run to $T=1.0$ and back to $T=2.0$. Alternating flux $ \widehat{\bE_h}=\bE_h^+, \; \widehat{\bB_h}=\bB_h^-$ for Maxwell's equations. The orders are measured with respect to $h_N,$ which is the size of the smallest mesh in each direction.  
	}
	\centering
	\begin{tabular}{|c|c| c c|c c|c c|}
		\hline
		
		\multirow{2}{*}{} & \multirow{2}{*}{$N$} &	\multicolumn{2}{|c|}{$ f$} &	\multicolumn{2}{|c|}{$ \bB$} &	\multicolumn{2}{|c|}{$ \bE$}   \\
		\cline{3-8}
		& & $L^2$ error & order & $L^2$ error &  order & $L^2$ error & order\\
		
		\hline 
		
		\multirow{5}{*}{$ k=1$}& 7	&	1.25E-01	&	& 7.67E-08	&	& 6.32E-08	& \\
		& 8	&	4.44E-02	&	1.49	& 2.96E-08	& 1.37	& 1.11E-08	& 2.51\\
		& 9	&	2.01E-02	&	1.14	& 5.83E-09	& 2.34	& 2.36E-09	& 2.23\\
		&10	&	5.74E-03	&	1.81	& 1.83E-09	& 1.67	& 3.53E-10	& 2.74\\
		
		\hline
		\multirow{5}{*}{$ k=2$}& 7	&	3.09E-02	&	& 1.38E-09	& 	& 8.84E-10	& \\
		& 8	&	6.65E-03	&	2.22	& 1.42E-10	& 3.28	& 8.24E-11	& 3.42\\
		& 9	&	2.56E-03	&	1.39	& 2.60E-11 	& 2.45	& 7.76E-12	& 3.41\\
		&10	&	4.07E-04	&	2.65	& 4.70E-12	& 2.47	& 1.07E-12	& 2.86 \\
		
		\hline
		\multirow{5}{*}{$ k=3$} & 7	&	8.47E-03	&	& 1.45E-11	&	& 7.78E-12	& \\
		& 8	&	9.37E-04	&	3.18	& 1.19E-12	& 3.61	& 3.45E-13	& 4.50\\
		& 9	&	9.48E-05	&	3.31	& 7.09E-14	& 4.07	& 3.87E-14	& 3.16\\
		&10	&	8.48E-06	&	3.48	& 9.49E-15	& 2.90	& 1.39E-15	& 4.80 \\
		\hline
		
	\end{tabular}
	\label{erroralt_set_1}
\end{table}

  We then perform   accuracy test for the adaptive sparse grid DG scheme. As in \cite{guo2017adaptive}, we run the simulations with a fixed maximum mesh level $N=7$ and different  $\varepsilon$ values, and report the $L^2$ errors and the number of active degrees of freedom at final time.  The following rates of convergence are calculated,  
\begin{align*}
\mbox{convergence rate with respect to the error threshold}  \quad &R_{\varepsilon}=\frac{\log(e_{l-1}/e_l)}{\log(\varepsilon_{l-1}/{\varepsilon_l})},\\
\mbox{convergence rate with respect to degrees of freedom } \quad & R_{\text{DOF}}=\frac{\log(e_{l-1}/e_l)}{\log(\text{DOF}_l/\text{DOF}_{l-1})}, 
\end{align*} 
where $e_l$ is the $L^2$ error defined in \eqref{eq:l2error} with refinement parameter $\varepsilon_l$, and $\text{DOF}_l$ is the associated number of active degrees of freedom at final time. The two widely used convergence rates above provide important measurement of accuracy and effectiveness  of adaptive schemes  \cite{bokanowski2013adaptive,guo2017adaptive}.  $R_{\epsilon}$ demonstrates how much the numerical error is reduced when one picks a smaller error threshold, while $R_{\text{DOF}}$ reveals the relation between the numerical error and the active degrees of freedom.
For comparison purposes, recall the DG schemes constructed by the tensor product full grid yields $R_{\epsilon}\approx 1$ and $R_{\text{DOF}}\approx \frac{k+1}{d}$.  Our numerical results are summarized in Tables \ref{table:ada_upwind_set_1} and \ref{table:ada_upwind_set_2} for the parameter choices 1 and 2 with the upwind flux for Maxwell's equations. The results computed by the alternating fluxes are very similar, and hence are omitted. We observe that for both test cases, rate $R_{\epsilon}$ is slightly smaller than 1, and  $R_{\text{DOF}}$ is much larger than $\frac{k+1}{d}$ but still smaller than $k+1$ for $f.$ This demonstrates the effectiveness as well as the computational savings of the adaptive scheme compared with the non-adaptive ones.  We also observe the error saturation for $\bB$ and $\bE$, i.e. when we reduce $\varepsilon$, unlike for $f$, the errors for $\bB$ and $\bE$ do not decrease. The error saturation can be ascribed to the fact that Maxwell's equations are solved on the fixed finest level mesh, and hence $\bB$ and $\bE$ are fully resolved by the scheme in this case. This is also evident from the magnitude of the errors for $\bB$ and $\bE$, which is much smaller than the $\epsilon$ parameters.

\begin{table}[htp]
	\caption{$L^2$ errors and orders of accuracy for the adaptive sparse grid DG method in example \ref{ex:SW} with parameter choice 1. Upwind flux for Maxwell's equations. Run to $T=1.0$ and back to $T=2.0$. $N=7$.  
	}
	\centering
	\begin{tabular}{|c|c| c c c c|c c|c c|}
		\hline
		
		\multirow{2}{*}{} & \multirow{2}{*}{$\varepsilon$} &	\multicolumn{4}{|c|}{$ f$} &	\multicolumn{2}{|c|}{$ \bB$} &	\multicolumn{2}{|c|}{$ \bE$}   \\
		\cline{3-10}
		&  &  DOF& $L^2$ error &  $R_{\text{DOF}}$& $R_\varepsilon$ & $L^2$ error &  $R_\varepsilon$& $L^2$ error & $R_\varepsilon$\\
		
		\hline 
		
		\multirow{4}{*}{$ k=1$} &5E-01	&	1216	&	6.06E-02	&	&	& 4.17E-07	& & 2.84E-06	&	\\
		& 1E-01	&	2272	&	3.35E-02	&	0.95	&	0.37	& 2.17E-07	& 0.41	& 1.44E-06	& 0.42\\
		& 5E-02	&	3744	&	9.27E-03	&	2.57	&	1.85	& 2.17E-07	& 0.00	& 1.44E-06	& 0.00\\
		& 1E-02	&	6656	&	4.10E-03	&	1.42	&	0.51	& 2.17E-07	& 0.00	& 1.44E-06	& 0.00\\
		
		\hline
		\multirow{4}{*}{$ k=2$} & 5E-02	&	648	    &	9.28E-01	&	&	& 1.45E-07	&	& 6.97E-07	&\\
		& 1E-02	&	10638	&	8.32E-04	&	2.51	&	4.36	& 9.71E-08	& 0.25	& 6.18E-07	& 0.07\\
		& 5E-03	&	12798	&	5.99E-04	&	1.78	&	0.47	& 9.71E-08	& 0.00	& 6.18E-07	& 0.00\\
		& 1E-03	&	20574	&	2.83E-04	&	1.58	&	0.47	& 7.24E-08	& 0.18	& 1.51E-07	& 0.88\\
		
		\hline
		\multirow{4}{*}{$ k=3$} & 1E-04	&	43904	&	7.44E-05	&	&	& 3.95E-08	&	& 1.85E-07	&\\
		& 5E-05	&	51584	&	4.95E-05	&	2.53	&	0.59	& 3.01E-08	& 0.39	& 5.20E-08	& 1.83\\
		& 1E-05	&	72704	&	1.08E-05	&	4.44	&	0.95	& 2.51E-08	& 0.11	& 2.23E-08	& 0.53\\
		& 5E-06	&	82176	&	7.22E-06	&	3.29	&	0.58	& 2.51E-08	& 0.00	& 2.23E-08	& 0.00\\
		\hline
		
	\end{tabular}
	\label{table:ada_upwind_set_1}
\end{table}

\begin{table}[htp]
	\caption{$L^2$ errors and orders of accuracy for the adaptive sparse grid DG method in example \ref{ex:SW} with parameter choice 2. Upwind flux for Maxwell's equations. Run to $T=1.0$ and back to $T=2.0$. $N=7$. 
	}
	\centering
	\begin{tabular}{|c|c| c c c c|c c|c c|}
		\hline
		
		\multirow{2}{*}{} & \multirow{2}{*}{$\varepsilon$} &	\multicolumn{4}{|c|}{$ f$} &	\multicolumn{2}{|c|}{$ \bB$} &	\multicolumn{2}{|c|}{$ \bE$}   \\
		\cline{3-10}
		&  &  DOF& $L^2$ error &  $R_{\text{DOF}}$& $R_\varepsilon$ & $L^2$ error &  $R_\varepsilon$& $L^2$ error & $R_\varepsilon$\\
		
		\hline 
		
		\multirow{4}{*}{$ k=1$} &5E-01	&	1088	&	7.71E-02	&	&	& 2.12E-07	& & 1.41E-06	&	\\
		& 1E-01	&	2672	&	1.87E-02	&	1.58	&	0.88	& 2.12E-07	& 0.00	& 1.41E-06	& 0.00\\
		& 5E-02	&	3648	&	1.08E-02	&	1.76	&	0.79	& 2.12E-07	& 0.00	& 1.41E-06	& 0.00\\
		& 1E-02	&	6224	&	4.90E-03	&	1.48	&	0.49	& 2.12E-07	& 0.00	& 1.41E-06	& 0.00\\
		
		\hline
		\multirow{4}{*}{$ k=2$} & 5E-02	&	5562 &	4.52E-03	&	&	& 9.66E-08	&	& 4.23E-07	&\\
		& 1E-02	&	9666	&	1.33E-03	&	2.21	&	0.76	& 9.66E-08	& 0.00	& 4.23E-07	& 0.00\\
		& 5E-03	&	12258	&	6.24E-04	&	3.19	&	1.09	& 9.66E-08	& 0.00	& 4.23E-07	& 0.00\\
		& 1E-03	&	19278	&	2.91E-04	&	1.68	&	0.47	& 7.13E-08	& 0.19	& 1.49E-08	& 0.65\\
		
		\hline
		\multirow{4}{*}{$ k=3$} & 1E-04	&	41728	&	5.81E-05	&	&	& 2.33E-08	&	& 1.06E-07	&\\
		& 5E-05	&	48384	&	2.83E-05	&	4.86	&	1.04	& 2.22E-08	& 0.07	& 1.02E-07	& 0.06\\
		& 1E-05	&	68736	&	1.02E-05	&	2.91	&	0.63	& 1.83E-08	& 0.12	& 1.60E-08	& 1.15\\
		& 5E-06	&	77568	&	7.05E-06	&	3.06	&	0.53	& 1.83E-08	& 0.00	& 1.59E-08	& 0.01\\
		\hline
		
	\end{tabular}
	\label{table:ada_upwind_set_2}
\end{table}

\begin{figure}[htp]
	\begin{center}
		\subfigure[k=1]{\includegraphics[width=.32\textwidth]{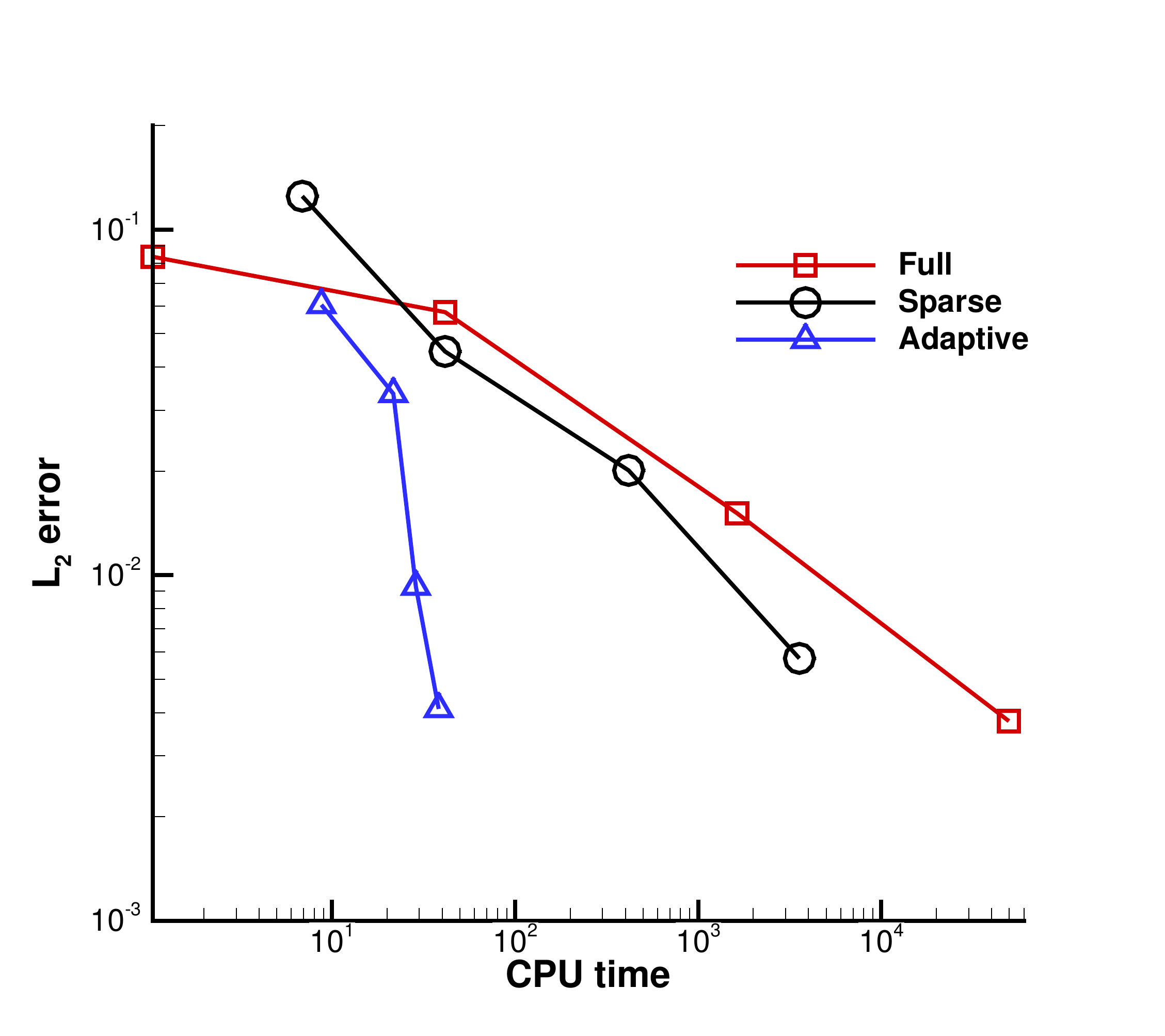}}
		\subfigure[k=2]{\includegraphics[width=.32\textwidth]{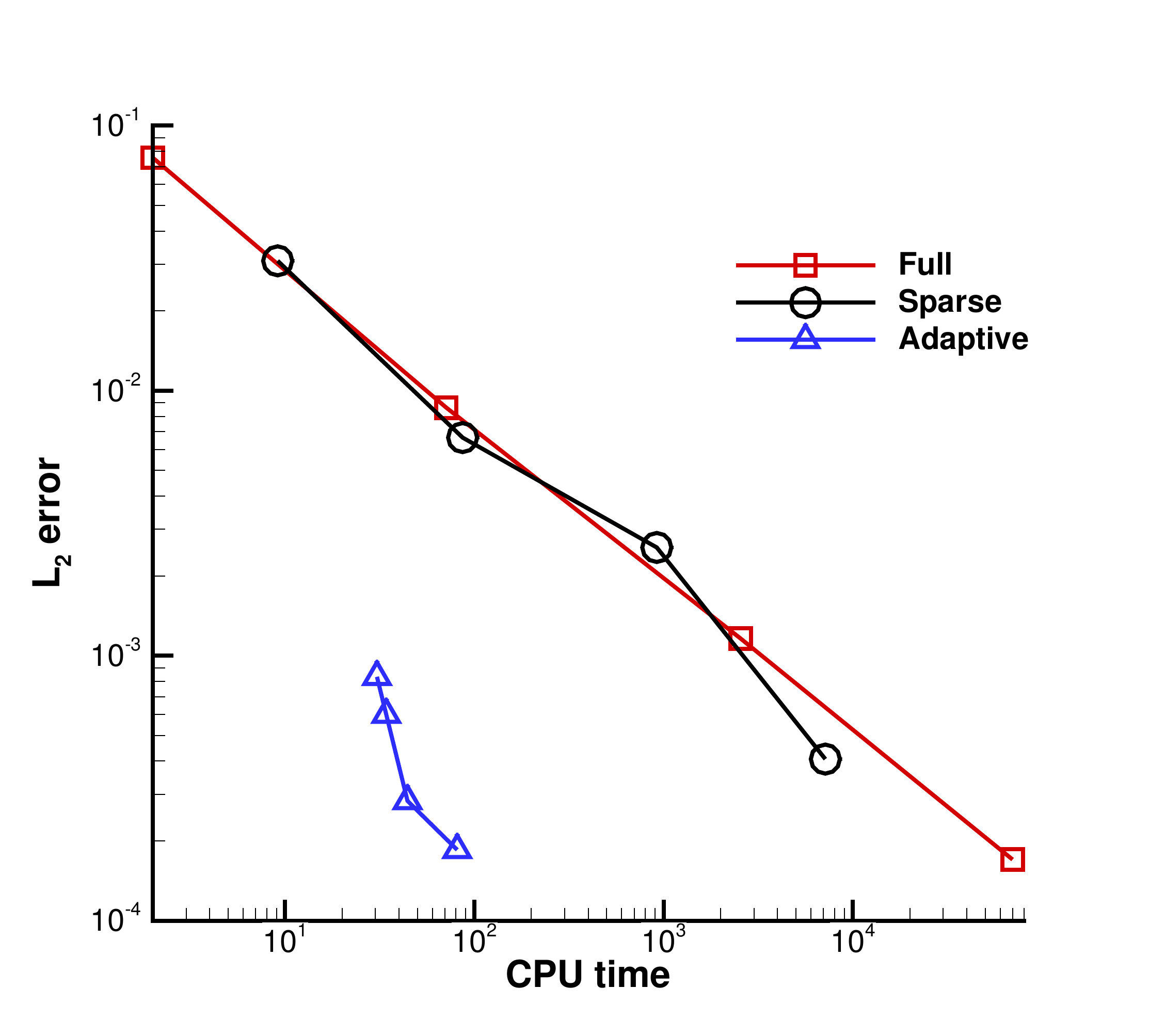}}
		\subfigure[k=3]{\includegraphics[width=.32\textwidth]{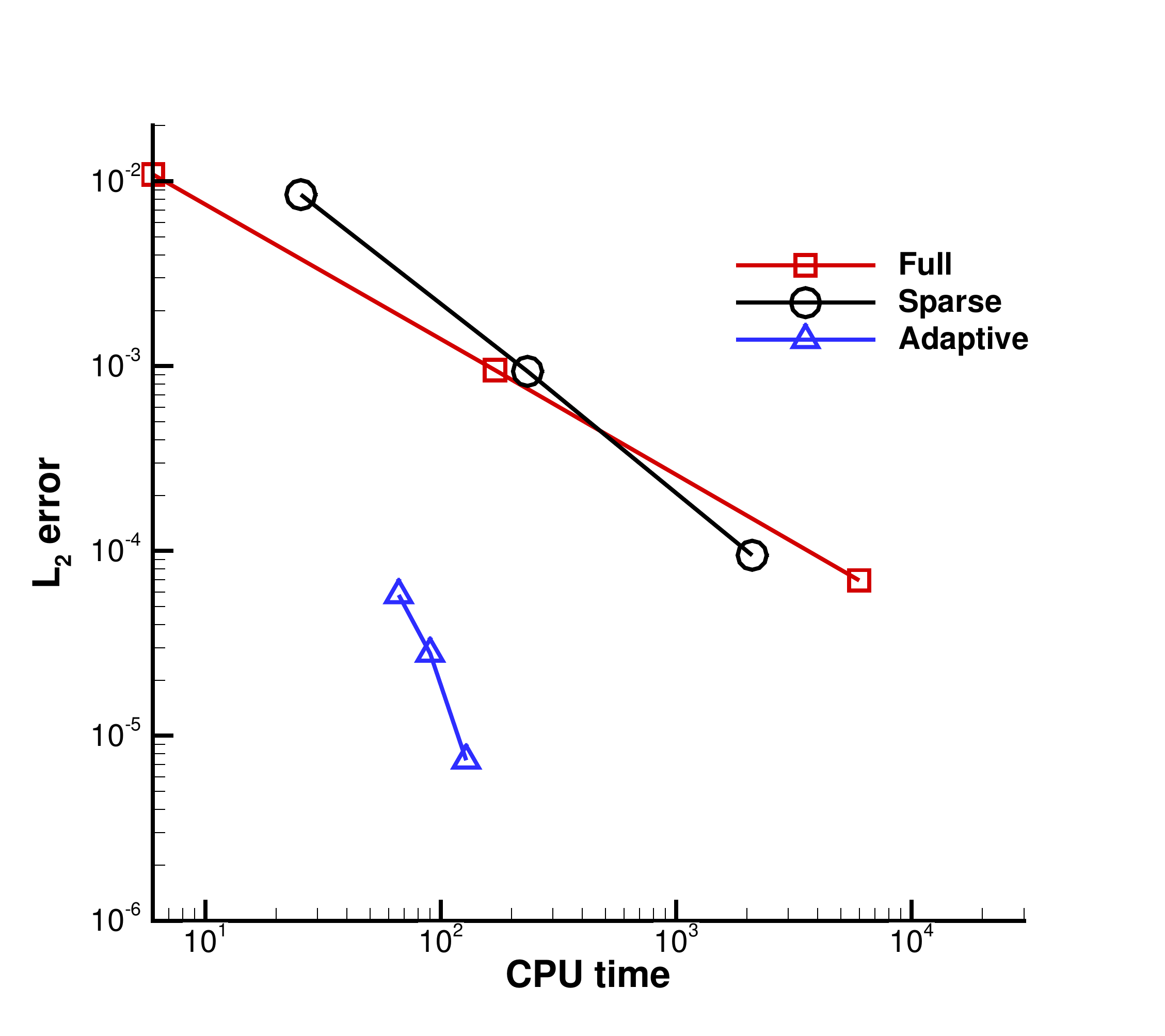}}
	\end{center}
	\caption{ $L^2$ errors of $f$ vs.  CPU time for full grid DG, sparse grid DG and, adaptive sparse grid DG methods with parameter choice 1 in Example 3.1. Run to $T=1.0$ and back to $T=2.0$. Upwind flux for Maxwell's equations.  (a) k=1, (b) k=2, (c) k=3.}
	\label{fig:cpu_2d}
\end{figure}

\begin{figure}[htp]
	\begin{center}
		\subfigure[k=1]{\includegraphics[width=.32\textwidth]{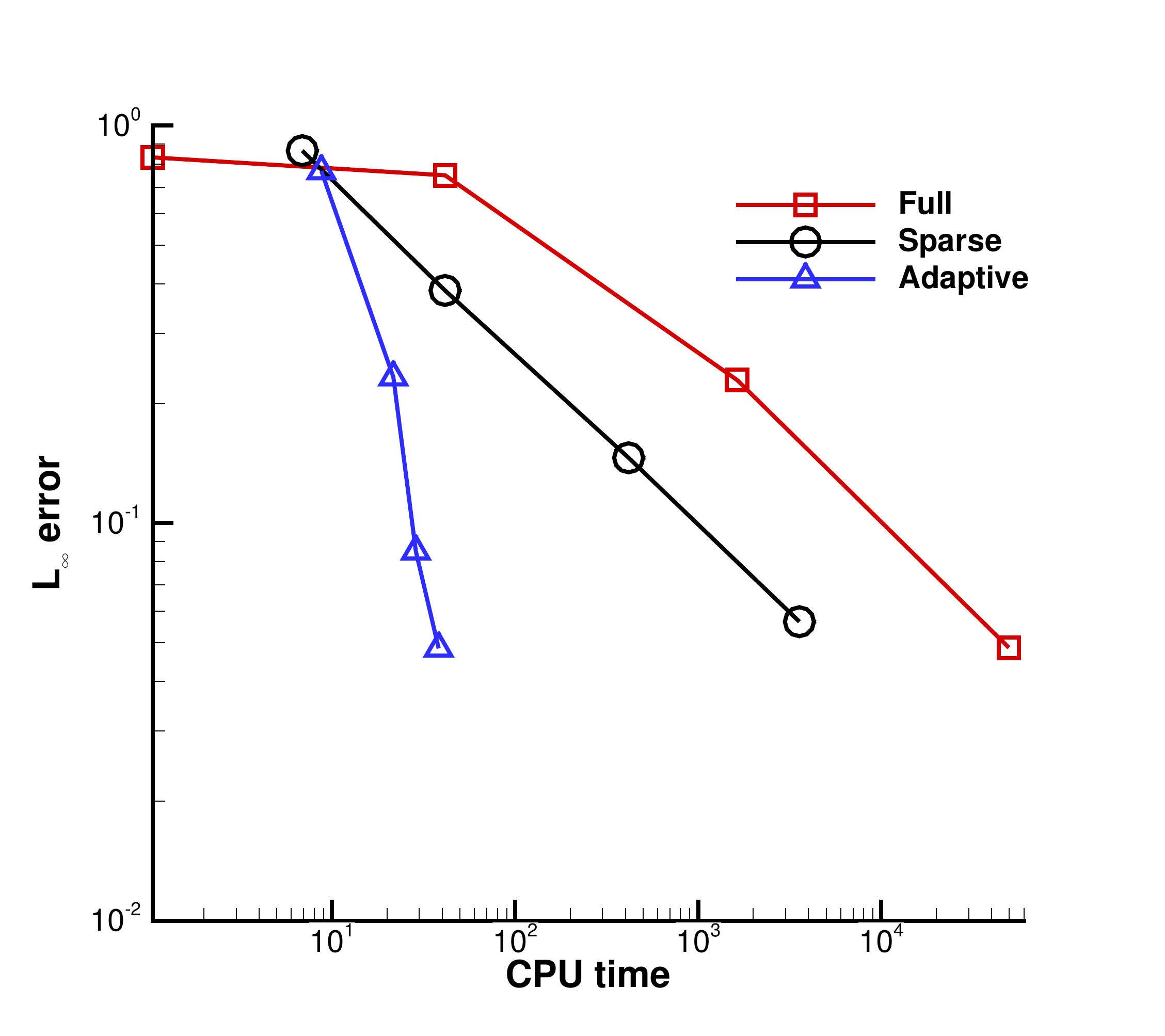}}
		\subfigure[k=2]{\includegraphics[width=.32\textwidth]{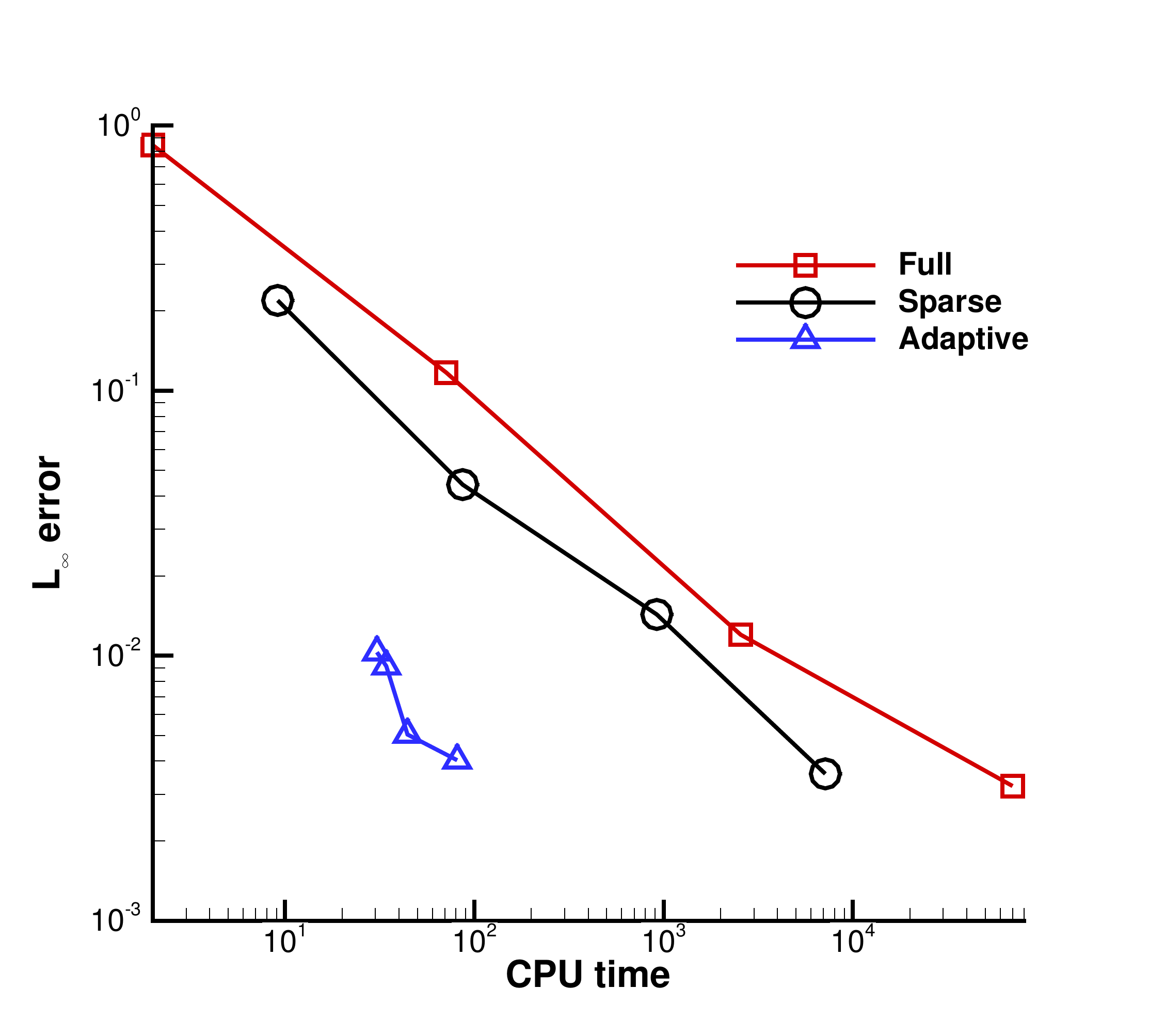}}
		\subfigure[k=3]{\includegraphics[width=.32\textwidth]{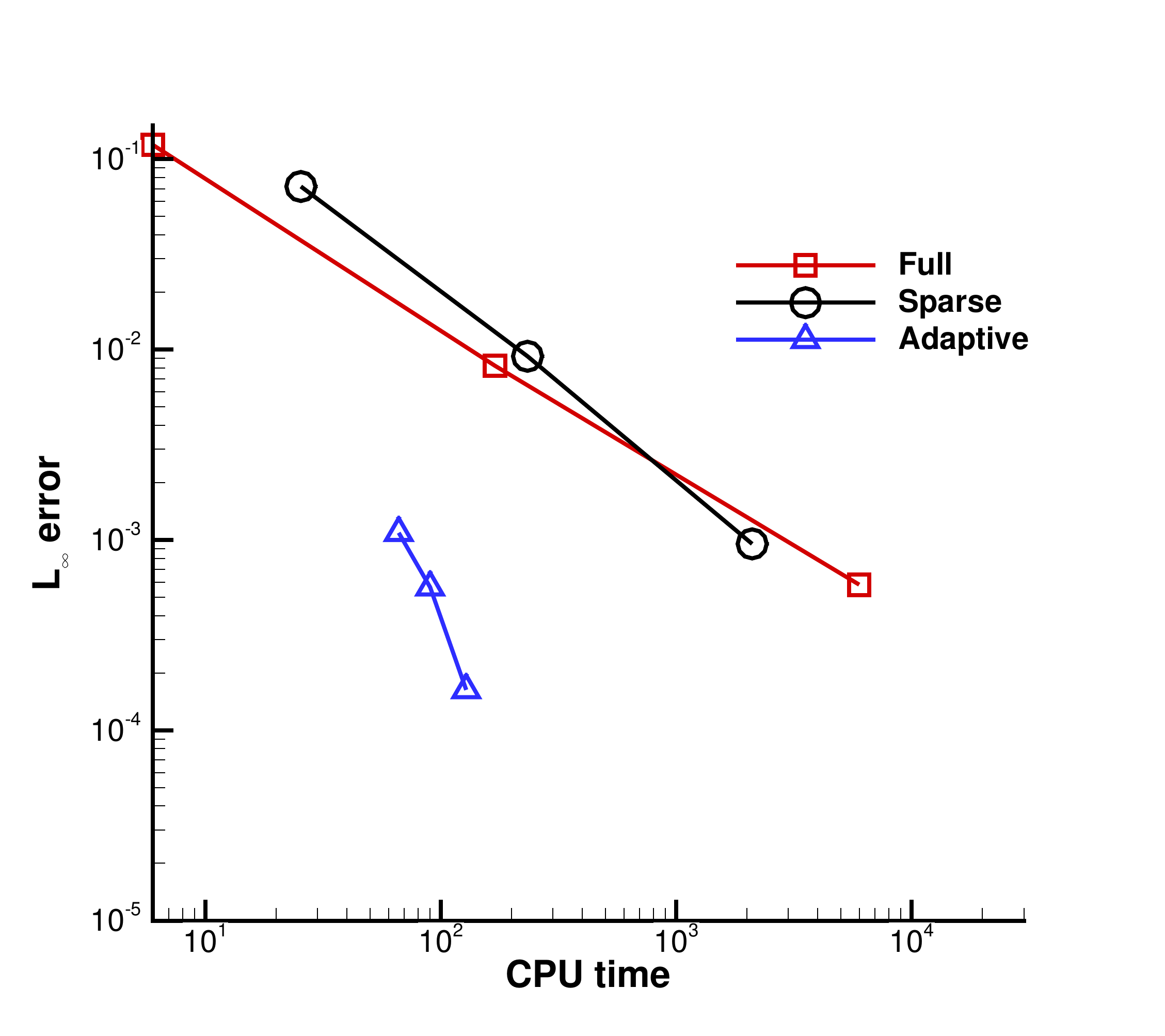}}
	\end{center}
	\caption{ $L^{\infty}$ errors of $f$ vs.  CPU time for full grid DG, sparse grid DG, and adaptive sparse grid DG methods with parameter choice 1 in Example 3.1. Run to $T=1.0$ and back to $T=2.0$. Upwind flux for Maxwell's equations.  (a) k=1, (b) k=2, (c) k=3.}
	\label{fig:cpu_3d}
\end{figure}

\begin{figure}[htp]
	\begin{center}
		\subfigure[$L^2$ error]{\includegraphics[width=.36\textwidth]{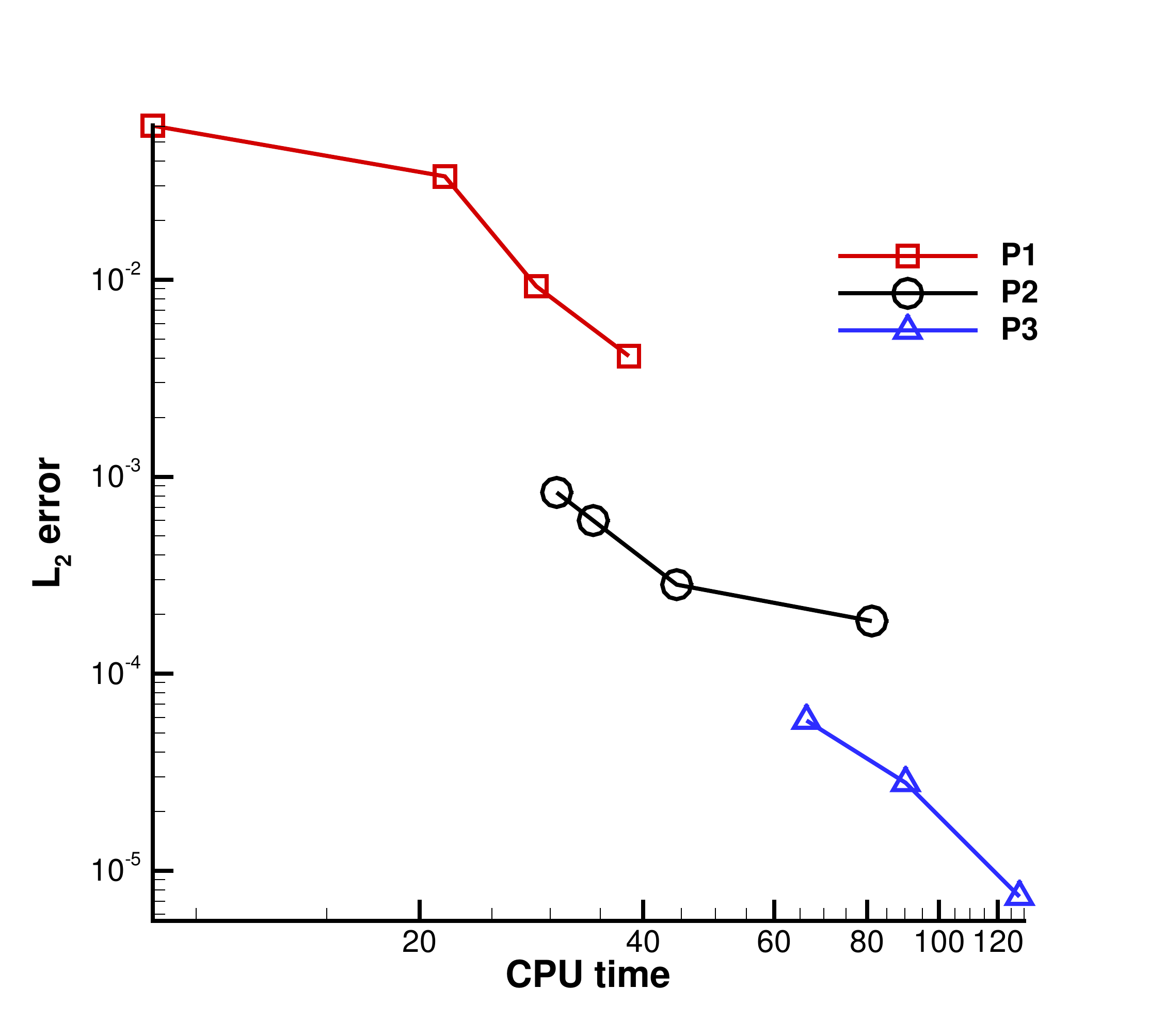}}
		\subfigure[$L^{\infty}$ error]{\includegraphics[width=.36\textwidth]{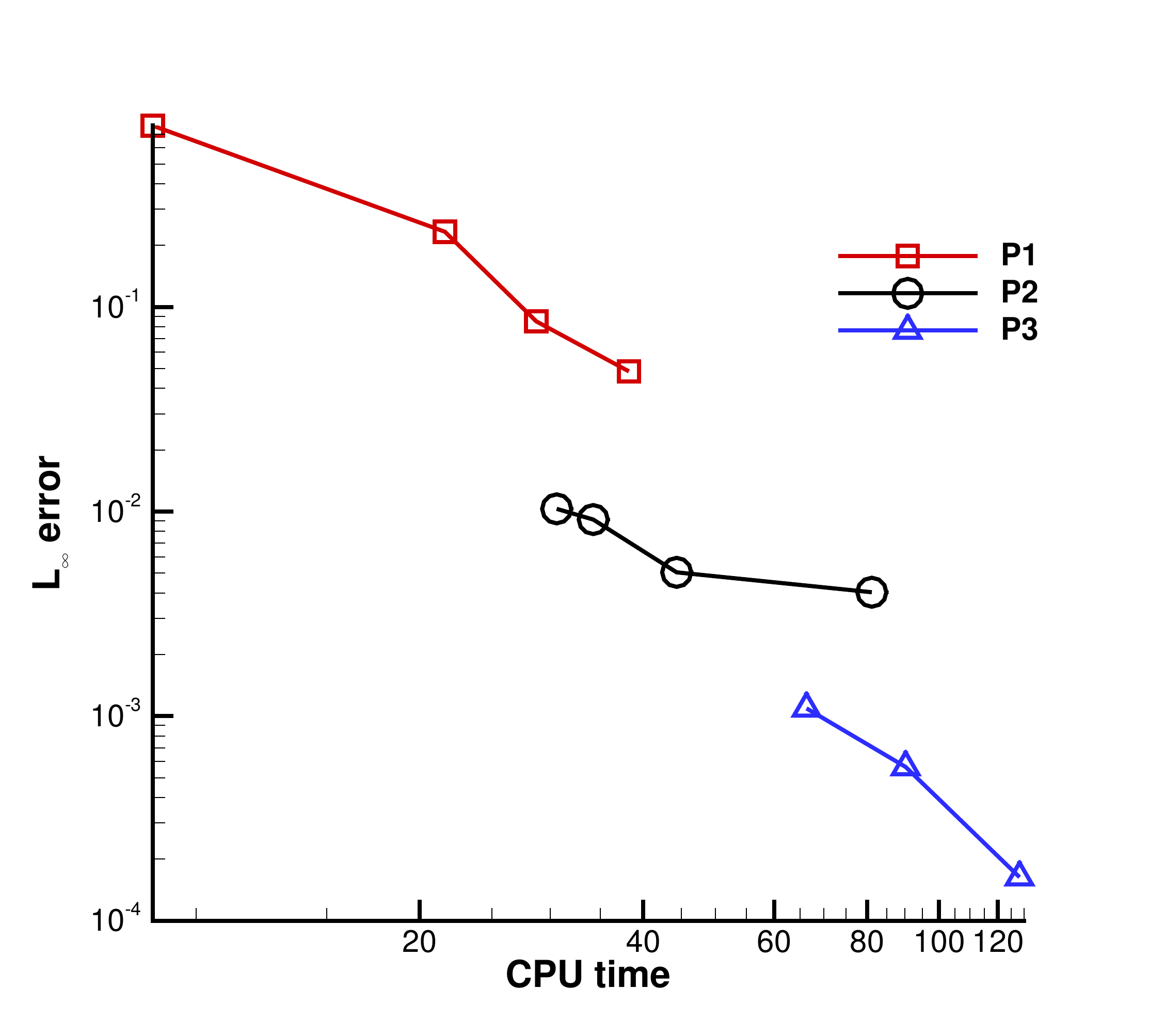}}
	\end{center}
	\caption{ $L^2$  and $L^{\infty}$ errors of $f$ vs.  CPU time for adaptive sparse grid DG methods with parameter choice 1 in Example 3.1. Run to $T=1.0$ and back to $T=2.0$. Upwind flux for Maxwell's equations.  (a) $L^2$ error, (b) $L^{\infty}$ error.}
	\label{fig:cpu_3d_ada}
\end{figure}

Now we use this example to compare the performance of the full grid DG, sparse grid DG, and adaptive sparse grid DG schemes. The efficiency will be evaluated based on the comparison of numerical error vs. CPU time. 
The computations in this example are implemented by  OpenMP codes  using computational resources from the Institute for Cyber-Enabled Research in  Michigan State University. \textcolor{blue}{}
Because the error in $f$ is more dominant than those of the electromagnetic field, we only plot the error  of $f$ in the $L^2$ and $L^\infty$ norms. Only upwind flux for Maxwell's equation is considered. In Figures \ref{fig:cpu_2d}-\ref{fig:cpu_3d_ada}, we run the SW instability to $T=1.0$ and back to $T=2.0,$ and plot the errors vs. CPU time for $k=1,2,3$ with the three numerical schemes. The time discretizations are TVD-RK3 method  for $k=1, 2,$ and RK4 method for $k=3.$ We use a fixed $CFL=0.1$ in all computations.
In the plots \ref{fig:cpu_2d}-\ref{fig:cpu_3d}, the mesh levels $N$ are taken as from 4 to 7, 4 to 7 and 4 to 6 for $k=1,2,3$, respectively, in full grid DG method and taken as from 7 to 10, 7 to 10 and 7 to 9 for $k=1,2,3$, respectively, in sparse grid DG method.
For the adaptive sparse grid method, we fix $N=7$ and take parameter $\varepsilon$ ranging from $0.5$ to $10^{-2}$, $10^{-2}$ to $10^{-4}$, and $10^{-4}$ to $10^{-5}$
for $k=1,2,3$, respectively.  From Figures \ref{fig:cpu_2d}-\ref{fig:cpu_3d}, it is evident that the adaptive sparse grid DG method is the most efficient among the three methods in all cases. The sparse grid DG method outperforms the full grid DG method when the resolution is finer, because we can see that the slope of the sparse grid DG method is steeper than that of the full grid DG method.
We note that the fast matrix-vector product as in \cite{shen2010efficient,shen2012efficient} can be used to further accelerate the sparse grid method, which will be investigated in a future work. In Figure \ref{fig:cpu_3d_ada}, we focus on the comparison of the adaptive sparse grid DG method with different polynomial orders. It is clear that the higher order adaptive scheme is more efficient. However, the slopes of three lines for $k=1,2,3$ seem comparable, which implies the convergence rate with respect to CPU time for different polynomial orders are similar for the adaptive scheme in this case.

In summary, the proposed sparse grid and adaptive sparse grid schemes are able to   simulate the VM system with smooth solutions with desired orders of accuracy. Both schemes require much less degrees of freedom compared with their  full grid counterpart in \cite{cheng2014discontinuous}.

\textbf{Numerical results for long time simulations.}
Here, we are concerned with the performance of the methods in long time, which are often required for plasma simulations. Most simulation results presented below are compared with the more expensive full grid DG method in \cite{cheng2014discontinuous}. The scaled macroscopic quantities under considerations are defined as follows.
$$K_1=\frac{1}{2 L_y} \int_\Omega f \xi_1^2\,d\xi_1 d\xi_2 dx_2,\quad K_2=\frac{1}{2 L_y} \int_\Omega f \xi_2^2\,d\xi_1 d\xi_2 dx_2,$$
$$ \mathcal{E}_1=\frac{1}{2 L_y} \int_\Ox E_1^2\, dx_2,\quad \mathcal{E}_2=\frac{1}{2 L_y} \int_\Ox E_2^2\, dx_2, \quad \mathcal{B}_3=\frac{1}{2 L_y} \int_\Ox B_3^2\, dx_2,$$
where $K_1, K_2$ are the scaled kinetic energies in each direction, $\mathcal{E}_1, \mathcal{E}_2$ are the scaled electric energies in each direction, and $\mathcal{B}_3$ is the scaled magnetic energy. Henceforth, the scaled kinetic and electric energies are the summation of the corresponding components in each direction, and the scaled total energy is defined as the summation of $K_1, K_2, \mathcal{E}_1, \mathcal{E}_2, \mathcal{B}_3.$
The scaled momentum  $P_1$ and $P_2$ are 
$$P_1=\frac{1}{L_y} \left(\int_\Omega  \xi_1 f \,d\xi_1 d\xi_2 dx_2 + \int_\Ox E_2B_3 \, dx_2 \right),\quad P_2=\frac{1}{L_y} \left( \int_\Omega  \xi_2 f \,d\xi_1 d\xi_2 dx_2 - \int_\Ox E_1B_3 \, dx_2 \right),$$
where the first term in each expression represents the momentum in particles while the second represents that of the electromagnetic field.

\begin{figure}[htb]
	\begin{center}
		\subfigure[ Sparse grid, relative error in mass]{\includegraphics[width=3in,angle=0]{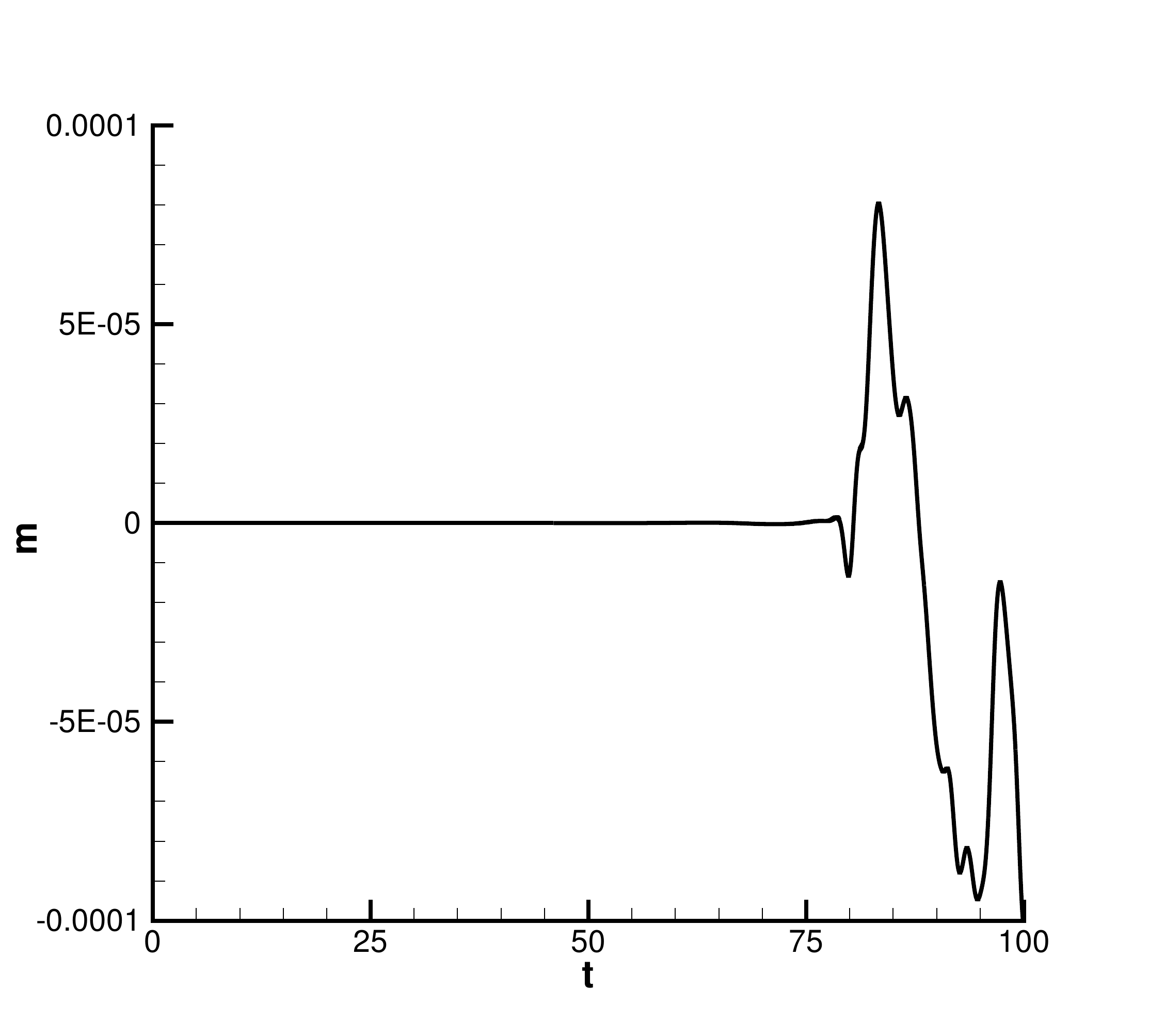}}
		\subfigure[Adaptive, relative error in mass]{\includegraphics[width=3in,angle=0]{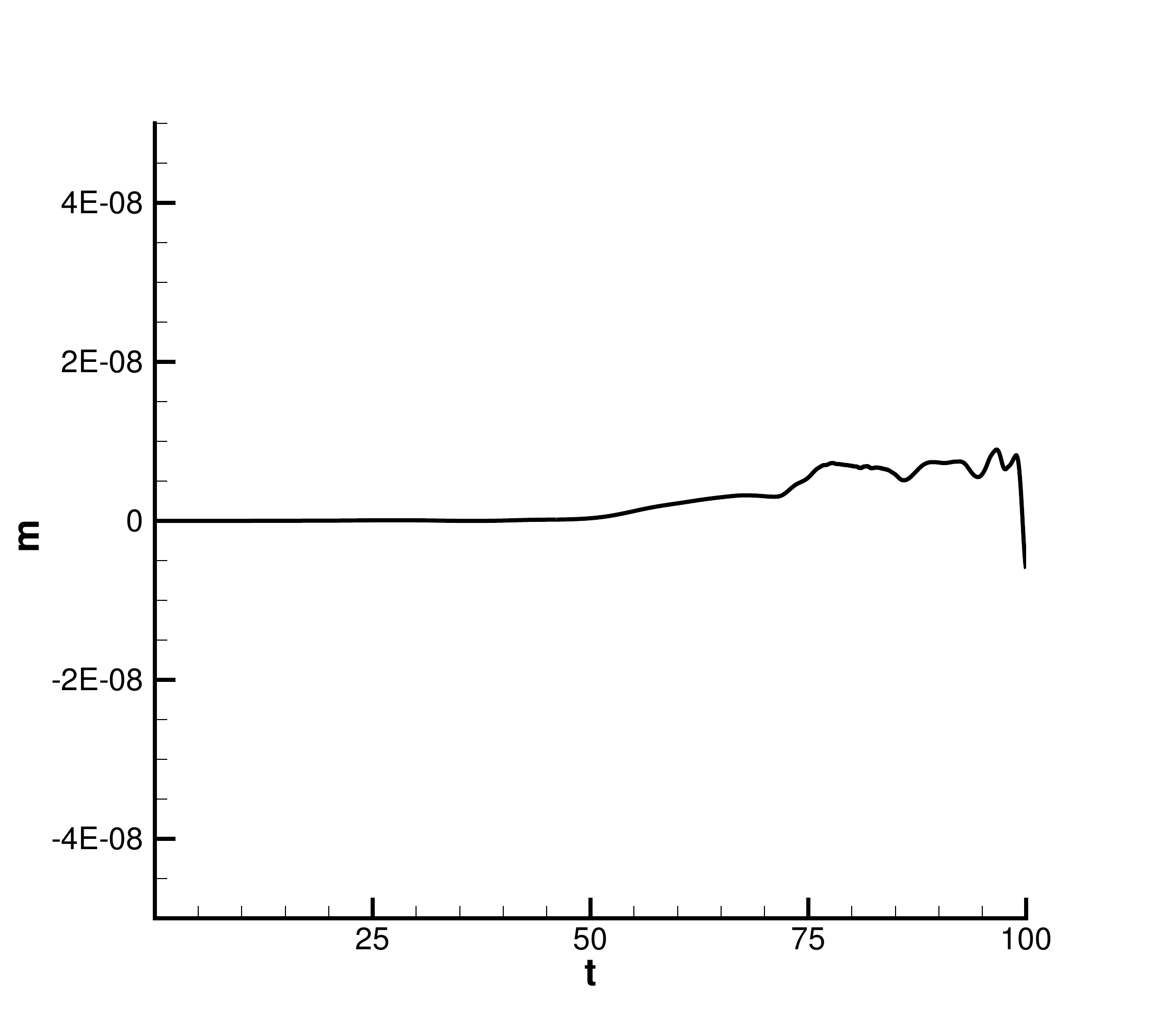}}
		\subfigure[Sparse grid, relative error in total energy]{\includegraphics[width=3in,angle=0]{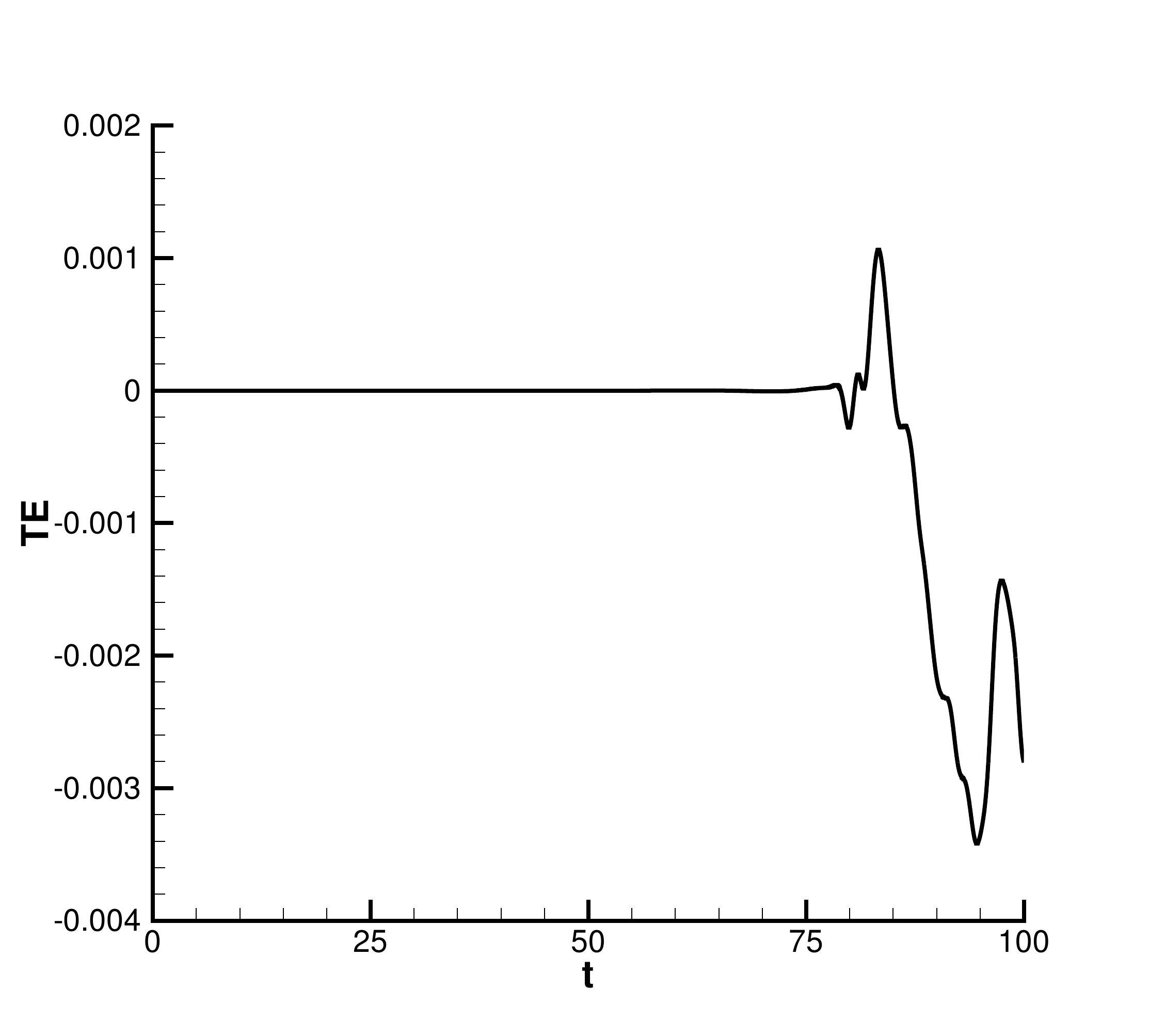}}
		\subfigure[Adaptive,  relative error in total energy]{\includegraphics[width=3in,angle=0]{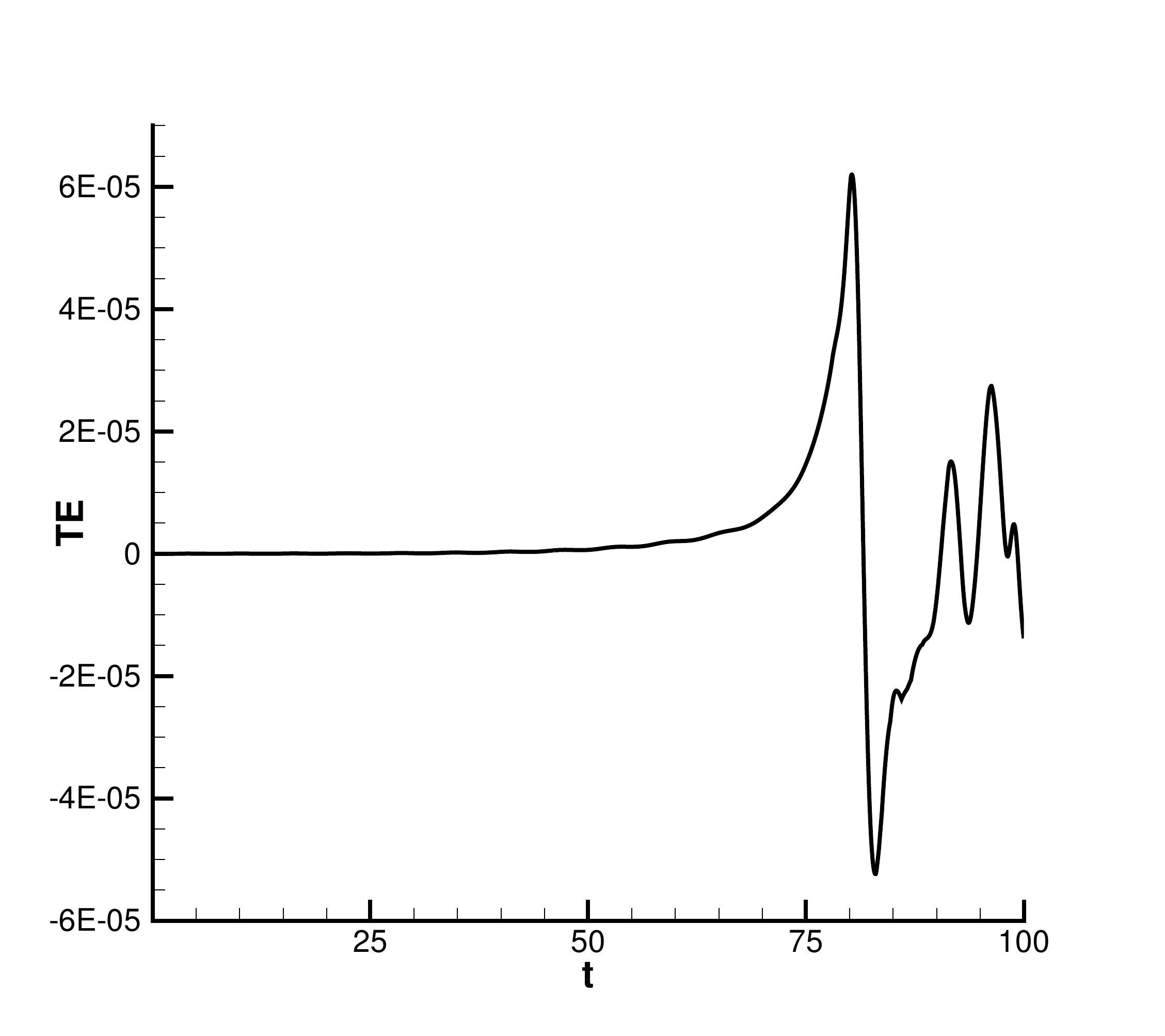}}
	\end{center}
	\caption{SW instability with parameter choice 1. Time evolution of relative error of mass and total energy by upwind flux for the Maxwell's equations. Sparse grid: $N=8$, $k=3$.  Adaptive sparse grid: $N=6$, $k=3$, $\eps=2\times10^{-7}$.}
	\label{masste1}
\end{figure}

\begin{figure}[htb]
	\begin{center}
		\subfigure[ Sparse grid, error in momentum $P1$]{\includegraphics[width=3in,angle=0]{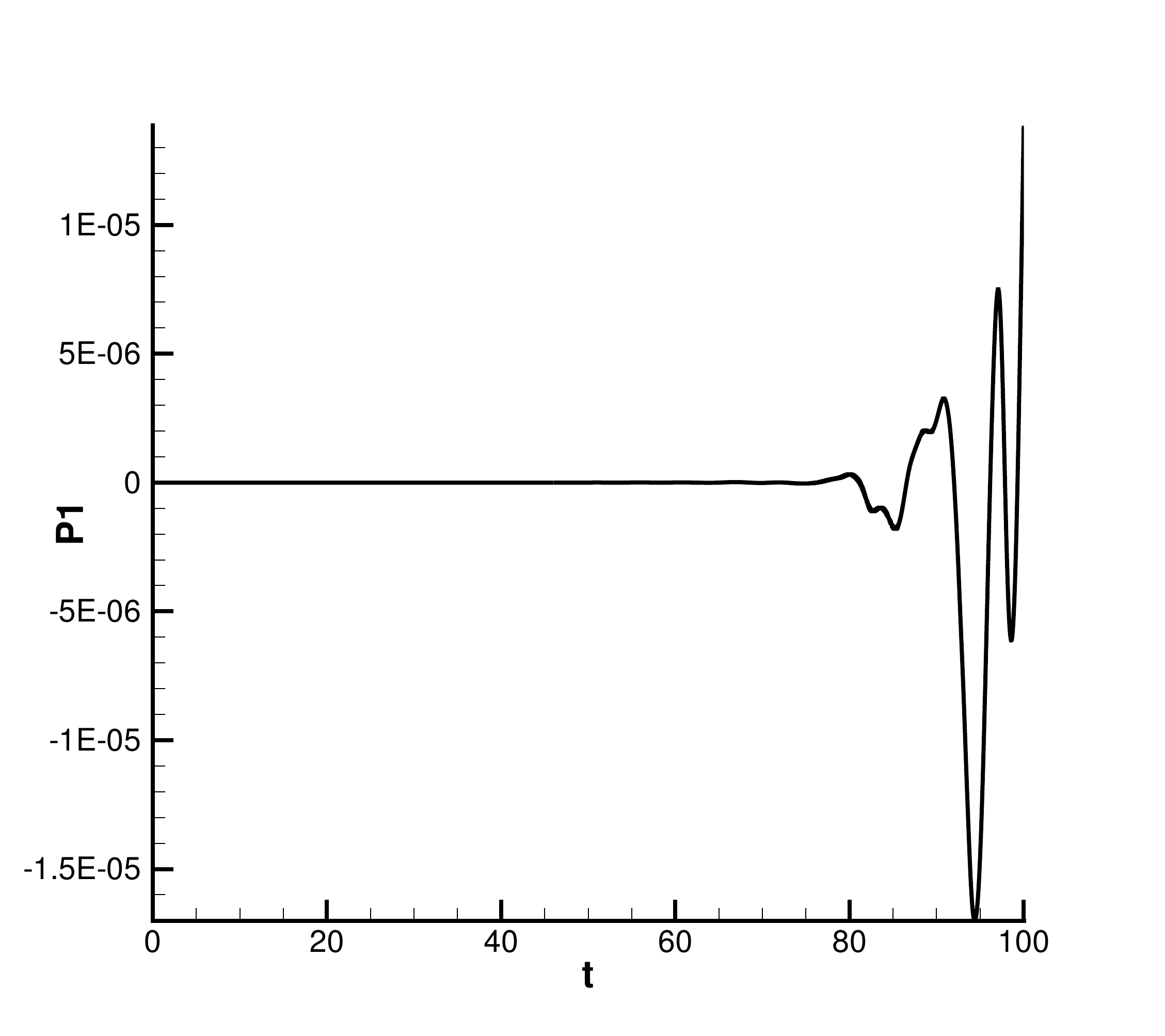}}
		\subfigure[Adaptive, error in momentum $P1$]{\includegraphics[width=3in,angle=0]{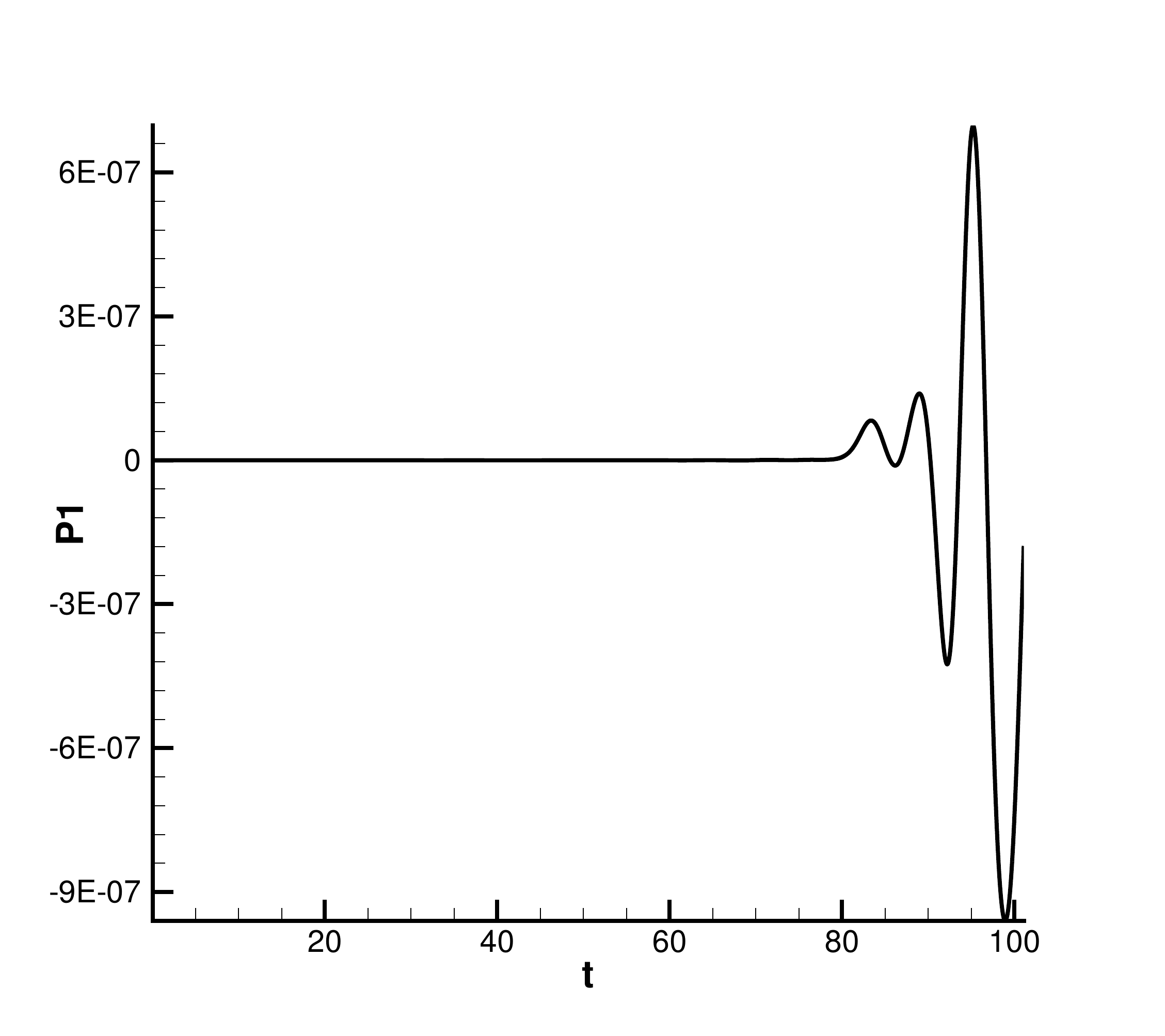}}
		\subfigure[Sparse grid, error in momentum $P2$]{\includegraphics[width=3in,angle=0]{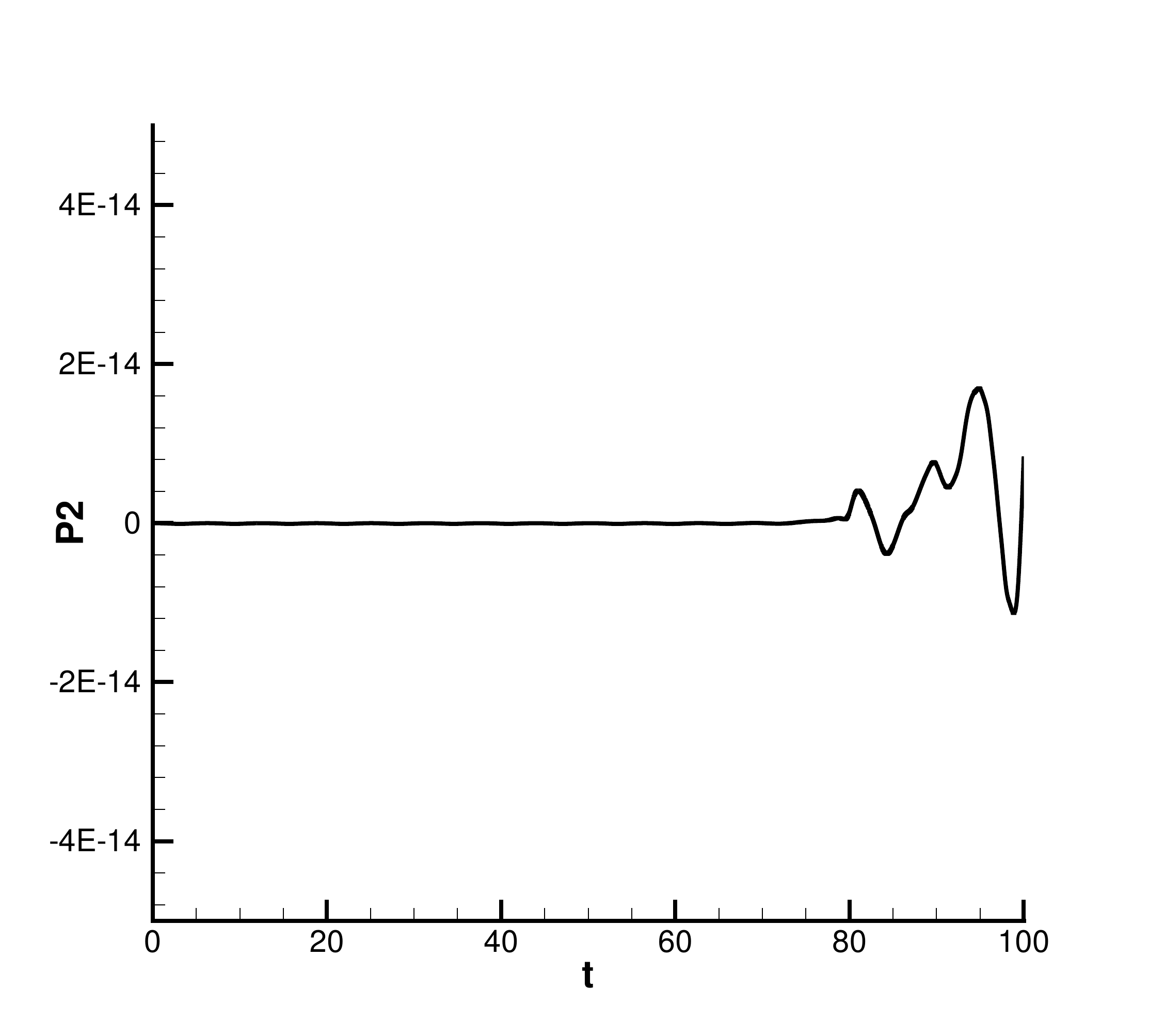}}
		\subfigure[Adaptive, error in momentum $P2$]{\includegraphics[width=3in,angle=0]{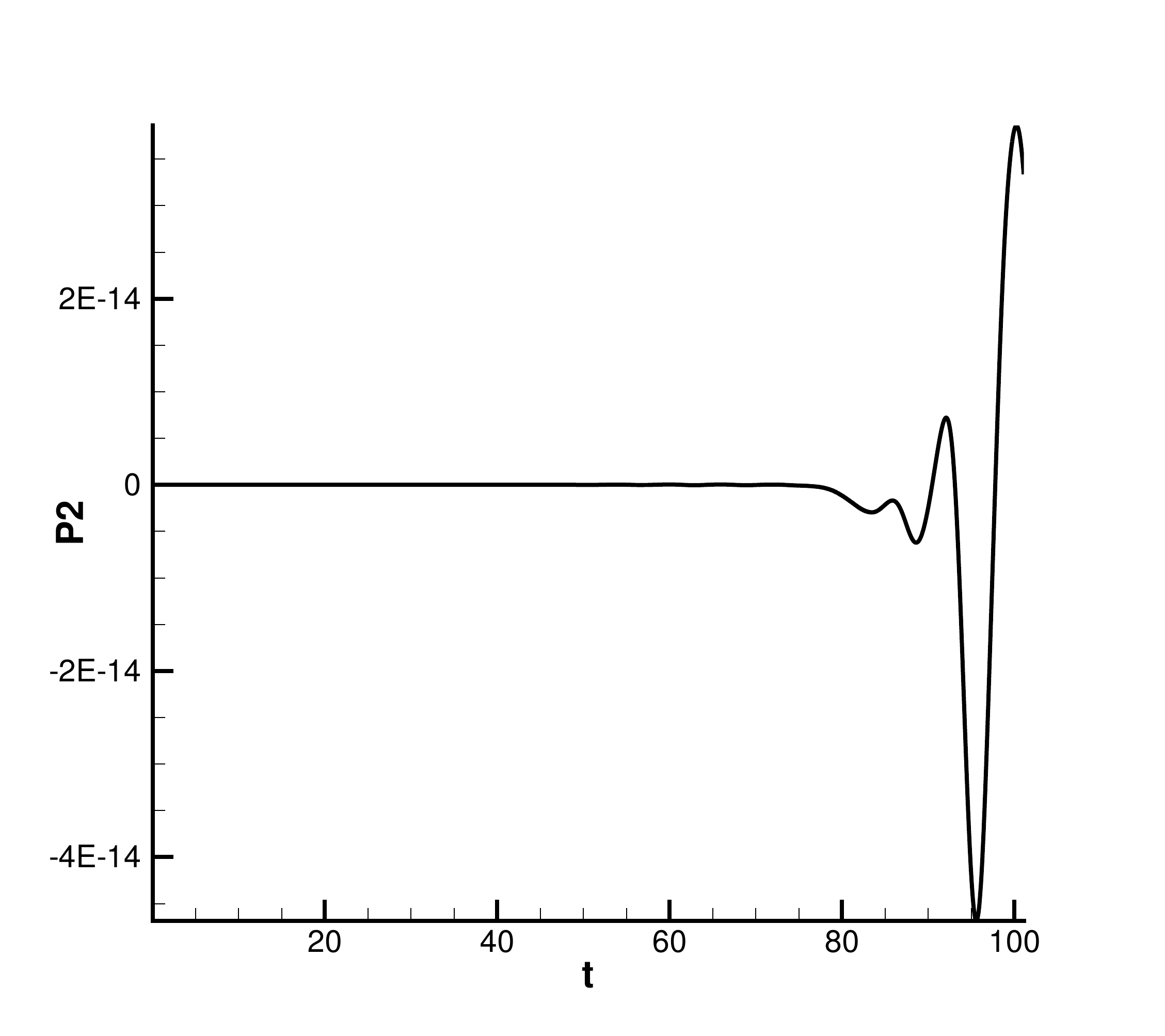}}
	\end{center}
	\caption{SW instability with parameter choice 1. Time evolution of error of momentum by upwind flux for the Maxwell's equations. Sparse grid: $N=8$, $k=3$.  Adaptive sparse grid: $N=6$, $k=3$, $\eps=2\times10^{-7}$.}
	\label{mom1}
\end{figure}

First, we consider the SW instability with parameter   choice 1 up to $t=100$. We take $N=8$, $k=3$ for the sparse grid scheme and  $N=6$, $k=3$, $\eps = 2\times10^{-7}$ for the adaptive sparse grid scheme. 
We  measure the errors in the conserved macroscopic quantities for both schemes. The time evolution of relative errors of mass (charge), total energy, and errors of momentum with the upwind flux are plotted in Figures \ref{masste1} and \ref{mom1}.
  When $t\leq75$, our scheme can conserve mass, total energy and momentum very well.
For longer time period, the errors start to accumulate, particularly for mass, total energy and the momentum $P_1.$ We notice that the largest relative error in mass is on the order of $10^{-4}$ and $10^{-8}$ for sparse grid and adaptive sparse grids respectively. Compared with  the results in \cite{cheng2014discontinuous} for traditional RKDG method, the error in conservation is much larger. 
According to the analysis performed in Section \ref{sec:method}, this is primarily due to the boundary effect at the velocity domain. Such errors can be reduced by enlarging $\Oxi$ or enhancing the resolution further. The adaptive sparse grid scheme offers better resolution (which is also validated by later plots in this section), therefore the contributions from the domain boundary is smaller, and the simulation retains better conservation in all quantities.
As for total energy,  the largest relative error  is on the order of $10^{-3}$ and $10^{-5}$ for sparse grid and adaptive sparse grids respectively. We also performed simulations based on alternating flux for the Maxwell part. The results show no   difference. This implies that the major contribution of error comes from the velocity domain boundary cutoff rather than the jump terms in the field, i.e. $\Theta_{h,2}(t)$ is much smaller than $\Theta_{h,3}(t)$ terms in Theorem \ref{lem:EneC}. In fact, we find that there is no significant   difference in performance for the two choices of fluxes for long time simulations. Therefore, for brevity, we only report numerical results for the scheme using the upwind numerical flux for the Maxwell DG solver from now on.





\begin{figure}[htb]
	\begin{center}
		\subfigure[ Sparse grid, kinetic energies]{\includegraphics[width=3in,angle=0]{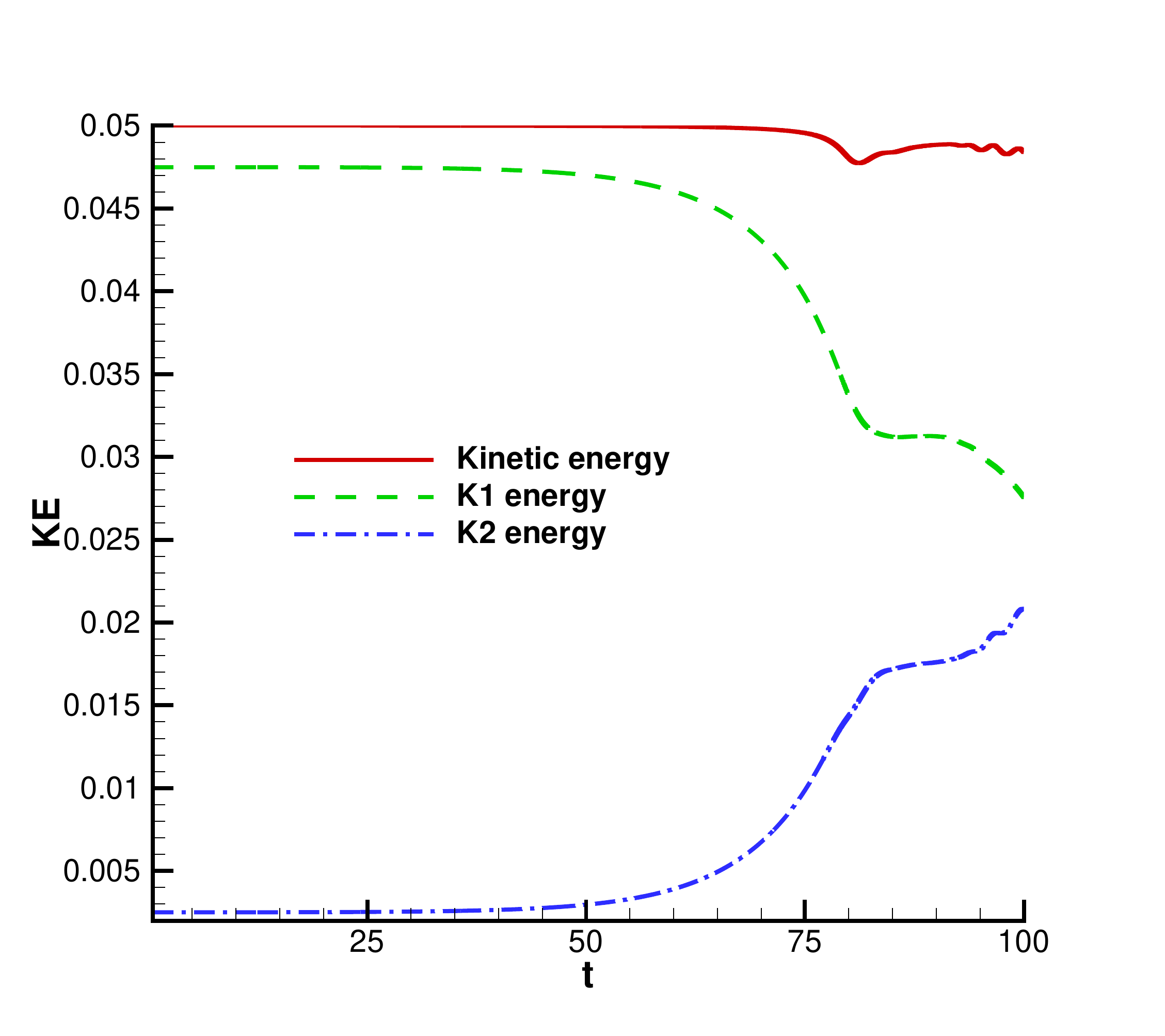}}
		\subfigure[ Adaptive, kinetic energies]{\includegraphics[width=3in,angle=0]{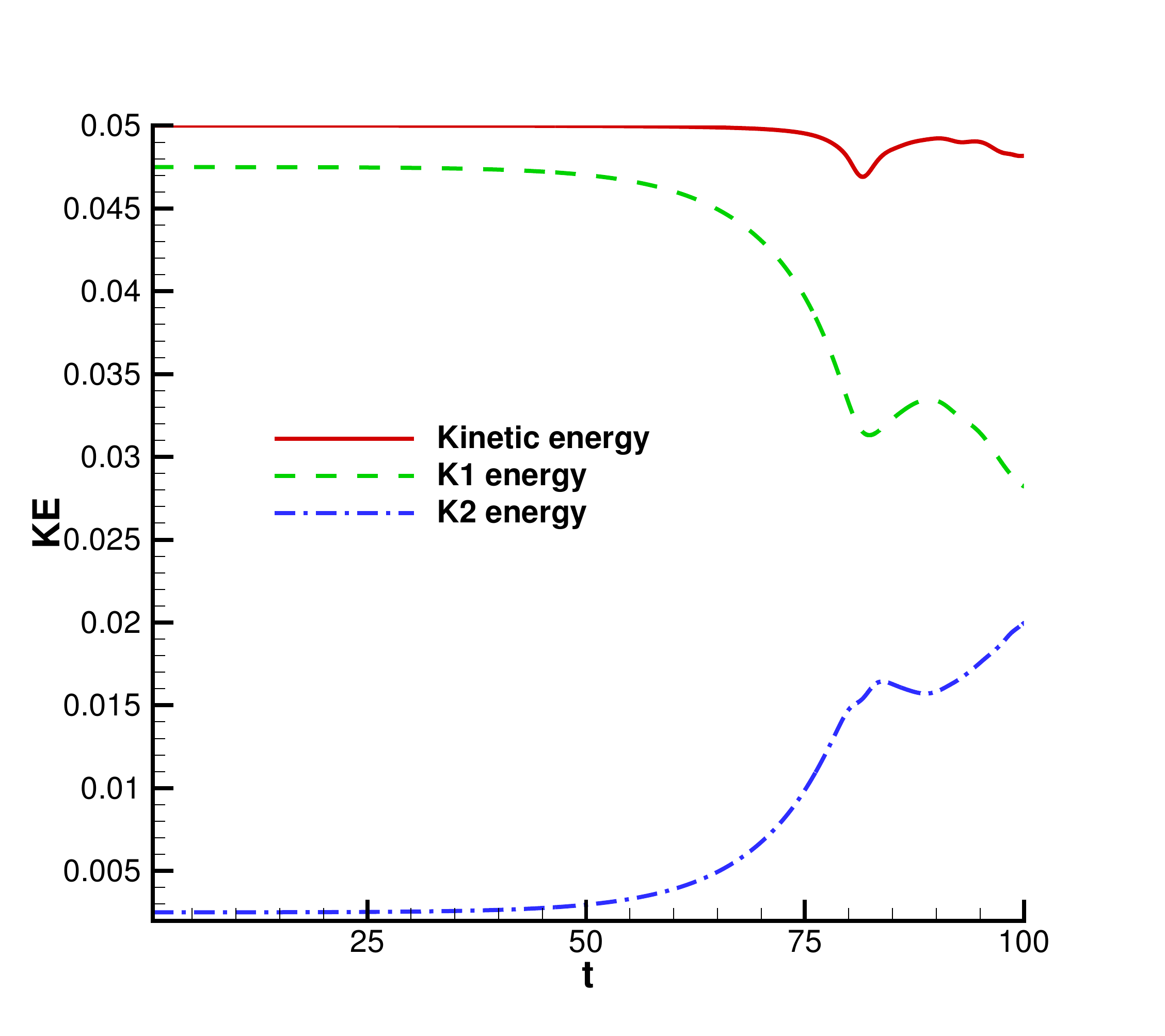}}
		\subfigure[Sparse grid, field energies]{\includegraphics[width=3in,angle=0]{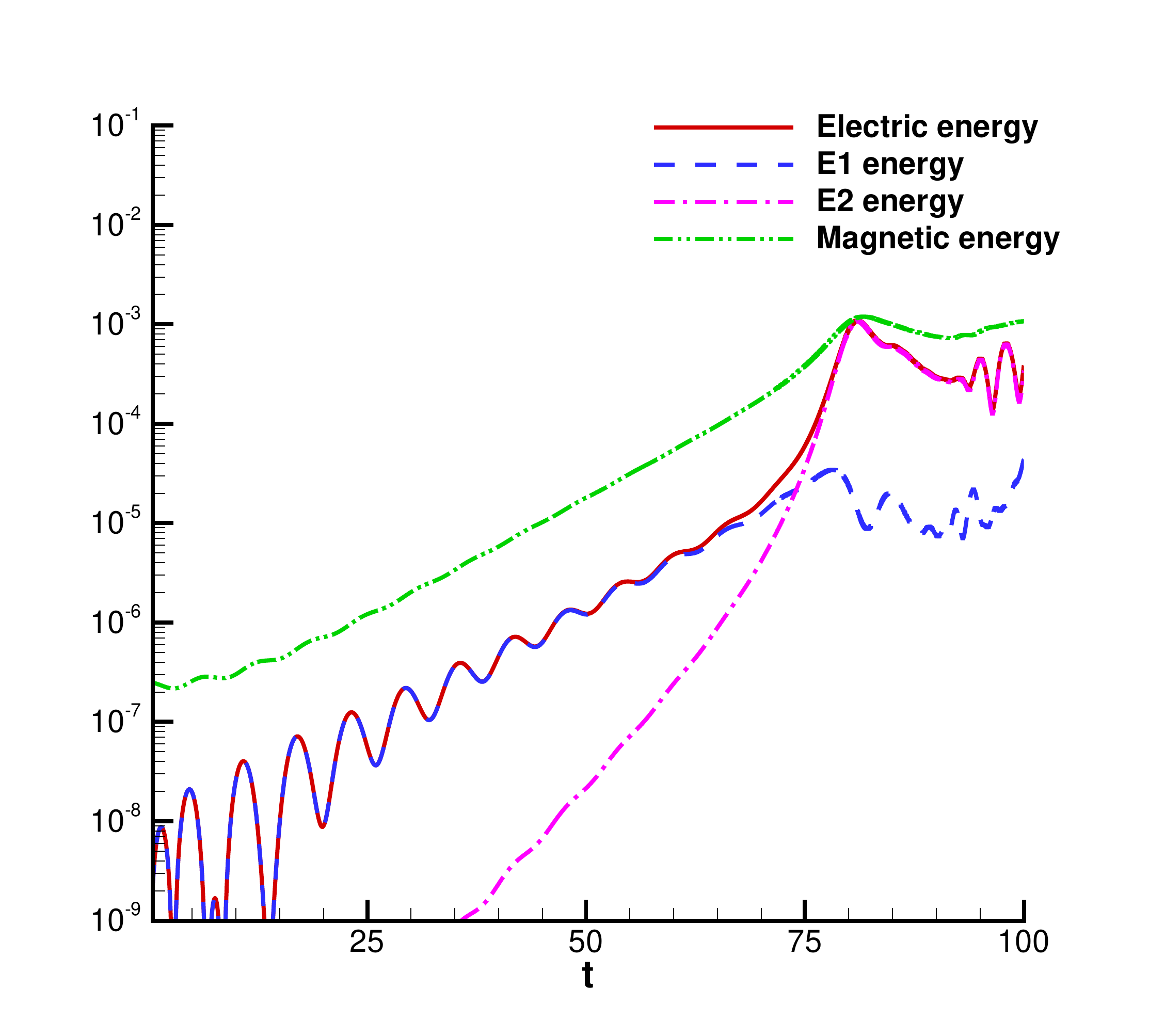}}
		\subfigure[Adaptive, field energies]{\includegraphics[width=3in,angle=0]{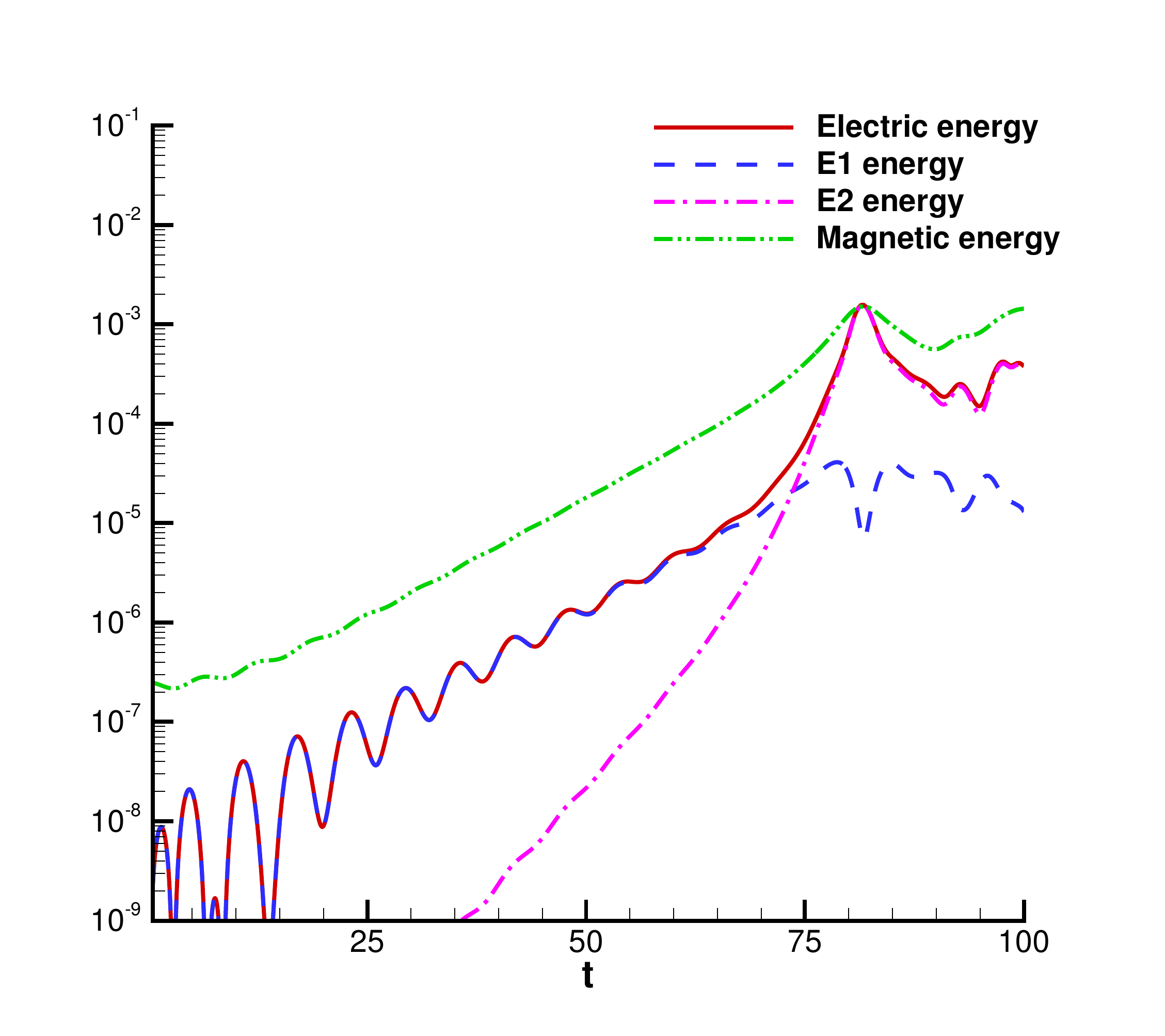}}
	\end{center}
	
	\caption{SW instability with parameter choice 1. Time evolution of kinetic, electric  and magnetic energies by upwind flux for the Maxwell's equations.  Sparse grid: $N=8$, $k=3$.  Adaptive sparse grid: $N=6$, $k=3$, $\eps=2\times10^{-7}$.}
	\label{keeme1}
\end{figure}

The time evolution of kinetic, electric  and magnetic energies are presented in Figure \ref{keeme1}.
The results from both simulations are very close to \cite{cheng2014discontinuous}. 
The transfer of kinetic energy from one component to the other is evident, and there is a small decay which is converted into the field energy. The magnetic and inductive electric fields grow initially at a linear growth rate. 
For $t\geq70$, kinetic effects come into play and the instability saturates.  The magnetic energy comes statistically constant, while the electric energy reaches its maximum value at saturation and then starts to decrease.     

\begin{figure}[htb]
	\begin{center}
		\subfigure[Sparse grid, Log Fourier modes of $E_1$]{\includegraphics[width=2.8in,angle=0]{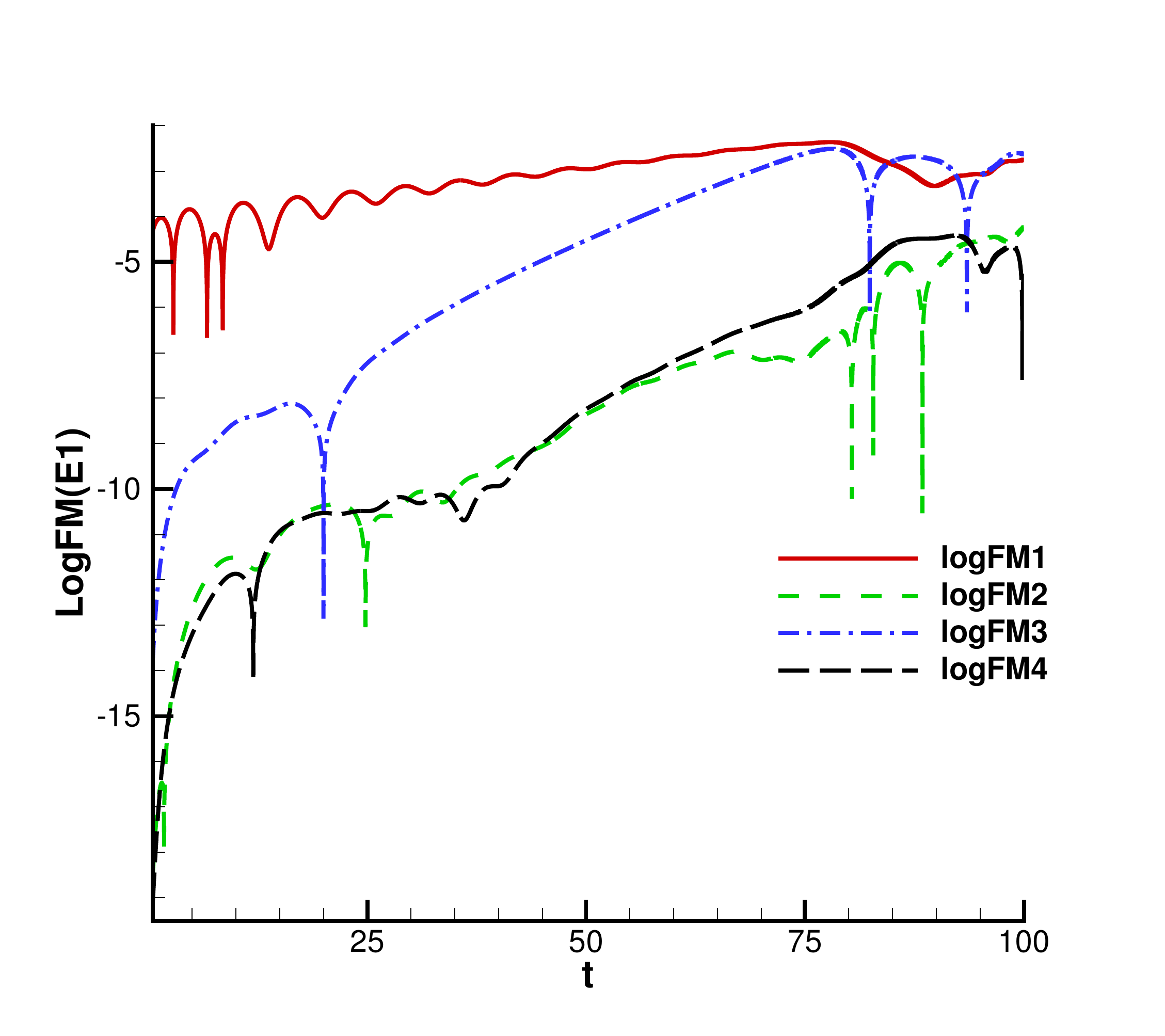}}
		\subfigure[Adaptive, Log Fourier modes of $E_1$]{\includegraphics[width=2.8in,angle=0]{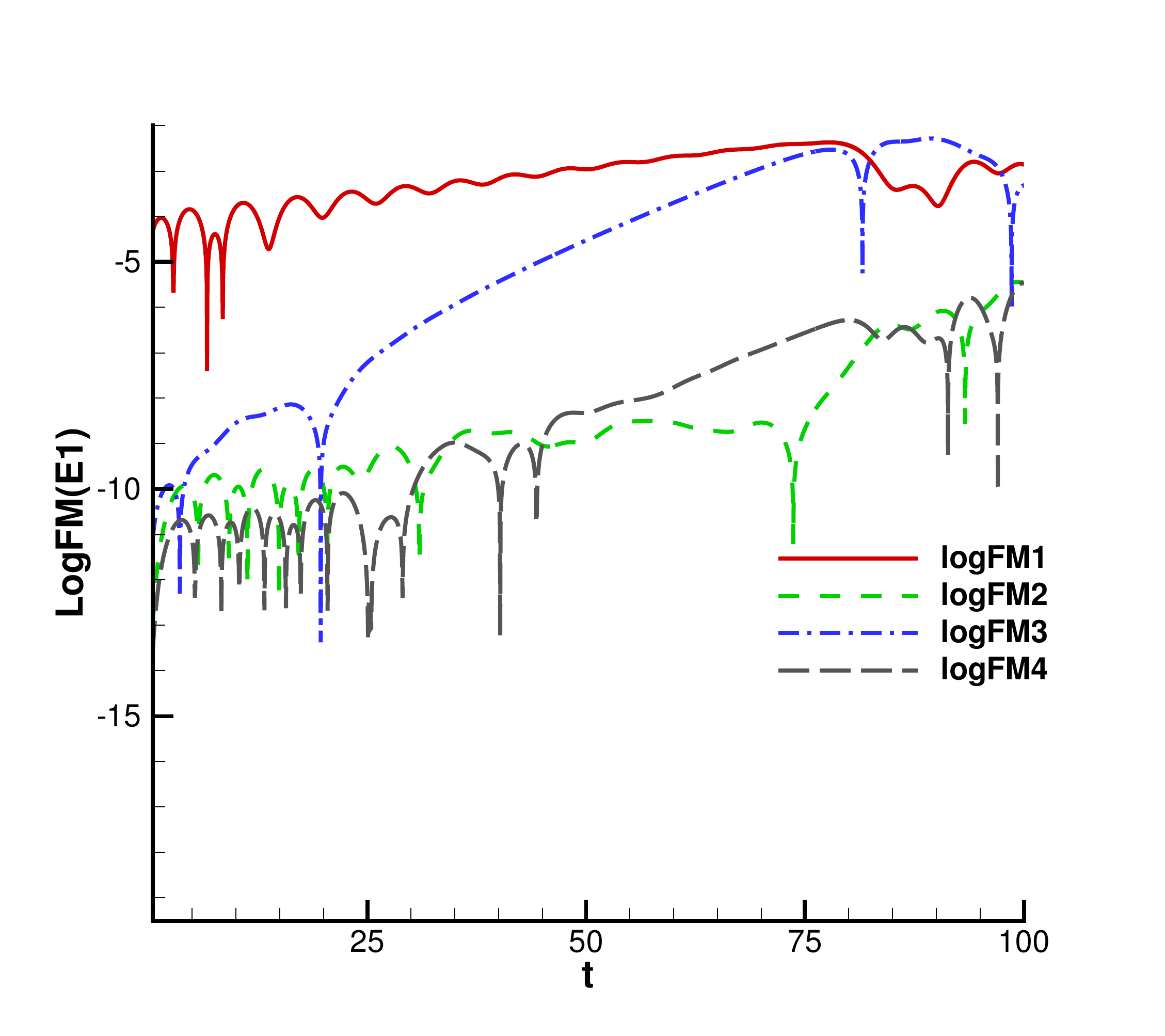}}
		\subfigure[Sparse grid, Log Fourier modes of $E_2$]{\includegraphics[width=2.8in,angle=0]{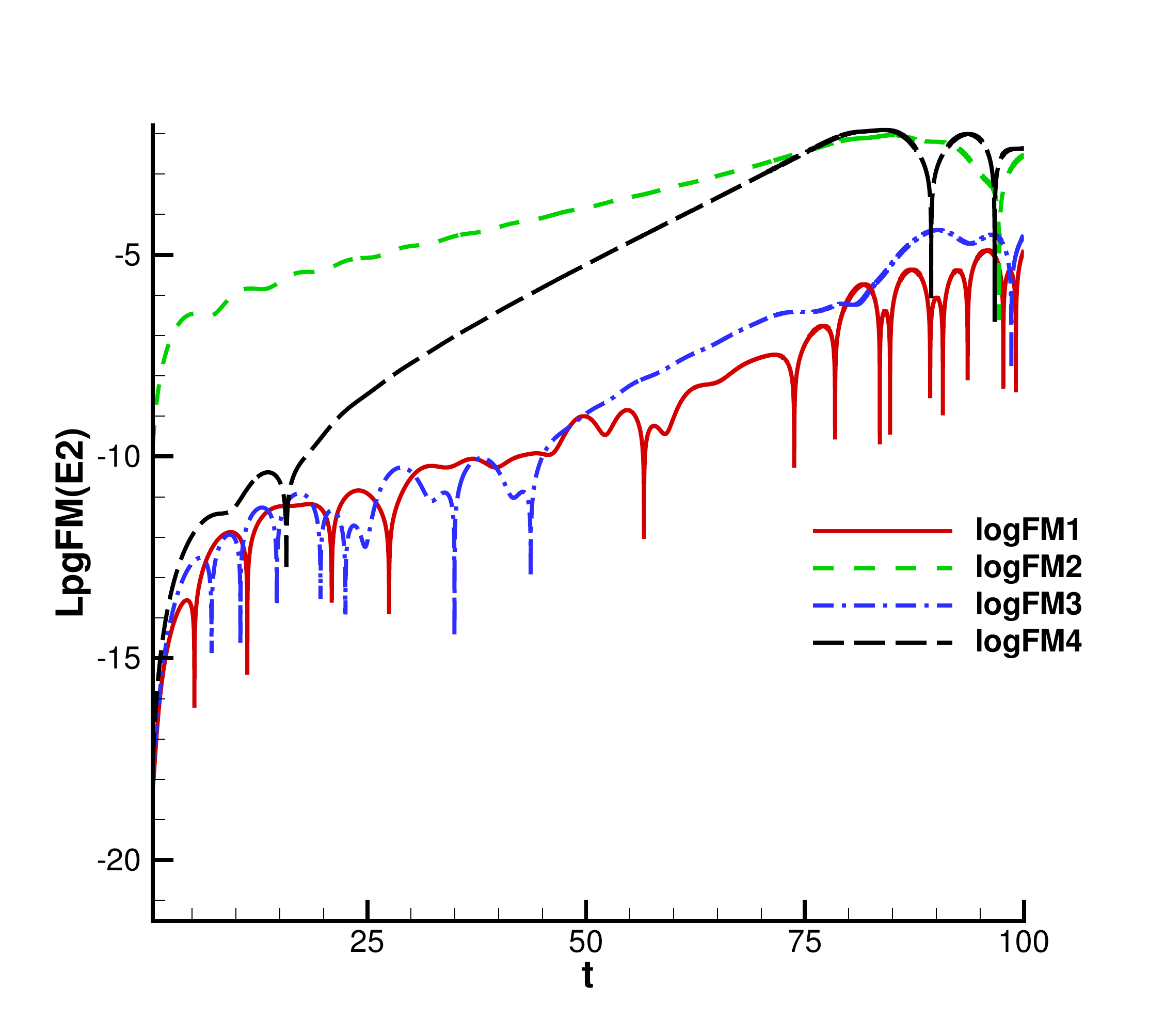}}
		\subfigure[Adaptive, Log Fourier modes of $E_2$]{\includegraphics[width=2.8in,angle=0]{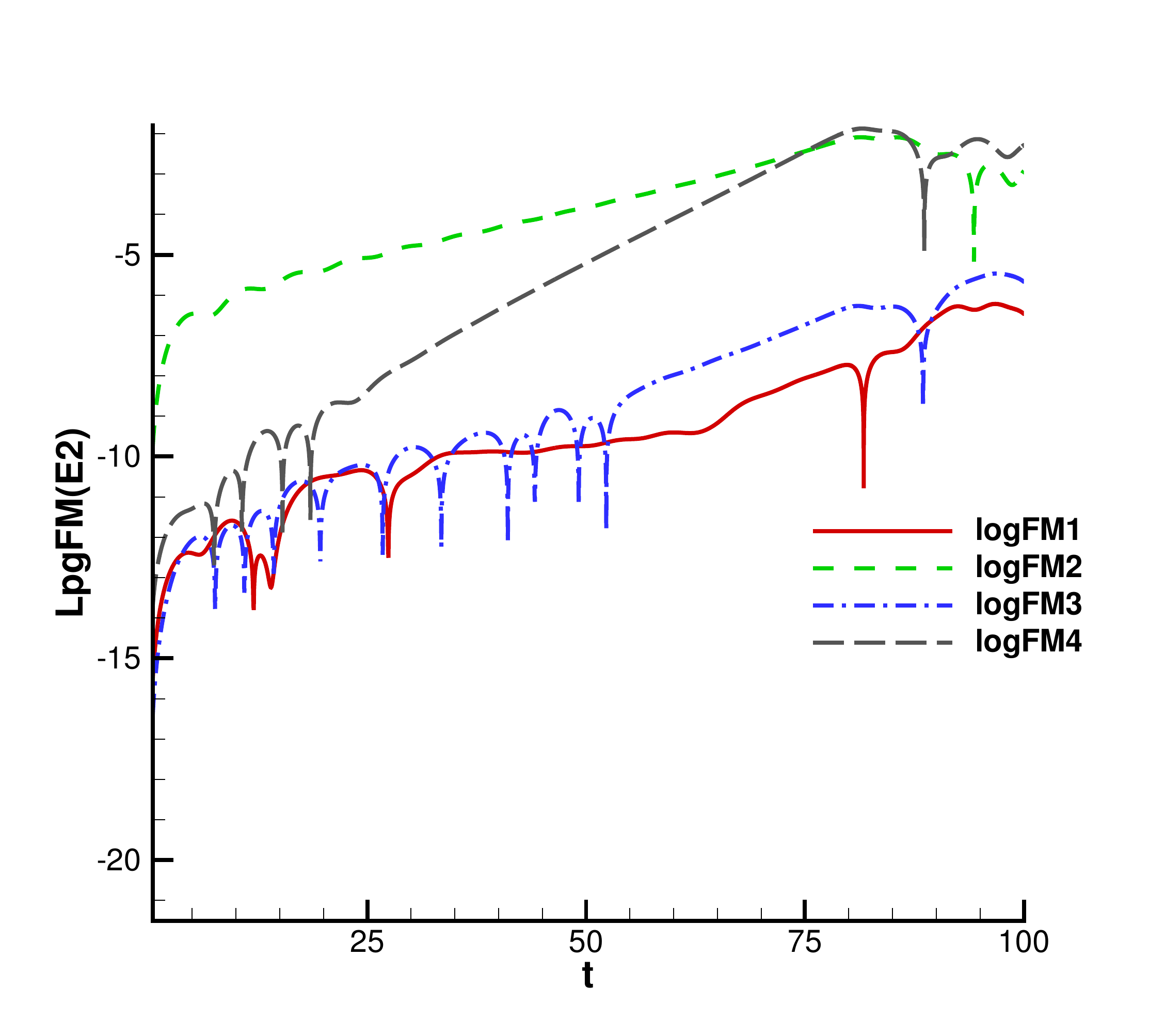}}
		\subfigure[Sparse grid, Log Fourier modes of $B_3$]{\includegraphics[width=2.8in,angle=0]{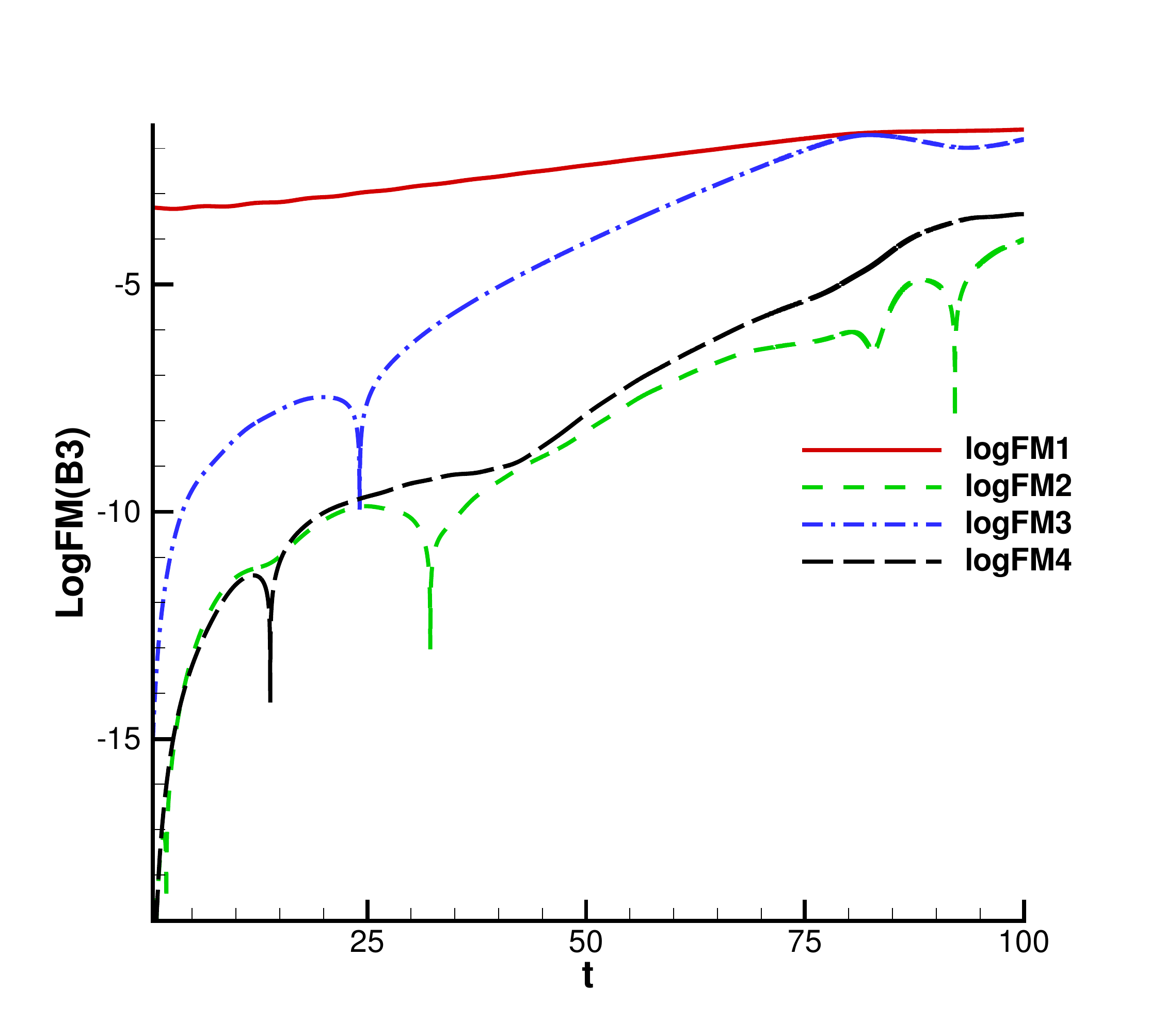}}
		\subfigure[Adaptive, Log Fourier modes of $B_3$]{\includegraphics[width=2.8in,angle=0]{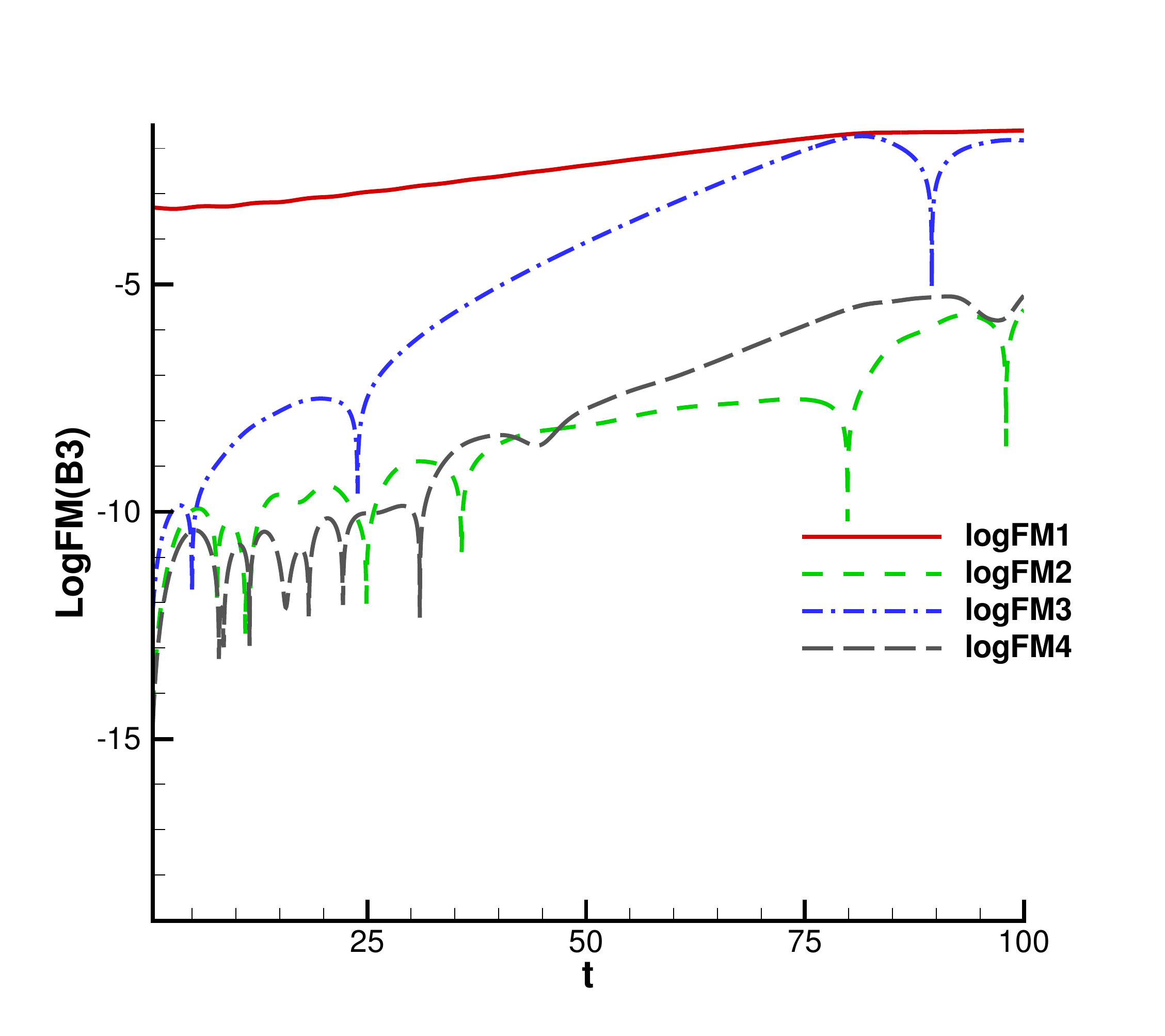}}
	\end{center}
	\caption{SW instability with parameter choice 1. The first four Log Fourier modes of $E_1$, $E_2$, $B_3$ by upwind flux for the Maxwell's equations. Sparse grid: $N=8$, $k=3$.  Adaptive sparse grid: $N=6$, $k=3$, $\eps=2\times10^{-7}$.}
	\label{logfm1}
\end{figure}

In Figure \ref{logfm1}, we present the first four Log Fourier modes of $E_1$, $E_2$, $B_3$.
Here, the $n$-th Log Fourier mode for a function $W(x,t)$ is defined as
$$\log FM_n(t)=\log_{10}\left(\frac1L\sqrt{\left| \int_0^L W(x,t)\,\sin(knx)\,dx\right|^2+\left| \int_0^L W(x,t)\,\cos(knx)\,dx\right|^2}\right).$$ 
Here, $k=k_0=0.2$. The first and third Log Fourier modes of $E_1$ and $B_3$ as well as the second and fourth Log Fourier modes of $E_2$ agree relatively well with \cite{cheng2014discontinuous}.    
For the remaining less significant modes, our schemes produce   larger magnitudes than those in \cite{cheng2014discontinuous}. 

\begin{figure}[htb]
	\begin{center}
		\subfigure[ Sparse grid, $x_2=0.05 \pi,  \, t=55.$]{\includegraphics[width=2.8in,angle=0]{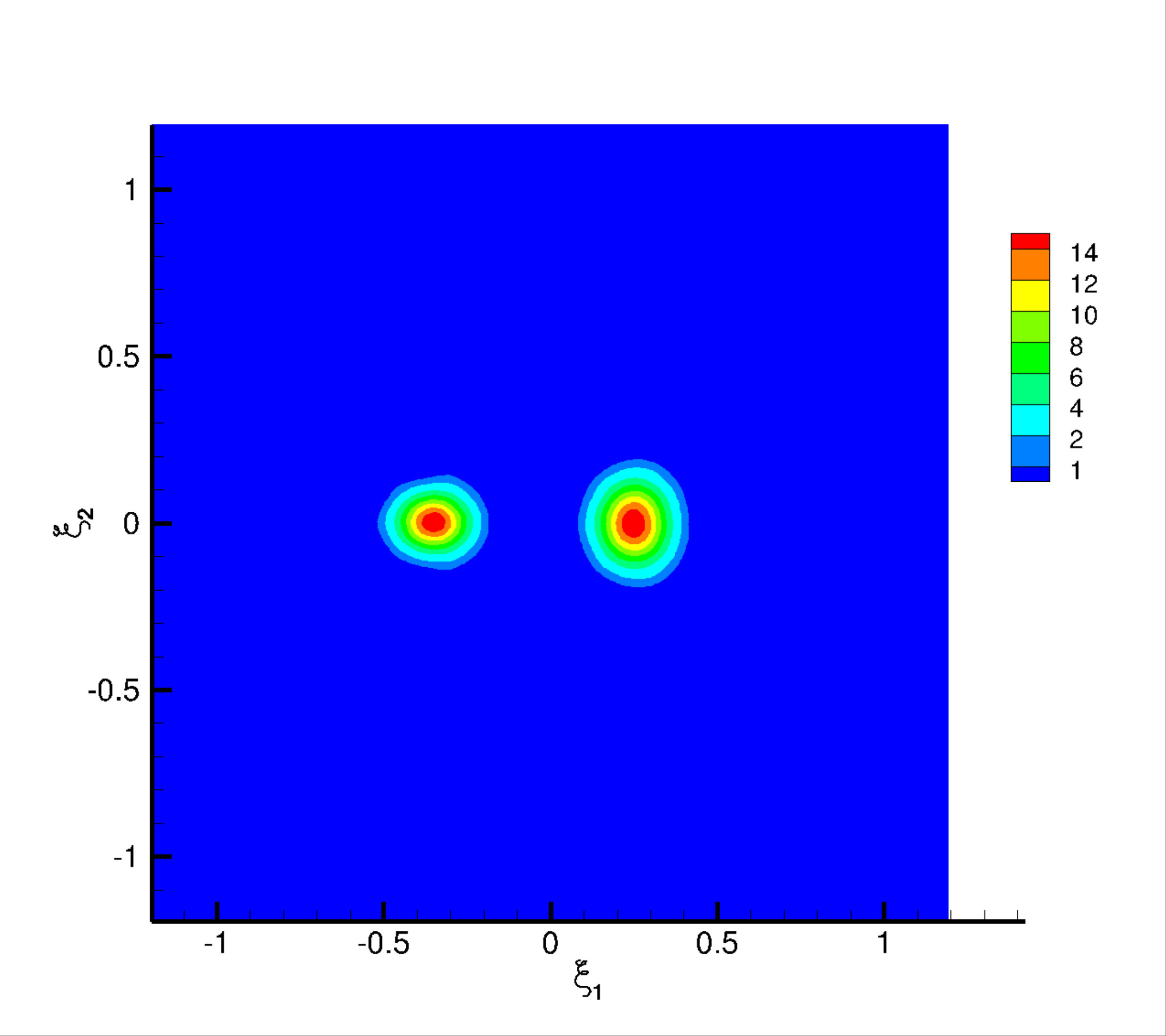}}
		\subfigure[ Adaptive, $x_2=0.05 \pi,  \, t=55.$]{\includegraphics[width=2.8in,angle=0]{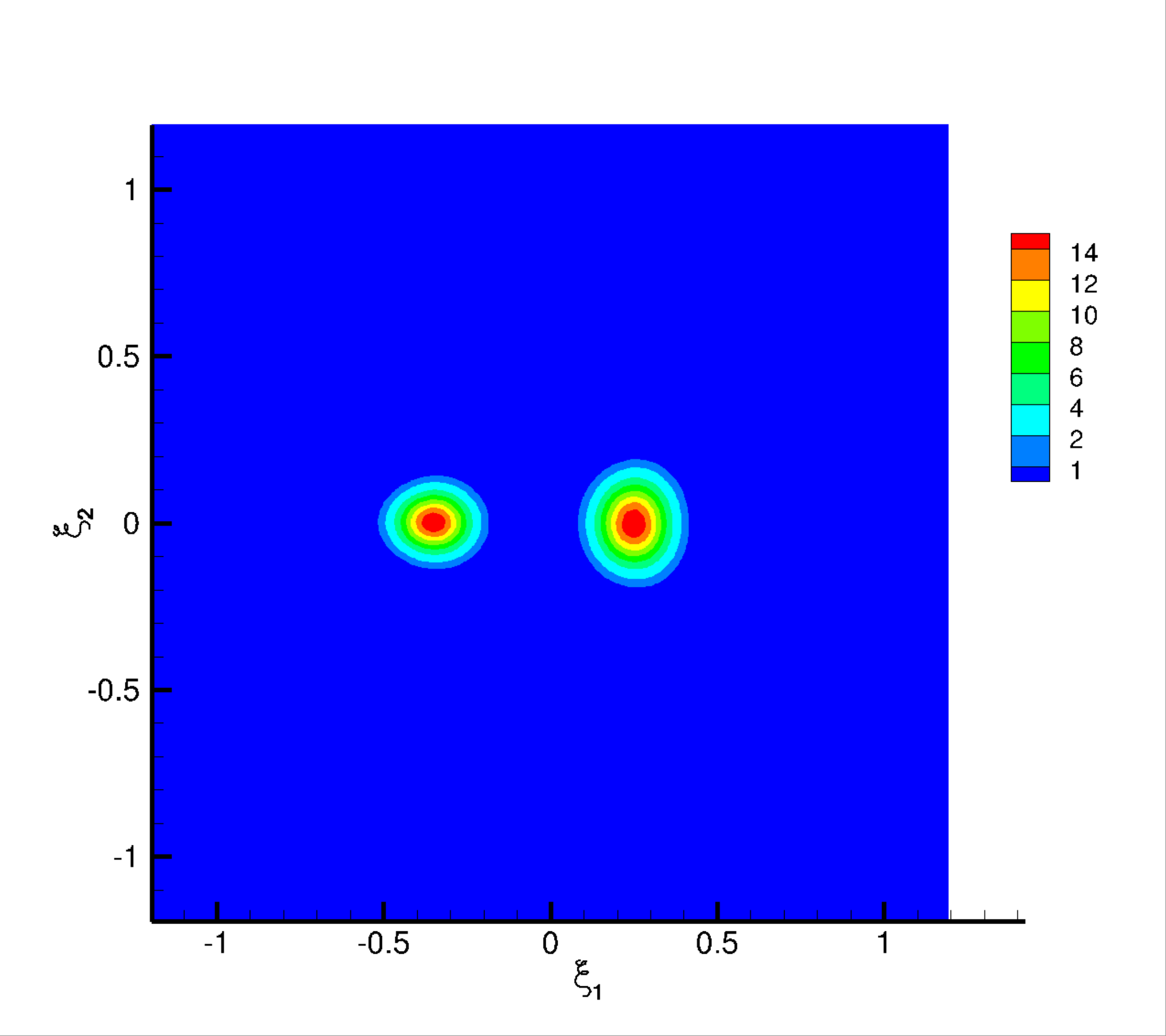}}
		\subfigure[Sparse grid,  $x_2=0.05 \pi, \,  t=82.$]{\includegraphics[width=2.8in,angle=0]{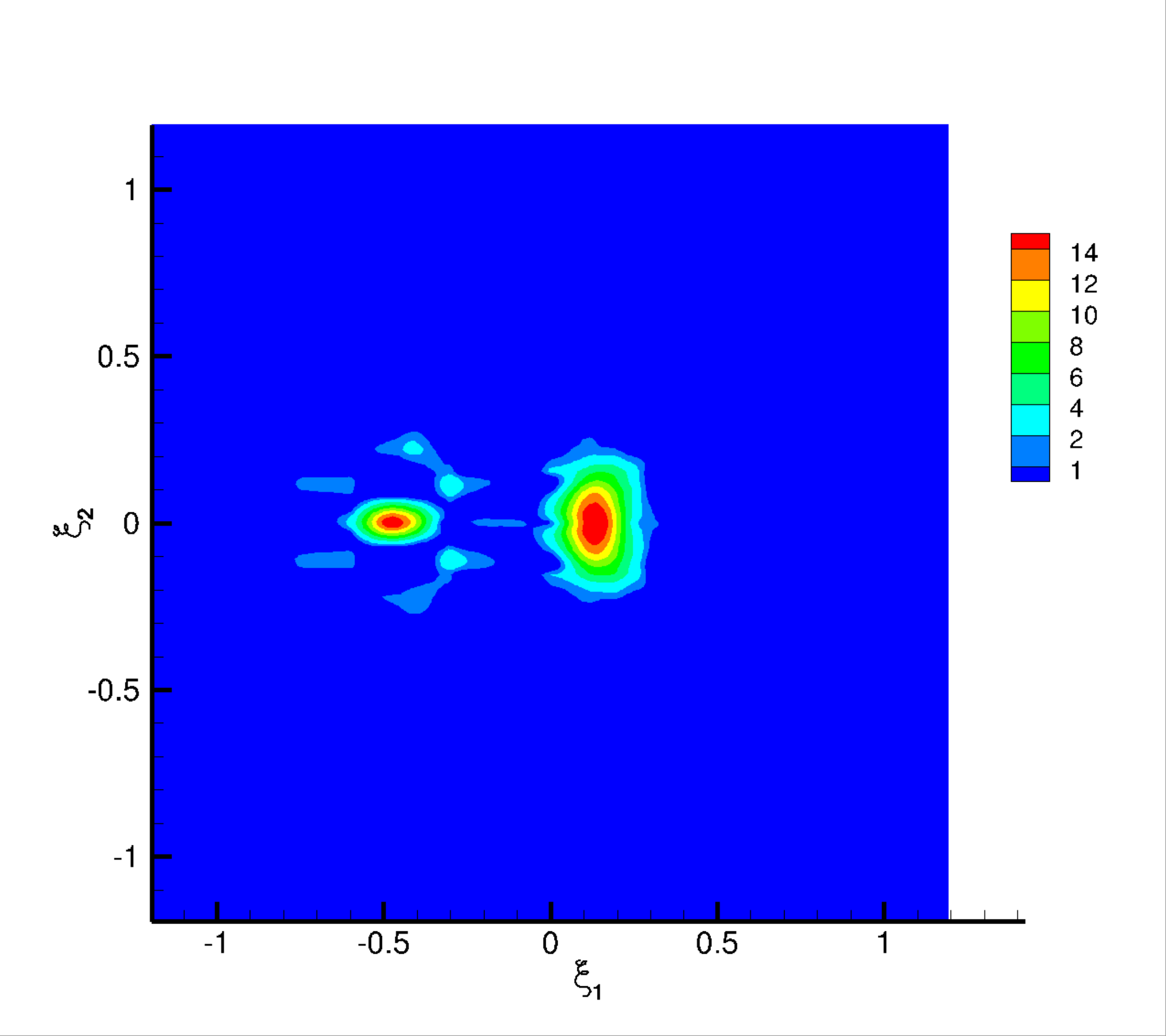}}
		\subfigure[Adaptive,  $x_2=0.05 \pi, \,  t=82.$]{\includegraphics[width=2.8in,angle=0]{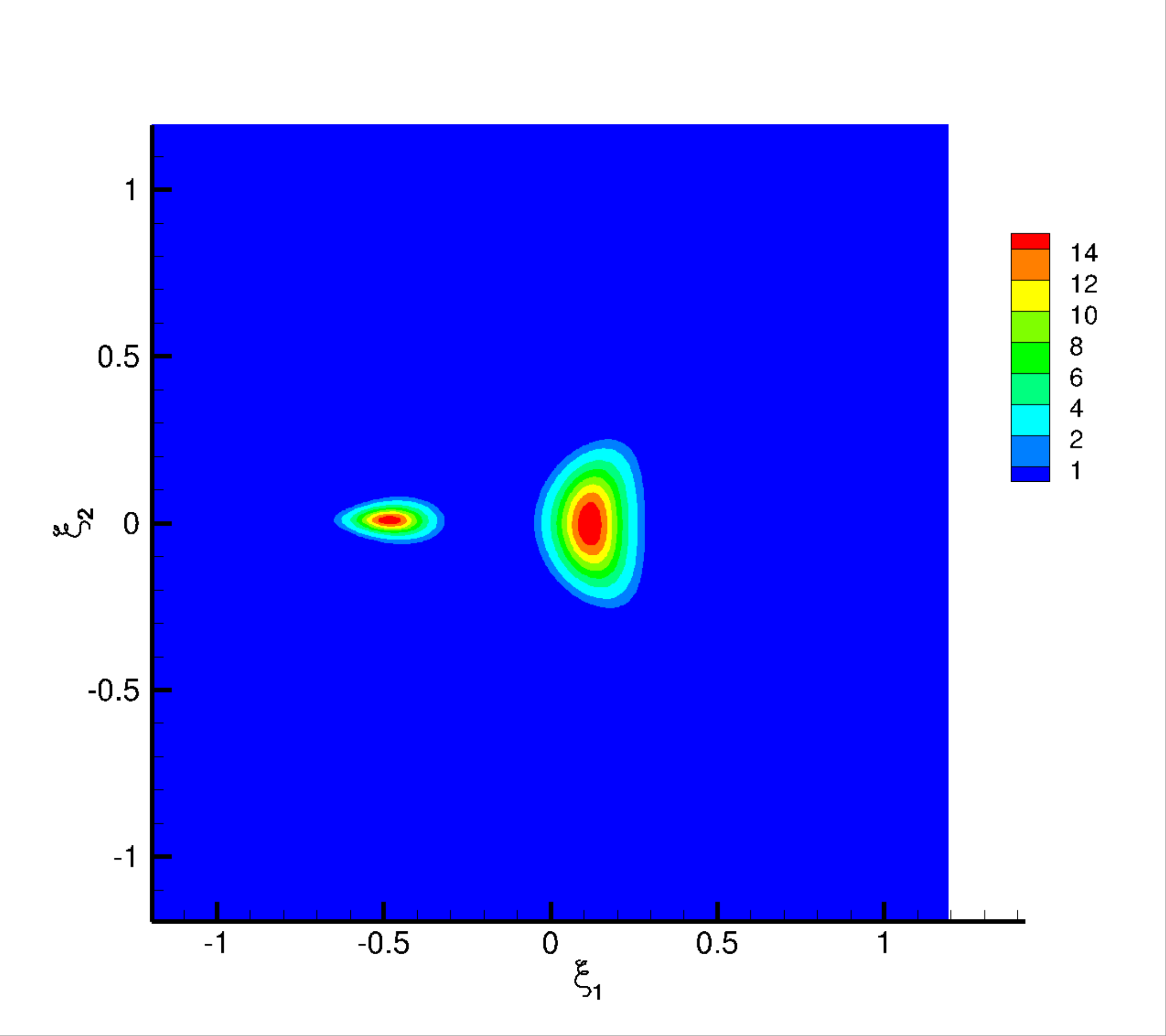}}
		\subfigure[ Sparse grid, $x_2=0.05 \pi,  \, t=100.$]{\includegraphics[width=2.8in,angle=0]{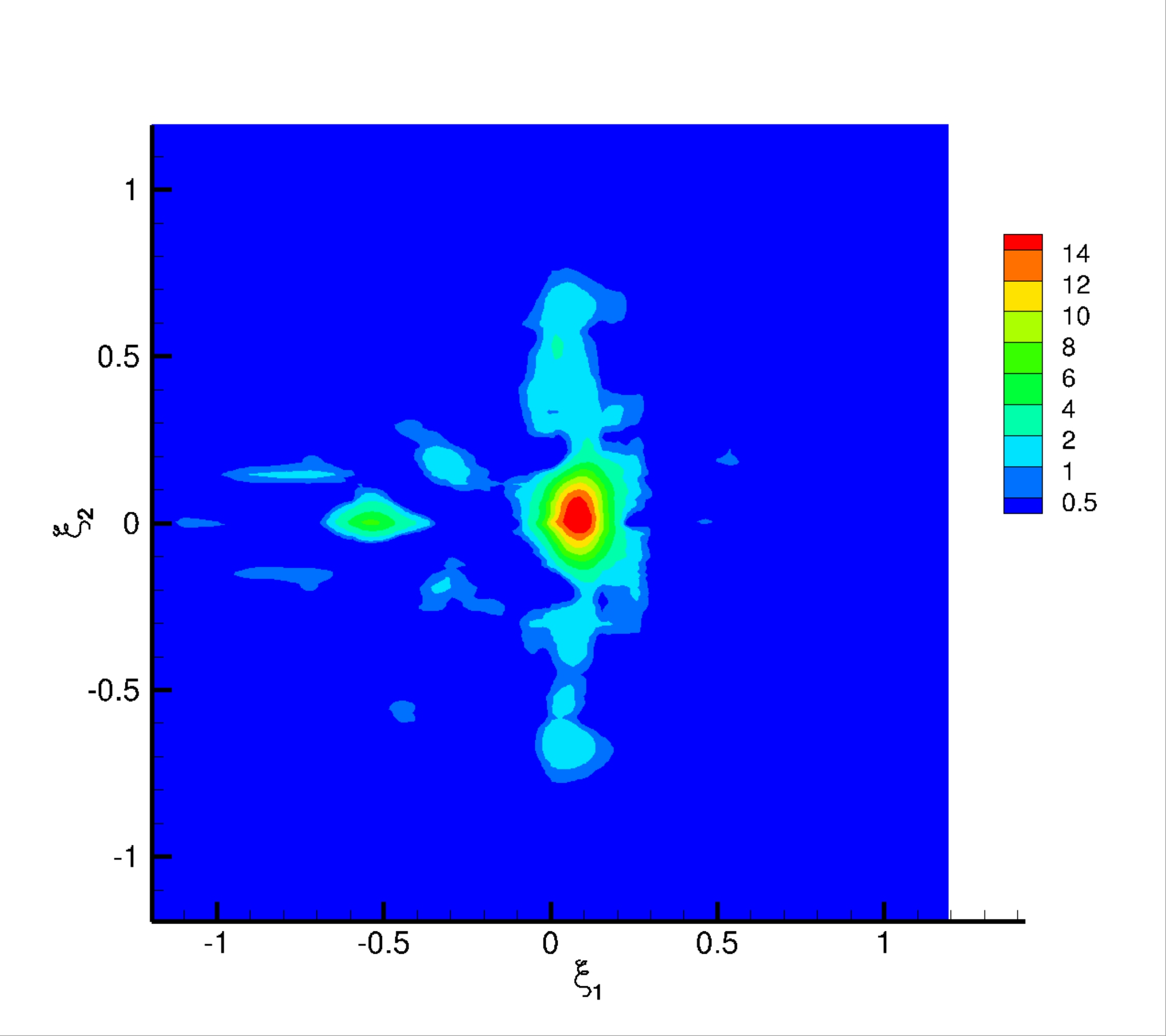}}
		\subfigure[ Adaptive, $x_2=0.05 \pi,  \, t=100.$]{\includegraphics[width=2.8in,angle=0]{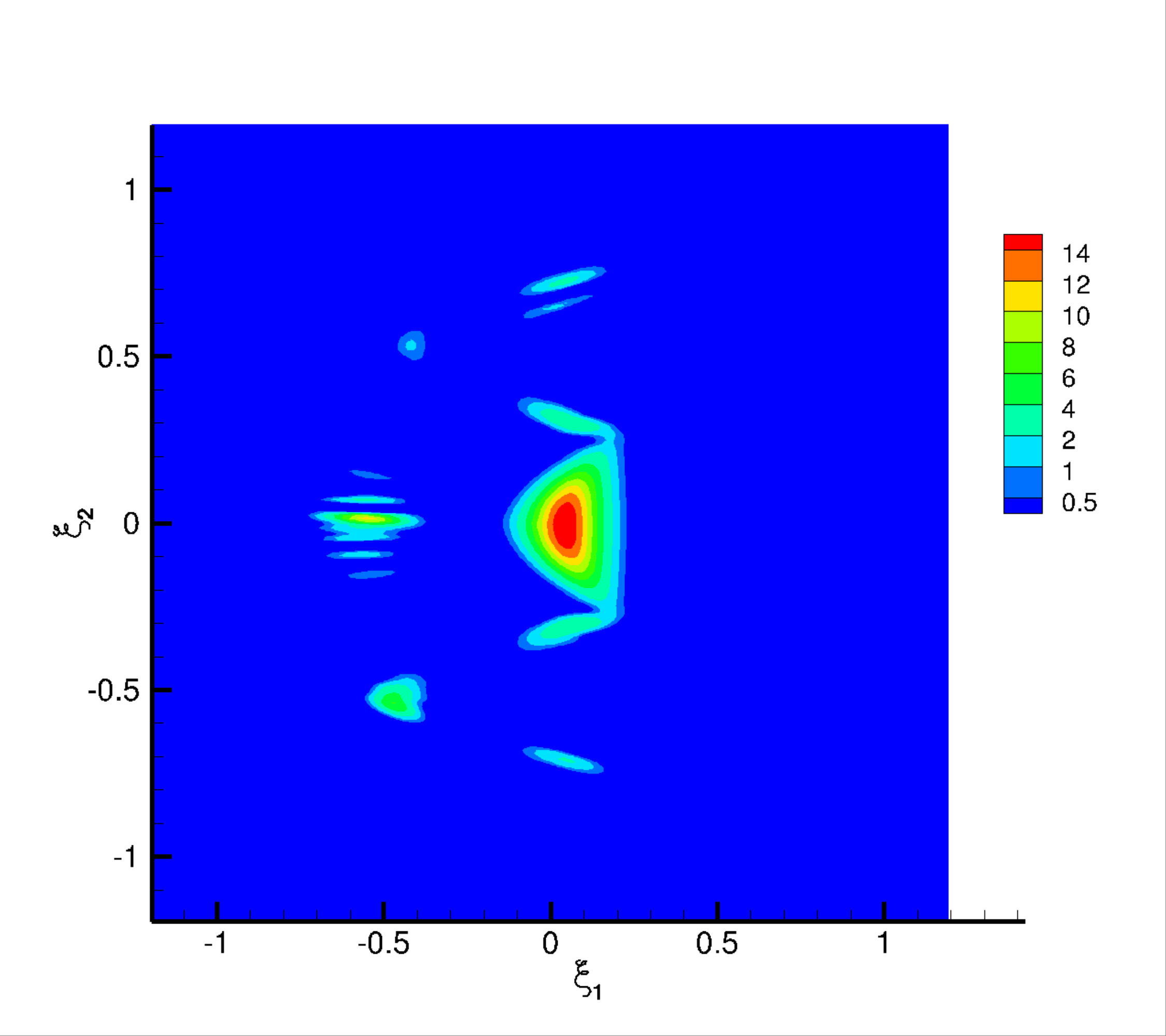}}
	\end{center}
	\caption{SW instability with parameter choice 1. 2D contour plots of the computed distribution function $f_h$ by upwind flux for the Maxwell's equations. Sparse grid: $N=8$, $k=3$.  Adaptive sparse grid: $N=6$, $k=3$, $\eps=2\times10^{-7}$.}
	\label{contour1}
\end{figure}
\begin{figure}[htb]
	\begin{center}
		\subfigure[ Sparse grid, $ t=55.$]{\includegraphics[width=2.8in,angle=0]{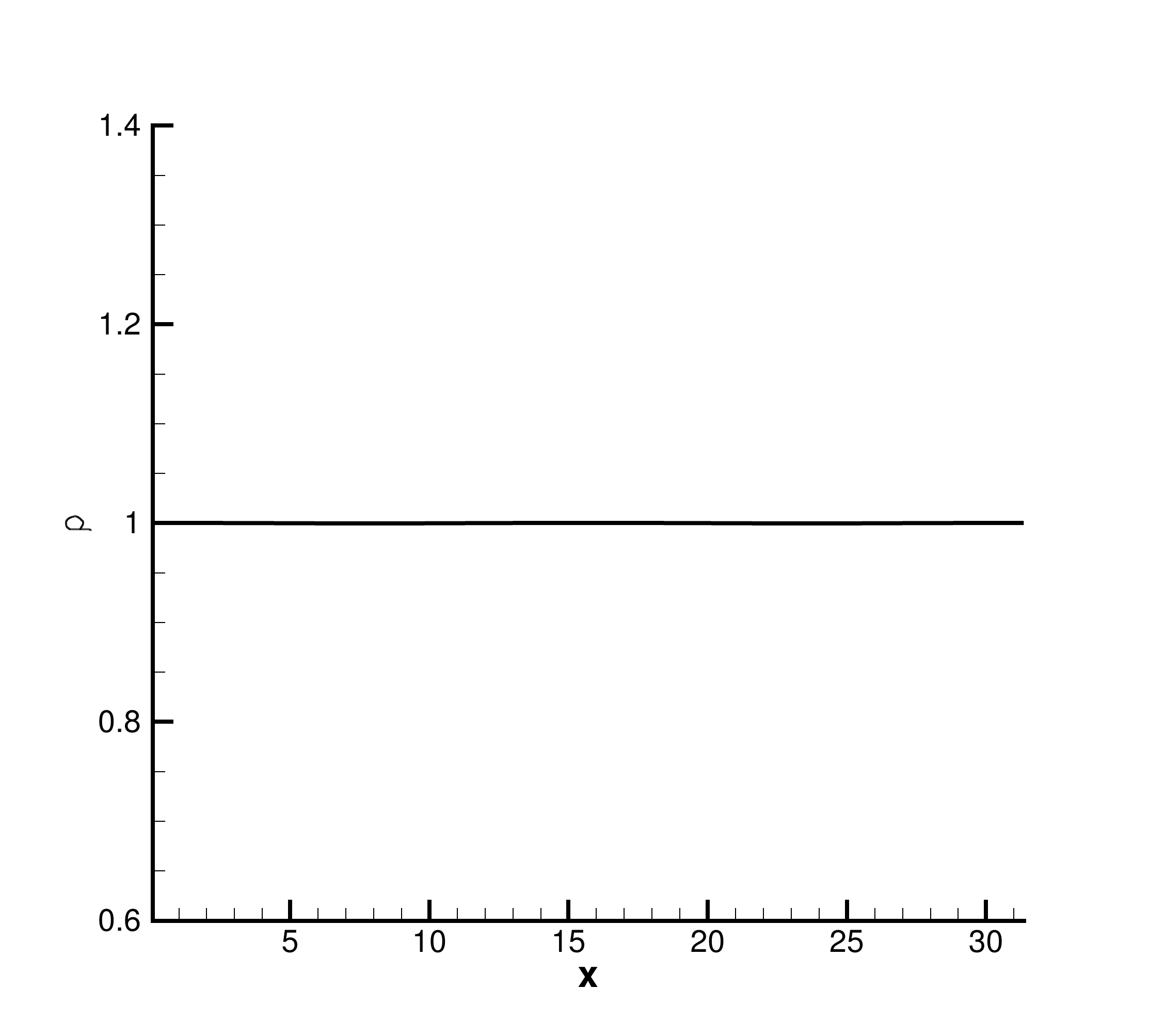}}
		\subfigure[ Adaptive, $ t=55.$]{\includegraphics[width=2.8in,angle=0]{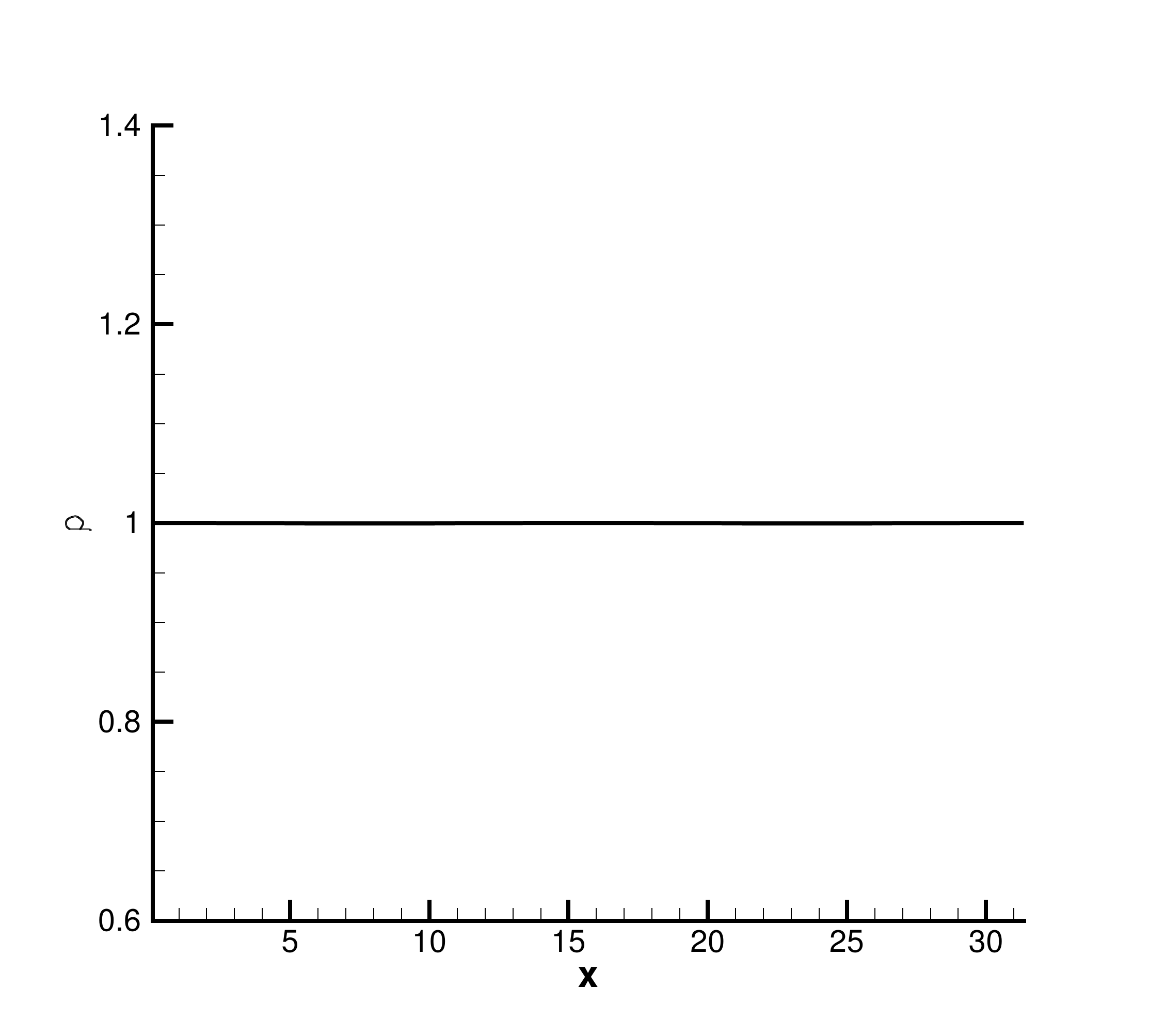}}
		\subfigure[ Sparse grid, $ t=82.$]{\includegraphics[width=2.8in,angle=0]{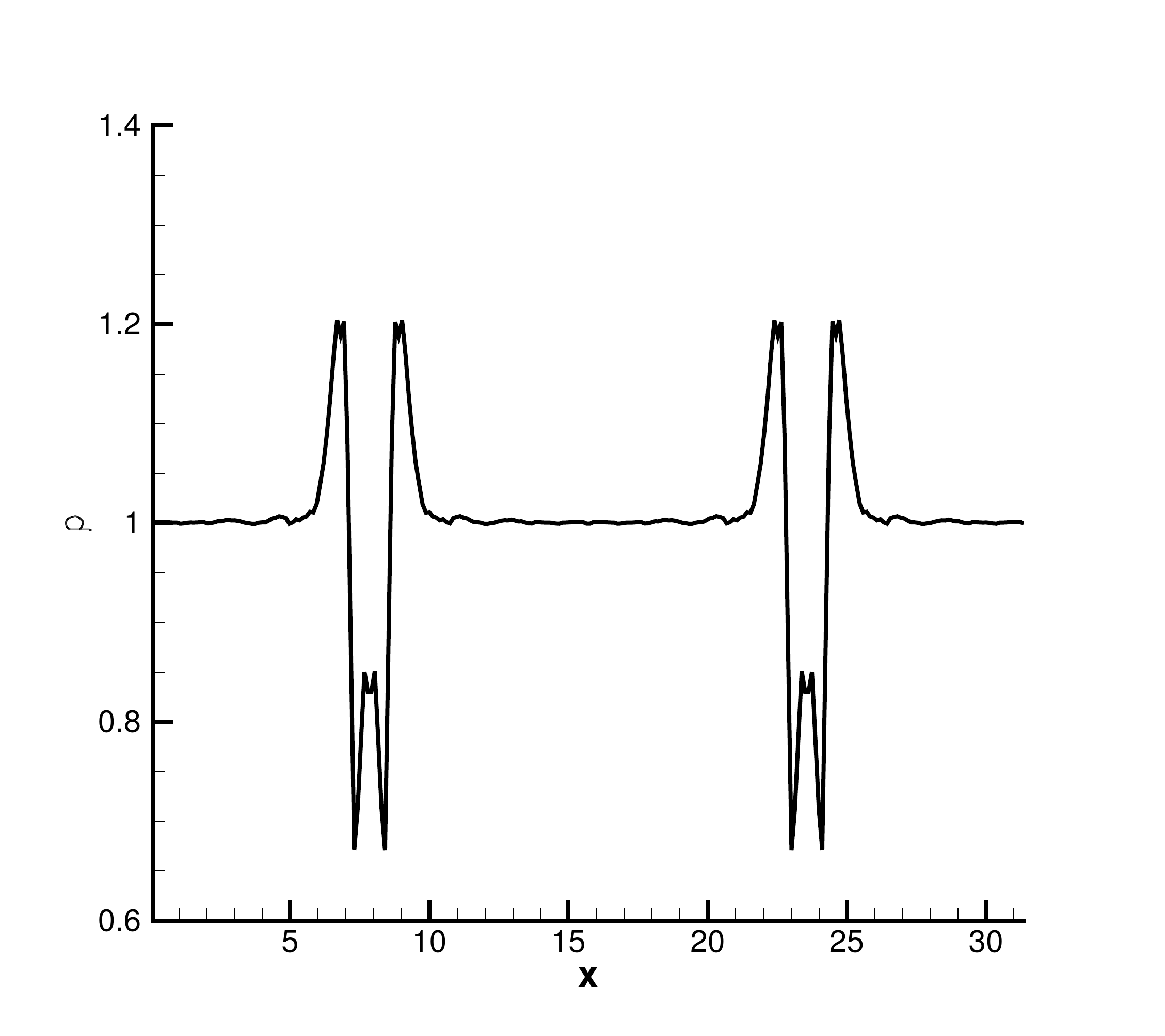}}
		\subfigure[ Adaptive, $ t=82.$]{\includegraphics[width=2.8in,angle=0]{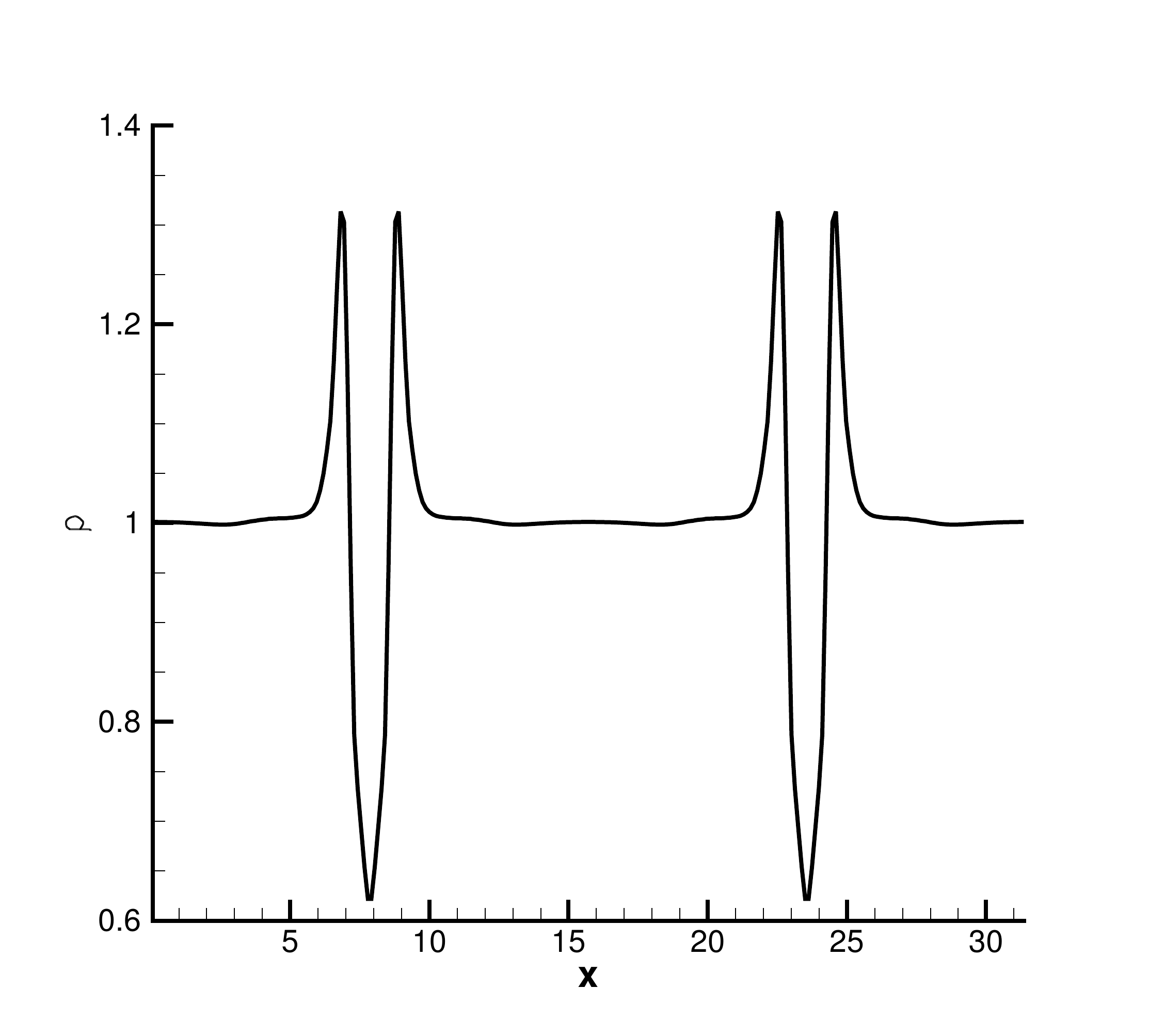}}
		\subfigure[ Sparse grid, $t=100.$]{\includegraphics[width=2.8in,angle=0]{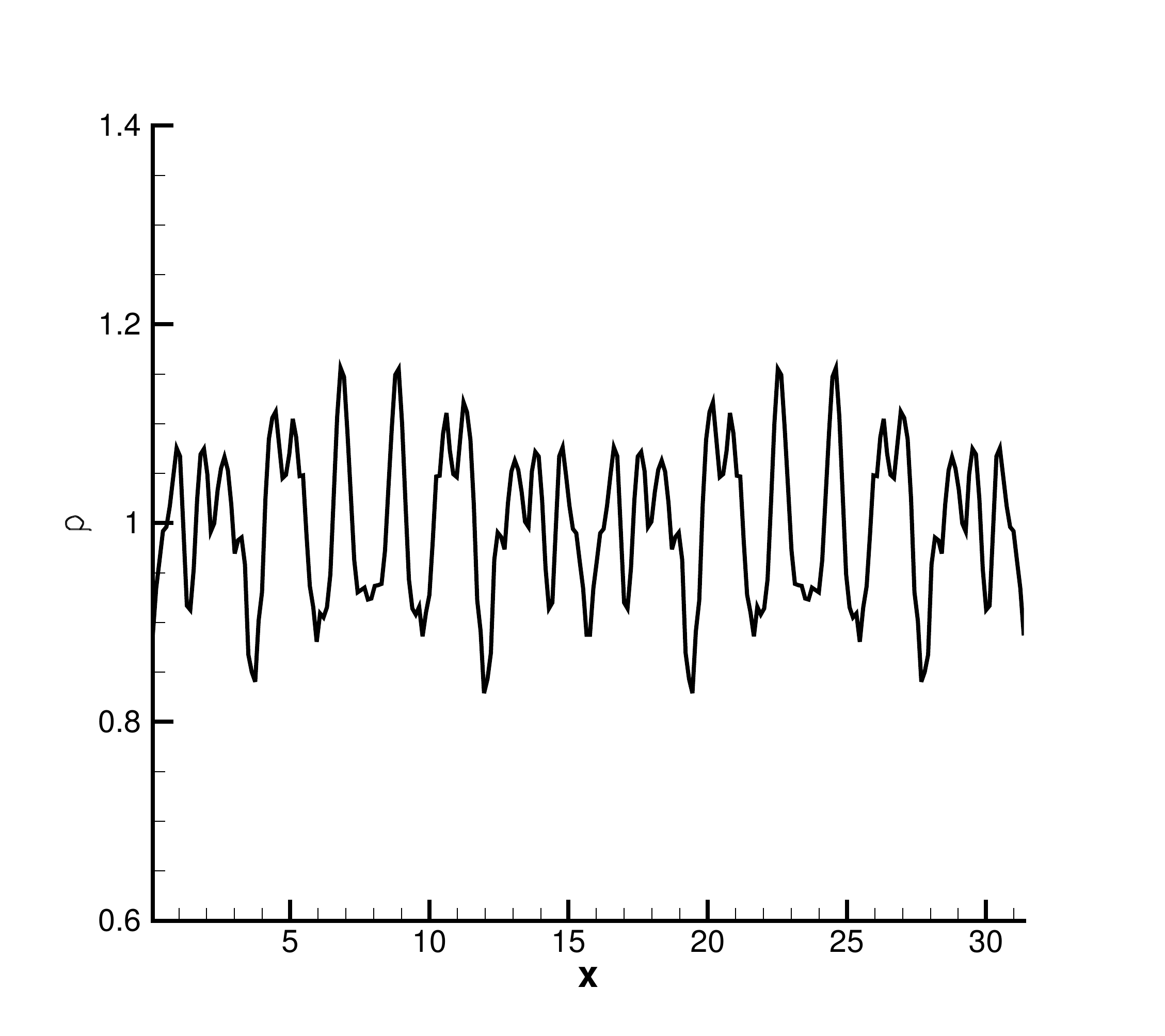}}
		\subfigure[ Adaptive, $t=100.$]{\includegraphics[width=2.8in,angle=0]{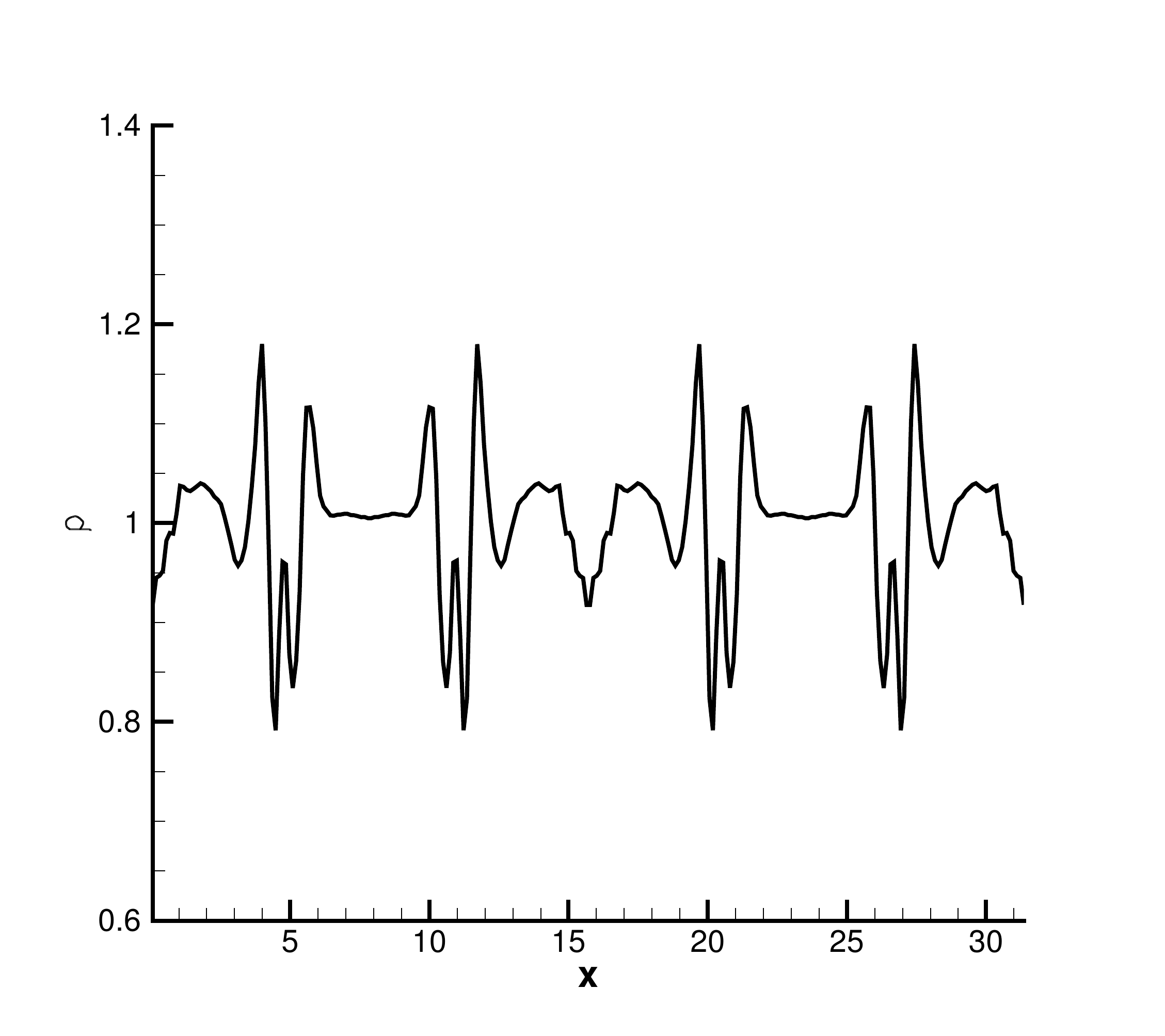}}
	\end{center}
	\caption{SW instability with parameter choice 1. Plots of the computed density function $\rho_h$ at selected  time $t$ by upwind flux for the Maxwell's equations.  Sparse grid: $N=8$, $k=3$.  Adaptive sparse grid: $N=6$, $k=3$, $\eps=2\times10^{-7}$.}
	\label{density1}
\end{figure}

In Figure \ref{contour1}, we further present the contour plots of distribution function $f$ at $x_2=0.05 \pi$ at several times which are consistent with those plotted for the density in Figure \ref{density1}. At later times,  
considerable fine structure generates in agreement with the Log Fourier plots.
We can also observe that the figures obtained by both schemes are in agreement with those in \cite{cheng2014discontinuous} at $t=55$ and $t=82$. At later times, e.g. $t=100,$ our scheme produces more filamented structures than those in \cite{cheng2014discontinuous}. This is due to the poorer resolution in our approximation spaces. At such time period,  the adaptive sparse grid scheme clearly offers better resolution.
In Figure \ref{density1}, we plot the density modulation, containing the expected spikes for the SW instability at selected  time $t$. At $t=55$ and $t=82$, the results agree overall with the calculations in \cite{cheng2014discontinuous}.  At $t=100$, the density shows visible   fluctuations comparing with those in \cite{cheng2014discontinuous}.

\begin{figure}[htb]
	\begin{center}
		
		\subfigure[ $t=0.$ Active elements: 0.73\%]{\includegraphics[width=2.9in,angle=0]{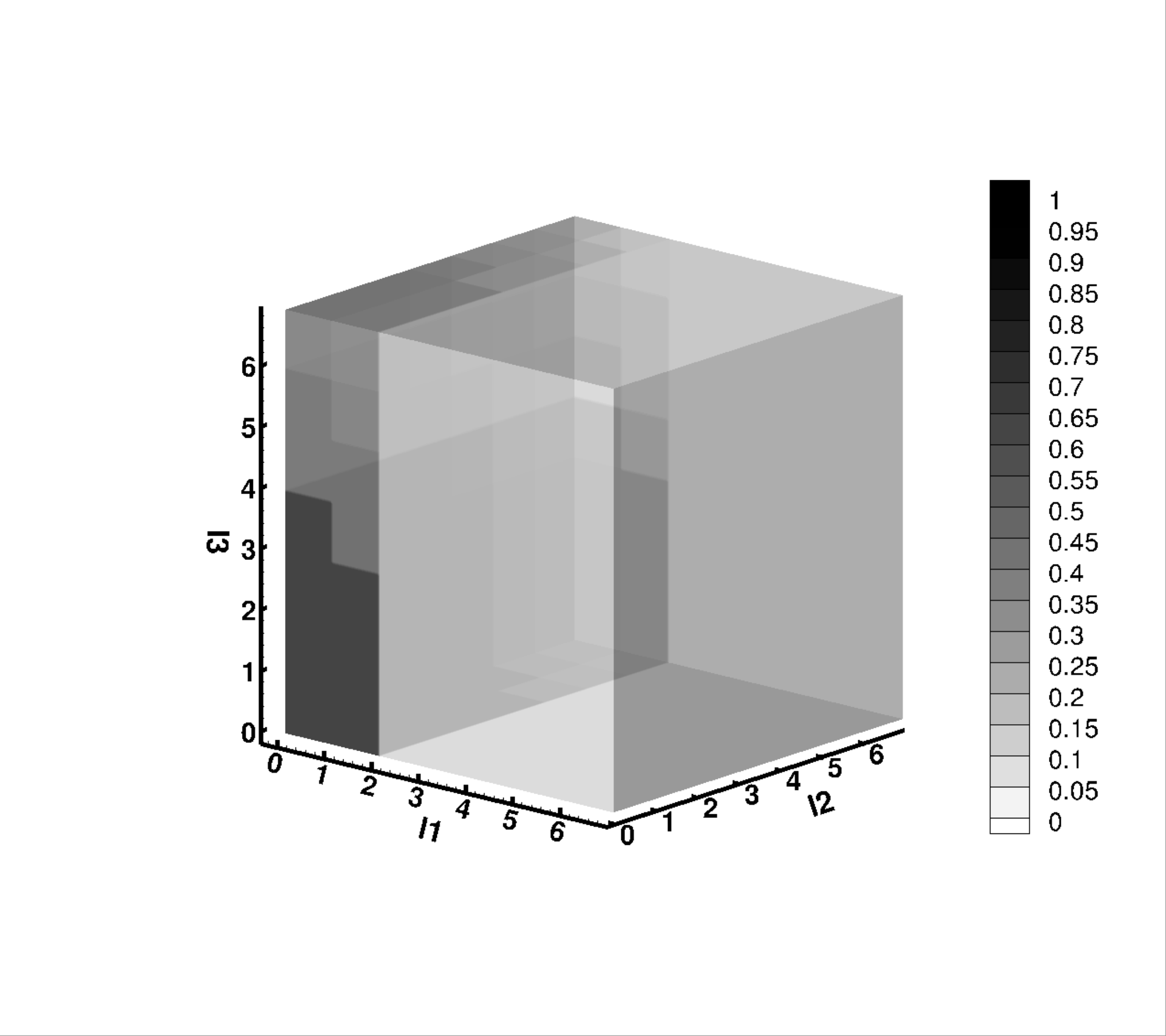}}	
		\subfigure[ $t=55.$  Active elements: 4.41\%]{\includegraphics[width=2.9in,angle=0]{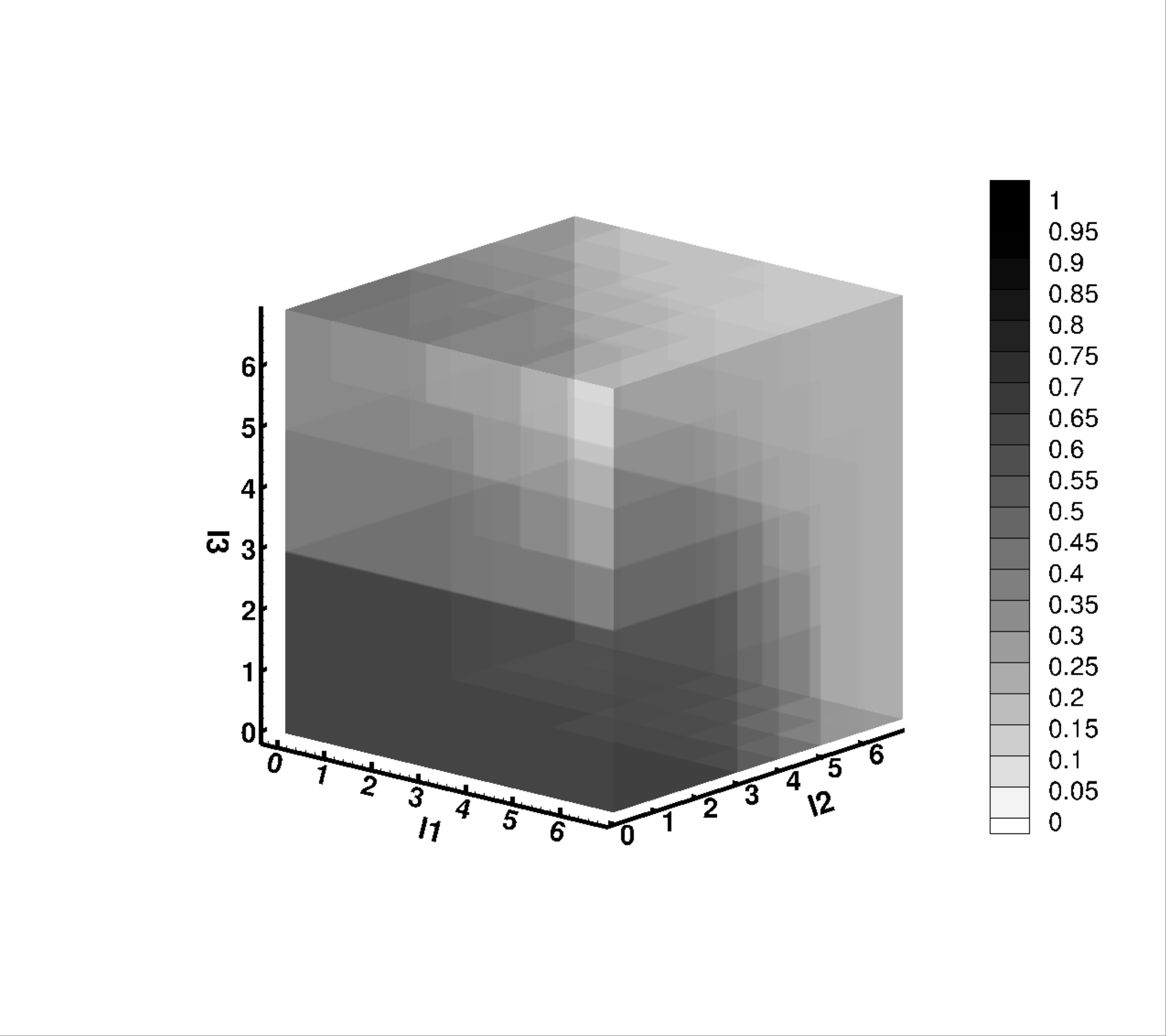}}
		\subfigure[ $t=82.$  Active elements: 26.61\%]{\includegraphics[width=2.9in,angle=0]{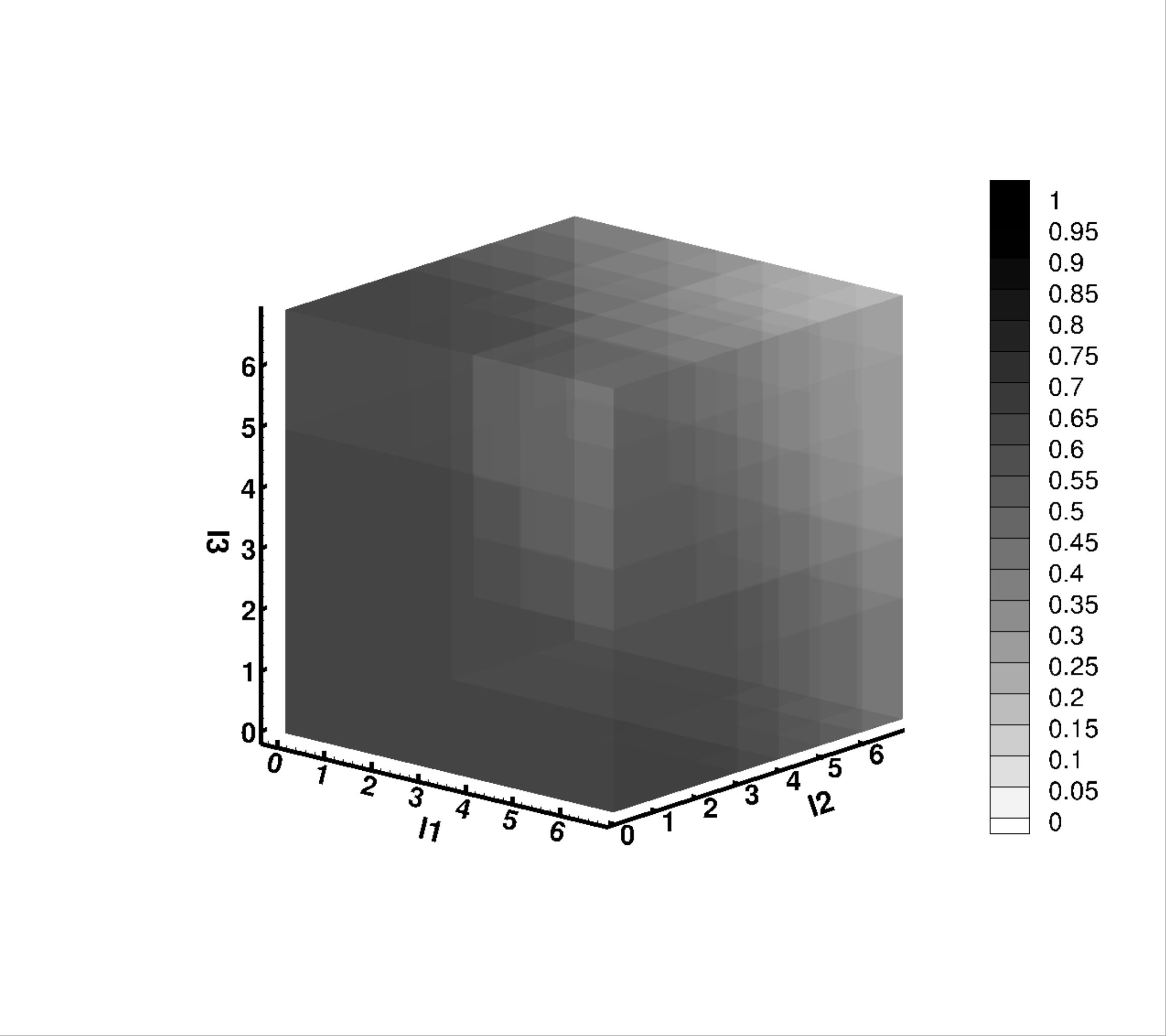}}
		\subfigure[ $t=100.$  Active elements: 52.41\%]{\includegraphics[width=2.9in,angle=0]{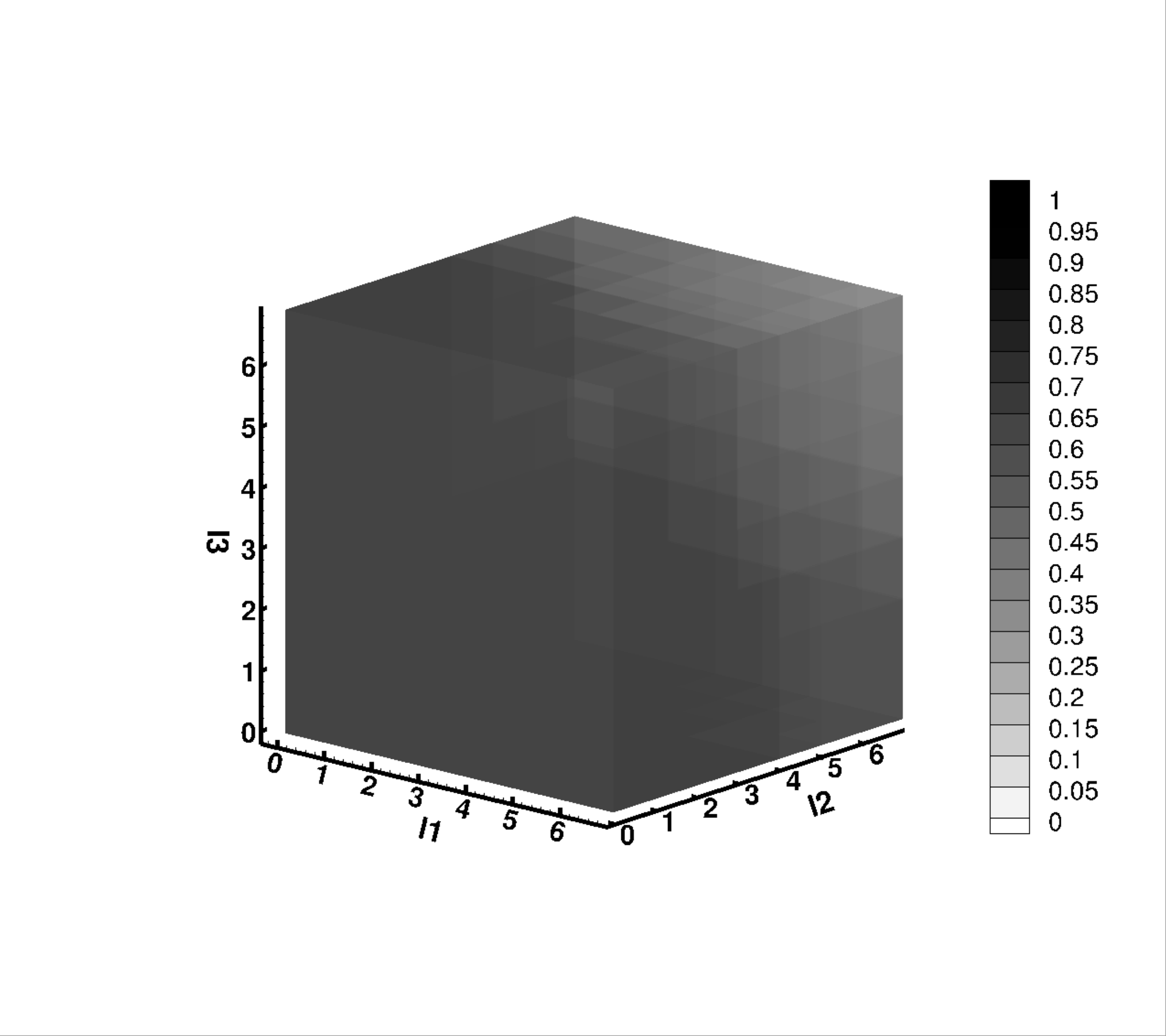}}
		
	\end{center}
	\caption{SW instability with parameter choice 1 and upwind flux for the Maxwell's equations. The percentage of active elements for each incremental space $\bW_\bl$, $\bl=(l_1,l_2,l_3)$ and $|\bl|_\infty\leq N$ in the Adaptive sparse grid. $N=6$, $k=3$,  $\eps=2\times10^{-7}$.}
	\label{percent1_ada}
\end{figure}

To verify the effectiveness of the adaptive scheme, in  Figure \ref{percent1_ada}, we   show the percentage of active elements for each incremental space $\bW_\bl$, $\bl=(l_1,l_2,l_3)$ in the adaptive scheme. If the percentage is 1, it means all elements on that level are active. If the percentage is 0, it means all elements on that level are deactivated. A full grid approximation corresponds to percentage being 1 on all levels, while the sparse grid approximation corresponds to percentage being 1 when $|\bl|_1 \le N,$ and $0$ otherwise. We observe that the number of active elements in the adaptive approximation space is increasing   with time. This is consistent with the fact that more fine structures are generated at late times due to filamentation. Therefore   more elements are needed to capture the fine structures than  previous times. That also explains why the adaptive scheme produces better numerical results than the standard sparse grid scheme. In particular, we count the total numbers of active elements at $t=0, 55, 82, 100$, which are $1924, 11556, 69750, 137406,$ respectively. The full grid space, on the other hand, has a total of $262144$ elements.  The   percentage of active elements, therefore grows from $0.73\%$ to $52.41\%$ from $t=0$ to $100.$  If $t$ becomes larger, the percentage will eventually reach $100\%,$ which amounts to a full grid calculation. This is intrinsically due to the properties of the VM solutions. After a certain time, it is known that any approximation with fixed mesh size (or grid resolution) will fail \cite{chengknorr_76, cheng_vp}. By using the adaptive scheme, the main advantage lies in the efficient calculation of the solution when $t$ is not so large. 

\begin{figure}[htb]
	\begin{center}
		\subfigure[ Sparse grid, kinetic energies]{\includegraphics[width=3in,angle=0]{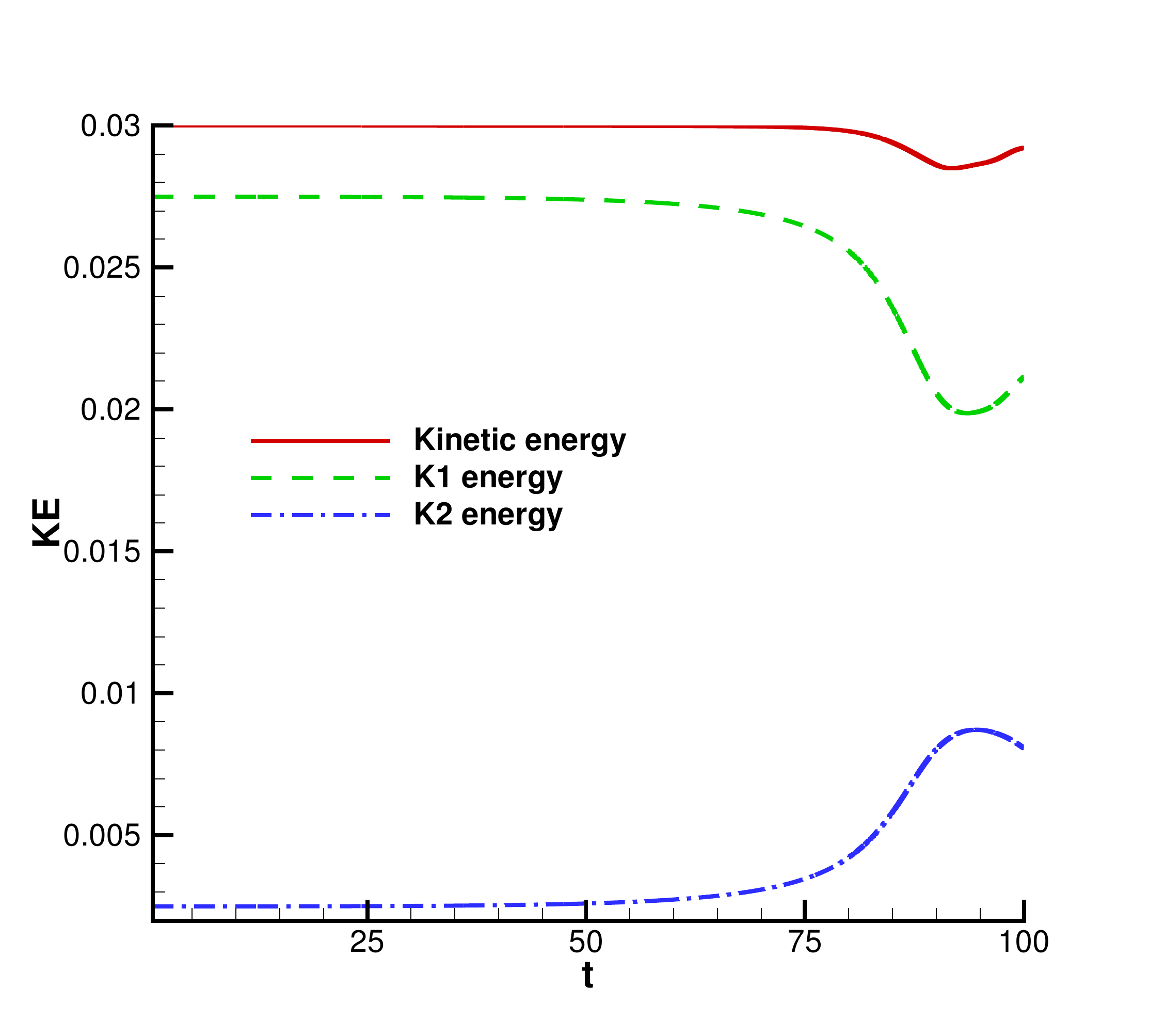}}
		\subfigure[ Adaptive, kinetic energies]{\includegraphics[width=3in,angle=0]{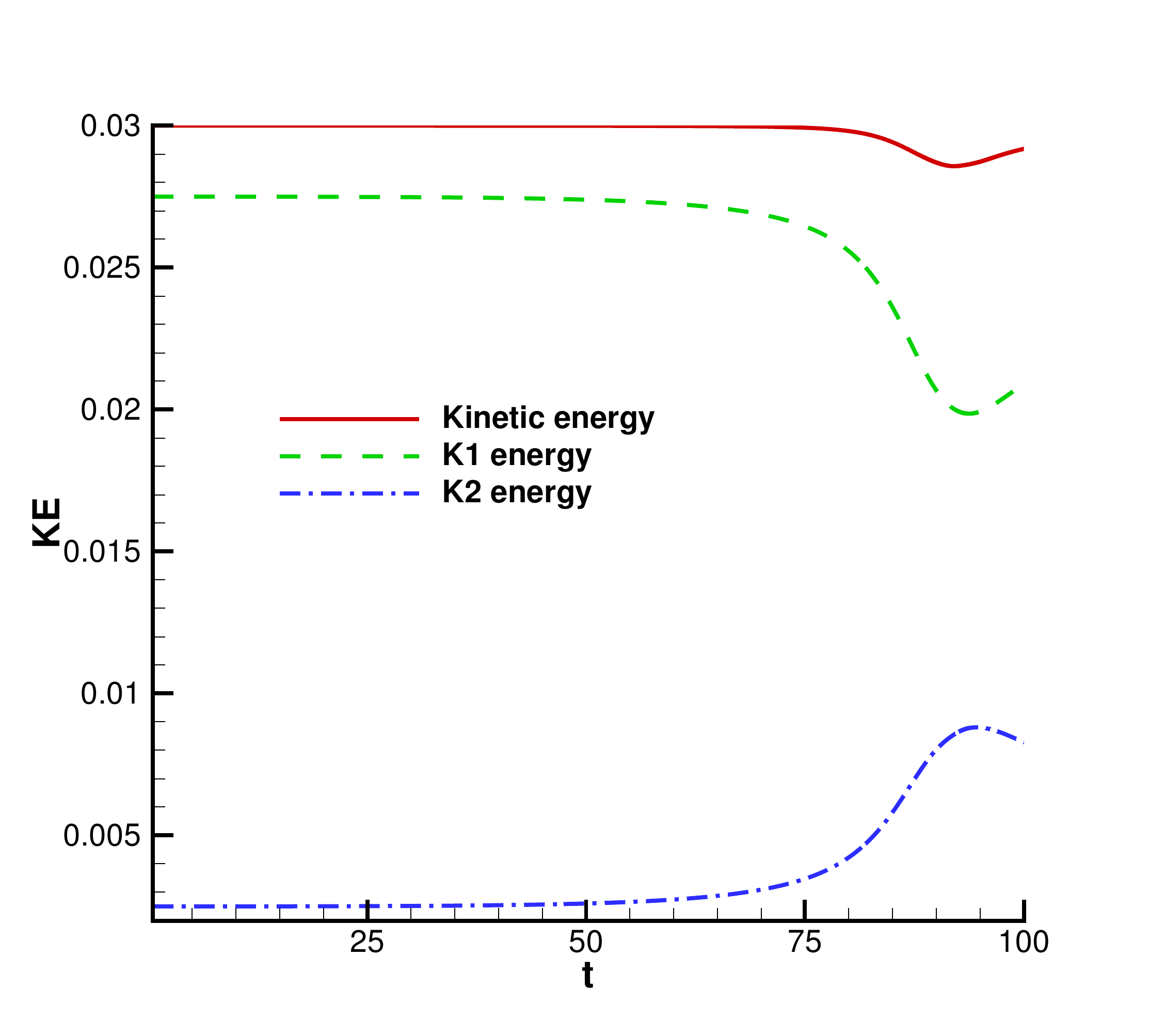}}
		\subfigure[Sparse grid, field energies]{\includegraphics[width=3in,angle=0]{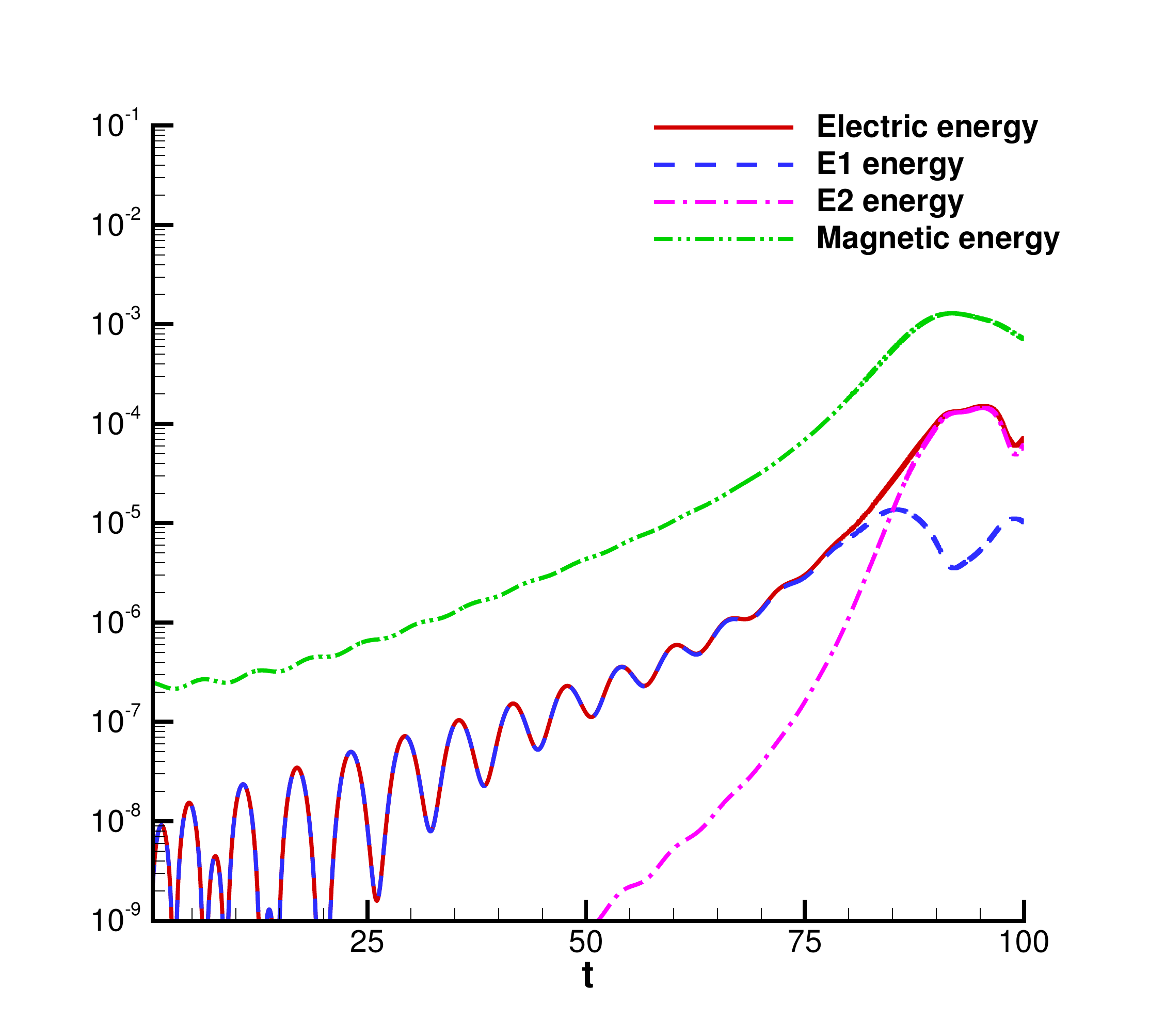}}
		\subfigure[Adaptive, field energies]{\includegraphics[width=3in,angle=0]{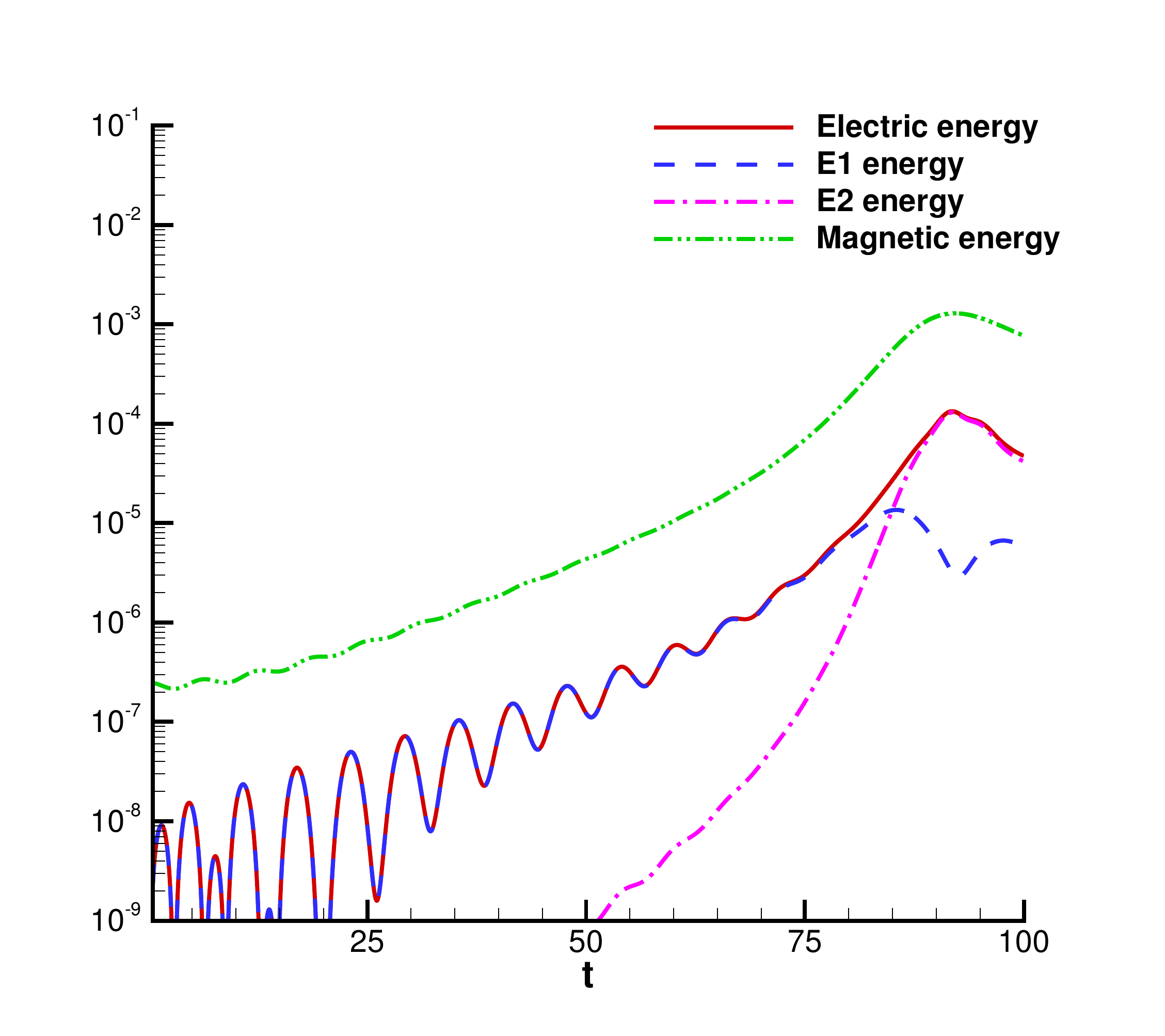}}
	\end{center}
	
	\caption{SW instability with parameter choice 2. Time evolution of kinetic, electric  and magnetic energies by upwind flux for the Maxwell's equations.  Sparse grid: $N=8$, $k=3$.  Adaptive sparse grid: $N=6$, $k=3$, $\eps=10^{-6}$.}
	\label{keeme2}
\end{figure}

\begin{figure}[htb]
	\begin{center}
		\subfigure[Sparse grid, Log Fourier modes of $E_1$]{\includegraphics[width=2.8in,angle=0]{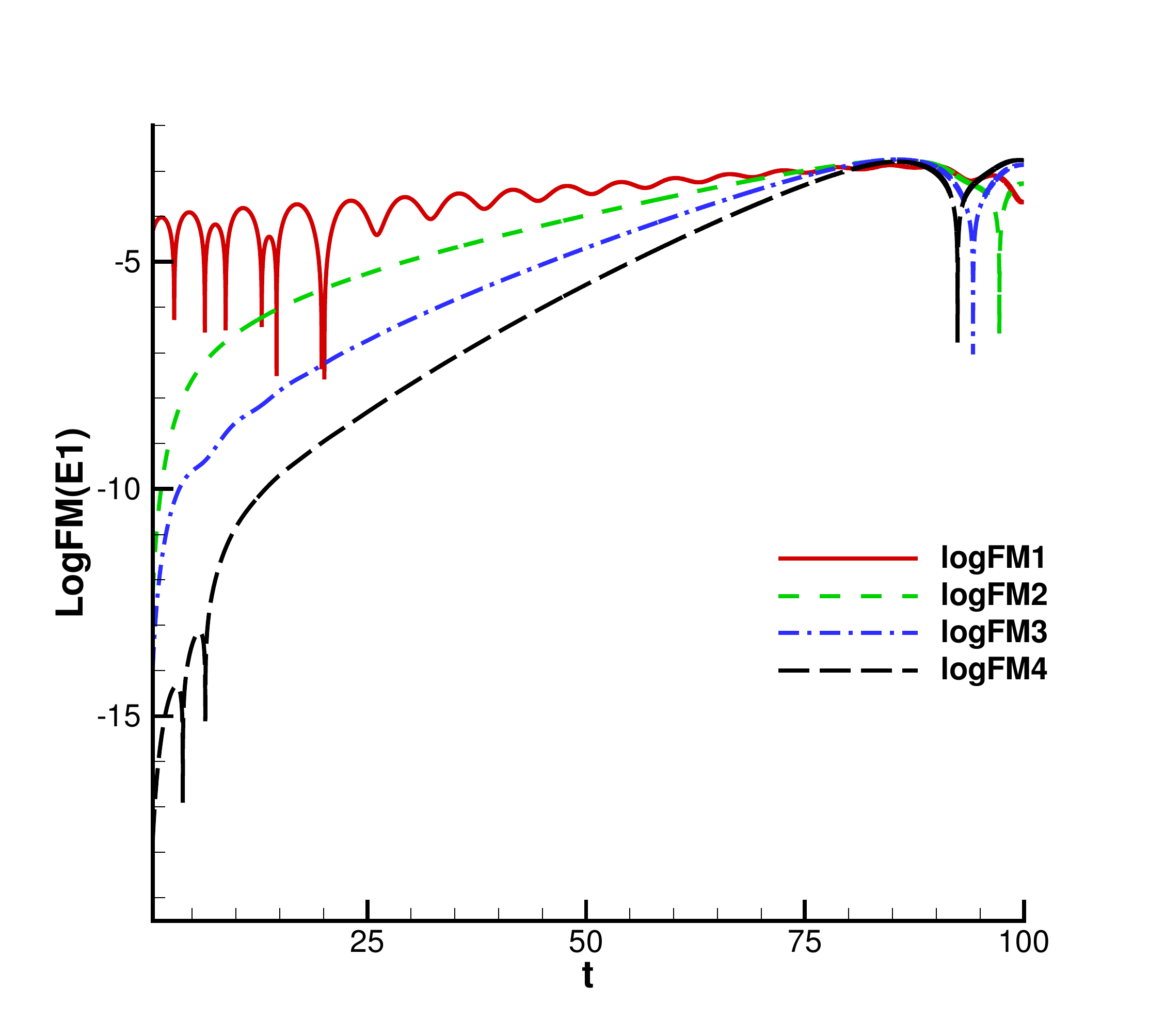}}
		\subfigure[Adaptive, Log Fourier modes of $E_1$]{\includegraphics[width=2.8in,angle=0]{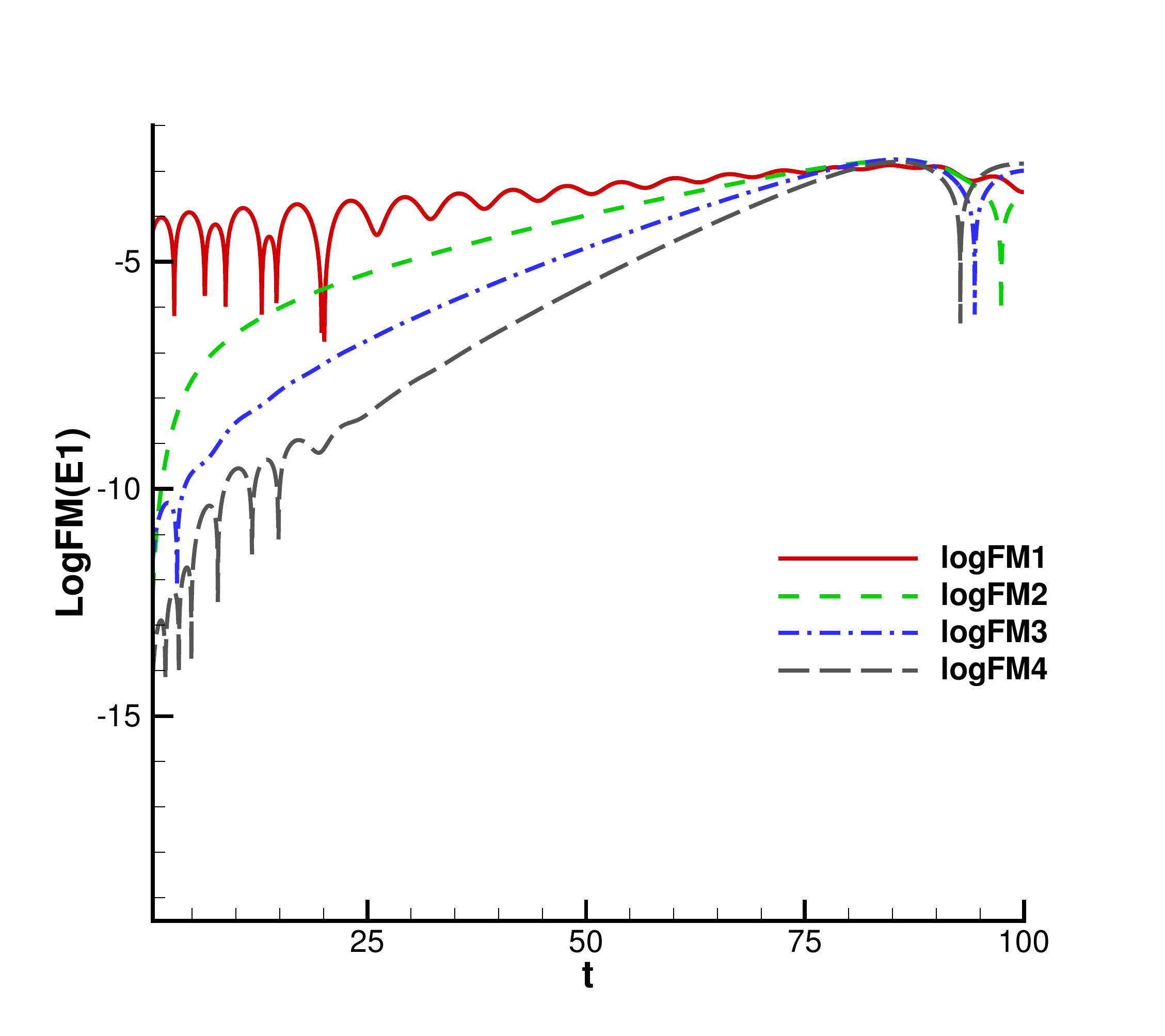}}
		\subfigure[Sparse grid, Log Fourier modes of $E_2$]{\includegraphics[width=2.8in,angle=0]{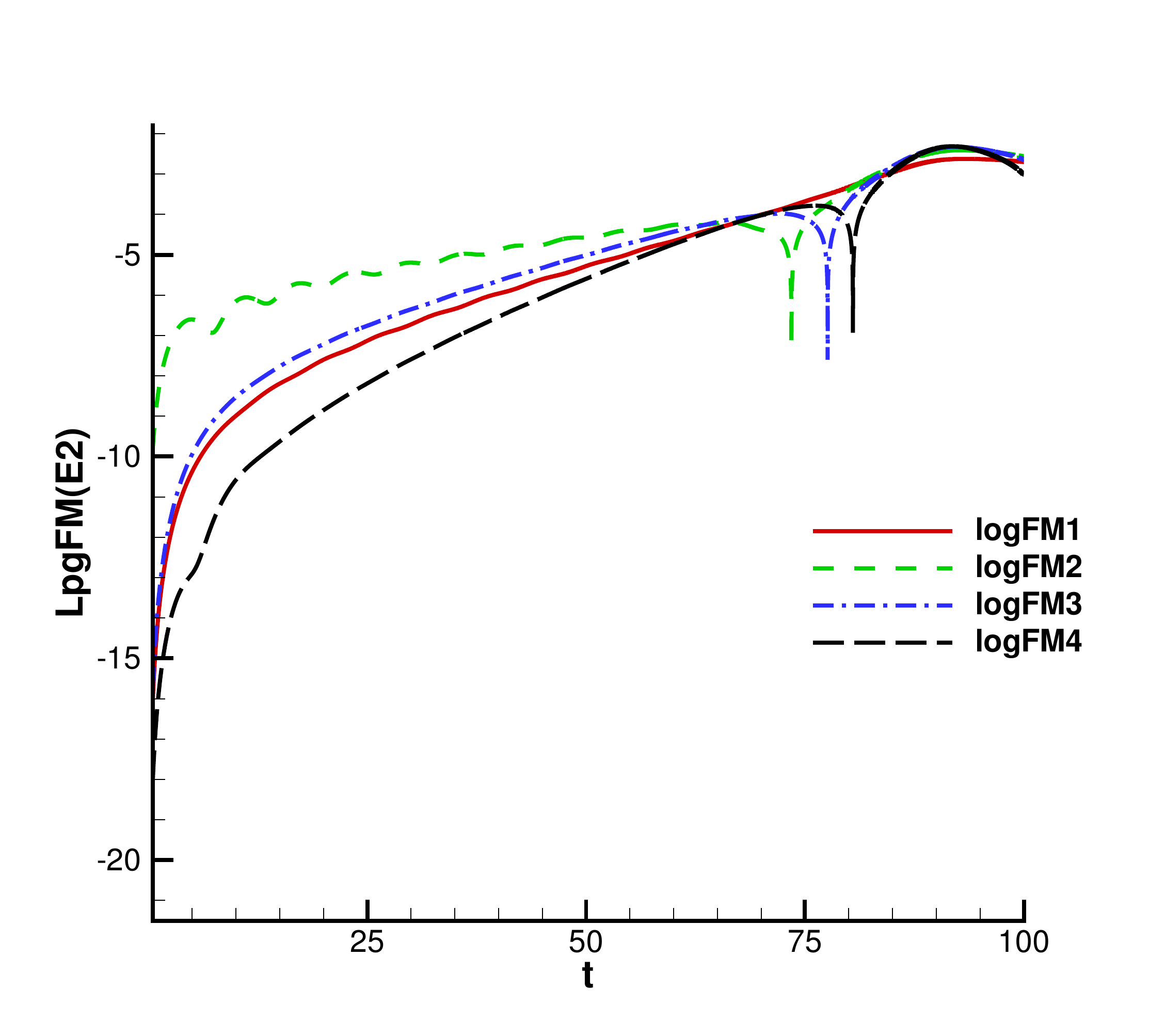}}
		\subfigure[Adaptive, Log Fourier modes of $E_2$]{\includegraphics[width=2.8in,angle=0]{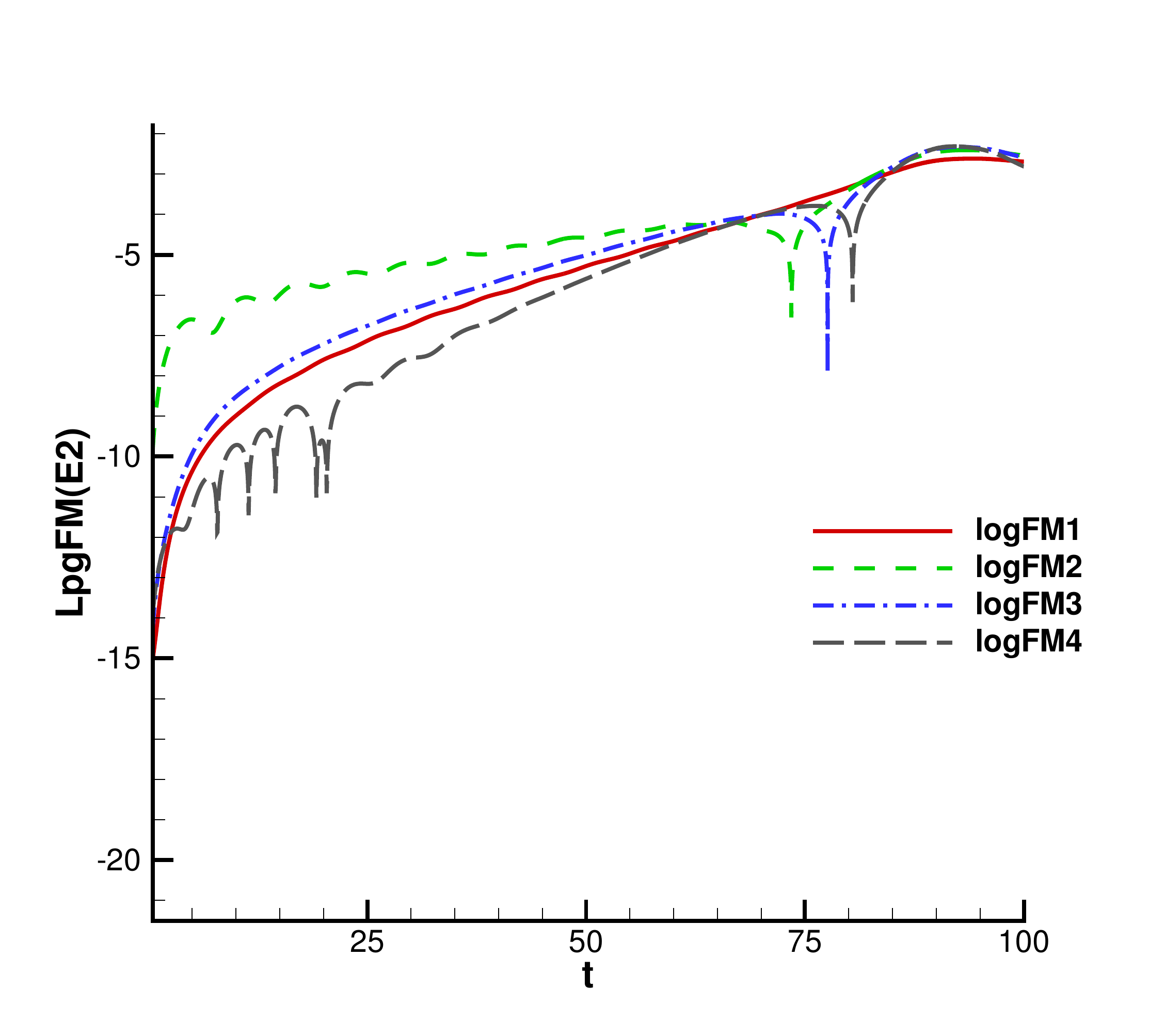}}
		\subfigure[Sparse grid, Log Fourier modes of $B_3$]{\includegraphics[width=2.8in,angle=0]{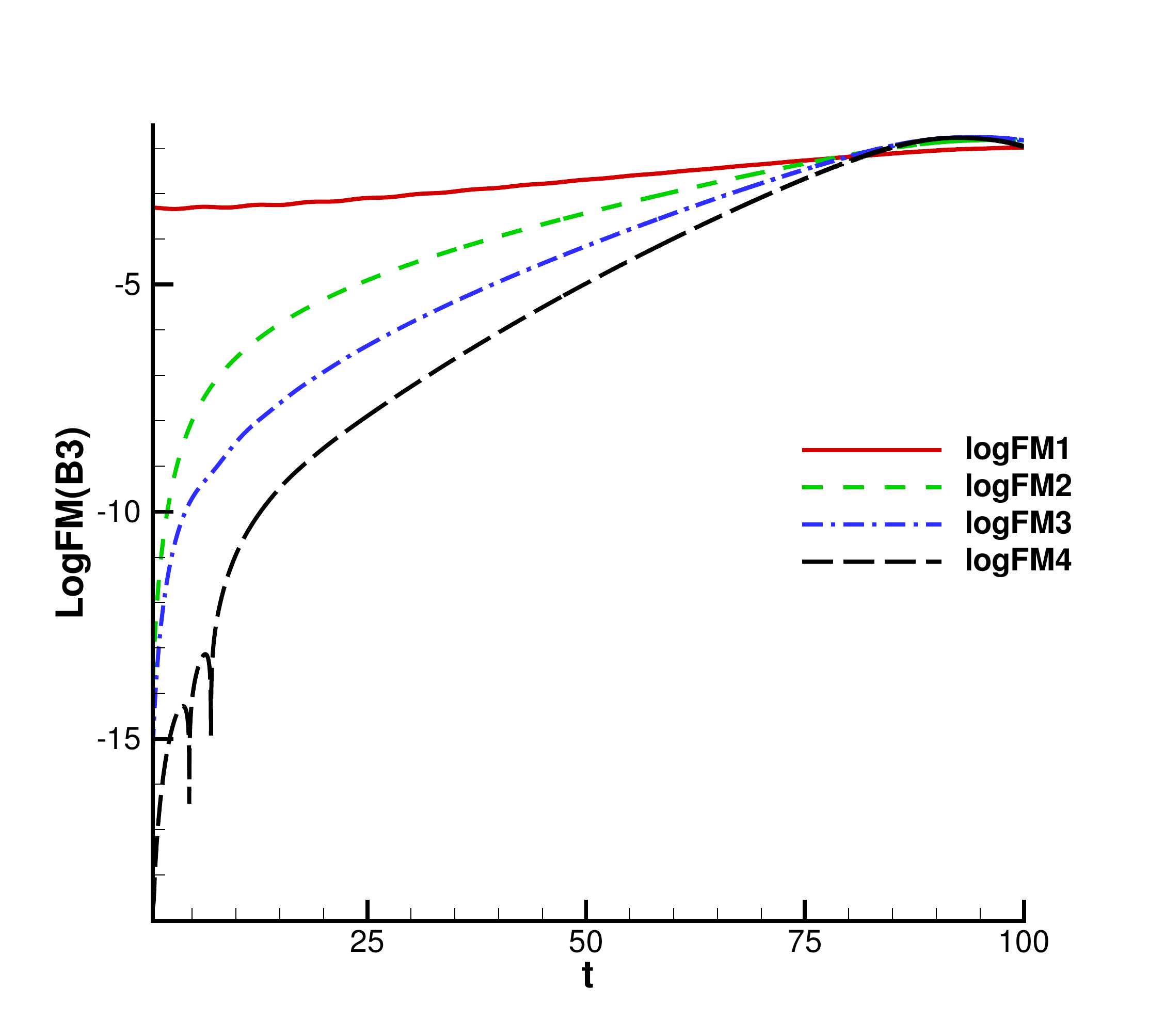}}
		\subfigure[Adaptive, Log Fourier modes of $B_3$]{\includegraphics[width=2.8in,angle=0]{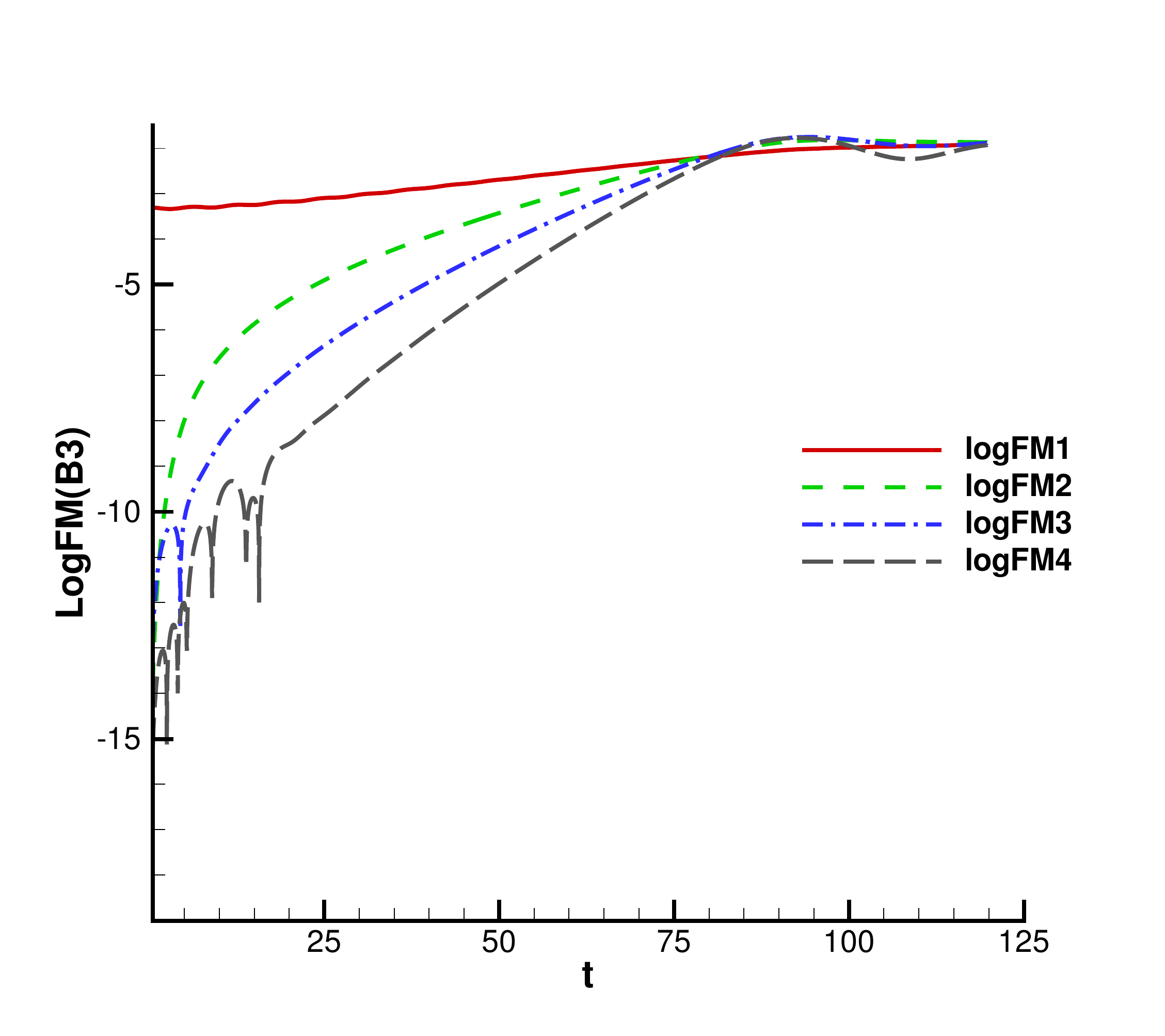}}
	\end{center}
	\caption{SW instability with parameter choice 2. The first four Log Fourier modes of $E_1$, $E_2$, $B_3$ by upwind flux for Maxwell's equations. Sparse grid: $N=8$, $k=3$.  Adaptive sparse grid: $N=6$, $k=3$, $\eps=10^{-6}$.}
	\label{logfm2}
\end{figure}

\begin{figure}[htb]
	\begin{center}
		\subfigure[ Sparse grid, $x_2=0.05 \pi,  \, t=55.$]{\includegraphics[width=2.8in,angle=0]{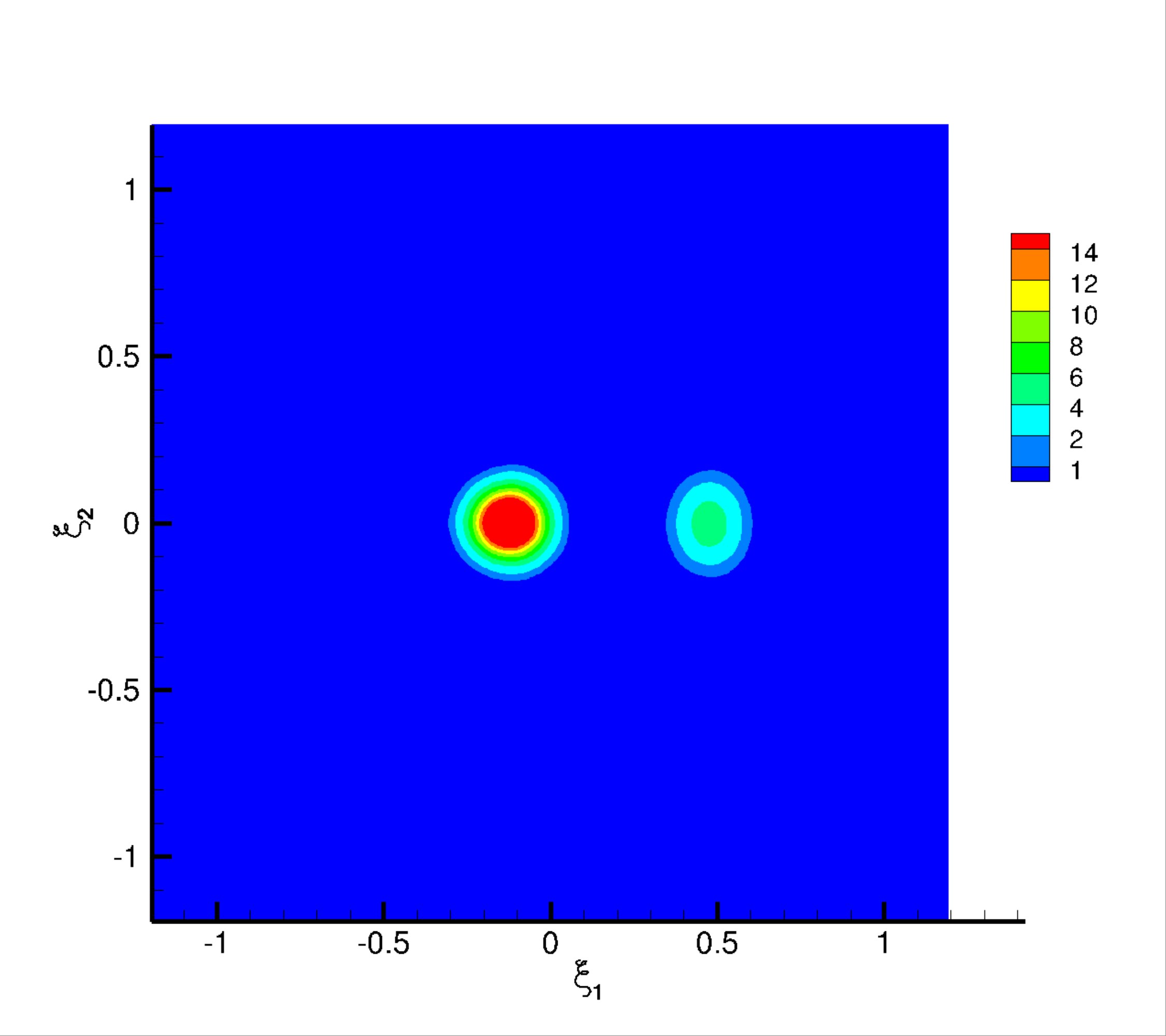}}
		\subfigure[ Adaptive, $x_2=0.05 \pi,  \, t=55.$]{\includegraphics[width=2.8in,angle=0]{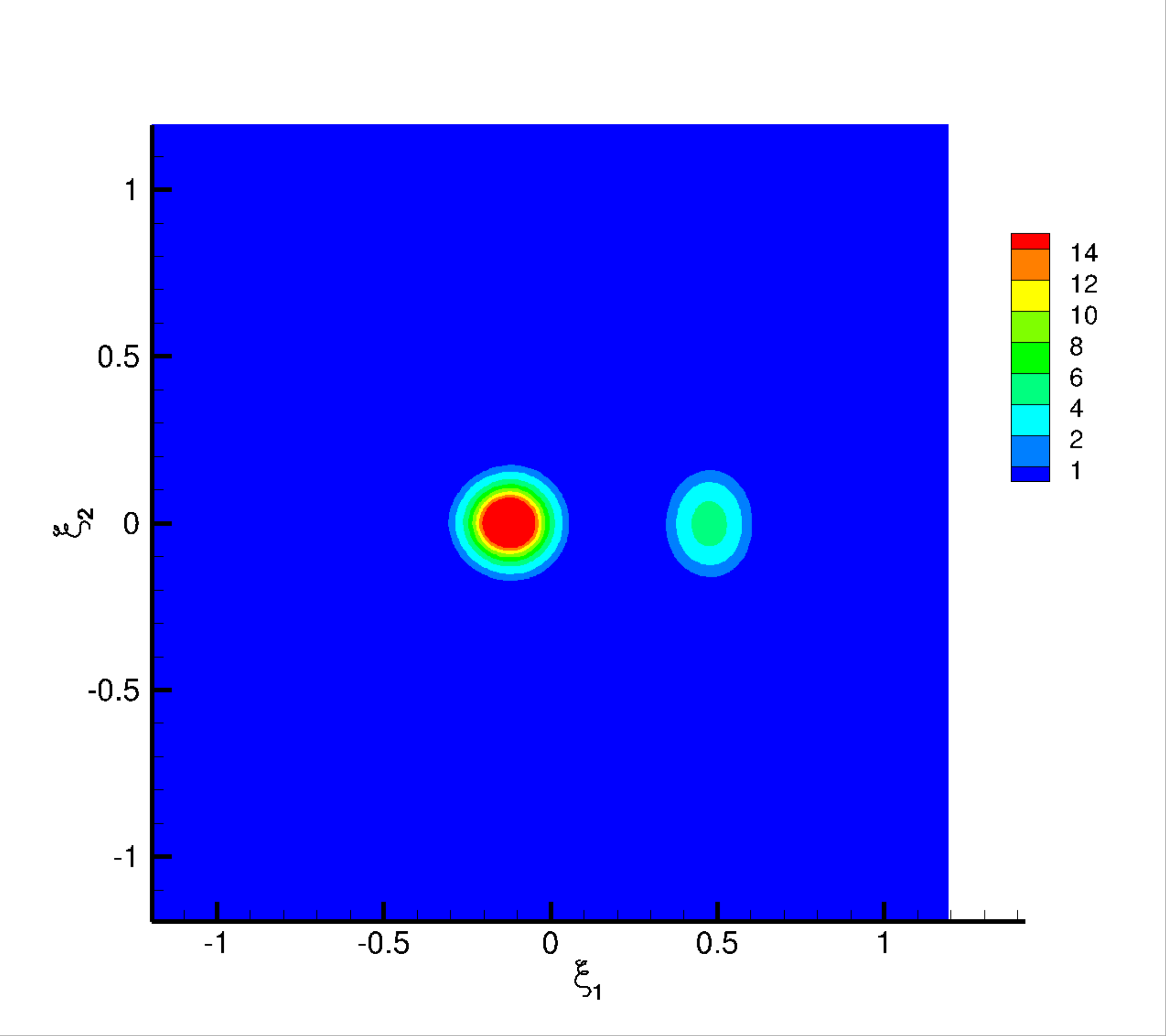}}
		\subfigure[Sparse grid,  $x_2=0.05 \pi, \,  t=82.$]{\includegraphics[width=2.8in,angle=0]{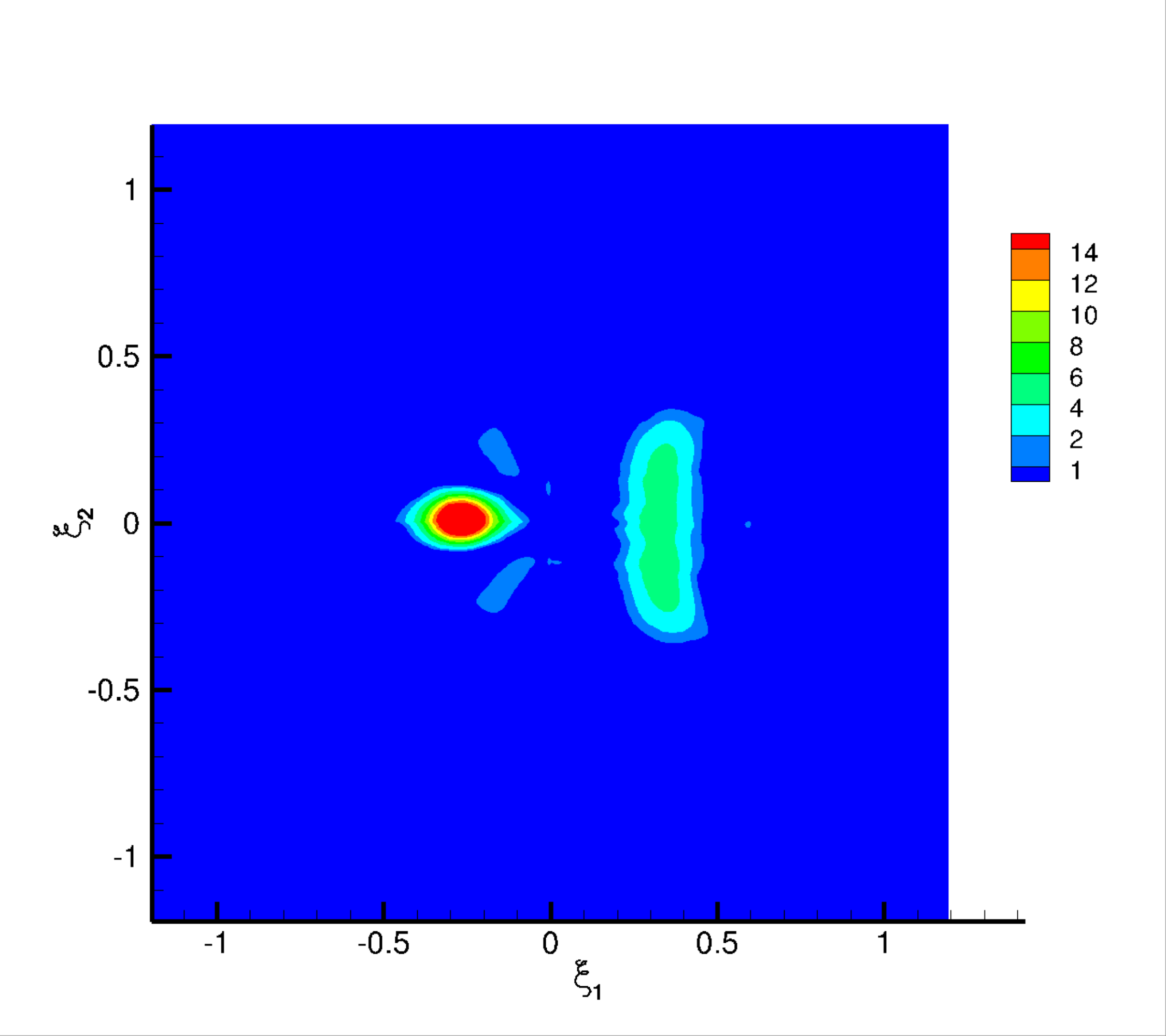}}
		\subfigure[Adaptive,  $x_2=0.05 \pi, \,  t=82.$]{\includegraphics[width=2.8in,angle=0]{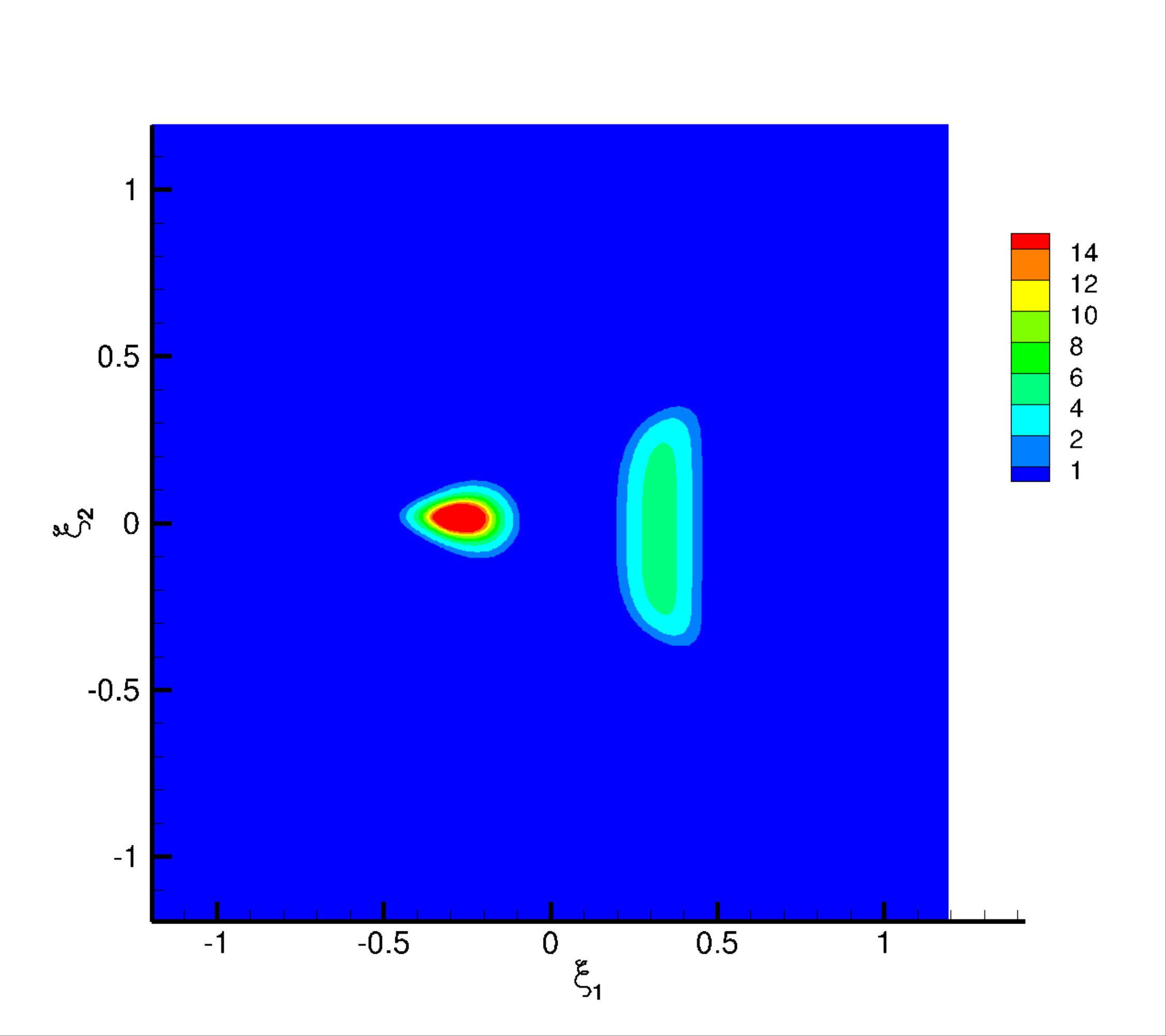}}
		\subfigure[ Sparse grid, $x_2=0.05 \pi,  \, t=100.$]{\includegraphics[width=2.8in,angle=0]{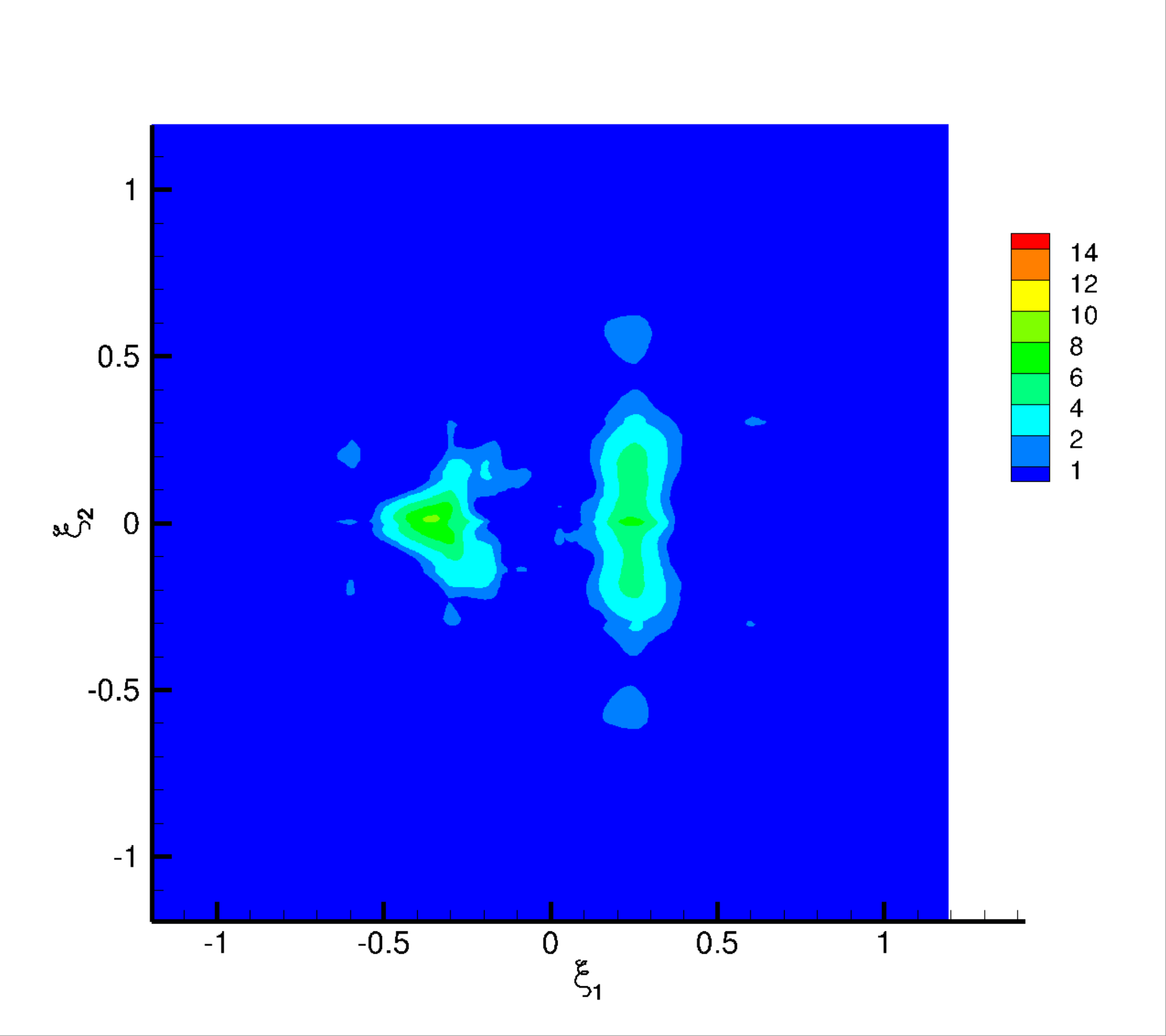}}
		\subfigure[Adaptive, $x_2=0.05 \pi,  \, t=100.$]{\includegraphics[width=2.8in,angle=0]{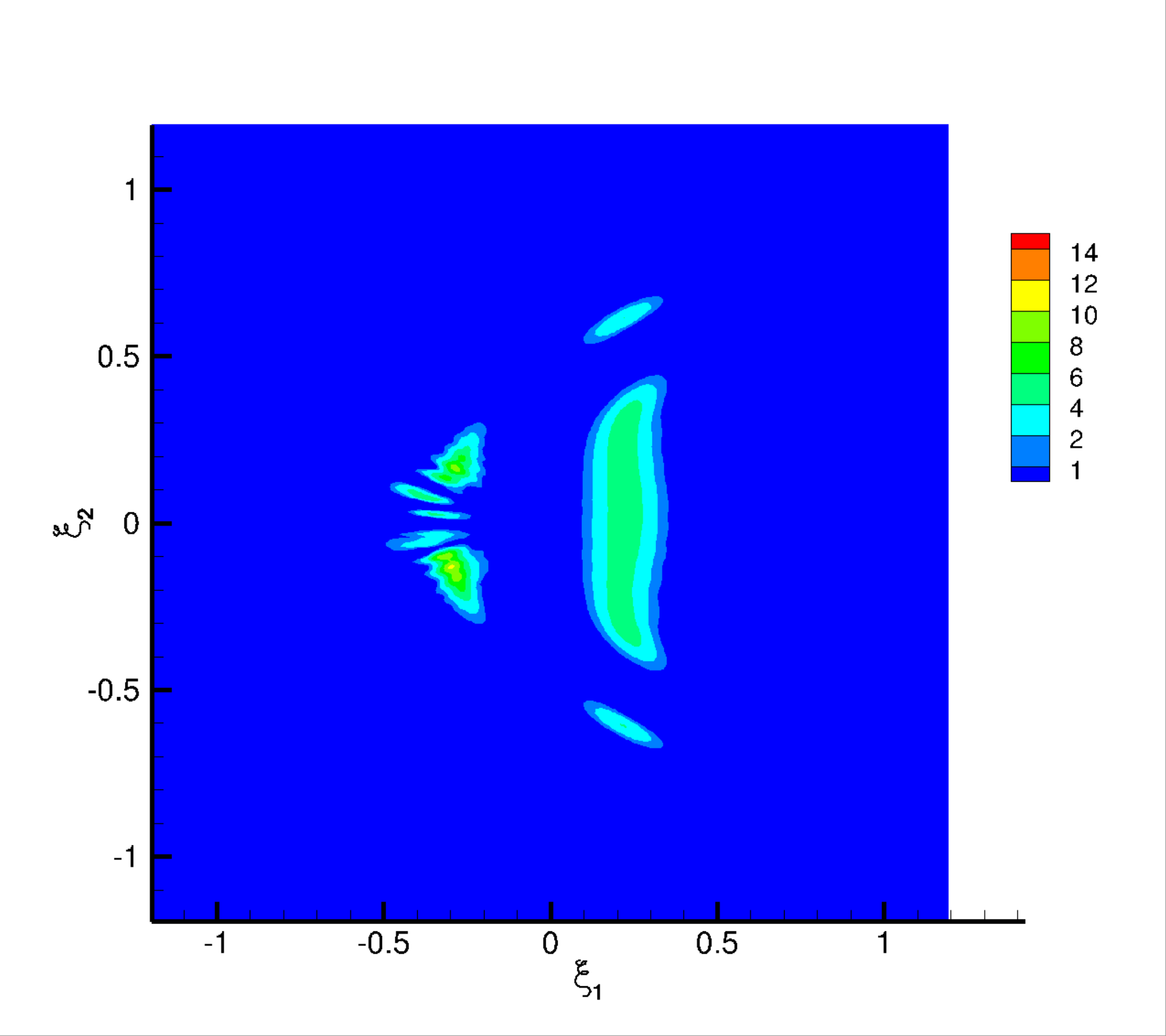}}
	\end{center}
	\caption{SW instability with parameter choice 2. 2D contour plots of the computed distribution function $f_h$ by upwind flux for Maxwell's equations.  Sparse grid: $N=8$, $k=3$.  Adaptive sparse grid: $N=6$, $k=3$, $\eps=10^{-6}$.}
	\label{contour2}
\end{figure}

\begin{figure}[htb]
	\begin{center}
		\subfigure[ Sparse grid, $ t=55.$]{\includegraphics[width=2.8in,angle=0]{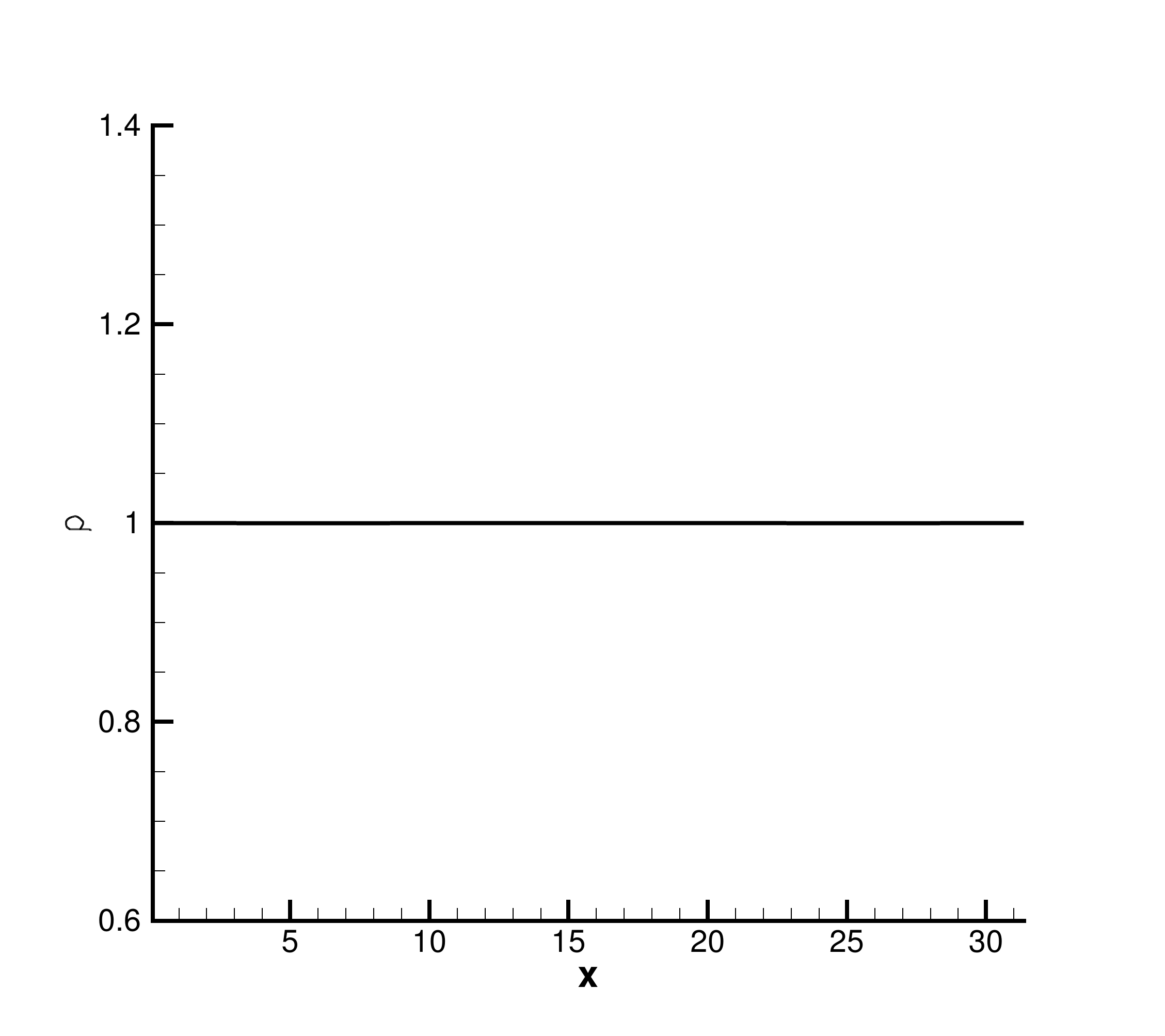}}
		\subfigure[ Adaptive, $ t=55.$]{\includegraphics[width=2.8in,angle=0]{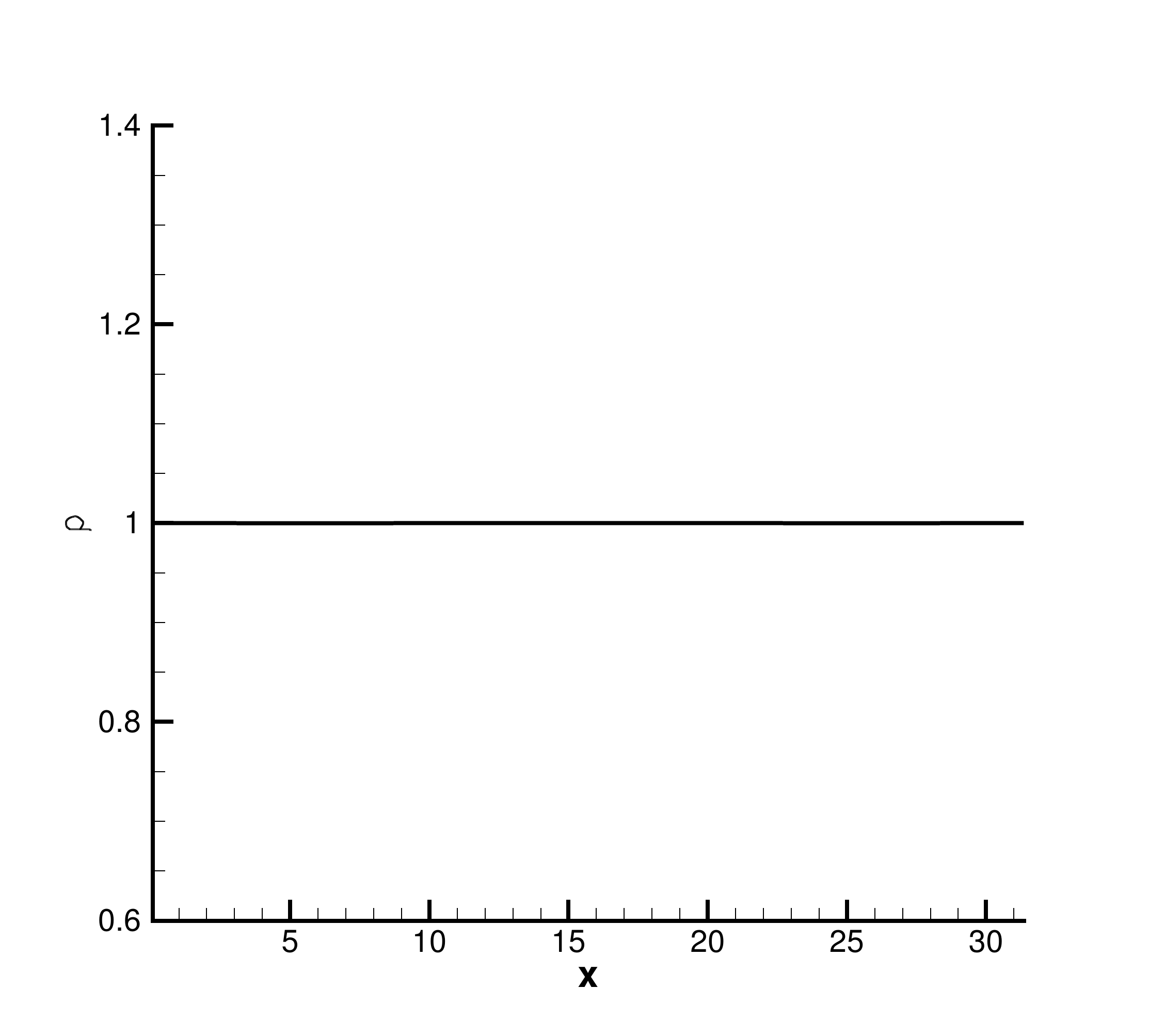}}
		\subfigure[ Sparse grid, $ t=82.$]{\includegraphics[width=2.8in,angle=0]{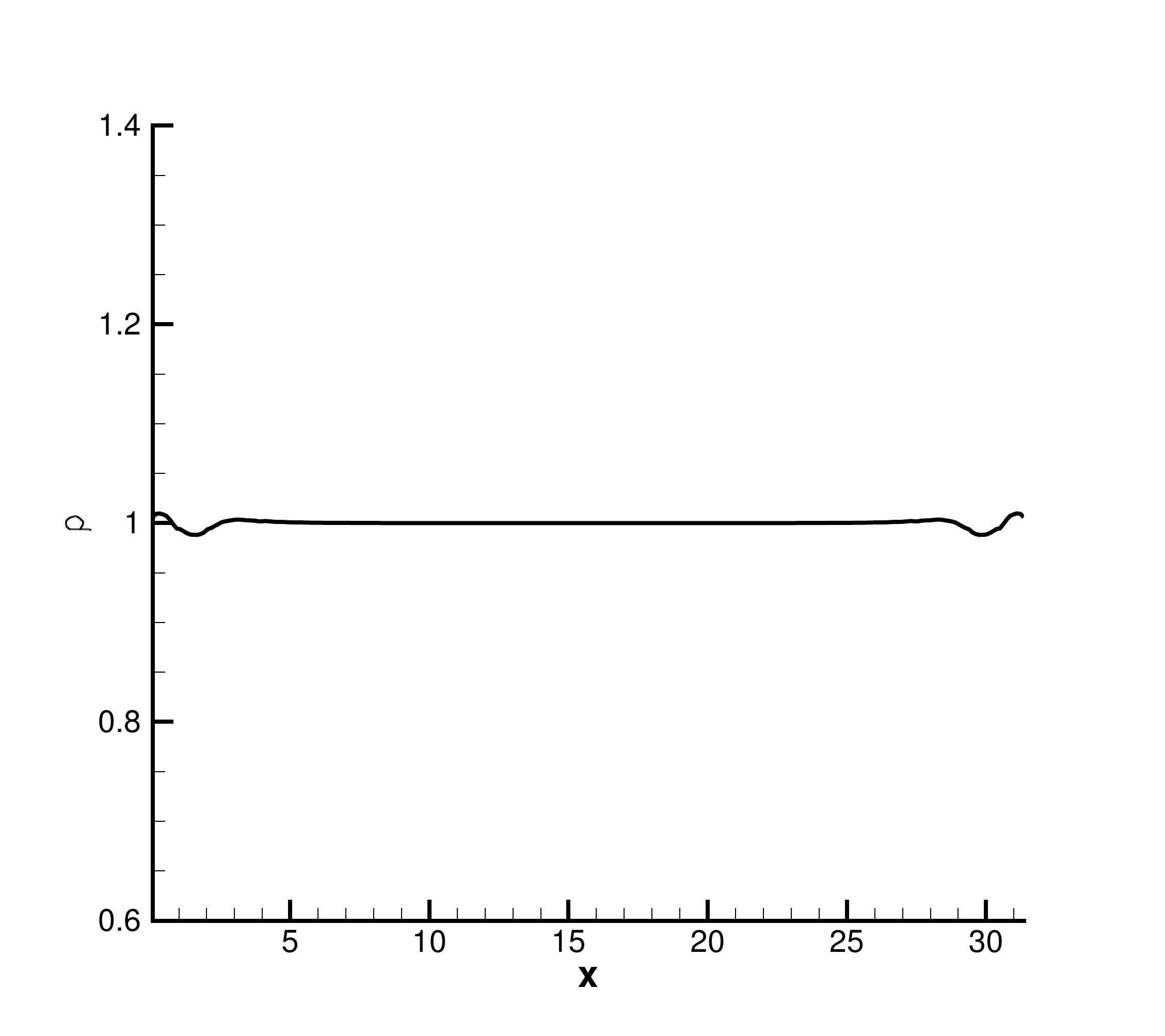}}
		\subfigure[ Adaptive, $ t=82.$]{\includegraphics[width=2.8in,angle=0]{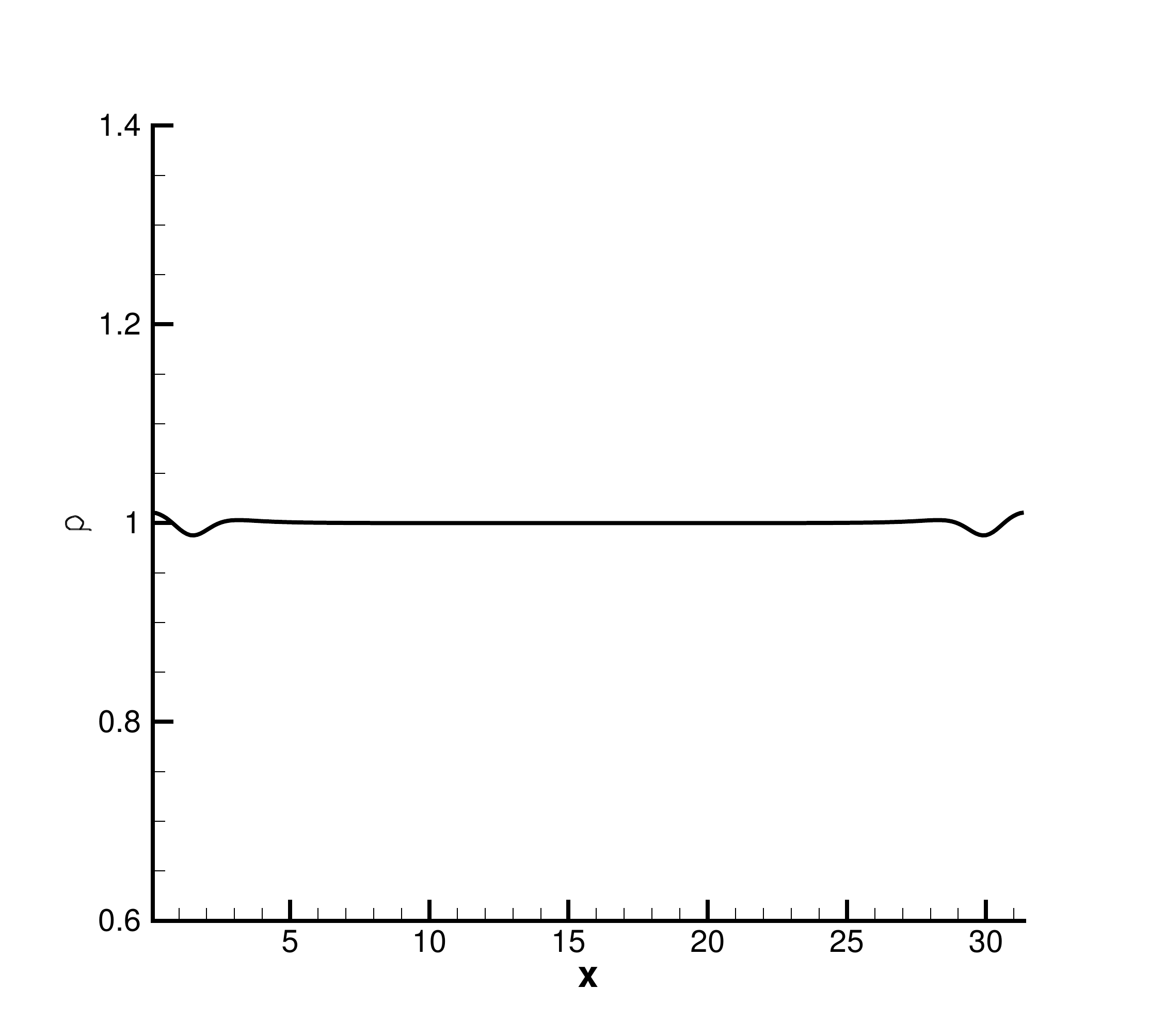}}
		\subfigure[ Sparse grid, $t=100.$]{\includegraphics[width=2.8in,angle=0]{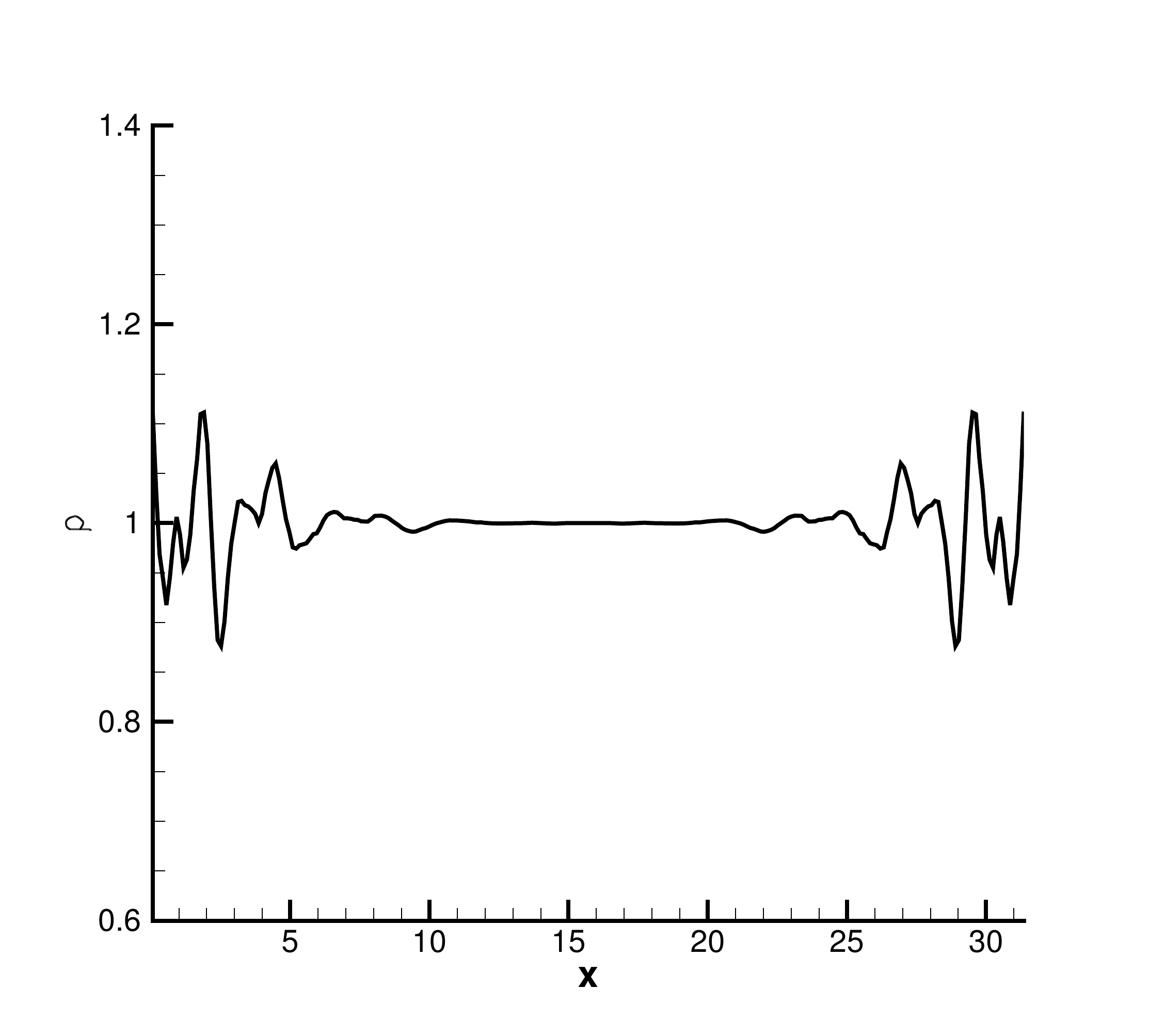}}
		\subfigure[ Adaptive, $t=100.$]{\includegraphics[width=2.8in,angle=0]{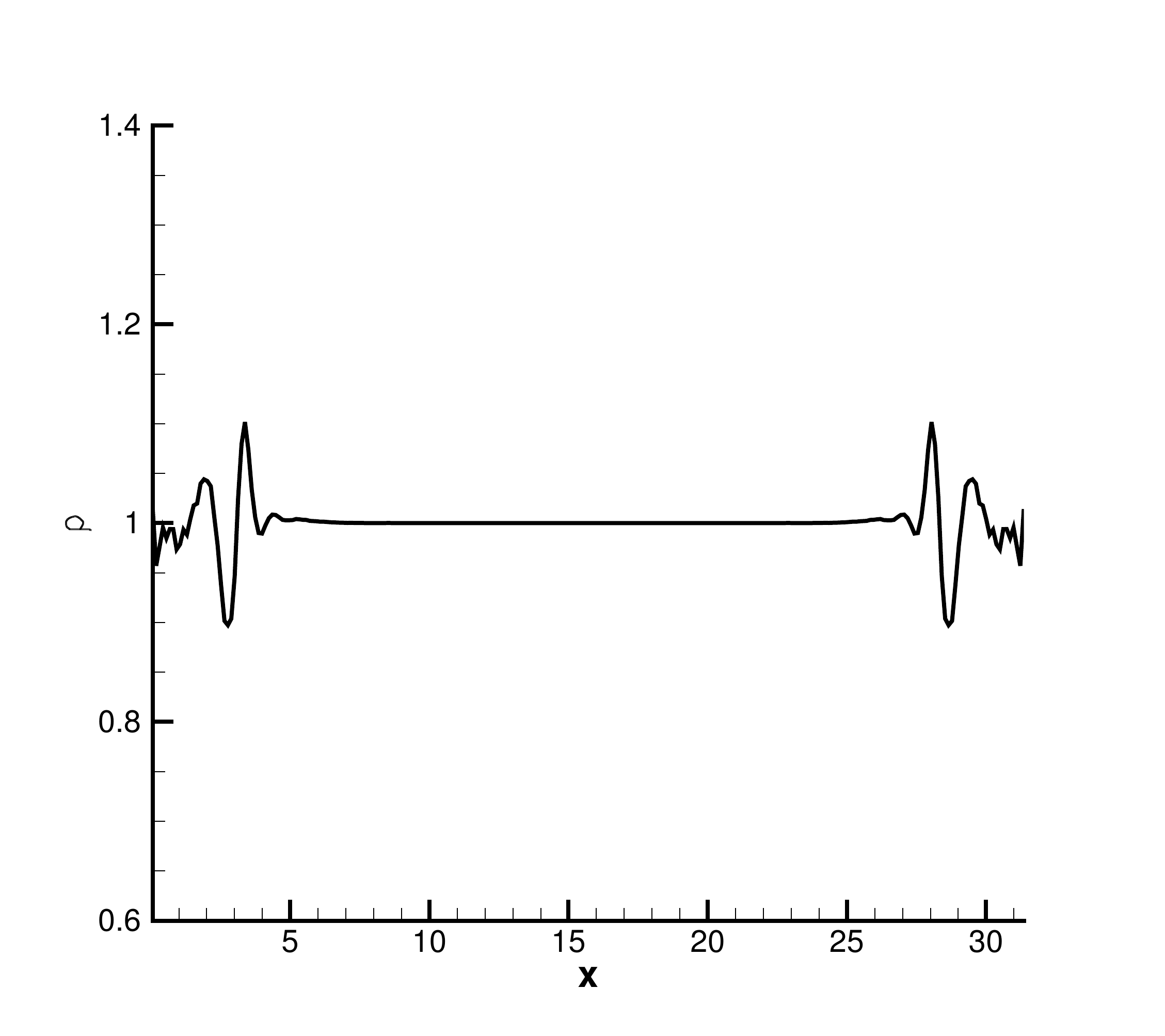}}
	\end{center}
	\caption{SW instability with parameter choice 2. Plots of the computed density function $\rho_h$ at selected  time $t$ by upwind flux for the Maxwell's equations.  Sparse grid: $N=8$, $k=3$.  Adaptive sparse grid: $N=6$, $k=3$, $\eps=10^{-6}$.}
	\label{density2}
\end{figure}

\begin{figure}[htb]
	\begin{center}
		
		\subfigure[ $t=0.$ Active elements: 0.64\% ]{\includegraphics[width=2.9in,angle=0]{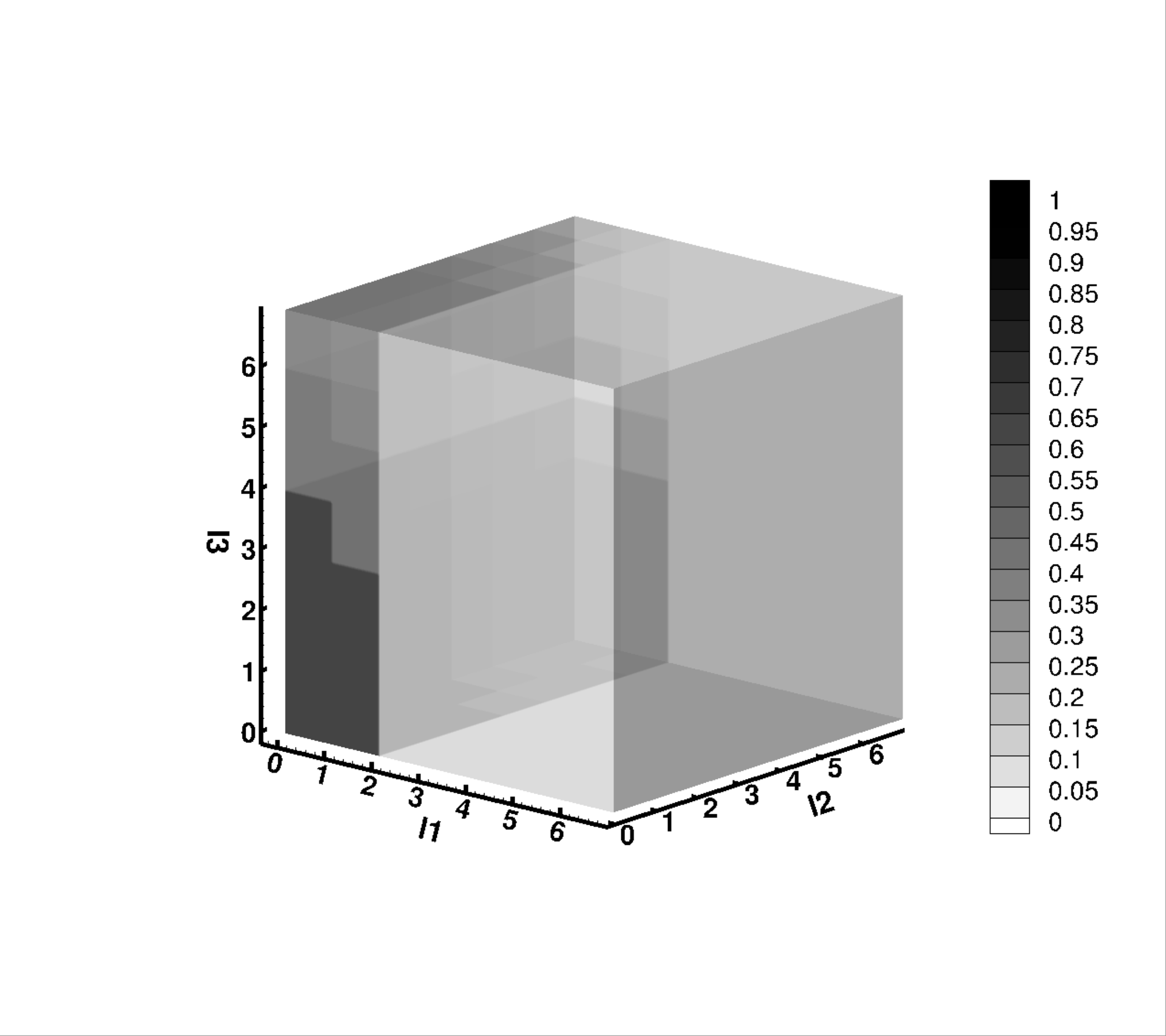}}	
		\subfigure[ $t=55.$ Active elements: 2.11\% ]{\includegraphics[width=2.9in,angle=0]{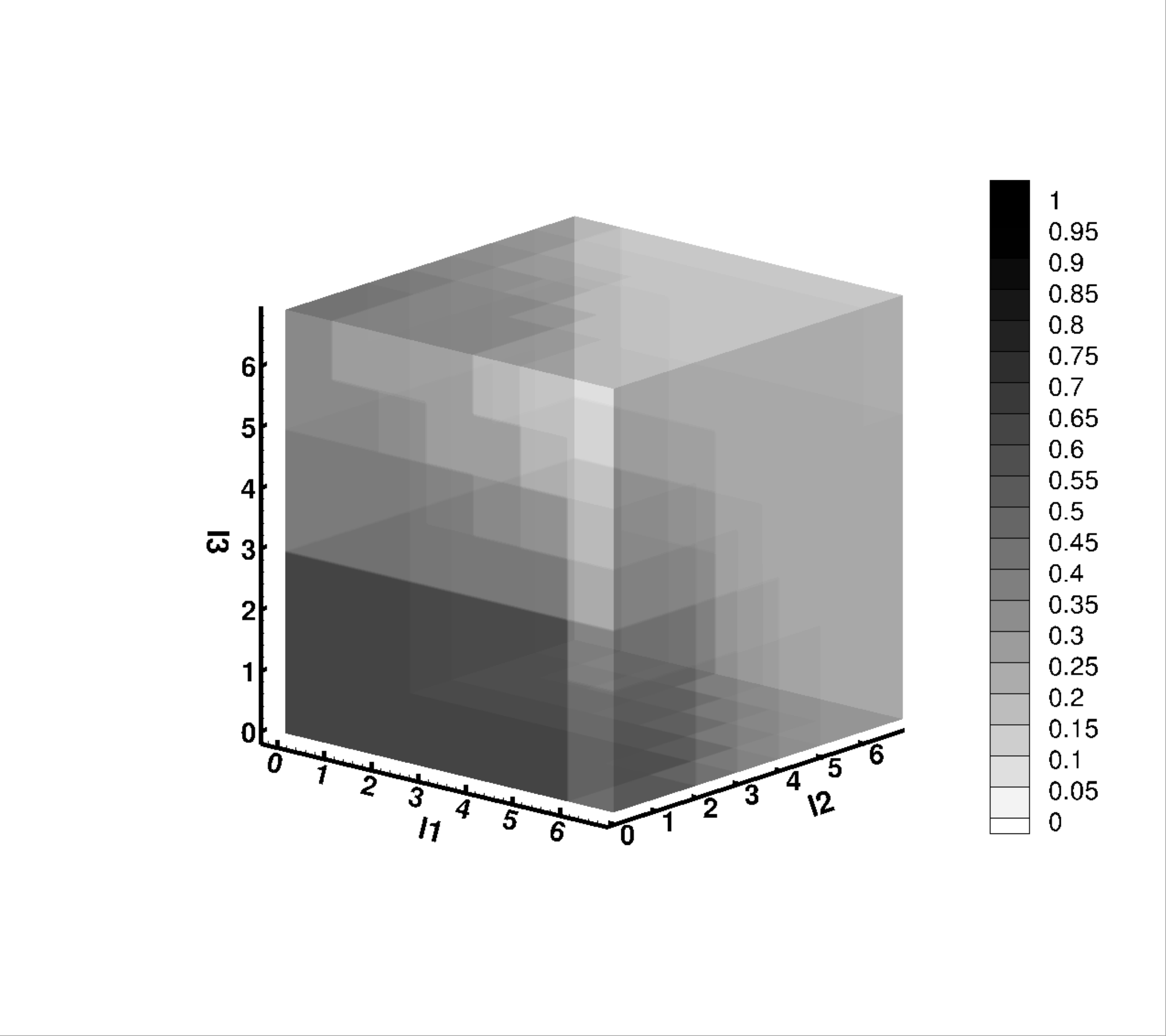}}
		\subfigure[  $t=82.$ Active elements: 5.88\%]{\includegraphics[width=2.9in,angle=0]{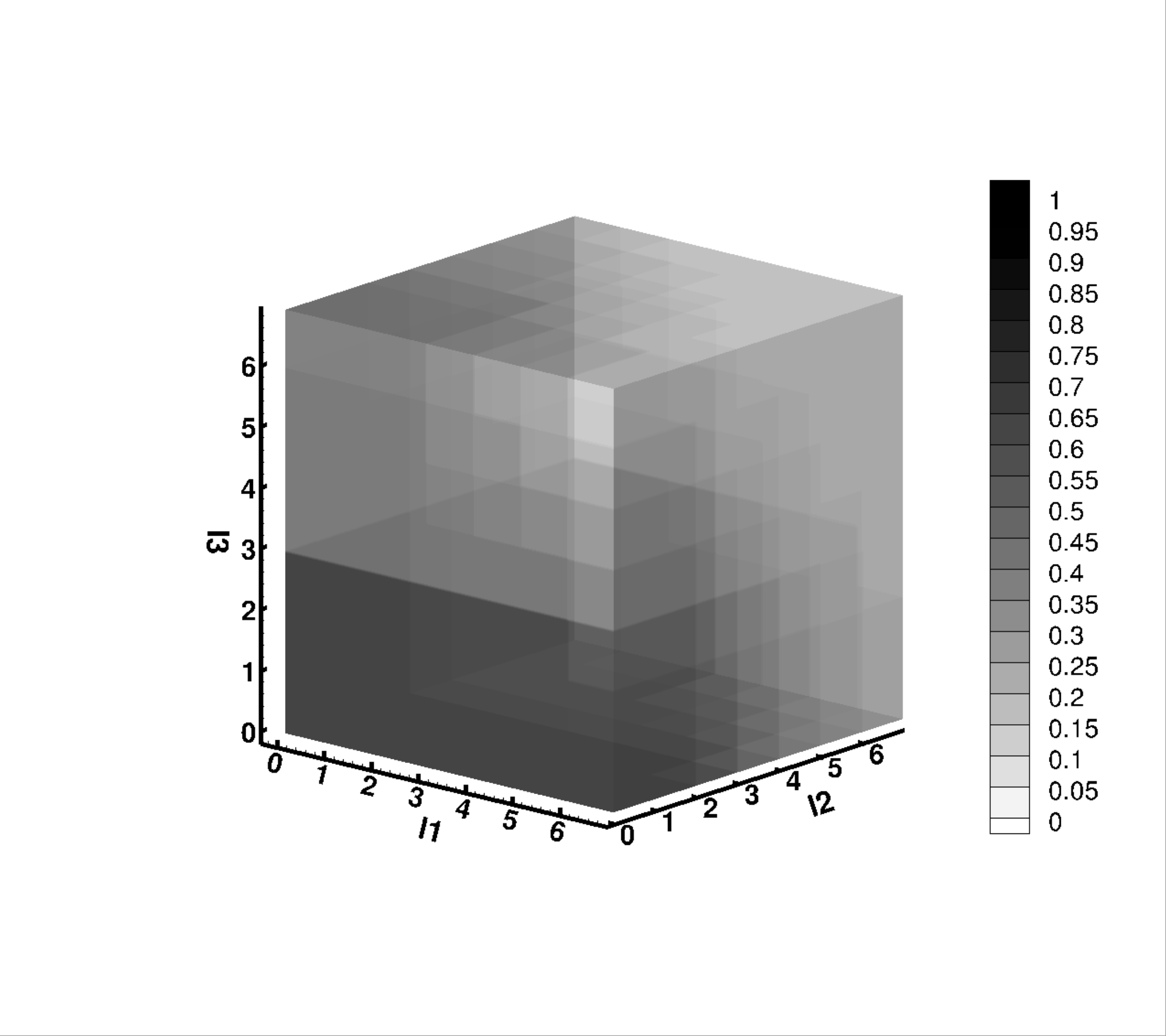}}
		\subfigure[ $t=100.$ Active elements: 22.15\%]{\includegraphics[width=2.9in,angle=0]{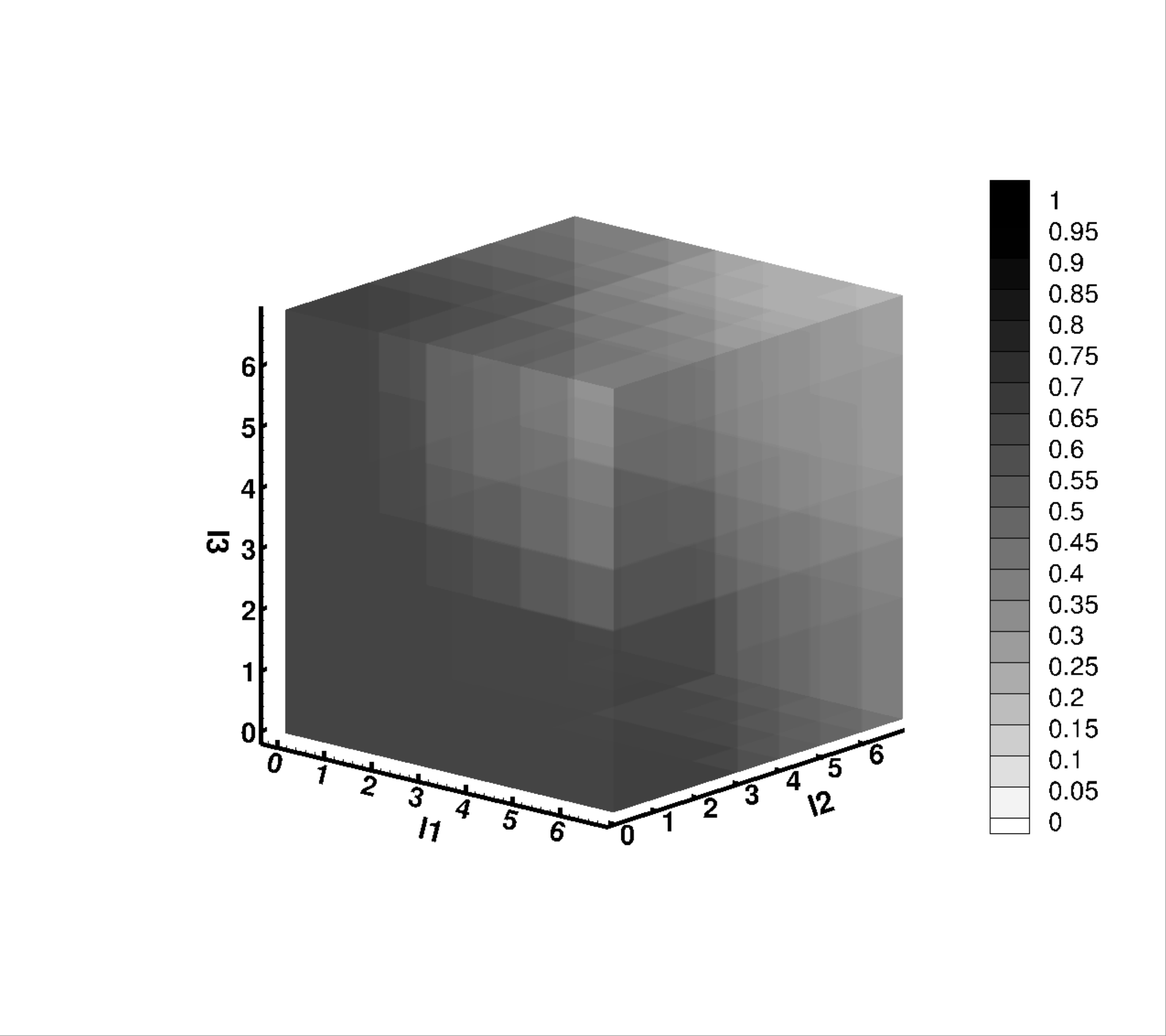}}
		
	\end{center}
	\caption{SW instability with parameter choice 2 and upwind flux for  Maxwell's equations. The percentage of active elements for each incremental space $\bW_\bl$, $\bl=(l_1,l_2,l_3)$ and $|\bl|_\infty\leq N$ in the adaptive sparse grid. $N=6$, $k=3$,  $\eps=10^{-6}$.}
	\label{percent2_ada}
\end{figure}

For choice 2, with the non-symmetric parameter set,  we take $N=8$, $k=3$ for the sparse grid scheme and  $N=6$, $k=3$, with a larger threshold $\eps = 10^{-6}$ for the adaptive sparse grid scheme.  Similar results in conservation errors are obtained as in choice 1, and the plots are omitted.
Figure \ref{keeme2} demonstrates energy transfers in this case. Compared with choice 1, the equipartition of the magnetic and electric energies at the peak is not achieved. The Log Fourier modes are plotted in Figure \ref{logfm2}. Compared with choice 1, for all four modes, the results from both schemes stay very close to the full grid approximations. The contour plots and the density plots are collected in Figures \ref{contour2}-\ref{density2}. The distribution of active elements for the adaptive scheme is presented in Figure \ref{percent2_ada}. The conclusions are similar to choice 1, i.e. the adaptive scheme offers better resolution by correctly tracking the filament structures, and the solutions are closer to the full grid approximations. Compared with choice 1, the distribution function is  less filamented, therefore less elements are needed. At $t=100,$ the percentage of active elements is $22.15\%$ compared with $52.41\%$ for choice 1. This also explains why the density plots in Figure \ref{density2} show less difference between the two methods when compared with Figure \ref{density1} for choice 1.

\subsection{2D2V Landau damping}

In this subsection, we study the   numerical performance of the proposed schemes   in the context of 2D2V simulations. Here, we consider the important example of   Landau damping, where the model is the zero magnetic limit of the VM system, the Vlasov-Amp\`{e}re (VA) system, i.e. we have $\bB=0.$ 
This is a classical test example that been previously computed by many algorithms \cite{FilbetS, einkemmer2019performance, kormann2015semi, kormann2016sparse}.  

%
%
%

\begin{exa}
	\label{ex:ld}
	We consider the 2D2V Landau damping with two spatial variables $x_1$ and $x_2$, and two velocity variables $\xi_1$ and $\xi_2$
	\begin{align}
	f_t + \xi_1 f_{x_1} &+ \xi_2 f_{x_2} + E_1 f_{\xi_1} + E_2 f_{\xi_2} = 0~, \\
	\frac{\df E_1}{\df t} &= -  j_1, \quad \frac{\df E_2}{\df t} = -  j_2~,
	\end{align}
	where
	\beq j_1=\int_{-\infty}^{\infty} \int_{-\infty}^{\infty}
	f(x_1, x_2, \xi_1, \xi_2, t) \xi_1 \,d\xi_1 d\xi_2,\quad j_2=\int_{-\infty}^{\infty} \int_{-\infty}^{\infty} f(x_1, x_2, \xi_1,
	\xi_2, t) \xi_2 \,d\xi_1 d\xi_2~.
	\eeq
	$f=f(x_1, x_2, \xi_1, \xi_2, t)$ is the distribution function, 
	$\textbf{E}=(E_1(x_1, x_2, t),E_2(x_1, x_2, t), 0)$ is a 2D electric field.  	
\end{exa}

The initial conditions are perturbed Maxwellian distributions
\begin{align}
f(x_1, x_2, \xi_1, \xi_2, 0)&=\frac{1}{2 \pi}
e^{- (\xi_1^2 + \xi_2^2) /2} [1+ \alpha (\cos(k_1 x_1)+\cos(k_2 x_2)) ],\\
E_1(x_1, x_2, 0)&=\frac{\alpha}{k_1} \sin(k_1 x_1), \quad   E_2(x_1, x_2, 0)=\frac{\alpha}{k_2} \sin(k_2 x_2) , 
\end{align}	
with $\alpha = 0.01$ and $\alpha = 0.5$ for weak and strong Landau damping, respectively. The velocity space is truncated to $\Omega_\xi=[-6, 6]^2$, the wave numbers are chosen as $k_1=k_2=0.5$, and the length of the period box in the physical space is $L_x=L_y=4 \pi$. The periodic boundary conditions and zero boundary condition are imposed on the physical and velocity spaces, respectively.

\textbf{Accuracy test and comparisons.}
Similar as in Section \ref{subsec:SW}, we can use the time reversibility of this system to test the accuracy. We compute the solutions $f(\bx, \xi, T),$ $\bE(\bx, T)$ to $T=0.5$ and reverse the velocity field, yielding $f(\bx, -\xi, 0.5),$ $\bE(\bx, 0.5)$, and evolve the system to $T=1$. Finally, we compare the numerical solutions with the initial condition  
$f(\bx, -\xi, 0),$ $\bE(\bx, 0)$. 
We test accuracy for the sparse grid DG method with $k=1,\,2,\,3$ on different levels of meshes. 
When $k=3$, we take $\displaystyle\Delta t=O( h_N^{4/3})$ to match the temporal and spatial orders in  the convergence study.
The  $L^2$ errors and orders of the numerical solutions for weak and strong Landau damping are shown in Tables \ref{table:ld_1} and \ref{table:ld_2}, respectively.  
We observe close to $(k+\frac12)$-th order accuracy for the weak Landau damping in Table \ref{table:ld_1}. There is a loss of accuracy for $\bE$ when $k=3, \, N=8$. This may be due to the domain cut-off in the velocity space, which causes truncation error of magnitude about $10^{-9}$.
For strong Landau damping, we observe about $(k+\frac12)$-th order accuracy with $k=1$. While for $k=2,3$, there is a reduction of accuracy for both $f$ and $\bE$. The reason is that for strong Landau damping, the solution becomes quickly filamented, and the sparse grid space in dimension $d=4$  can not provide enough resolution at a coarse mesh level. 

\begin{table}
	\caption{$L^2$ errors and orders of accuracy for the sparse grid DG method in Example \ref{ex:ld} with $\alpha=0.01$. Run to $T=0.5$ and back to $T=1.0$. The orders are measured with respect to $h_N,$ which is the size of the smallest mesh in each direction.
	}
	\centering
	\begin{tabular}{|c|c| c c|c c|}
		\hline
		
		\multirow{2}{*}{} & \multirow{2}{*}{$N$} & 	\multicolumn{2}{|c|}{$ f$} &		\multicolumn{2}{|c|}{$ \bE$}   \\
		\cline{3-6}
		&  & $L^2$ error & order & $L^2$ error & order\\
		
		\hline 
		
		\multirow{4}{*}{$ k=1$} & 5	& 	1.01E-03&	& 4.71E-05	&	\\
		& 6		&	3.72E-04	&	1.44	& 8.77E-06	& 2.43\\
		& 7		&	1.24E-04	&	1.58	& 2.32E-06	& 1.92\\
		& 8		&	4.05E-05	&	1.61	& 5.07E-07	& 2.19\\
		
		\hline
		\multirow{4}{*}{$ k=2$} & 5 &	1.61E-04	&	& 1.20E-06	& \\
		& 6		&	4.42E-05	&	1.86	& 2.36E-07	& 2.35\\
		& 7		&	8.54E-06	&	2.37	& 3.30E-08	& 2.84\\
		& 8		&	1.22E-06	&	2.81	& 4.32E-09	& 2.93\\
		
		\hline
		\multirow{4}{*}{$ k=3$} & 5		&	2.76E-05	&		& 2.32E-08	& \\
		& 6		&	8.75E-07	&	4.98	& 1.41E-09	& 4.04\\
		& 7		&	1.85E-07	&	2.24	& 1.18E-10	& 3.58\\
		& 8		&	1.42E-08	&	3.70	& 7.90E-11	& 0.58\\ 
		\hline
		
	\end{tabular}
	\label{table:ld_1}
\end{table}

\begin{table}
	\caption{$L^2$ errors and orders of accuracy for the sparse grid DG method in Example \ref{ex:ld} with $\alpha=0.5$. Run to $T=0.5$ and back to $T=1.0$. The orders are measured with respect to $h_N,$ which is the size of the smallest mesh in each direction.}
	\centering
	\begin{tabular}{|c|c| c c|c c|}
		\hline
		
		\multirow{2}{*}{} & \multirow{2}{*}{$N$} &	\multicolumn{2}{|c|}{$ f$} &		\multicolumn{2}{|c|}{$ \bE$}   \\
		\cline{3-6}
		&  & $L^2$ error & order & $L^2$ error & order\\
		
		\hline 
		
		\multirow{4}{*}{$ k=1$} & 6 &	1.62E-03&	& 1.68E-03	&	\\
		& 7		&	7.75E-04	&	1.06	& 8.54E-04	& 0.98\\
		& 8 	&	3.30E-04	&	1.23	& 3.52E-04	& 1.28\\
		& 9 	&	1.32E-04	&	1.32	& 1.34E-04	& 1.39\\
		
		\hline
		\multirow{4}{*}{$ k=2$} & 6	 &	1.60E-04	&	& 1.06E-04	& \\
		& 7 	&	4.91E-05	&	1.70	& 3.03E-05	& 1.81\\
		& 8		&	1.30E-05	&	1.92	& 1.01E-05	& 1.58\\
		& 9		&	2.71E-06	&	2.26	& 3.41E-06	& 1.57\\
		
		\hline
		\multirow{4}{*}{$ k=3$} & 5	& 5.57E-05	&	& 2.38E-05	&\\
		& 6		&	9.59E-06	&	2.54	& 6.53E-06	& 1.87\\
		& 7		&	1.86E-06	&	2.37	& 1.30E-06	& 2.33\\
		& 8		&	3.45E-07	&	2.43	& 2.37E-07	& 2.46\\
		\hline
		
	\end{tabular}
	\label{table:ld_2}
\end{table}

\begin{figure}[htp]
	\begin{center}
		\subfigure[k=1]{\includegraphics[width=.32\textwidth]{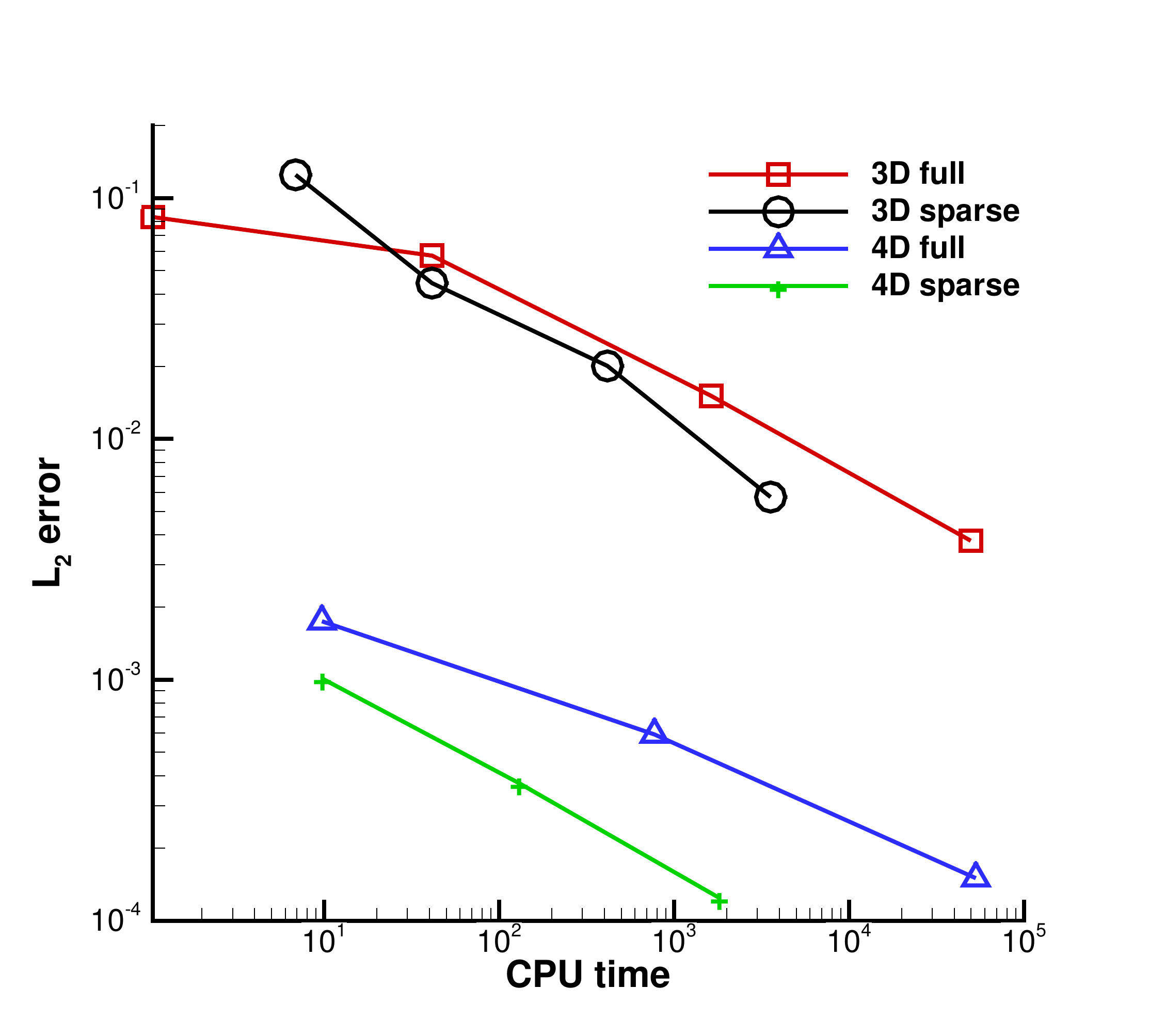}}
		\subfigure[k=2]{\includegraphics[width=.32\textwidth]{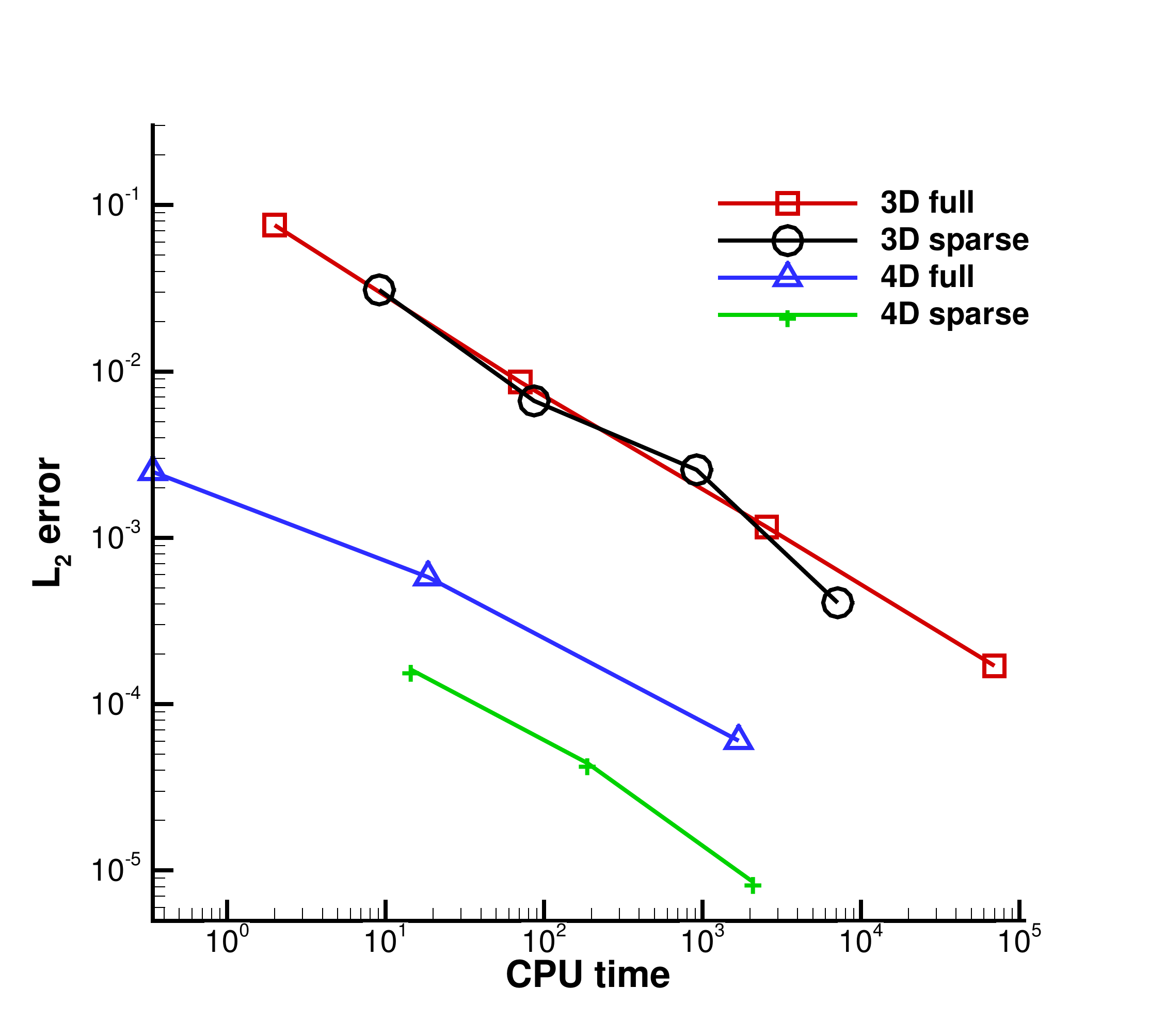}}
		\subfigure[k=3]{\includegraphics[width=.32\textwidth]{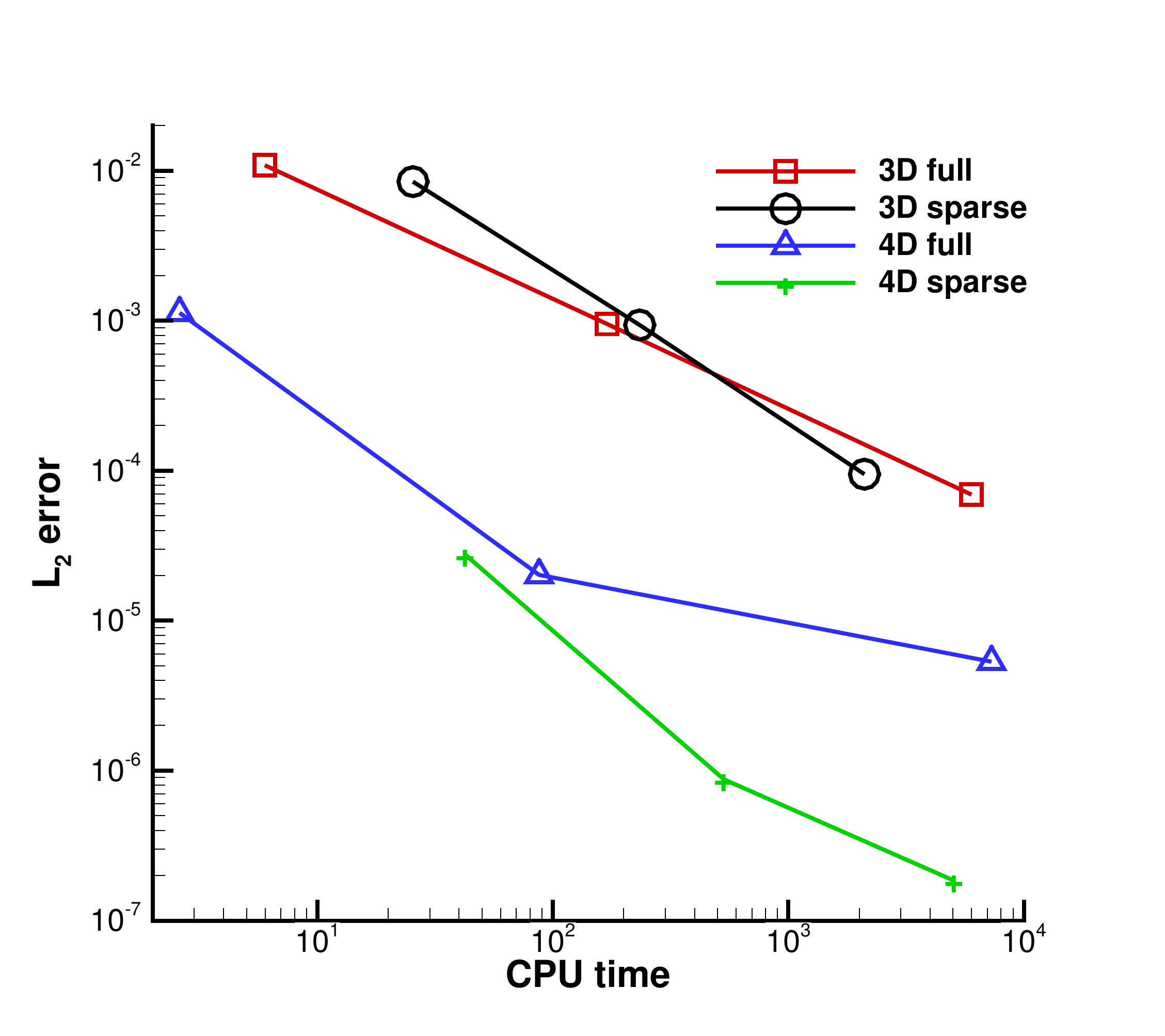}}
	\end{center}
	\caption{ $L^2$ errors of $f$ vs.  CPU time for full grid DG and sparse grid DG methods for 3D and 4D simulations. 3D: Example \ref{ex:SW}. Parameter choice 1.   Upwind flux for Maxwell's equations. 4D: Example \ref{ex:ld}. $\alpha=0.01.$  (a) k=1, (b) k=2, (c) k=3.}
	\label{fig:cpu_3d_4d_L2}
\end{figure}

\begin{figure}[htp]
	\begin{center}
		\subfigure[k=1]{\includegraphics[width=.32\textwidth]{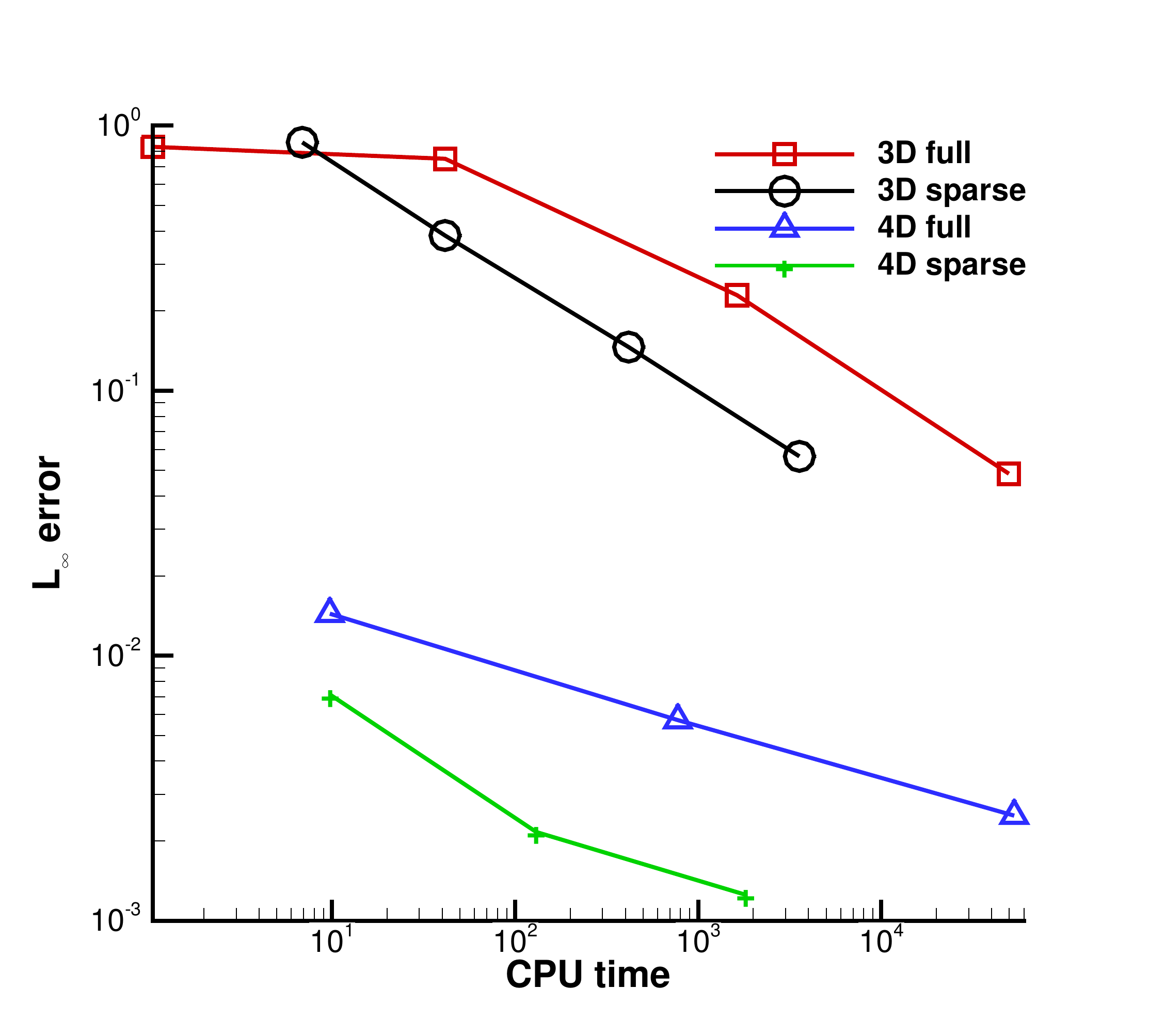}}
		\subfigure[k=2]{\includegraphics[width=.32\textwidth]{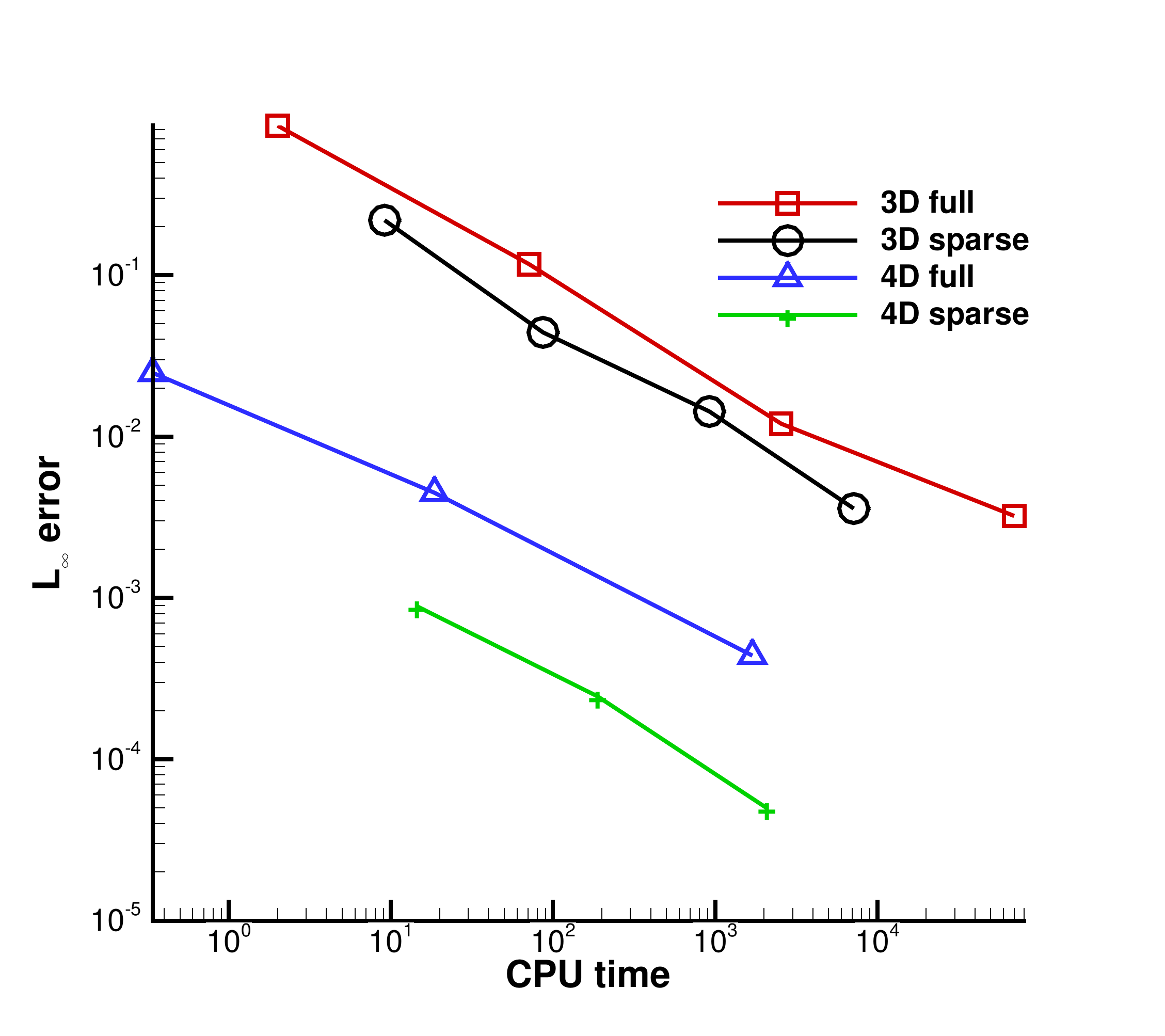}}
		\subfigure[k=3]{\includegraphics[width=.32\textwidth]{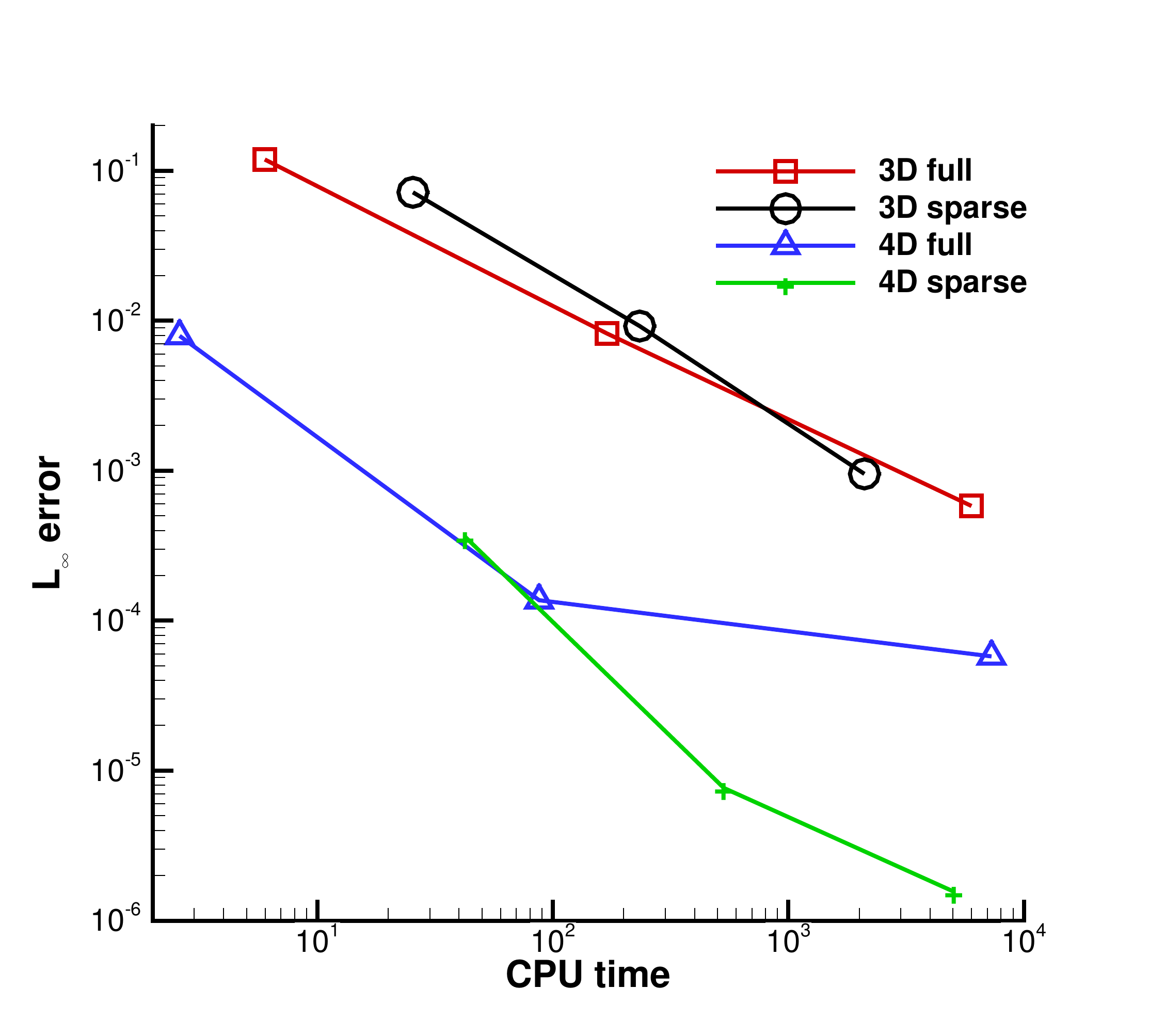}}
	\end{center}
	\caption{ $L^{\infty}$ errors of $f$ vs.  CPU time for full grid DG and sparse grid DG methods for 3D and 4D simulations. 3D: Example \ref{ex:SW}. Parameter choice 1.  Upwind flux for Maxwell's equations. 4D: Example \ref{ex:ld}. $\alpha=0.01.$  (a) k=1, (b) k=2, (c) k=3.}
	\label{fig:cpu_3d_4d_L8}
\end{figure}

We then use this example to compare the performance of the full grid DG vs. sparse grid DG schemes in 4D vs. 3D. Notice that though the setup in streaming Weibel instability and weak Landau damping is different,  we hope this can still illustrate the main difference of the schemes in different dimensions by focusing on the rate of change of errors vs. CPU time.
 In Figures \ref{fig:cpu_3d_4d_L2}-\ref{fig:cpu_3d_4d_L8}, we  plot the errors vs. CPU time with $k=1,2,3$   for full grid and sparse grid DG schemes for Examples \ref{ex:SW}-\ref{ex:ld}. The time discretizations are TVD-RK3 method  for $k=1, 2,$ and RK4 method for $k=3.$ We use a fixed $CFL=0.1$ in all computations.
For the 4D case, the mesh levels $N$ are taken as from 3 to 5, 2 to 4 and 2 to 4 for $k=1,2,3$, respectively, in the full grid DG method and taken as from 5 to 7 for $k=1,2,3$, respectively, in the sparse grid DG method. The mesh levels $N$ for the 3D case is the same as those in Figures \ref{fig:cpu_2d}-\ref{fig:cpu_3d}.
From Figures  \ref{fig:cpu_3d_4d_L2}-\ref{fig:cpu_3d_4d_L8}, it is evident that the full grid method suffers from the curse of dimensionality, i.e. when $d$ increases from 3 to 4, the absolute value of the slope of log-log plot of errors vs. CPU times becomes smaller. This is not as pronounced for the sparse grid method. We can see the slope stays roughly the same for the 3D and 4D computations. This shows  effectiveness of the sparse grid approach in dimension reduction.

\textbf{Numerical results for long time simulations.}
We show the numerical results of the methods in long time simulation for the VA system.
The macroscopic quantities are defined as
$$K_1=\frac{1}{2 L_x L_y} \int_\Omega f \xi_1^2\,d\xi_1 d\xi_2 dx_1 dx_2,\quad K_2=\frac{1}{2 L_x L_y} \int_\Omega f \xi_2^2\,d\xi_1 d\xi_2 dx_1 dx_2,$$
$$ \mathcal{E}_1=\frac{1}{2 L_x L_y} \int_\Ox E_1^2 \, dx_1 dx_2, \quad \mathcal{E}_2=\frac{1}{2 L_x L_y} \int_\Ox E_2^2\, dx_1 dx_2,$$
where $K_1, K_2$ are the scaled kinetic energies in each direction, $\mathcal{E}_1, \mathcal{E}_2$ are the scaled electric energies in each direction.
Therefore, the scaled kinetic and electric energies are the summation of the corresponding components in each direction, and the scaled total energy is defined as the summation of $K_1, K_2, \mathcal{E}_1, \mathcal{E}_2.$
The scaled momentum  $P_1$, $P_2$ and enstrophy are 
$$P_1=\frac{1}{L_x L_y} \int_\Omega  \xi_1 f \,d\xi_1 d\xi_2 dx_1 dx_2 ,\quad P_2=\frac{1}{L_x L_y}  \int_\Omega  \xi_2 f \,d\xi_1 d\xi_2 dx_1 dx_2, $$ 
$$ En=\frac{1}{L_x L_y} \int_\Omega  f^2 \,d\xi_1 d\xi_2 dx_1 dx_2.$$

\begin{figure}[htb]
	\begin{center}
		\subfigure[Electric energy]{\includegraphics[width=3in,angle=0]{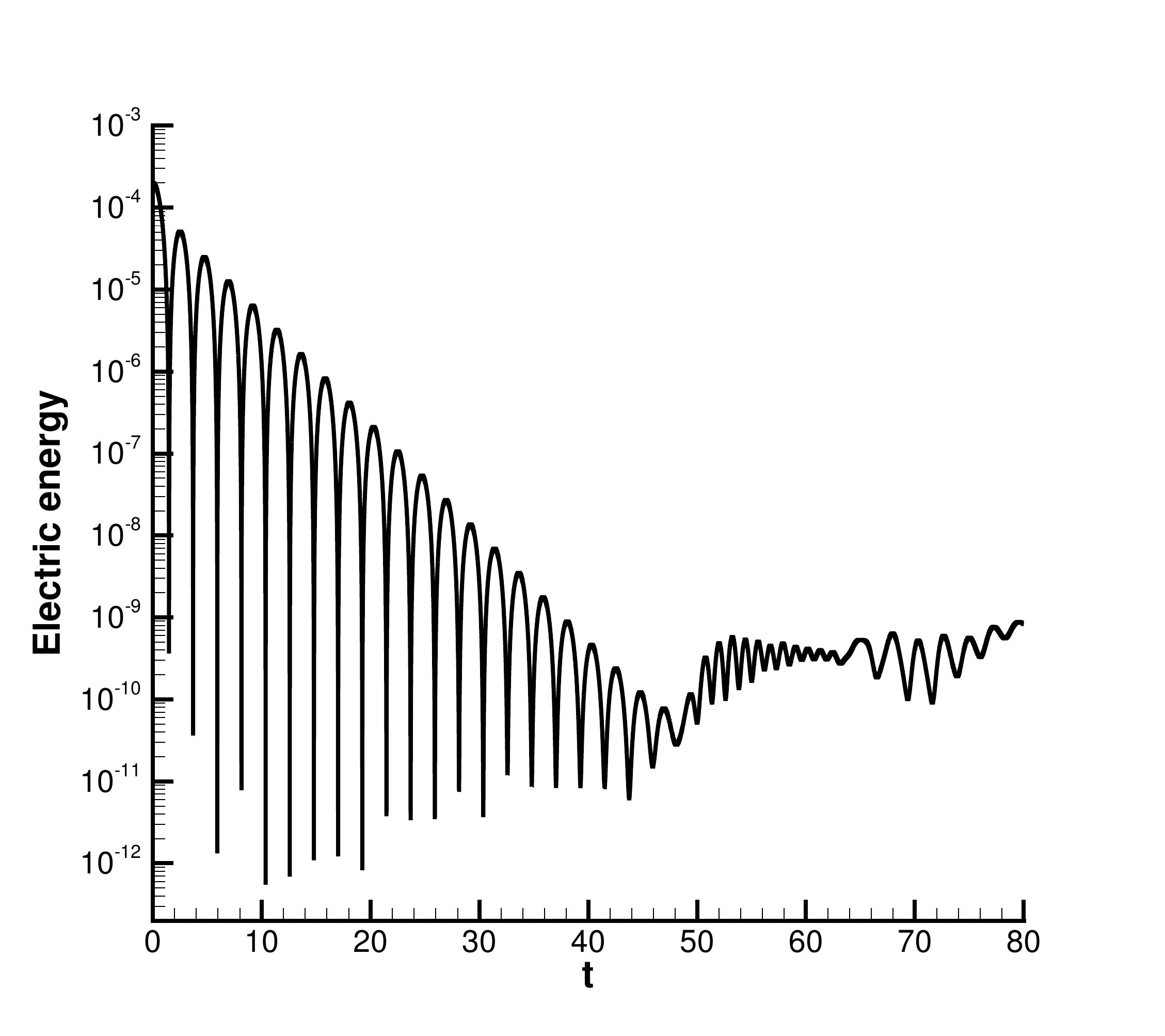}}
		\subfigure[Relative error in mass]{\includegraphics[width=3in,angle=0]{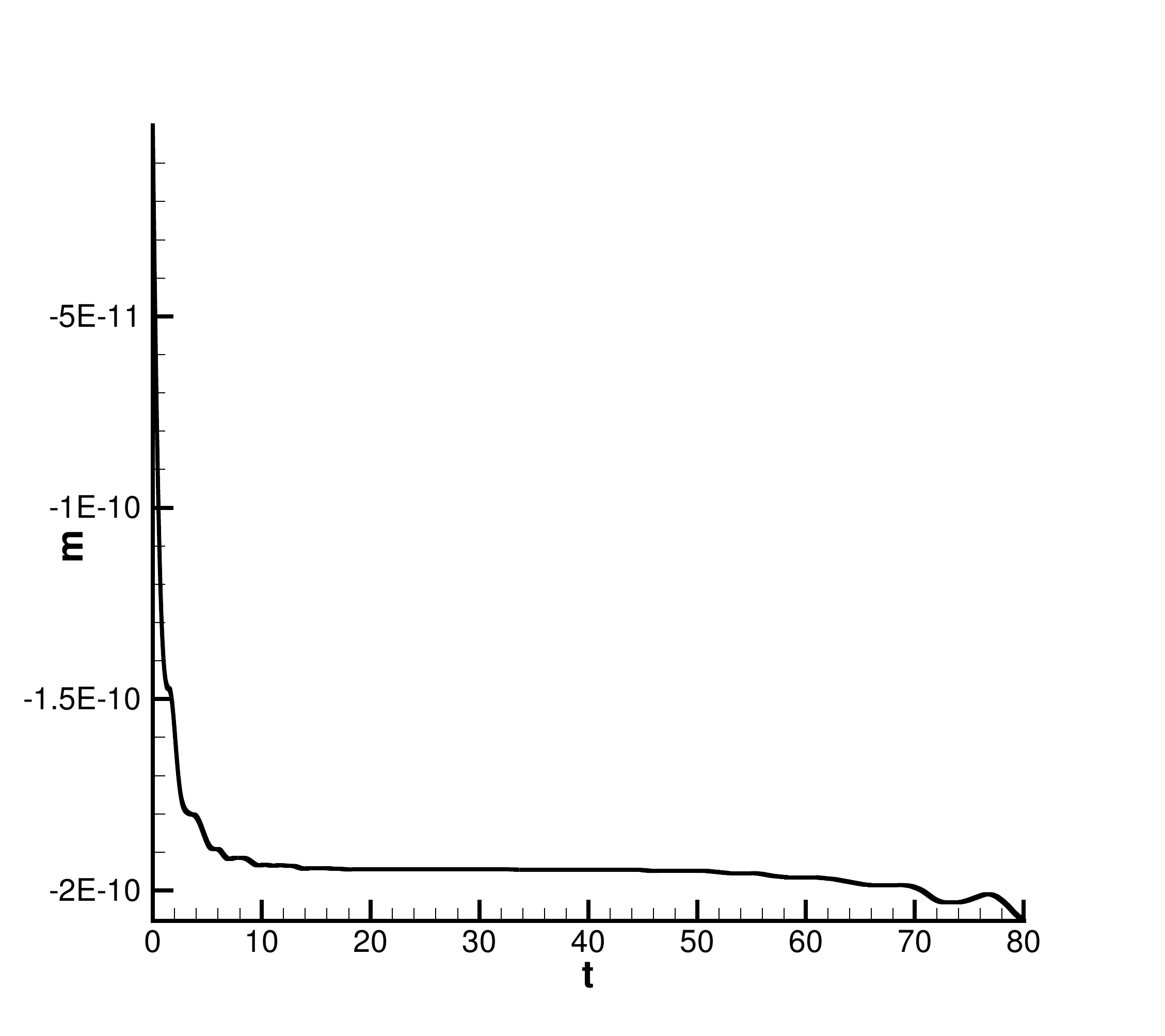}}
		\subfigure[Relative error in total energy]{\includegraphics[width=3in,angle=0]{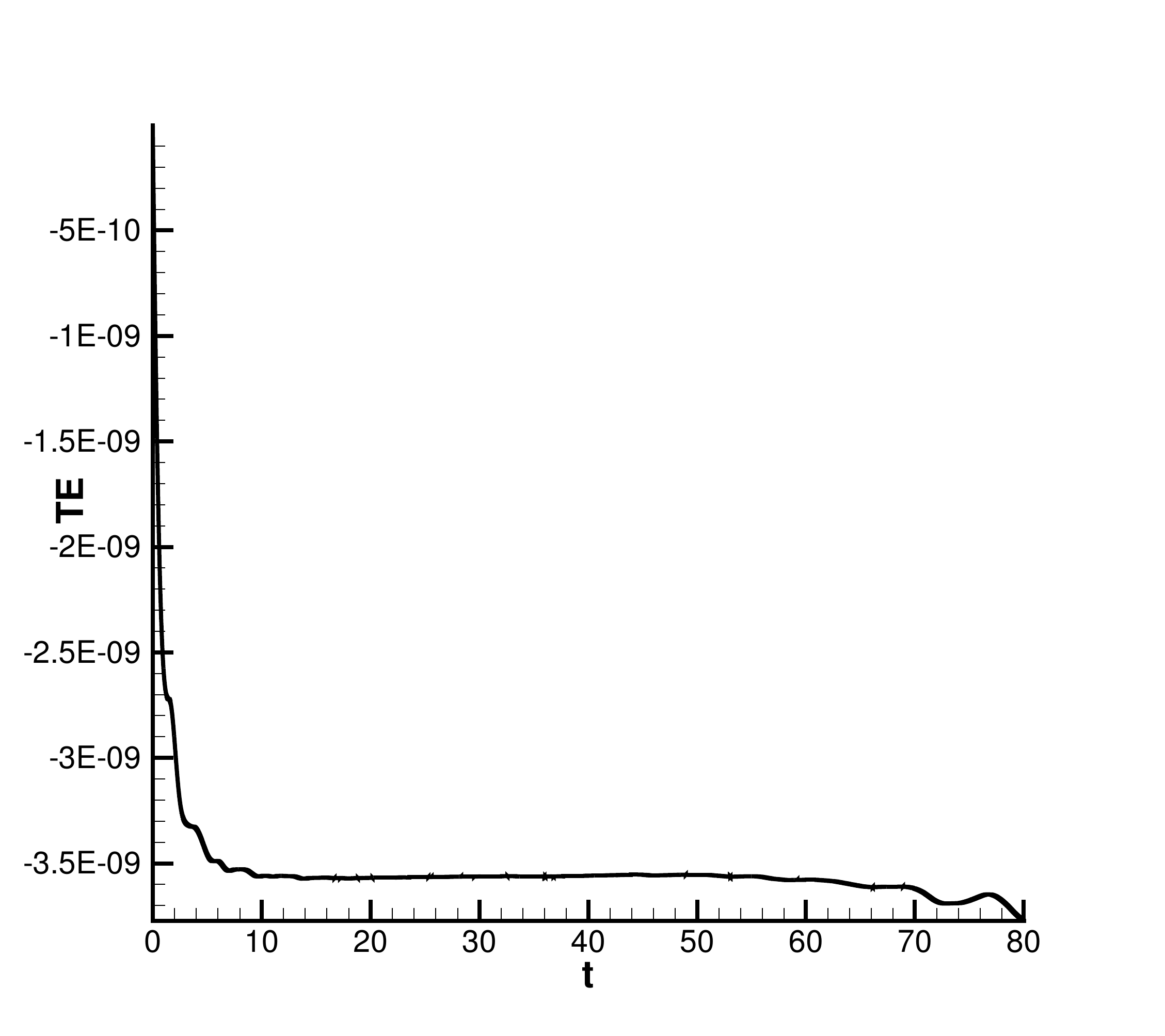}}
		\subfigure[error in momentum P1]{\includegraphics[width=3in,angle=0]{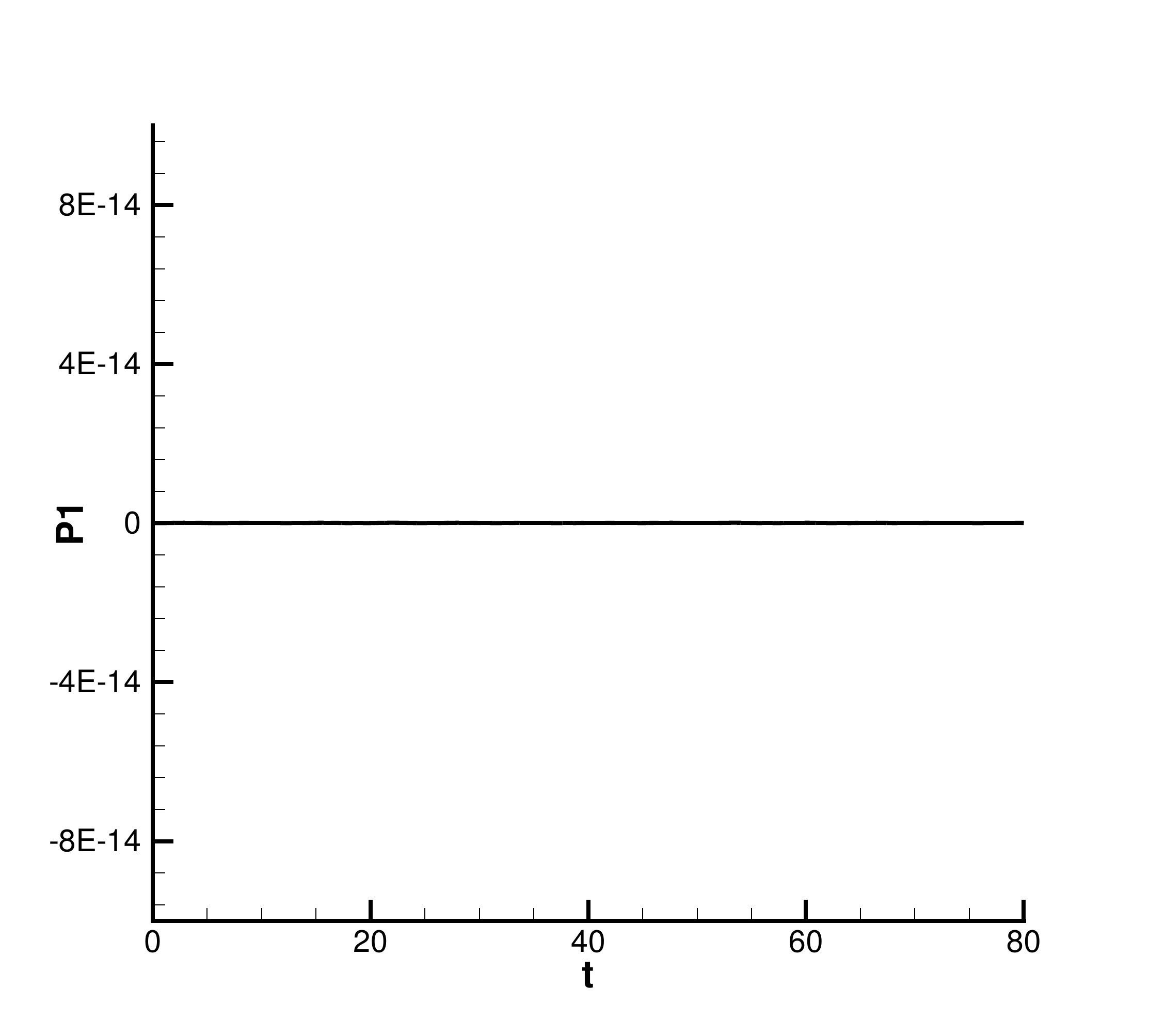}}
		\subfigure[error in momentum P2]{\includegraphics[width=3in,angle=0]{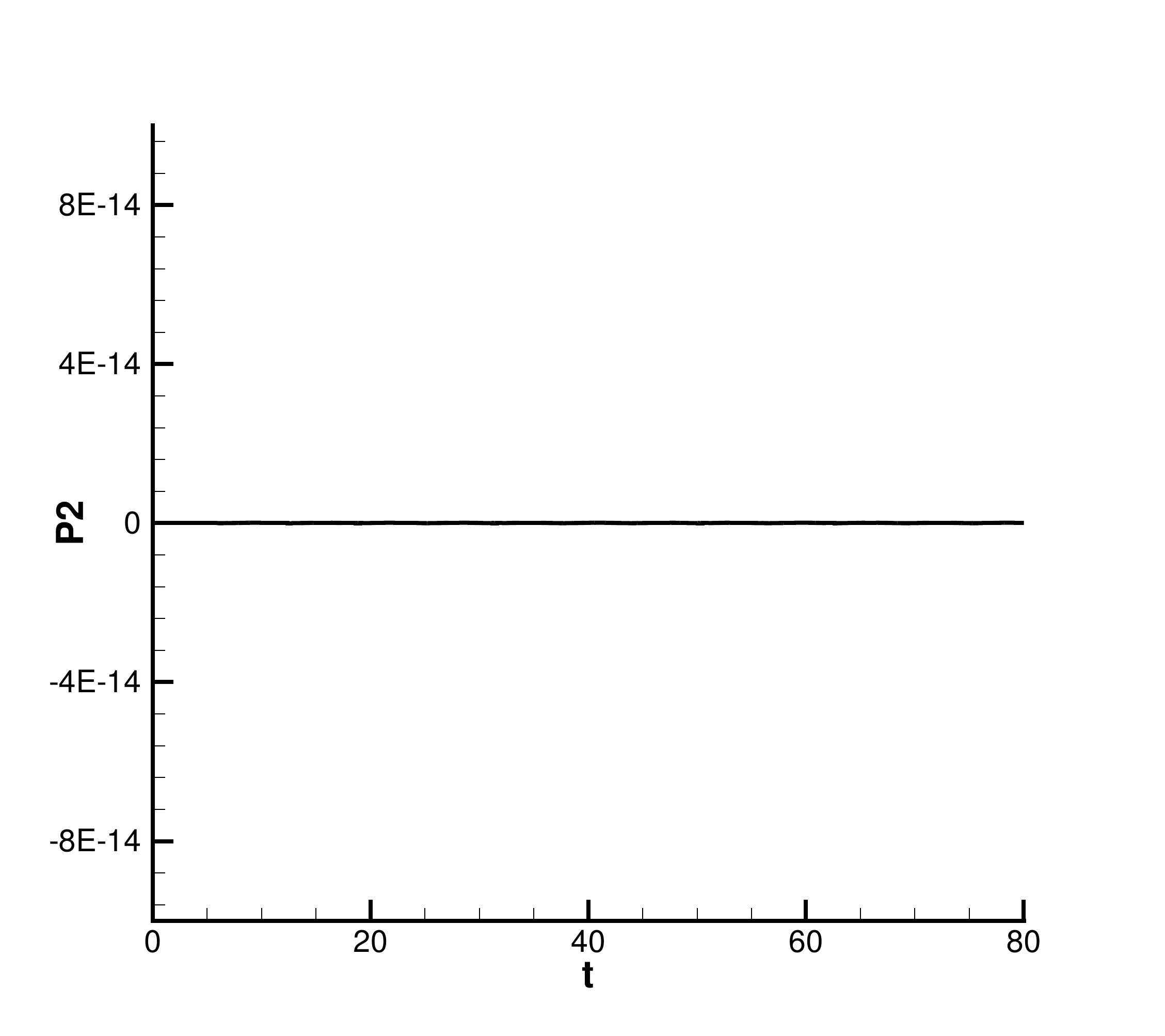}}
		\subfigure[Enstrophy]{\includegraphics[width=3in,angle=0]{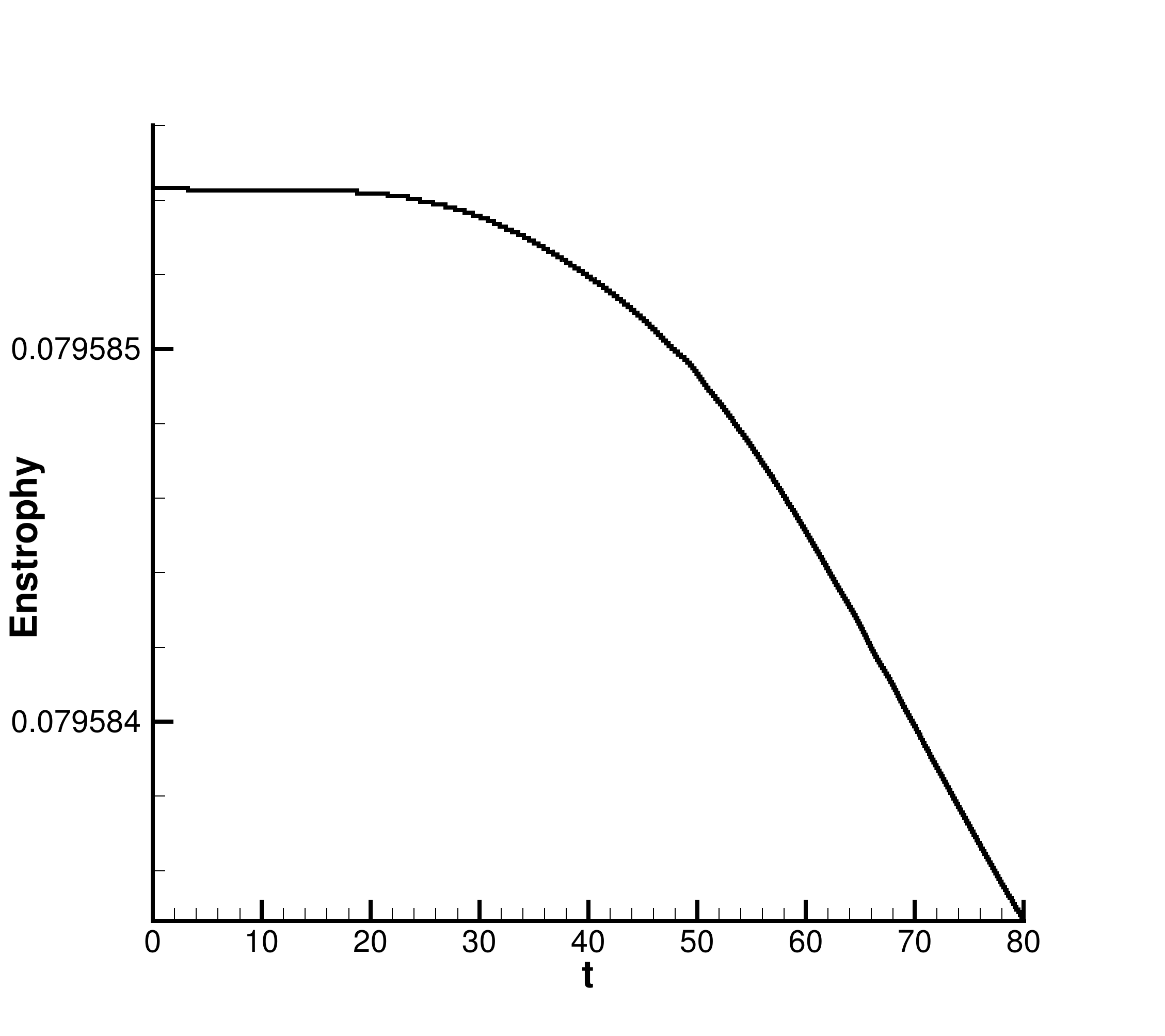}}
	\end{center}
	\caption{2D2V weak Landau damping with $\alpha=0.01$. Time evolution of electric energy, relative error of mass and total energy, error of momentum and enstrophy. Sparse grid: $N=7$, $k=3$. }
	\label{2d2v_ener_mass_te1}
\end{figure}

For weak Landau damping, we take $N=7, k=3$ for the sparse grid method. The solutions are computed up to $T=80$.  
We present time evolution of the electric energy, relative error of mass and total energy, error in momentum and enstrophy in Figure \ref{2d2v_ener_mass_te1}.
We measure the errors of the conserved macroscopic quantities: mass, total energy and momentum. From Figure \ref{2d2v_ener_mass_te1}, we observe that the sparse method can conserve these quantities very well. The largest relative errors in mass and total energy are on the order of $10^{-10}$ and $10^{-9}$. The largest error in enstrophy is on the order of $10^{-6}$. The damping in enstrophy comes from the numerical dissipation enforced by the Lax-Friedrichs flux.

\begin{figure}[htb]
	\begin{center}
		\subfigure[Electric energy]{\includegraphics[width=3in,angle=0]{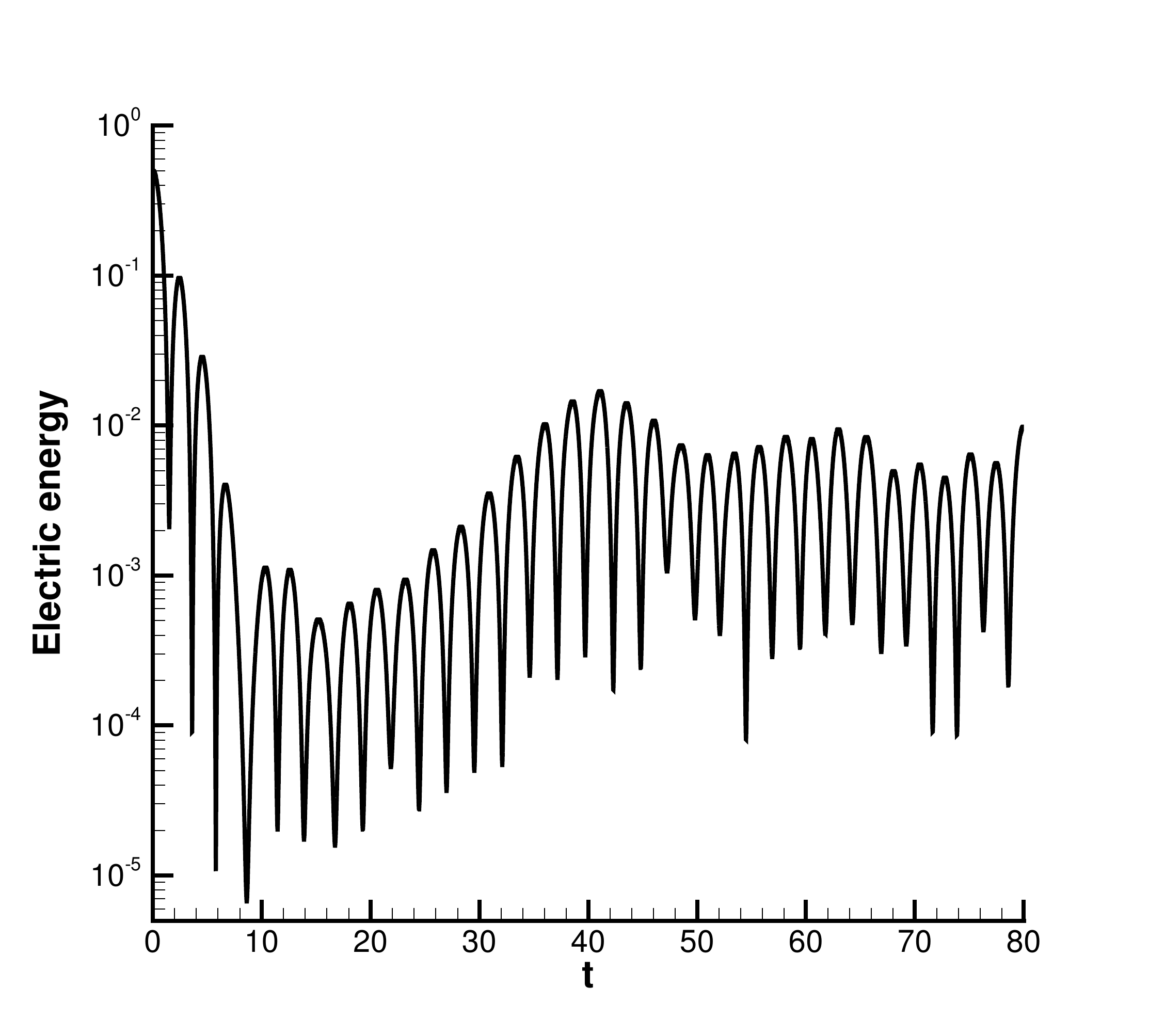}}
		\subfigure[Relative error in mass]{\includegraphics[width=3in,angle=0]{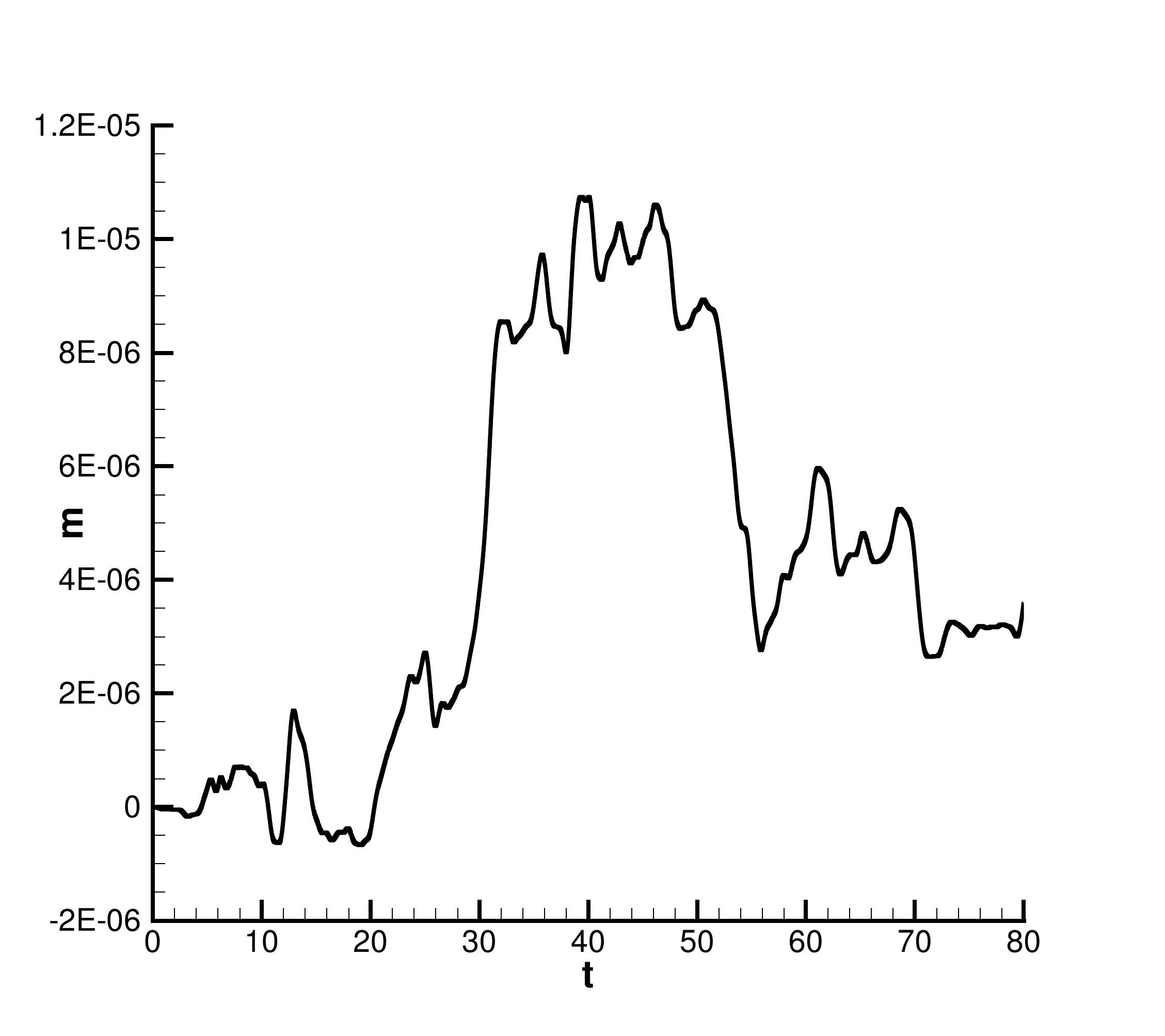}}
		\subfigure[Relative error in total energy]{\includegraphics[width=3in,angle=0]{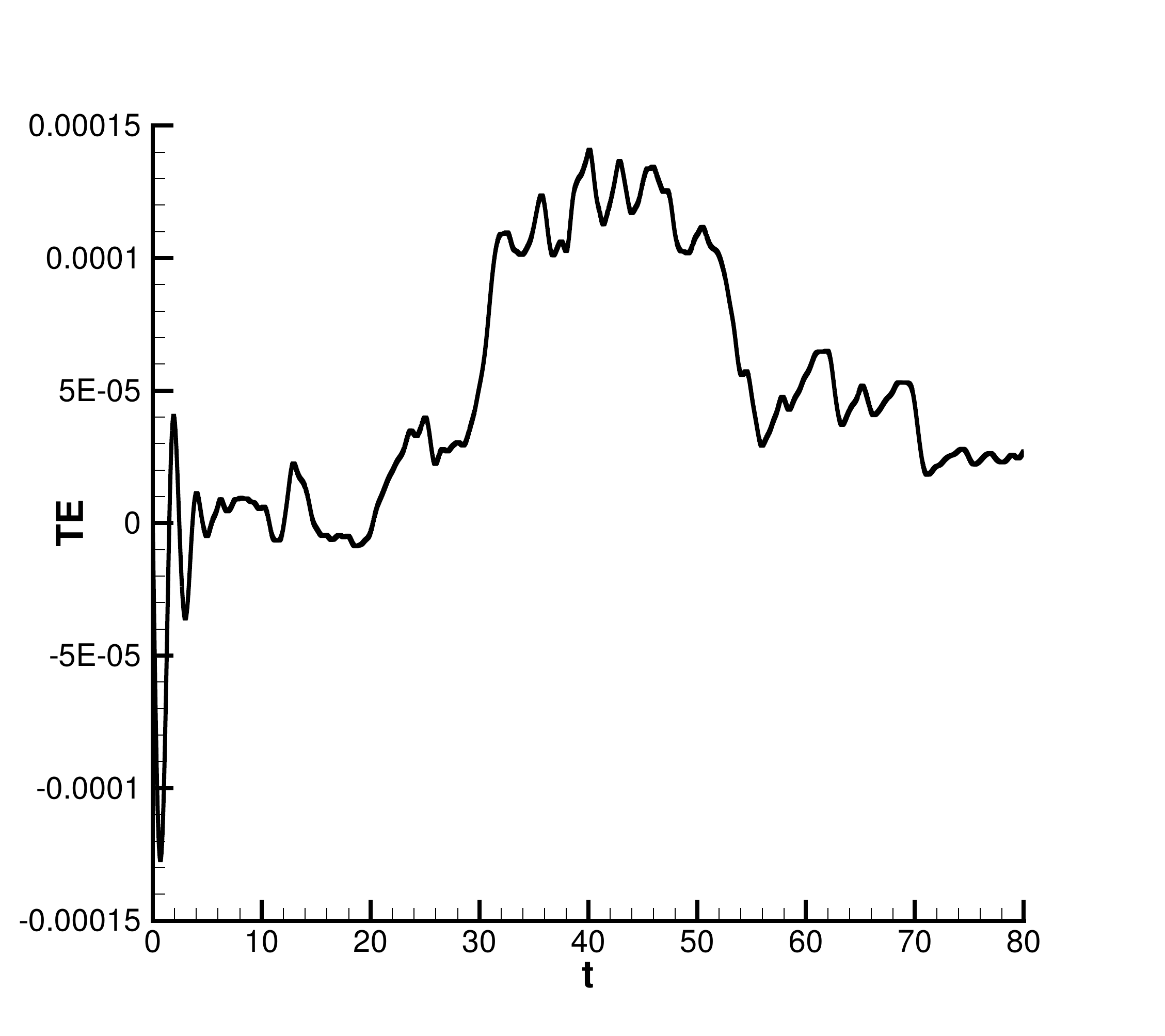}}
		\subfigure[error in momentum P1]{\includegraphics[width=3in,angle=0]{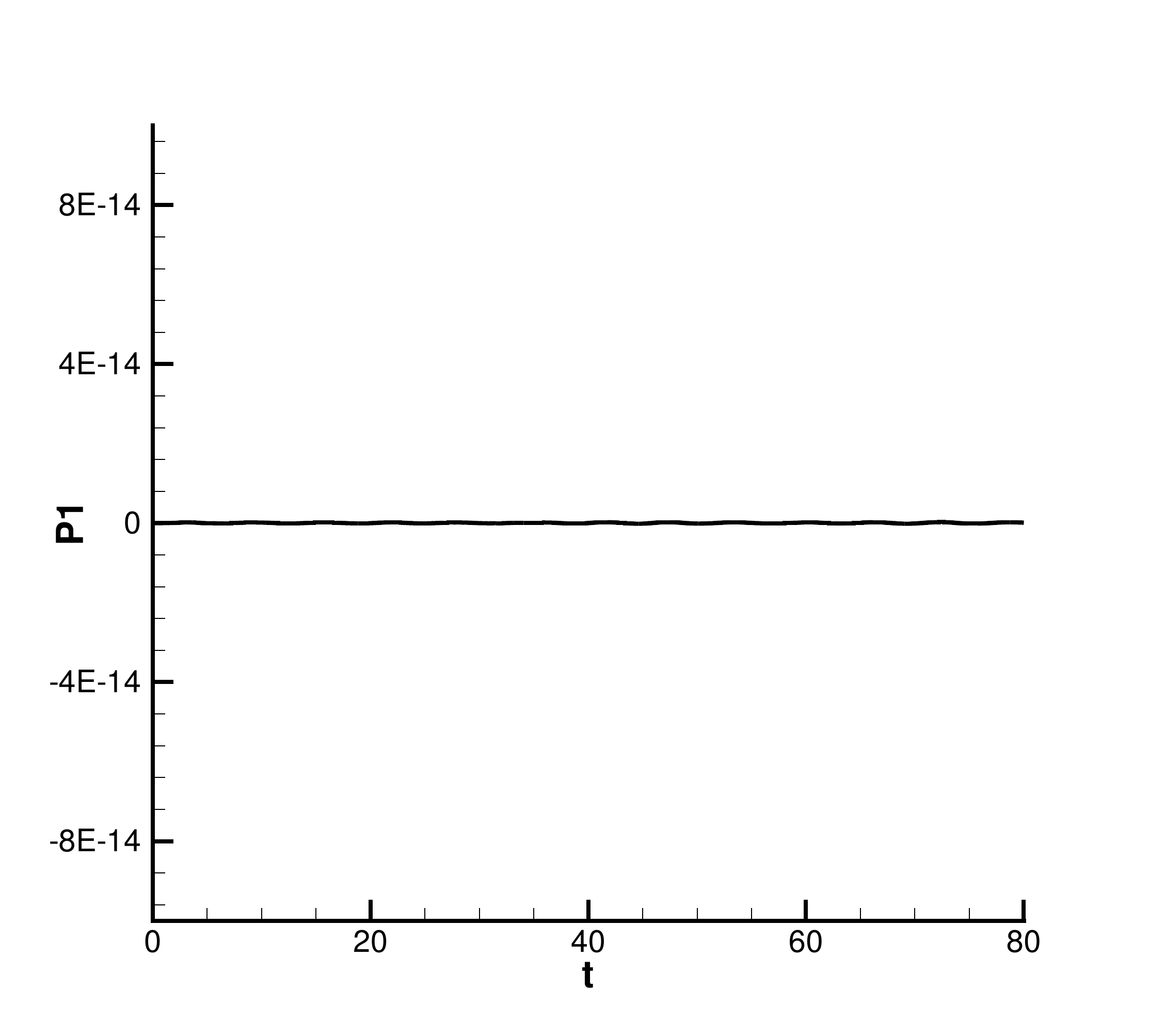}}
		\subfigure[error in momentum P2]{\includegraphics[width=3in,angle=0]{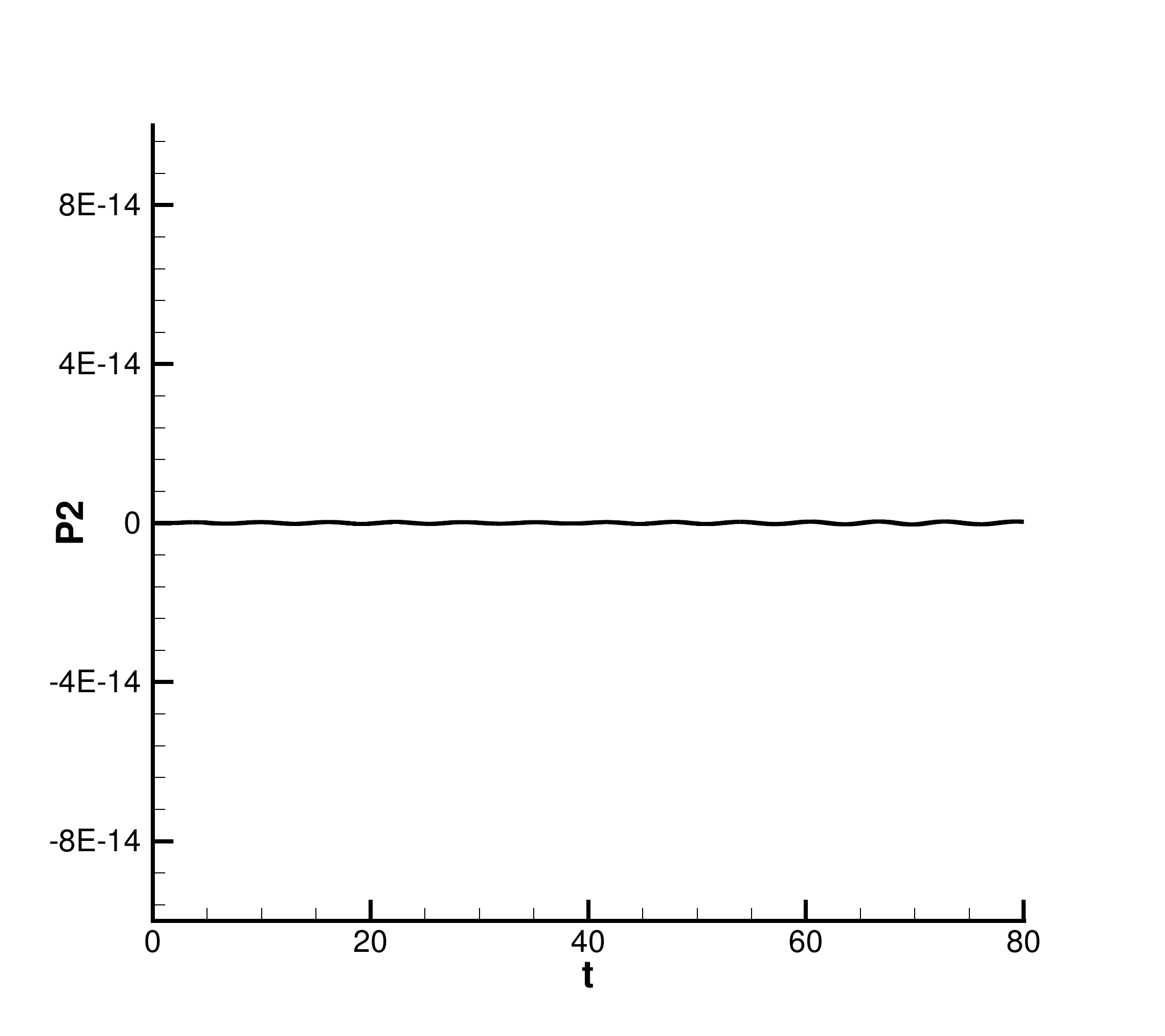}}
		\subfigure[Enstrophy]{\includegraphics[width=3in,angle=0]{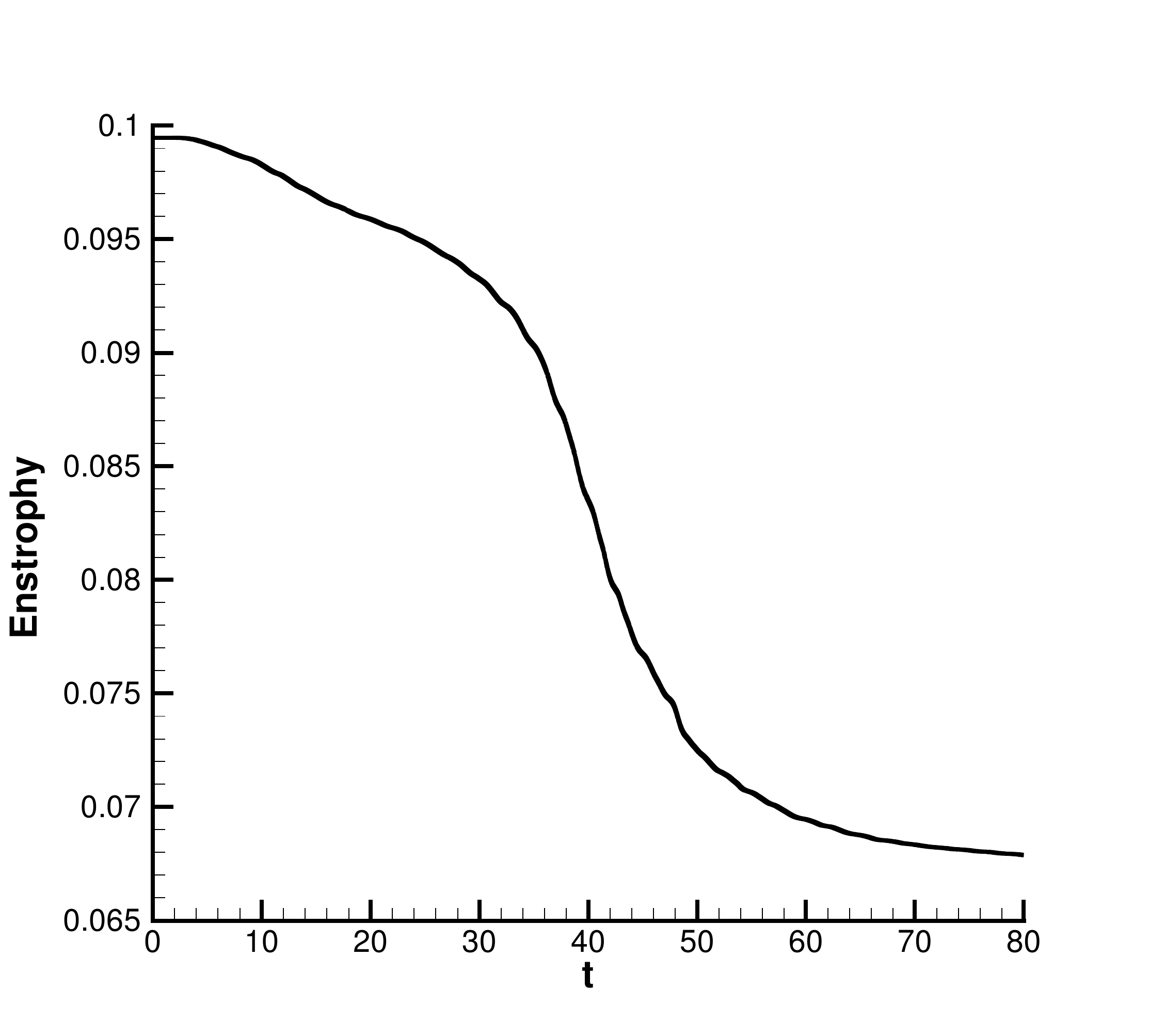}}
	\end{center}
	\caption{2D2V strong Landau damping with $\alpha=0.5$. Time evolution of electric energy, relative error of mass and total energy, error of momentum and enstrophy. Sparse grid: $N=7$, $k=3$.  }
	\label{2d2v_ener_mass_te2}
\end{figure}
\begin{figure}[htb]
	\begin{center}
		\subfigure[$x_2=2 \pi, \, \xi_2=0, \,  t=1.$]{\includegraphics[width=3in,angle=0]{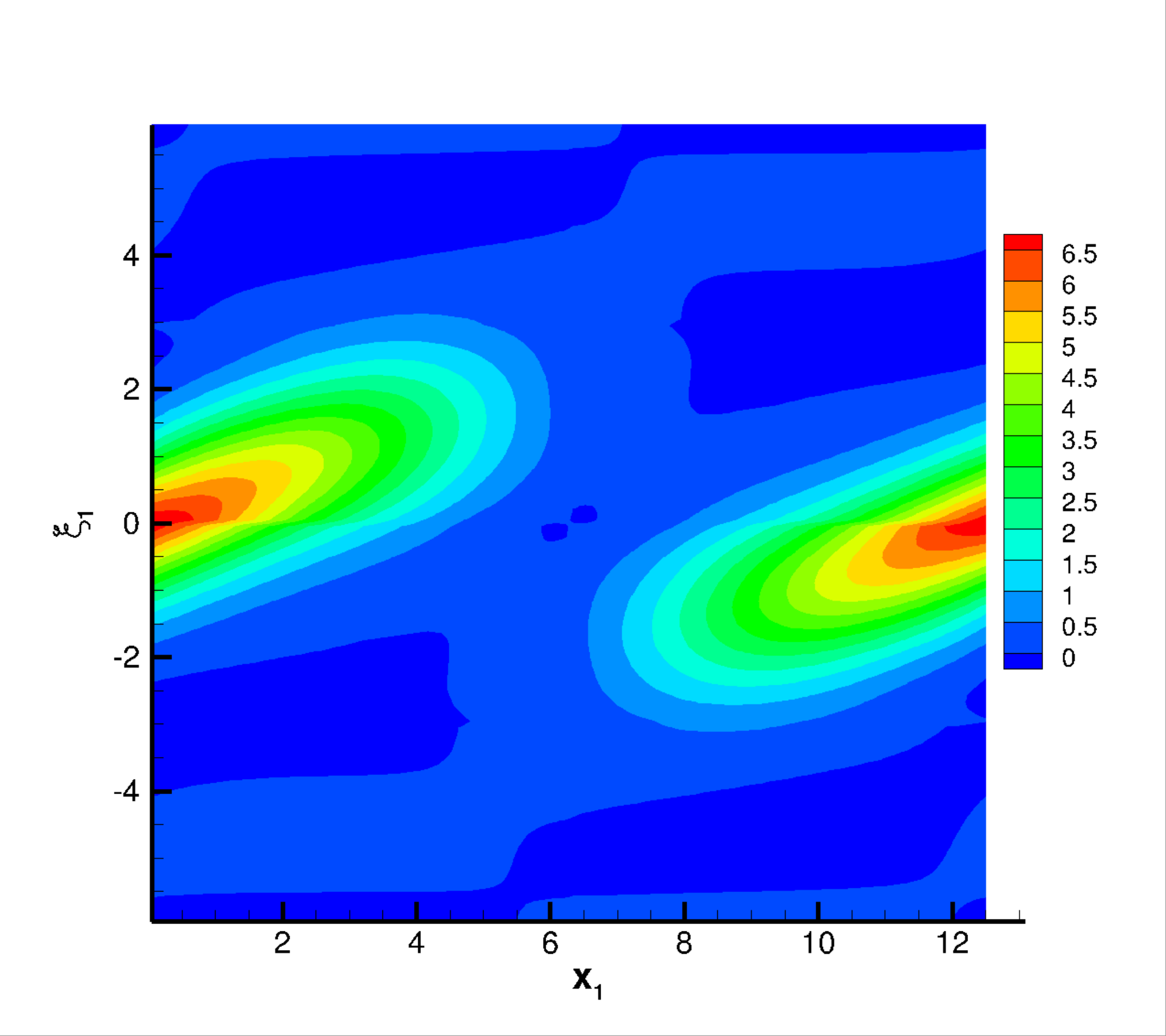}}
		\subfigure[$x_2=2 \pi, \, \xi_2=0, \,  t=5.$]{\includegraphics[width=3in,angle=0]{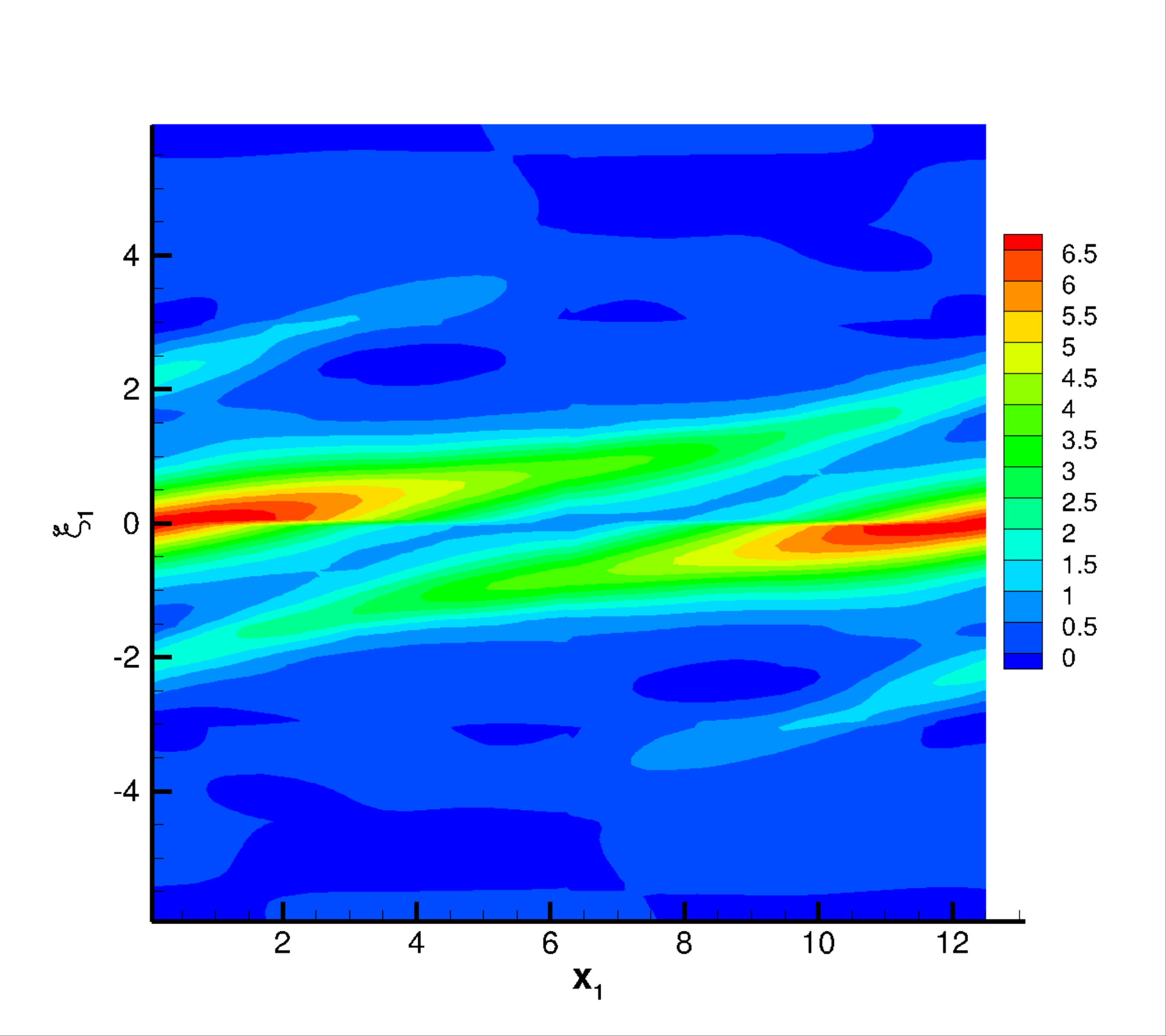}}
		\subfigure[$x_2=2 \pi, \, \xi_2=0, \,  t=10.$]{\includegraphics[width=3in,angle=0]{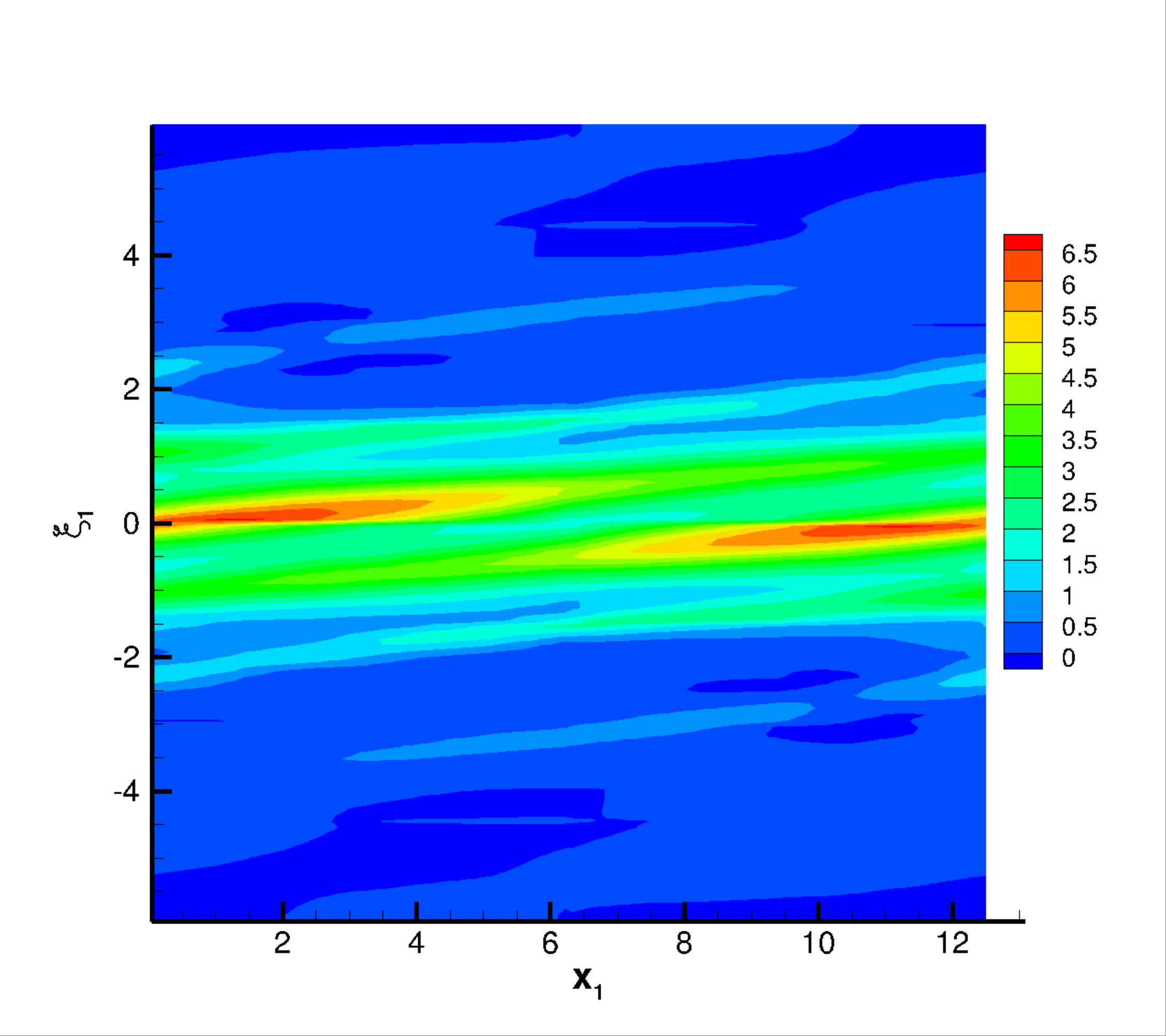}}
		\subfigure[$x_2=2 \pi, \, \xi_2=0, \,  t=20.$]{\includegraphics[width=3in,angle=0]{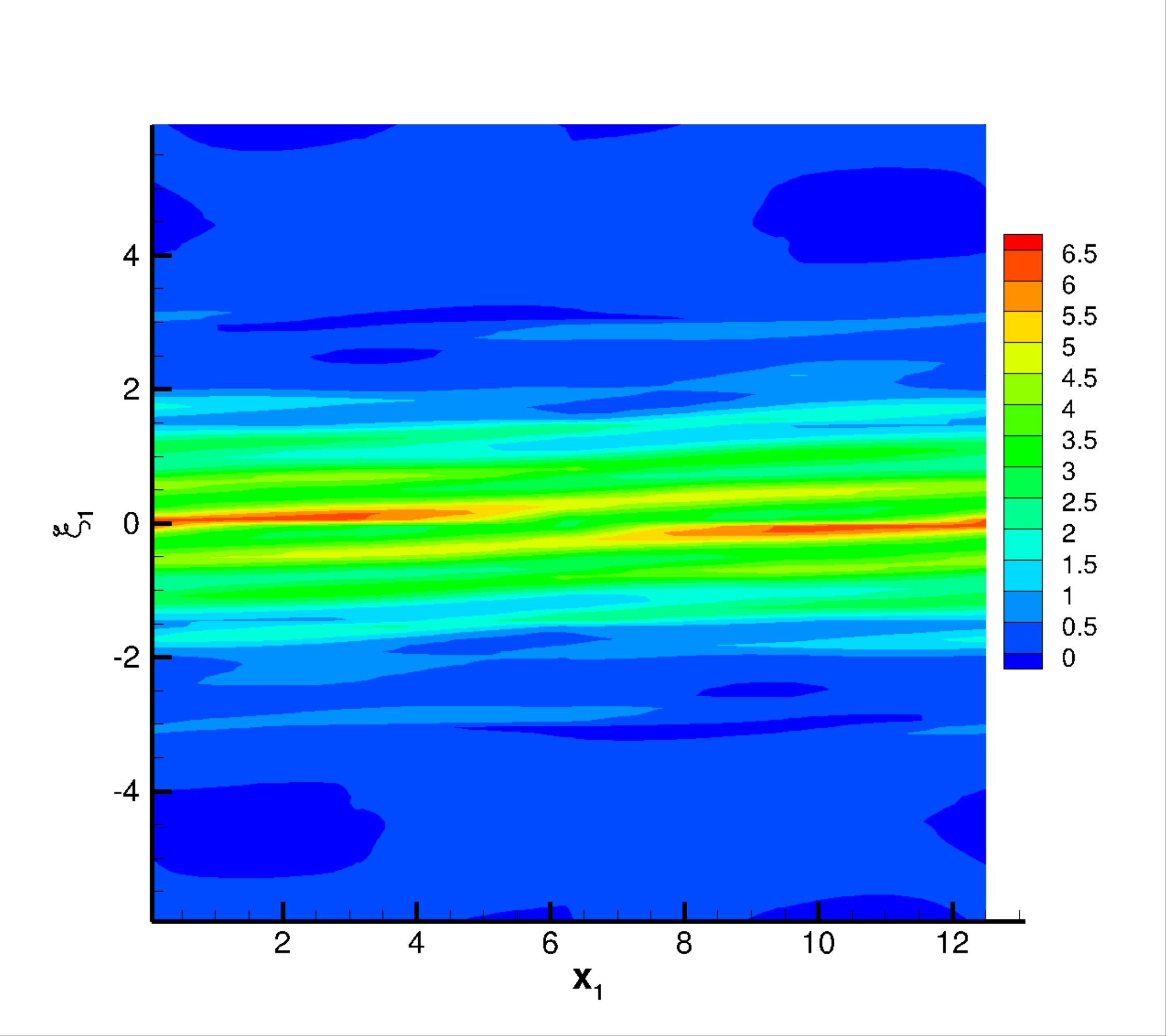}}
	\end{center}
	\caption{2D2V strong Landau damping with $\alpha=0.5$. 2D contour plots of the computed distribution function $f_h$ at selected time $t$. Sparse grid: $N=7$, $k=3$. }
	\label{2d2v_contour2}
\end{figure}

If the value of $\alpha$ increases, nonlinear effects start to dominate, and as such linear theory is no longer a good approximation. 
We show the electric energy, relative error in mass and total energy, error in momentum and enstrophy in Figure  \ref{2d2v_ener_mass_te2} up to $T=80$ for strong Landau damping. The same parameters $N=7, k=3$ are used as those in the weak Landau damping simulations. 
In this example, a larger error in conserved quantities is observed. The largest relative error is on the order of $10^{-5}$ and $10^{-4}$ in mass and total energy, respectively. A decrease in the  enstrophy is on the order of $10^{-2}.$
In Figure \ref{2d2v_contour2}, we further present the contour plots of distribution function $f$ at $x_2=2 \pi, \, \xi_2=0$ at several instances of time. From the plots of $f$, we can observe that fine structures are generated at a later time. We remark that the computational results can be improved if we use the adaptive method or a more refined mesh level in the velocity direction as in 
\cite{kormann2016sparse}.

%
%

\section{Conclusions and future work}
\label{sec:conclusion}
In this paper, we developed sparse grid and adaptive sparse grid DG schemes for simulating VM equations.  
The key ingredients of the scheme are the multiwavelet tensorized finite element approximation space and its hierarchical basis representation. Based on this construction, it is possible to reduce the storage and computational cost in high-dimensional simulations. 
Numerical tests show that our schemes achieve the desired of accuracy and the adaptive scheme is very efficient in capturing the solution structure. In the future, we will extend the scheme to other high-dimensional kinetic simulations.

\bibliographystyle{abbrv}
\bibliography{ref_sparse_VM,refer,ref_cheng_plasma_2,ref_cheng,ref_cheng_2}

\end{document}